\documentclass[twoside]{article}

\usepackage{tocbibind}
\usepackage{fancyhdr,bm,mathtools,graphicx}
\usepackage{amssymb,amsmath}
\usepackage{marvosym}
\usepackage{mathstyle}
\usepackage{hieroglf}
\usepackage[T1]{fontenc}
\usepackage{eufrak}
\usepackage{romanbar}
\usepackage[dvipsnames]{xcolor}
\usepackage{ifthen}
\usepackage{tikz-cd}
\tikzcdset{
arrow style=tikz,
diagrams={>={Stealth[scale=1]}}
}
\usepackage[T1]{fontenc}
\usepackage{textcomp}

\usepackage{stmaryrd}
\usepackage[mathscr]{euscript}
\usepackage{chngcntr}
\usepackage[style]{fncychap}
\usepackage{latexsym}

\def\Ref#1{[{equation~\ref{#1}}]}

\newcounter{chapterr}
\counterwithin{section}{chapterr}

\newcounter{definition}
\counterwithin{definition}{chapterr}
\def\definition
{\refstepcounter{definition}\vskip0.4\baselineskip\noindent{{\textbf{Definition~}}}{$\bf\arabic{chapterr}.$}{$\bf\arabic{definition}$}:\hskip 0.4\parindent}
\newcounter{theorem}
\counterwithin{theorem}{chapterr}
\def\theorem
{\refstepcounter{theorem}\vskip0.4\baselineskip\noindent{{\textbf{Theorem~}}}{$\bf\arabic{chapterr}.$}{$\bf\arabic{theorem}$}:\hskip 0.4\parindent}
\newcounter{symbol}
\counterwithin{symbol}{chapterr}

\newcounter{QMaxiom}
\counterwithin{QMaxiom}{chapterr}

\newcounter{corollary}
\counterwithin{corollary}{chapterr}
\def\corollary
{\refstepcounter{corollary}\vskip0.4\baselineskip\noindent{{\textbf{Corollary~}}}{$\bf\arabic{chapterr}.$}{$\bf\arabic{corollary}$}:\hskip 0.4\parindent}
\newcounter{lemma}
\counterwithin{lemma}{chapterr}
\def\lemma
{\refstepcounter{lemma}\vskip0.4\baselineskip\noindent{{\textbf{Lemma~}}}{$\bf\arabic{chapterr}.$}{$\bf\arabic{lemma}$}:\hskip 0.4\parindent}
\newcounter{proposition}
\counterwithin{proposition}{chapterr}
\def\proposition
{\refstepcounter{proposition}\vskip0.4\baselineskip\noindent{{\textbf{Proposition~}}}{$\bf\arabic{chapterr}.$}{$\bf\arabic{proposition}$}:\hskip 0.4\parindent}
\newcounter{remark}
\counterwithin{remark}{chapterr}
\def\remark
{\refstepcounter{remark}\vskip0.4\baselineskip\noindent{{\textbf{Remark~}}}{$\bf\arabic{chapterr}.$}{$\bf\arabic{remark}$}:\hskip 0.4\parindent}
\def\proof
{\vskip0.3\baselineskip\noindent{\textbf{proof}}\hskip.08\baselineskip:\hskip 0.4\parindent}

\newcounter{fixed}
\counterwithin{fixed}{chapterr}
\def\fixed
{\refstepcounter{fixed}\vskip0.4\baselineskip\noindent{{\textbf{Fixed Objects~}}}{$\bf\arabic{chapterr}.$}{$\bf\arabic{fixed}$}:\hskip 0.4\parindent}
\def\caution{\textasteriskcentered\hskip0.25\baselineskip}
\def\V{V}
\def\C{\mathbb{C}}

\def\c{c}

\def\R{\mathbb{R}}

\def\Z{\mathbb{Z}}

\newcommand{\p}[1]{{#1^{\prime}}}

\def\F{\mathscr{F}}
\def\C{\mathbb{C}}
\def\R{\mathbb{R}}

\def\card{{\textsf{card}}}
\def\endef{{\flushright{\noindent$\blacksquare$\\}}\noindent}
\def\endthm{{\flushright{\noindent$\square$\\}}\noindent}
\def\endcor{{\flushright{\noindent$\square$\\}}\noindent}

\def\endlem{{\flushright{\noindent$\square$\\}}\noindent}
\def\endpro{{\flushright{\noindent$\square$\\}}\noindent}

\def\endlem{{\flushright{\noindent$\square$\\}}\noindent}
\def\endfixed{{\flushright{\noindent$\Diamond$\\}}\noindent}
\def\endremark{\vskip0.5\baselineskip}

\def\then{\Rightarrow}
\def\thenn{\Leftrightarrow}

\def\({\left(}
\def\){\right)}
\def\[{\left[}
\def\]{\right]}

\newcommand{\CSs}[1]{{\mathcal{P}}\(#1\)}

\newcommand{\Card}[1]{\left| #1\right|}
\newcommand{\CarD}[1]{\card\(#1\)}
\def\cardeq{\overset{\underset{\mathrm{card}}{}}{=}}

\newcommand{\union}[1]{\bigcup#1}
\newcommand{\intersection}[1]{\bigcap#1}

\newcommand{\Union}[3]{\bigcup_{#1\in#2}#3}

\newcommand{\Dproduct}[3]{\prod_{{#1}\in{#2}}{#3}}
\newcommand{\cmp}[2]{#1\circ#2}
\newcommand{\Func}[2]{{\textsf{F}}\bpair{#1}{#2}}
\newcommand{\surFunc}[2]{{{\mathfrak{s}}\textsf{F}}\opair{#1}{#2}}
\newcommand{\IF}[2]{{\mathsf{B\hphantom{}F}}\opair{#1}{#2}}

\newcommand{\domain}[1]{{\mathsf{dom}}\bbsingle{#1}}
\newcommand{\codomain}[1]{{\mathsf{codom}}\bbsingle{#1}}
\newcommand{\funcimage}[1]{{\mathsf{img}}\bbsingle{#1}}
\newcommand{\EqR}[1]{{\mathsf{EqR}}\bbsingle{#1}}
\newcommand{\EqClass}[2]{{#1}\left/{#2}\right.}
\newcommand{\pEqclass}[2]{\[{#1}\]_{#2}}
\newcommand{\PEqclass}[2]{{\mathsf{EqC}}\bbpair{#1}{#2}}
\newcommand{\Cprod}[2]{{#1}\times{#2}}

\newcommand{\OR}[2]{#1\thinspace\lor\thinspace#2}
\newcommand{\AND}[2]{#1\thinspace\land\thinspace#2}
\newcommand*{\suchthat}{\;\ifnum\currentgrouptype=16 \middle\fi|\;}
\newcommand{\Foreach}[2]{\forall\thinspace#1\in#2\negthinspace:\thickspace}
\newcommand{\Exists}[2]{\exists\thinspace#1\in#2\negthinspace:\thickspace}

\newcommand{\defset}[3]{\left\{#1\in#2\thickspace:\thickspace#3\right\}}
\newcommand{\defsets}[3]{\left\{#1\subseteq#2\thickspace:\thickspace#3\right\}}
\newcommand{\defSet}[2]{\left\{#1\thickspace|\thickspace#2\right\}}
\newcommand{\negation}[1]{\neg{#1}}
\def\eqdef{\overset{\underset{\mathrm{def}}{}}{=}}
\def\indef{\overset{\underset{\mathrm{def}}{}}{\in}}
\newcommand{\resd}[1]{{\mathfrak{res}}{\mathsf{D}}_{#1}}
\newcommand{\rescd}[1]{{\mathfrak{res}}{\mathsf{C\negthinspace D}}_{#1}}
\newcommand{\res}[1]{{\mathfrak{res}}_{#1}}
\newcommand{\finv}[1]{{#1}^{-1}}
\newcommand{\Res}[2]{{\textsf{res}}\left({#1};{#2}\right)}

\def\index{{\mathscr{I}}}
\def\Zp{\Z^{+}}

\def\empty{\varnothing}
\def\Xt{{\mathbb{X}}}
\def\Yt{{\mathbb{Y}}}
\newcommand{\opair}[2]{\(#1,\thinspace #2\)}

\newcommand{\bpair}[2]{\negthinspace\left(#1,\thinspace #2\right)}
\newcommand{\bbpair}[2]{\negthinspace\left(#1;\thinspace #2\right)}
\newcommand{\bbsingle}[1]{\negthinspace\left({#1}\right)}
\newcommand{\cpair}[2]{\[#1,\thinspace #2\]}

\newcommand{\alltopologies}[1]{{\mathbf{Top}}\bbsingle{#1}}

\newcommand{\closedsets}[1]{{\mathsf{Closed}}\bbsingle{#1}\negthinspace}
\newcommand{\compl}[2]{#1\setminus#2}

\newcommand{\seta}[1]{\left\{#1\right\}}
\newcommand{\func}[2]{#1\(#2\)}

\def\x{x}

\def\U{U}

\def\asubset{A}

\def\point{p}

\newcommand{\nei}[1]{{\mathsf{Nei}}_{#1}}

\newcommand{\Cl}[1]{{\mathsf{Cl}}_{#1}}

\def\cf{f}
\def\cg{g}

\newcommand{\image}[1]{#1^{\rightarrow}}
\newcommand{\pimage}[1]{#1^{\leftarrow}}
\newcommand{\CF}[2]{{\mathsf{CF}}\bpair{#1}{#2}}
\newcommand{\HOF}[2]{{\mathsf{HF}}\bpair{#1}{#2}}

\newcommand{\function}[3]{#1:\thinspace#2\to#3}

\newcommand{\binary}[2]{#1,\thinspace#2}

\newcommand{\connecteds}[1]{{\mathsf{Cnd}}\bbsingle{#1}}
\newcommand{\maxcon}[1]{{\mathsf{MaxCnd}}\bbsingle{#1}}

\newcommand{\identity}[1]{{\rm{Id}}_{#1}}

\newcommand{\VVS}[1]{\VS_{#1}}

\newcommand{\WW}[1]{W_{#1}}

\newcommand{\subvec}[2]{{\mathsf{SubVec}}\bbpair{#1}{#2}}
\newcommand{\Vspan}[1]{{\mathsf{span}}_{#1}}
\def\VS{{\mathbb{V}}}
\newcommand{\NVS}[1]{{\mathscr{V}}_{#1}}

\newcommand{\vsum}[1]{+_{#1}}
\newcommand{\vsprod}[1]{\times_{#1}}

\newcommand{\Lin}[2]{{\mathsf{L}}\bpair{#1}{#2}}
\newcommand{\Linisom}[2]{{\mathsf{LIsom}}\bpair{#1}{#2}}

\newcommand{\VLin}[2]{{\mathbb{L}}\bpair{#1}{#2}}

\newcommand{\ovecbasis}[1]{{\mathsf{O\negthinspace VBases}}\bbsingle{#1}}

\newcommand{\Det}[1]{{\mathrm{det}}_{#1}}
\newcommand{\triple}[3]{\opair{#1}{\binary{#2}{#3}}}

\newcommand{\mtuple}[2]{\(\suc{#1}{#2}\)}

\newcommand{\norm}[2]{\left\lVert{#1}\right\rVert_{#2}}

\newcommand{\vv}[1]{v_{#1}}
\newcommand{\ww}[1]{w_{#1}}
\newcommand{\uu}[1]{u_{#1}}

\newcommand{\abs}[1]{\left|{#1}\right|}

\newcommand{\suc}[2]{{#1},\ldots,{#2}}

\newcommand{\Times}[2]{{#1}\times\ldots\times{#2}}

\def\RR{\tilde{\R}}
\newcommand{\Ropenman}[2]{{\tilde{#1}}_{\RR^{#2}}}
\newcommand{\topR}[1]{\widehat{#1}}

\newcommand{\funcprod}[2]{{#1}{\overset{\underset{\mathrm{\mathfrak{f}}}{}}{\times}}{#2}}
\newcommand{\Injection}[2]{{\mathrm{Inj}}_{{#1}\to{#2}}}
\newcommand{\mantop}[1]{{{\mathsf{Top}}^{\mathfrak{m}}}\bbsingle{#1}}
\newcommand{\mantops}[1]{\widehat{#1}}
\newcommand{\NVLin}[2]{{\overline{\mathbb{L}}}\bpair{#1}{#2}}
\newcommand{\Mat}[3]{{\mathsf{MAT}}\bbpair{#1}{\binary{#2}{#3}}}
\newcommand{\BMat}[3]{{\mathbb{MAT}}\bbpair{#1}{\binary{#2}{#3}}}
\newcommand{\TMat}[3]{{\widehat{\mathbb{MAT}}}\bbpair{#1}{\binary{#2}{#3}}}

\newcommand{\Man}[1]{\mathscr{M}_{#1}}
\newcommand{\Emsubman}[1]{{\mathsf{EmSubMan}}\bbsingle{#1}}
\newcommand{\emsubman}[2]{{\mathsf{emsubman}}\bbpair{#1}{#2}}

\newcommand{\charttransfer}[3]{{\mathsf{chart}}^{\star}_{\triple{#1}{#2}{#3}}}
\newcommand{\manprod}[2]{{#1}{\overset{\underset{\mathrm{\mathfrak{m}}}{}}{\times}}{#2}}
\newcommand{\topprod}[2]{{#1}{\overset{\underset{\mathrm{\mathfrak{t}}}{}}{\times}}{#2}}
\newcommand{\topologyofspace}[1]{{\mathcal{T}}\bbsingle{#1}}
\newcommand{\maxatlases}[3]{{\mathsf{maxAtl}}^{\(#1\)}\bpair{#2}{#3}}
\newcommand{\atlases}[3]{{\mathsf{Atl}}^{\(#1\)}\bpair{#2}{#3}}

\newcommand{\maxatlasgen}[3]{{\mathfrak{maxAtl}}^{\(#1\)}_{\opair{#2}{#3}}}
\newcommand{\M}[1]{\mathrm{M}_{#1}}
\newcommand{\maxatlas}[1]{{\boldsymbol{\mathscr{A}}}_{#1}}
\newcommand{\atlas}[1]{\mathscr{A}_{#1}}
\newcommand{\difclass}[1]{\mathsf{C}^{#1}}
\newcommand{\vecf}[2]{{\mathsf{V\negthinspace F}}^{#2}\bbsingle{#1}}
\newcommand{\Vecf}[2]{{\mathbb{V\negthinspace F}}^{#2}\bbsingle{#1}}

\newcommand{\avecf}[1]{{\mathscr{X}}_{#1}}

\newcommand{\tanspace}[2]{{\mathsf{T}}_{#1}{#2}}
\newcommand{\Tanspace}[2]{{\mathbb{T}}_{#1}{#2}}
\newcommand{\avec}[1]{{\boldsymbol{v}}_{#1}}
\newcommand{\tanbun}[1]{{\mathsf{T}}#1}
\newcommand{\Tanbun}[1]{{\mathbf{T}}#1}
\newcommand{\basep}[1]{\pi_{#1}}
\newcommand{\mapdifclass}[3]{{\mathsf{C}}^{#1}\bpair{#2}{#3}}
\newcommand{\Lmapdifclass}[3]{{\mathbf{C}}^{#1}\bpair{#2}{#3}}
\newcommand{\Diffeo}[3]{{\mathsf{Diff}}^{#1}\bpair{#2}{#3}}
\newcommand{\Diff}[2]{{\mathsf{Diff}}^{#1}\bbsingle{#2}}
\newcommand{\GDiff}[2]{{\mathbf{Diff}}^{#1}\bbsingle{#2}}
\newcommand{\banachmapdifclass}[5]{{\mathsf{C}}^{#1}_{\opair{#2}{#3}}\bpair{#4}{#5}}
\newcommand{\Derivation}[2]{{\mathsf{Der}}^{#2}\bbsingle{#1}}
\newcommand{\LDerivation}[2]{{\mathbb{D}}^{#2}\bbsingle{#1}}
\newcommand{\aderivation}[1]{\Delta_{#1}}
\newcommand{\map}[3]{{#1}\negthinspace:{#2}\to{#2}}
\def\rdot{\cdot}
\newcommand{\tanchart}[2]{{\mathscr T}^{#1}_{#2}}
\newcommand{\tanatlas}[1]{{\mathfrak{TA}}_{#1}}
\newcommand{\tanspaceiso}[3]{{\boldsymbol{\theta}}_{#1}^{\opair{#2}{#3}}}
\newcommand{\der}[3]{{\mathsf{D}}^{\opair{#2}{#3}}{#1}}
\newcommand{\derop}[2]{{\mathsf{D}}^{\opair{#1}{#2}}}
\newcommand{\derr}[1]{{\mathsf{D}}{#1}}
\newcommand{\banachder}[3]{{\mathscr{D}}^{\opair{#2}{#3}}{#1}}
\newcommand{\banachderr}[1]{{\mathscr{D}}{#1}}
\newcommand{\Rder}[2]{{\mathsf{d}}_{#2}{#1}}
\newcommand{\Rderop}[1]{{\mathsf{d}}_{#1}}
\newcommand{\projection}[2]{{\boldsymbol{\mathscr{P}}}^{\(#1\)}_{#2}}

\newcommand{\Eucbase}[2]{{\boldsymbol{e}}^{\(#1\)}_{#2}}
\newcommand{\deltaf}[2]{\delta_{{#1},{#2}}}

\newcommand{\zerovec}[1]{{\boldsymbol{0}}_{#1}}
\newcommand{\mandim}[1]{{\mathrm{dim}}_{\mathrm{man}}\bbsingle{#1}}

\newcommand{\immersion}[3]{{\mathsf{Immersion}}^{#1}\bpair{#2}{#3}}
\newcommand{\submersion}[3]{{\mathsf{Submersion}}^{#1}\bpair{#2}{#3}}
\newcommand{\embedding}[3]{{\mathsf{Embedding}}^{#1}\bpair{#2}{#3}}
\newcommand{\leftparinj}[3]{{\mathfrak{LI}}^{\opair{#1}{#2}}_{#3}}
\newcommand{\rightparinj}[3]{{\mathfrak{RI}}^{\opair{#1}{#2}}_{#3}}
\newcommand{\prodmantan}[4]{{\boldsymbol{\Theta}}^{\opair{#1}{#2}}_{\opair{#3}{#4}}}

\newcommand{\Group}[1]{{\mathbb{G}}_{#1}}

\newcommand{\IndTop}[1]{{\mathsf{IndTop}}_{#1}}

\newcommand{\GHom}[2]{{\mathsf{GHom}}\bpair{#1}{#2}}
\newcommand{\GIsom}[2]{{\mathsf{GIsom}}\bpair{#1}{#2}}
\newcommand{\Subgroups}[1]{{\mathsf{SubGr}}\bbsingle{#1}}

\newcommand{\asubgroup}[1]{{\mathrm{H}}_{#1}}
\newcommand{\LCoset}[1]{{\mathsf{Lcoset}}_{#1}}

\newcommand{\smoothring}[2]{{\mathscr{R}}\difclass{#1}\bbsingle{#2}}

\newcommand{\OO}[1]{O_{#1}}
\newcommand{\TFB}[4]{{\mathsf{TF}}^{\infty}_{\opair{#1}{#2}}\bbpair{#3}{#4}}
\newcommand{\VTFB}[4]{{\mathbb{TF}}^{\infty}_{\opair{#1}{#2}}\bbpair{#3}{#4}}
\newcommand{\DVTFB}[2]{\overset{\underset{\mathrm{TF}}{}}{\bigotimes}\bbpair{#1}{#2}}
\newcommand{\tftensor}[1]{\overset{\underset{\mathrm{#1}}{}}{\otimes}}
\newcommand{\BanachDiff}[5]{{\mathsf{Diff}}^{#1}_{\opair{#2}{#3}}\bpair{#4}{#5}}
\newcommand{\TFpullback}[3]{{#1}^{\star}_{\opair{#2}{#3}}}
\newcommand{\TFPullback}[1]{{#1}^{\star}}
\newcommand{\TFpullbackcov}[2]{{#1}^{\star}_{\(#2\)}}

\newcommand{\CRing}[1]{{\mathbb{K}}_{#1}}
\newcommand{\vecs}[1]{{\mathbb{W}}_{#1}}

\newcommand{\Cmodule}[1]{{\mathrm{M}}_{#1}}
\newcommand{\CModule}[1]{{\mathbb{M}}_{#1}}

\newcommand{\mlinmap}[1]{T_{#1}}
\newcommand{\multiprod}[2]{\underbrace{\Times{#1}{#1}}_{#2}}
\newcommand{\Vdual}[1]{{#1}^{\star}}
\newcommand{\Vddual}[1]{{#1}^{\star\negthinspace\star}}
\newcommand{\Tensors}[3]{{\mathbf{T}}_{#1}^{#2}\bbsingle{#3}}
\newcommand{\MTensors}[3]{{\mathbb{T}}_{#1}^{#2}\bbsingle{#3}}
\newcommand{\Dsum}[1]{\bigoplus#1}
\newcommand{\directsum}[2]{{#1}\oplus{#2}}
\newcommand{\DTensors}[1]{\bigotimes#1}
\newcommand{\tensor}[1]{\overset{\underset{\mathrm{#1}}{}}{\otimes}}
\newcommand{\vsbase}[1]{{\mathbf{e}}_{#1}}
\newcommand{\Bdual}[2]{{\overline{#1}}^{#2}}
\newcommand{\Bddual}[2]{\overline{\overline{#1}}_{#2}}
\newcommand{\tensorvsbase}[3]{{\mathbf{E}}^{\opair{#1}{#2}}_{#3}}
\newcommand{\dualpb}[1]{{#1}^{\ast}}
\newcommand{\Vpullback}[3]{{#1}^{\star}_{\opair{#2}{#3}}}
\newcommand{\VPullback}[1]{{#1}^{\star}}
\newcommand{\Vpullbackcov}[2]{{#1}^{\star}_{\(#2\)}}

\newcommand{\quadruple}[4]{\opair{#1}{\binary{#2}{\binary{#3}{#4}}}}
\newcommand{\fbundle}[1]{{\mathbb{F}}_{#1}}
\newcommand{\fbtotal}[1]{{\mathscr{E}}_{#1}}
\newcommand{\fbbase}[1]{{\mathscr{B}}_{#1}}
\newcommand{\fbfiber}[1]{{\mathscr{F}}_{#1}}
\newcommand{\fbprojection}[1]{{\boldsymbol{\pi}}_{#1}}
\newcommand{\Tot}[1]{{\mathrm{E}}_{#1}}
\newcommand{\B}[1]{{\mathrm{B}}_{#1}}
\newcommand{\Fib}[1]{{\mathrm{F}}_{#1}}
\newcommand{\subman}[2]{\left. #1\right|_{#2}}
\newcommand{\proj}[3]{{\mathrm{pr}}^{\opair{#1}{#2}}_{#3}}
\newcommand{\reS}[2]{\left. #1\right|_{#2}}
\newcommand{\fbatlas}[1]{{\mathsf{FB}}{\mathscr{A}}\bbsingle{#1}}
\newcommand{\plt}[2]{\Lambda_{#1}{#2}}
\newcommand{\pltfib}[2]{\reS{\widetilde{#1}}{#2}}
\newcommand{\transition}[3]{\overset{\underset{\mathrm{#1}}{}}{\boldsymbol{\chi}}\cpair{#2}{#3}}
\newcommand{\fbsec}[1]{\sigma_{#1}}
\newcommand{\fbsections}[1]{{\boldsymbol{\Gamma}}^{\infty}\bbsingle{#1}}
\newcommand{\fbsectionsl}[2]{{\boldsymbol{\Gamma}}^{\infty}\bbpair{#1}{#2}}

\newcommand{\fbmorb}[1]{\widehat{#1}}
\newcommand{\fbmorphisms}[2]{{\mathsf{FBMor}}\bpair{#1}{#2}}
\newcommand{\fbisomorphisms}[2]{{\mathsf{VBIsom}}\bpair{#1}{#2}}
\newcommand{\quintuple}[5]{\opair{#1}{\binary{#2}{\binary{#3}{\binary{#4}{#5}}}}}
\newcommand{\vbundle}[1]{{\mathbb{V}}_{#1}}
\newcommand{\vbtotal}[1]{{\mathscr{E}}_{#1}}
\newcommand{\vbbase}[1]{{\mathscr{B}}_{#1}}
\newcommand{\vbfiber}[1]{{\mathbb{X}}_{#1}}
\newcommand{\vbprojection}[1]{{\boldsymbol{\pi}}_{#1}}
\newcommand{\vTot}[1]{{\mathrm{E}}_{#1}}
\newcommand{\vB}[1]{{\mathrm{B}}_{#1}}

\newcommand{\vbatlas}[1]{{\mathscr{A}}_{#1}}
\newcommand{\Fvbundle}[1]{\widetilde{#1}}
\newcommand{\algfield}[1]{\Bbbk_{#1}}
\newcommand{\vecsmanifold}[1]{\widetilde{#1}}
\newcommand{\GL}[2]{{\mathsf{GL}}_{#2}\bbsingle{#1}}
\newcommand{\fibervecs}[2]{{\mathbf{Fib}}_{#1}\bbsingle{#2}}
\newcommand{\vbsec}[1]{\sigma_{#1}}
\newcommand{\vbsections}[1]{{\boldsymbol{\Gamma}}^{\infty}\bbsingle{#1}}
\newcommand{\vbsectionsl}[2]{{\boldsymbol{\Gamma}}^{\infty}\bbpair{#1}{#2}}
\newcommand{\vbsadd}[1]{\overset{\underset{\mathrm{#1}}{}}{+}}
\newcommand{\vbsprod}[1]{\overset{\underset{\mathrm{#1}}{}}{\times}}
\newcommand{\vbsscalprod}[1]{\overset{\underset{\mathrm{#1}}{}}{\cdot}}
\newcommand{\lframe}[2]{{\boldsymbol{\Phi}}\bbpair{#1}{#2}}
\newcommand{\vbTensors}[4]{{\mathcal{T}}_{#1}^{#2}\bbpair{#4}{#3}}
\newcommand{\VBTensors}[3]{{\mathcal{T}}_{#1}^{#2}\bbsingle{#3}}
\newcommand{\VBTensorsMan}[3]{{\boldsymbol{\mathscr{T}}}_{#1}^{#2}\bbsingle{#3}}
\newcommand{\vbchartlocalt}[1]{{\boldsymbol{\Omega}}\bbsingle{#1}}
\newcommand{\btriple}[3]{\bbpair{#1}{\opair{#2}{#3}}}
\newcommand{\vbtensorchart}[4]{{\Xi}^{\opair{#1}{#2}}_{\opair{#3}{#4}}}
\newcommand{\vbtensorlocalt}[3]{{\Delta}^{\opair{#1}{#2}}_{#3}}
\newcommand{\vbtensormaxatlas}[1]{{\mathsf{max}}\maxatlas{#1}}
\newcommand{\vbtensorprojection}[3]{{\boldsymbol{\pi}}_{#1}^{\opair{#2}{#3}}}
\newcommand{\vbtensoratlas}[3]{\boldsymbol{{\mathscr{A}}}_{#3}^{\opair{#1}{#2}}}
\newcommand{\vbtensorbundle}[3]{{\mathbb{TB}}^{\opair{#1}{#2}}\bbsingle{#3}}
\newcommand{\TF}[3]{{\mathsf{TF}}^{\infty}_{\opair{#1}{#2}}\bbsingle{#3}}
\newcommand{\MTF}[3]{{\mathscr{M}}{\mathsf{TF}}^{\infty}_{\opair{#1}{#2}}\bbsingle{#3}}
\newcommand{\VTF}[3]{{\mathbb{TF}}^{\infty}_{\opair{#1}{#2}}\bbsingle{#3}}
\newcommand{\DVTF}[1]{\overset{\underset{\mathrm{VTF}}{}}{\bigotimes}#1}

\newcommand{\vbtensor}[1]{\overset{\underset{\mathrm{#1}}{}}{\otimes}}
\newcommand{\atf}[1]{{\mathscr{X}}_{#1}}
\newcommand{\atff}[1]{{\mathscr{Y}}_{#1}}
\newcommand{\vbmorphisms}[2]{{\mathsf{VBMor}}\bpair{#1}{#2}}
\newcommand{\vbisomorphisms}[2]{{\mathsf{VBIsom}}\bpair{#1}{#2}}
\newcommand{\VBpullback}[3]{{#1}^{\star}_{\opair{#2}{#3}}}
\newcommand{\VBPullback}[1]{{#1}^{\star}}
\newcommand{\VBpullbackcov}[2]{{#1}^{\star}_{\(#2\)}}

\newcommand{\subvbatlas}[3]{{\mathsf{sub}}{\boldsymbol{\mathscr{A}}}\bbpair{#1}{\binary{#2}{#3}}}
\newcommand{\vbHom}[3]{{\mathcal{H}}\btriple{#3}{#1}{#2}}
\newcommand{\VBHom}[2]{{\mathcal{H}}\opair{#1}{#2}}
\newcommand{\vbHomchart}[3]{{\mathbf{H}}_{\triple{#1}{#2}{#3}}}
\newcommand{\Endo}[1]{{\mathsf{End}}\bbsingle{#1}}
\newcommand{\vbchartlocaltt}[2]{\overline{\boldsymbol{\Omega}}\opair{#1}{#2}}
\newcommand{\vbHomlocalt}[3]{{\overline{\bf{H}}}_{\triple{#1}{#2}{#3}}}
\newcommand{\vbHommaxatlas}[2]{{\mathscr{H}}{\mathsf{max}}\maxatlas{\opair{#1}{#2}}}
\newcommand{\vbHomMan}[2]{{\boldsymbol{\mathscr{H}}}\opair{#1}{#2}}
\newcommand{\vbHomprojection}[2]{{\boldsymbol{\pi}}_{\opair{#1}{#2}}}
\newcommand{\vbHomatlas}[2]{\boldsymbol{{\mathscr{A}}}^{\mathscr{H}}_{\opair{#1}{#2}}}
\newcommand{\vbHombundle}[2]{{\mathbb{HOM}}\opair{#1}{#2}}
\newcommand{\vbDualbundle}[1]{{#1}^{\star}}
\newcommand{\lframedual}[2]{{\overline{#1}}^{#2}}

\newcommand{\vbSum}[3]{{\mathcal{W}}\btriple{#3}{#1}{#2}}
\newcommand{\VBSum}[2]{{\mathcal{W}}\opair{#1}{#2}}
\newcommand{\vbSumchart}[3]{{\mathbf{w}}_{\triple{#1}{#2}{#3}}}
\newcommand{\vbSumlocalt}[3]{{\overline{\bf{w}}}_{\triple{#1}{#2}{#3}}}
\newcommand{\vbSummaxatlas}[2]{{\mathscr{W}}{\mathsf{max}}\maxatlas{\opair{#1}{#2}}}
\newcommand{\vbSumMan}[2]{{\boldsymbol{\mathscr{W}}}\opair{#1}{#2}}
\newcommand{\vbSumprojection}[2]{{\boldsymbol{\pi}}_{\opair{#1}{#2}}}
\newcommand{\vbSumatlas}[2]{\boldsymbol{{\mathscr{A}}}^{\mathscr{W}}_{\opair{#1}{#2}}}
\newcommand{\vbSumbundle}[2]{#1\oplus#2}
\newcommand{\vbPro}[4]{{\mathcal{P}}\btriple{\binary{#3}{#4}}{#1}{#2}}
\newcommand{\VBPro}[2]{{\mathcal{P}}\opair{#1}{#2}}
\newcommand{\provbchartlocaltt}[2]{\overline{\overline{\boldsymbol{\Omega}}}\opair{#1}{#2}}
\newcommand{\vbProlocalt}[4]{{\overline{\bf{p}}}_{\quadruple{#1}{#2}{#3}{#4}}}
\newcommand{\vbProchart}[4]{{\mathbf{P}}_{\quadruple{#1}{#2}{#3}{#4}}}
\newcommand{\vbPromaxatlas}[2]{{\mathscr{P}}{\mathsf{max}}\maxatlas{\opair{#1}{#2}}}
\newcommand{\vbProMan}[2]{{\boldsymbol{\mathscr{P}}}\opair{#1}{#2}}
\newcommand{\vbProprojection}[2]{{\boldsymbol{\pi}}_{\opair{#1}{#2}}}
\newcommand{\vbProatlas}[2]{\boldsymbol{{\mathscr{A}}}^{\mathscr{P}}_{\opair{#1}{#2}}}
\newcommand{\vbProbundle}[2]{{#1}\times{#2}}
\newcommand{\inducedbundle}[3]{{#1}^{\star}\bbpair{#2}{#3}}
\newcommand{\inducedvbprojection}[2]{{#1}^{\star}_{#2}}
\newcommand{\inducedvblocalt}[2]{\boldsymbol{\mathscr{I}}\opair{#1}{#2}}
\newcommand{\vbIndatlas}[2]{\boldsymbol{{\mathscr{A}}}^{\mathrm{IND}}_{\bbpair{#1}{#2}}}
\newcommand{\vbIndbundle}[3]{{\mathbb{IND}}\bbpair{\binary{#1}{#2}}{#3}}
\newcommand{\vbResatlas}[2]{\boldsymbol{{\mathscr{A}}}^{\mathrm{RES}}_{\bbpair{#1}{#2}}}
\newcommand{\vbResbundle}[2]{{\mathbb{RES}}\bbpair{#1}{#2}}
\newcommand{\myitem}[1]{${\fontsize{6.65}{7}\selectfont{\textbf{#1}}}$}
\def\varfill{\dotfill}

\def\toclevel@section{1}\def\toclevel@subection{2}
\addtocontents{toc}{\string\let\string\l@section\string\l@subsection}
\def\toclevel@subsection{2}\def\toclevel@subsubection{3}
\addtocontents{toc}{\string\let\string\l@subsection\string\l@subsubsection}
\usepackage{tocloft}
\newcommand{\chapteR}[1]{\cleardoublepage
{\refstepcounter{chapterr}\vskip\baselineskip\centering{\fontsize{21}{21}\selectfont${\bf\Roman{chapterr}}$
\vskip0.6\baselineskip{\fontsize{21}{21}\selectfont{\textbf{#1}}}}
\vskip 5.9\baselineskip}
\addcontentsline{toc}{0}
{\protect\vskip0.5\baselineskip\noindent\bf{\Roman{chapterr}\hskip0.5\baselineskip#1\hspace{\fill}}}\par
\fancyhead[LO]{\ifthenelse{\value{chapterr}=0}{#1}{$\bf\Roman{chapterr}$\hskip0.7\baselineskip#1}}
}
\newcommand{\Bibliography}[1]{\vskip0.5\baselineskip\centering{\huge\bf{References}}\vskip \baselineskip
\addcontentsline{toc}{0}
{\protect\vskip0.5\baselineskip\noindent\bf{References}\hspace{\fill}}\par
\fancyhead[LO,RE]{\bf{References}#1}
}

\def\refthm#1{[{theorem~\ref{#1}}]}
\def\reflem#1{[{lemma~\ref{#1}}]}
\def\refdef#1{[{definition~\ref{#1}}]}
\def\refcor#1{[{corollary~\ref{#1}}]}
\def\refpro#1{[{proposition~\ref{#1}}]}

\newcommand{\mathleft}{\@fleqntrue\@mathmargin0pt}

\renewcommand{\sectionmark}[1]{\ifthenelse{\value{section}=0}{\markright{#1}{}}
{\markright{${\arabic{section}}$ #1}{}}}
\makeatletter
\def\cleardoublepage{\clearpage\if@twoside \ifodd\c@page\else
\hbox{}
\vspace*{\fill}
\vspace{\fill}
\thispagestyle{empty}
\newpage
\if@twocolumn\hbox{}\newpage\fi\fi\fi}
\makeatother
\fancypagestyle{plain}{
\fancyhf{}
\fancyhead{}
\fancyhead[R]{\thepage}
\fancyhead[L]{\leftmark}
}
\newcommand{\newsymp}[1]{{#1}\equiv}
\newcommand{\newsymb}[1]{{#1}\equiv}

\newcommand{\SET}[1]{{\mathscr{S}}_{#1}}
\newcommand{\FUNCTION}[1]{{{f}}_{#1}}

\newcommand{\prop}[1]{{\mathfrak{p}}\llparenthesis{#1}\rrparenthesis}
\def\Prop{\mathfrak{p}}
\newcommand{\propos}[1]{{\mathbf{\mathfrak{p}}}_{#1}}
\def\dummy{\centerdot}
\usepackage[hang,flushmargin]{footmisc}
\renewcommand{\footnoterule}{%
  \kern 20pt
  \hrule width \textwidth height 0.5pt
  \kern 5pt
}
\def\quotl{``}
\def\quotr{"}
\usepackage[inner=3.8cm,outer=3cm,bottom=3.9cm,twoside]{geometry}
\begin{document}
\thispagestyle{empty}
%\rule[0pt]{1.5pt}{200pt}
%\section{}
\noindent
{\\ \\\textbf{\fontsize{40}{40}\selectfont
{\textsf{Smooth Vector Bundles}}}}
%\\[0.8\baselineskip]{\textbf{\fontsize{20}{20}\selectfont
%{\textsf{of Finite Rank}}}}
\\[8\baselineskip]
\noindent
{\fontsize{21}{21}\selectfont
{\textsf{Farzad Shahi}}
}
\\[11\baselineskip]
\noindent
{\fontsize{11}{11}\selectfont
{\textrm{Version:} 1.00}
}
\vfill\hfill
$\underline{\Huge{\textsf{\bf F}}\negthickspace\negthickspace\negthinspace{{\rotatebox{90}{\textsf{\bf S}}}}}$
%{{\Huge{\Bat}}}
%\setlength{\unitlength}{2mm}
%\begin{picture}(30,20)
%\linethickness{0.075mm}
%\multiput(-100,0)(1,0){1}%
%{\line(0,1){20}}
%\end{picture}
%%%%%%%%%%%%%%%%%%%%%%%%%%%%%%%%%%%%%%%%%
\newpage
\thispagestyle{empty}
\noindent
{\fontsize{9.4}{9.4}\selectfont
{\underline{{\bf\textsf{Title:}} Smooth Vector Bundles}}}\\
{\fontsize{9.4}{9.4}\selectfont
{\underline{{\bf\textsf{Author:}} Farzad Shahi}}}\\
{{\bf\textsf{email}}}:~{\texttt{shahi.farzad@gmail.com}}
\vskip 4\baselineskip
\noindent
{\fontsize{9.4}{9.4}\selectfont
{\underline{{\bf\textsf{Version:}} 1.00}}
}\\
\noindent
{\fontsize{9.4}{9.4}\selectfont
{\bf\textsf{2022}}}
\vskip 4\baselineskip
\noindent
{\fontsize{9.4}{9.4}\selectfont
{\bf\textsf{Typesetting:}} By the author, using \TeX}
\vskip 5\baselineskip
\noindent
{\fontsize{9.4}{9.4}\selectfont
{\bf\textsf{Abstract:}}
The current version serves as an introduction to the basics of the theory of real and complex smooth vector bundles with finite rank,
and the general concept of a tensor field on a smooth vetor bundle.}
%%%%%%%%%%%%%%%%%%%%%%%%%%%%%%%%%%%%%%%%%%
\newpage
\thispagestyle{empty}
\noindent
{\fontsize{20}{20}\selectfont
{\bf Preface}}
\vskip1.5\baselineskip
\noindent
The practice of assigning a vector from $\R^3$ to each point of the space-time or a portion of it is a common exercise of
classical physics, famously known as the idea of a field. Among the well known examples are the notions of
the electric and magnetic fields in the theory of electromagnetism, the gravitational field in the classical theory of gravity,
and the velocity field in the theory of fluid dynamics. The framework for treating such problems can be conceived with
some further degree of subtlety; each point of the space-time or the portion under consideration can
be thought to be equipped with a copy of $\R^3$. So any field can be considered as a rule that assigns
to every point of space-time a vector from the corresponding vector space of that point. But taking such subtlety
into account in dealing with fields defined over the space-time with values in a fixed vector space, is absolutely
unnecessary. On the other hand, such a consideration is inherently required when the nature of the vector space of
interest depends on the point, and hence there is no fixed vector space among which the vectors of a field
can be chosen. Such problems arise, for instance, when motion is constrained to a bent surface in space,
and thus the tangent plane to every point of the surface well suits to the vector space of interest at that point.
By generalizing the scope of the problem to an abstract manifold, the corresponded generalization of this framework
meets the notion of the tangent bundle of an abstract manifold.\\
A plain vector bundle consists of a set called the base set, to each point of which a vector space is assigned among
a collection of isomorphic finite-dimensional vector spaces, and the vector space assigned to each point is
identified with that point called the fiber of the vector bundle at that point.
Hence the corresponding vector spaces of any pair of distinct points do not intersect.
The union of all vector spaces corresponded to points of the base set is called the total set of the vector bundle.
Moreover, the map sending each vector in the total set to its corresponding point in the base set is called the
projection map of the vector bundle. No further structure is assumed for a plain vector bundle.\\
There is a class of vector bundles that are endowed with topological structures in a naturally consistent manner.
The objects of this class are famously referred to as $\quotl$topological vector bundles$\quotr$.
A topological vector bundle is basically a plain vector bundle whose base and total sets are equipped with
suitable topological structures that make the projection map continuous, and locally trivialize the whole construction
by means of topological behavior. Precisely, when it is said that a topological vector bundle is locally trivialized
everywhere, it means that there is an open cover of the base space, the product of each element of which with a fixed
vector space isomorphic to all fibers is homeomorphic to the pre-image of that element in the total space under the projection map.
Intuitively, the naturality of the construction of a topological vector bundle lies mainly in the fact that
it provides a setting in which the fibers are arranged in a continuous way, and so that locally all fibers
can be replaced by a fixed vector space without altering the topological behavior. An structure-preserving map
between any pair of topological vector bundles is considered to be a continuous map between the total spaces of them
that sends fibers into fibers linearly. Topological vector bundles together with these structure-preserving maps between them
form a category.\\
There is a specific class of topological vector bundles that prepare the setting for dealing with differentiability
aspects of relevant problems. The objects of this class are widely known as smooth vector bundles.
A smooth vector bundle is basically a plain vector bundle having its base and total sets endowd with
suitable differentiable structures that make the projection map infinitely differentiable, and locally trivialize the whole construction
by means of differentiable behavior. Similar to the case of topological vector bundles, locally trivializing property of
a smooth vector bundle refers to the existence of an open cover of the base space, the product of each element of which with a fixed
vector space (endowed with its canonical differentiable structure) isomorphic to all fibers is diffeomorphic to the pre-image
of that element in the total space under the projection map. The common dimension of the linear structure of the fibers
is called the rank of the smooth vector bundle. The general purpose behind such constructions may be
the preparation of an environment for smoothly arranging the fibers, in a way that replacement of all fibers with a fixed isomorphic
vector space would not effect locally the structure up to diffeomorphism. An structure-preserving map
between any pair of smooth vector bundles is considered to be an infinitely differentiable map between the total spaces of them
that sends fibers into fibers linearly. Smooth vector bundles together with the structure-preserving maps between them
form a category, which is actually a sub-category of the category of topological vector bundles. The category of smooth vector bundles
is the main focus of this book.\\
The category of smooth vector bundles is a specific sub-category of a general category, named the category of smooth fiber bundles.
The structure of a smooth fiber bundle is similar to that of smooth vector bundles, with the exception that
the fibers are considered to be general smooth manifolds rather than vector spaces. So clearly ignoring the
linear structures of fibers of a smooth vector bundle, it has an intrinsic smooth fiber bundle structure. Thus, the study of
the structure of a general smooth fiber bundle can be viewd as the abstraction of the features of smooth vector bundles
not related to linear structures of fibers, and hence crucial to their study.\\
Generalizing the framewok for considering the idea of a field, from space-time to a general smooth manifold, the notion
of the tangent bundle of a smooth manifold will come into play naturally, as mentioned earlier. Moreover, the concept of a
field on a manifold is not regarded as just an arbitrary assignment of a vector to each point of the manifold
from the correponded tangent space of that point. Actually, the main goal of constructing manifolds is to prepare the setting
for dealing with properties related to smoothness. So the fields of interest on a smooth manifold are taken to be such rules of
assigning vectors from the tangent bundle to the points of the manifold that possess the further property of being smooth.
Further generalization of the framework from the tangent bundle of a manifold to a general smooth vector bundle,
a field on a smooth vector bundle would naturally be regarded as a smooth map from the base space to the total space,
assigning to each point a vector from the corresponding fiber of that point, which are called the sections of the smooth vector bundle.\\
Given a finite-dimensional vector space, there are a variety of ways to construct new vector spaces from it.
One of the most important constructions from a given vector space is the space of all tensors or multilinear maps
on it. The space of all tensors on a vector space is a very large space encompassing the spaces of covariant, contravariant, and
mixed tensors as its sub-spaces. This construction in the category of finite-dimensional vector spaces can be carried
to the category of smooth vector bundles in a canonical way. Unifying the spaces of tensors on all fibers of a smooth vector bundle
into a single space, and endowing it with a differentiable structure canonically constructed based on the
structure of the initial smooth vector bundle, a new smooth vector bundle arises which is called the
$\quotl$tensor bundle of the initial smooth vector bundle$\quotr$. Any section of the tensor bundle of a smooth vector bundle
is called a tensor field on that smooth vector bundle. So a tensor field on a smooth vector bundle is a smooth map assigning to
each point of the base space a covariant, contravariant, or mixed multilinear map on the fiber corresponded to that point.
The consideration of tensor fields is mostly limited to the special case of the tangent bundle of a smooth manifold in the literature,
but this text takes into account the most general case of a tensor field on a smooth vector bundle.\\
The are other constructions from initially given finite-dimensional vector spaces that can be simply extended to
the category of smooth vector bundles. The dual of a vector space, the space of all linear maps between a pair of
vector spaces, and the direct product of a pair of vector spaces are such instances, the counterpart of which are called
the dual vector bundle, the Hom bundle, and the Whitney sum of vector bundles, respectively, in the category of
smooth vector bundles.\\
Finally, it is important to pay attention to a delicate technical detail. We consider two types of smooth vector bundles
in this text with regard to the underlying field of the fibers, namely the real and complex vector bundles.
In both cases, the total space wich is formed by combining together all of the fibers, is actually a real manifold,
that is a differentiable structure modeled on a real Euclidean space. In addition, there are some occurrences of
model spaces when constructing the total spaces of a new vector bundle, such as the tensor bundle of a complex vector bundle,
that are expressed as a combination of real and complex vector spaces. In such cases, the realification of the complex part
of the model space is taken into account and hence is regarded as a real vector space.
%%%%%%%%%%%%%
\vskip0.5\baselineskip
\noindent
{\bf{A brief explanation about chapters:}}
The first chapter is devoted to the review of basic differential geometry necessary to the
development of the theory of smooth vector bundles.
The second chapter is an introduction to the theory of multilinear algebra, starting
with the most general case of tensors on modules and then continuing within the category of finite dimensional vector spaces.
The third chapter briefly introduces the basic concepts of the theory of smooth fiber bundles.
The fourth chapter is the core part of this text that technically deals with the problems mentioned above.
\vskip0.5\baselineskip
\noindent
{\bf{Prerequisites:}}
A considerable knowledge of linear algebra, differential calculus, topology, and differential geometry,
and a basic knowledge of the theory of modules is presumed.
\hfill
{\textsf{Farzad Shahi}}
%%%%%%%%%%%%%%%%%%%%%%%%%%%%
%%%%%%%%%%%%%%%%%%%%%%%%%%%%
%%%%%%%%%%%%%%%%%%%%%%%%%%%%%%%%%%%%%%%%%%
%\tableofcontents
\newpage
\thispagestyle{empty}
%\Thecontents{}
\section*{\fontsize{21}{21}\selectfont\bf{Contents}}
\addtocontents{toc}{\protect\setcounter{tocdepth}{-1}}
\tableofcontents
\addtocontents{toc}{\protect\setcounter{tocdepth}{3}}
%\tableofcontents
\newpage
%%%%%%%%%%%%%%%%%%%%%%%%%%%%%%%%%%%%%%%%%%%%%%%%%%%%%%%%%%
%%%%%%%%%%%%%%%%%%%%%%%%%%%%%%%%%%%%%%%%%%%%%%%%%%%%%%%%%%
\newpage
\chapteR{
Mathematical Symbols
}
\thispagestyle{fancy}
\section*{
Set-theory
}
%%%%%%%%%%%%%%%%%%%%%%%%%%%%%%%%%%%%%%%%%%%%%%%%%%%%%%%%%%%%%%%%%%%%%%%%%%%%%%%%%%%%%%%%
$\newsymb{\empty}$
empty-set
\varfill
$\empty$\\
%%%%%%%%%%%%%%%%%%%%%%%%%%%%%%%%%%%%%%%%%%%%%%%%%%%%%%%%%%%%%%%%%%%%%%%%%%%%%%%%%%%%%%%%
$\newsymp{\SET{1}=\SET{2}}$
$\SET{1}$
equals
$\SET{2}$.
\varfill
$=$\\
%%%%%%%%%%%%%%%%%%%%%%%%%%%%%%%%%%%%%%%%%%%%%%%%%%%%%%%%%%%%%%%%%%%%%%%%%%%%%%%%%%%%%%%%
%%%%%%%%%%%%%%%%%%%%%%%%%%%%%%%%%%%%%%%%%%%%%%%%%%%%%%%%%%%%%%%%%%%%%%%%%%%%%%%%%%%%%%%%
$\newsymp{\SET{1}\in\SET{2}}$
$\SET{1}$
is an element of
$\SET{2}$.
\varfill
$\in$\\
%%%%%%%%%%%%%%%%%%%%%%%%%%%%%%%%%%%%%%%%%%%%%%%%%%%%%%%%%%%%%%%%%%%%%%%%%%%%%%%%%%%%%%%%
$\newsymp{\SET{1}\ni\SET{2}}$
$\SET{1}$
contains
$\SET{2}$.
\varfill
$\ni$\\
%%%%%%%%%%%%%%%%%%%%%%%%%%%%%%%%%%%%%%%%%%%%%%%%%%%%%%%%%%%%%%%%%%%%%%%%%%%%%%%%%%%%%%%%
$\newsymb{\seta{\binary{\SET{1}}{\SET{2}}}}$
the set composed of
$\SET{1}$
and
$\SET{2}$
\varfill
$\seta{\binary{\dummy}{\dummy}}$\\
%%%%%%%%%%%%%%%%%%%%%%%%%%%%%%%%%%%%%%%%%%%%%%%%%%%%%%%%%%%%%%%%%%%%%%%%%%%%%%%%%%%%%%%%
$\newsymb{\defset{\SET{1}}{\SET{2}}{\prop{\SET{1}}}}$
all elements of
$\SET{2}$
having the property
$\Prop$
\varfill
$\defset{\dummy}{\dummy}{\dummy}$\\
%%%%%%%%%%%%%%%%%%%%%%%%%%%%%%%%%%%%%%%%%%%%%%%%%%%%%%%%%%%%%%%%%%%%%%%%%%%%%%%%%%%%%%%%
$\newsymb{\union{\SET{}}}$
union of all elements of
$\SET{}$
\varfill
$\bigcup$\\
%%%%%%%%%%%%%%%%%%%%%%%%%%%%%%%%%%%%%%%%%%%%%%%%%%%%%%%%%%%%%%%%%%%%%%%%%%%%%%%%%%%%%%%%
$\newsymb{\intersection{\SET{}}}$
intersection of all elements of
$\SET{}$
\varfill
$\bigcap$\\
%%%%%%%%%%%%%%%%%%%%%%%%%%%%%%%%%%%%%%%%%%%%%%%%%%%%%%%%%%%%%%%%%%%%%%%%%%%%%%%%%%%%%%%%
$\newsymb{\CSs{\SET{}}}$
power-set of
$\SET{}$
\varfill
$\CSs{\dummy}$\\
%%%%%%%%%%%%%%%%%%%%%%%%%%%%%%%%%%%%%%%%%%%%%%%%%%%%%%%%%%%%%%%%%%%%%%%%%%%%%%%%%%%%%%%%
$\newsymb{\Dproduct{\alpha}{\index}{\SET{\alpha}}}$
Cartesian-product of the collection of indexed sets
${\seta{\SET{\alpha}}}_{\alpha\in\index}$
\varfill
$\prod$\\
%%%%%%%%%%%%%%%%%%%%%%%%%%%%%%%%%%%%%%%%%%%%%%%%%%%%%%%%%%%%%%%%%%%%%%%%%%%%%%%%%%%%%%%%
$\newsymp{\SET{1}\subseteq\SET{2}}$
$\SET{1}$
is a subset of
$\SET{2}$.
\varfill
$\subseteq$\\
%%%%%%%%%%%%%%%%%%%%%%%%%%%%%%%%%%%%%%%%%%%%%%%%%%%%%%%%%%%%%%%%%%%%%%%%%%%%%%%%%%%%%%%%
$\newsymp{\SET{1}\supseteq\SET{2}}$
$\SET{1}$
includes
$\SET{2}$.
\varfill
$\supseteq$\\
%%%%%%%%%%%%%%%%%%%%%%%%%%%%%%%%%%%%%%%%%%%%%%%%%%%%%%%%%%%%%%%%%%%%%%%%%%%%%%%%%%%%%%%%
$\newsymp{\SET{1}\subset\SET{2}}$
$\SET{1}$
is a proper subset of
$\SET{2}$.
\varfill
$\subset$\\
%%%%%%%%%%%%%%%%%%%%%%%%%%%%%%%%%%%%%%%%%%%%%%%%%%%%%%%%%%%%%%%%%%%%%%%%%%%%%%%%%%%%%%%%
$\newsymp{\SET{1}\supset\SET{2}}$
$\SET{1}$
properly includes
$\SET{2}$.
\varfill
$\supset$\\
%%%%%%%%%%%%%%%%%%%%%%%%%%%%%%%%%%%%%%%%%%%%%%%%%%%%%%%%%%%%%%%%%%%%%%%%%%%%%%%%%%%%%%%%
$\newsymb{\SET{1}\cup\SET{2}}$
union of
$\SET{1}$
and
$\SET{2}$
\varfill
$\cup$\\
%%%%%%%%%%%%%%%%%%%%%%%%%%%%%%%%%%%%%%%%%%%%%%%%%%%%%%%%%%%%%%%%%%%%%%%%%%%%%%%%%%%%%%%%
$\newsymb{\SET{1}\cap\SET{2}}$
intersection of
$\SET{1}$
and
$\SET{2}$
\varfill
$\cap$\\
%%%%%%%%%%%%%%%%%%%%%%%%%%%%%%%%%%%%%%%%%%%%%%%%%%%%%%%%%%%%%%%%%%%%%%%%%%%%%%%%%%%%%%%%
$\newsymb{\SET{1}\times\SET{2}}$
Cartesian-product of
$\SET{1}$
and
$\SET{2}$
\varfill
$\times$\\
%%%%%%%%%%%%%%%%%%%%%%%%%%%%%%%%%%%%%%%%%%%%%%%%%%%%%%%%%%%%%%%%%%%%%%%%%%%%%%%%%%%%%%%%
$\newsymb{\compl{\SET{1}}{\SET{2}}}$
the relative complement of
$\SET{2}$
with respect to
$\SET{1}$
\varfill
$\setminus$\\
%%%%%%%%%%%%%%%%%%%%%%%%%%%%%%%%%%%%%%%%%%%%%%%%%%%%%%%%%%%%%%%%%%%%%%%%%%%%%%%%%%%%%%%%
$\newsymb{\func{\FUNCTION{}}{\SET{}}}$
value of the function
$\FUNCTION{}$
at
$\SET{}$
\varfill
$\func{\dummy}{\dummy}$\\
%%%%%%%%%%%%%%%%%%%%%%%%%%%%%%%%%%%%%%%%%%%%%%%%%%%%%%%%%%%%%%%%%%%%%%%%%%%%%%%%%%%%%%%%
$\newsymb{\domain{\FUNCTION{}}}$
domain of the function
$\FUNCTION{}$
\varfill
$\domain{\dummy}$\\
%%%%%%%%%%%%%%%%%%%%%%%%%%%%%%%%%%%%%%%%%%%%%%%%%%%%%%%%%%%%%%%%%%%%%%%%%%%%%%%%%%%%%%%%
$\newsymb{\codomain{\FUNCTION{}}}$
codomain of the function
$\FUNCTION{}$
\varfill
$\codomain{\dummy}$\\
%%%%%%%%%%%%%%%%%%%%%%%%%%%%%%%%%%%%%%%%%%%%%%%%%%%%%%%%%%%%%%%%%%%%%%%%%%%%%%%%%%%%%%%%
$\newsymb{\funcimage{\FUNCTION{}}}$
image of the function
$\FUNCTION{}$
\varfill
$\funcimage{\dummy}$\\
%%%%%%%%%%%%%%%%%%%%%%%%%%%%%%%%%%%%%%%%%%%%%%%%%%%%%%%%%%%%%%%%%%%%%%%%%%%%%%%%%%%%%%%%
$\newsymb{\image{\FUNCTION{}}}$
image-map of the function
$\FUNCTION{}$
\varfill
$\image{\dummy}$\\
%%%%%%%%%%%%%%%%%%%%%%%%%%%%%%%%%%%%%%%%%%%%%%%%%%%%%%%%%%%%%%%%%%%%%%%%%%%%%%%%%%%%%%%%
$\newsymb{\pimage{\FUNCTION{}}}$
inverse-image-map of the function
$\FUNCTION{}$
\varfill
$\pimage{\dummy}$\\
%%%%%%%%%%%%%%%%%%%%%%%%%%%%%%%%%%%%%%%%%%%%%%%%%%%%%%%%%%%%%%%%%%%%%%%%%%%%%%%%%%%%%%%%
$\newsymb{\resd{\FUNCTION{}}}$
domain-restriction-map of the function
$\FUNCTION{}$
\varfill
$\resd{\dummy}$\\
%%%%%%%%%%%%%%%%%%%%%%%%%%%%%%%%%%%%%%%%%%%%%%%%%%%%%%%%%%%%%%%%%%%%%%%%%%%%%%%%%%%%%%%%
$\newsymb{\rescd{\FUNCTION{}}}$
codomain-restriction-map of the function
$\FUNCTION{}$
\varfill
$\rescd{\dummy}$\\
%%%%%%%%%%%%%%%%%%%%%%%%%%%%%%%%%%%%%%%%%%%%%%%%%%%%%%%%%%%%%%%%%%%%%%%%%%%%%%%%%%%%%%%%
$\newsymb{\func{\res{\FUNCTION{}}}{\SET{}}}$
domain-restriction and codomain-restriction of\\ the function
$\FUNCTION{}$ to $\SET{}$ and $\func{\image{\FUNCTION{}}}{\SET{}}$, respectively
\varfill
$\res{\dummy}$\\
%%%%%%%%%%%%%%%%%%%%%%%%%%%%%%%%%%%%%%%%%%%%%%%%%%%%%%%%%%%%%%%%%%%%%%%%%%%%%%%%%%%%%%%%
$\newsymb{\Func{\SET{1}}{\SET{2}}}$
the set of all maps from
$\SET{1}$
to
$\SET{2}$
\varfill
$\Func{\dummy}{\dummy}$\\
%%%%%%%%%%%%%%%%%%%%%%%%%%%%%%%%%%%%%%%%%%%%%%%%%%%%%%%%%%%%%%%%%%%%%%%%%%%%%%%%%%%%%%%%
$\newsymb{\IF{\SET{1}}{\SET{2}}}$
the set of all bijective functions from
$\SET{1}$
to
$\SET{2}$
\varfill
$\IF{\dummy}{\dummy}$\\
%%%%%%%%%%%%%%%%%%%%%%%%%%%%%%%%%%%%%%%%%%%%%%%%%%%%%%%%%%%%%%%%%%%%%%%%%%%%%%%%%%%%%%%%
$\newsymb{\finv{\FUNCTION{}}}$
the inverse mapping of the bijective function $\FUNCTION{}$
\varfill
$\finv{\dummy}$\\
%%%%%%%%%%%%%%%%%%%%%%%%%%%%%%%%%%%%%%%%%%%%%%%%%%%%%%%%%%%%%%%%%%%%%%%%%%%%%%%%%%%%%%%%
$\newsymb{\surFunc{\SET{1}}{\SET{2}}}$
the set of all surjective functions from
$\SET{1}$
to
$\SET{2}$
\varfill
$\surFunc{\dummy}{\dummy}$\\
%%%%%%%%%%%%%%%%%%%%%%%%%%%%%%%%%%%%%%%%%%%%%%%%%%%%%%%%%%%%%%%%%%%%%%%%%%%%%%%%%%%%%%%%
$\newsymb{\cmp{\FUNCTION{1}}{\FUNCTION{2}}}$
composition of the function
$\FUNCTION{1}$
with the function
$\FUNCTION{2}$
\varfill
$\cmp{}{}$\\
%%%%%%%%%%%%%%%%%%%%%%%%%%%%%%%%%%%%%%%%%%%%%%%%%%%%%%%%%%%%%%%%%%%%%%%%%%%%%%%%%%%%%%%%
$\newsymb{\Injection{\SET{1}}{\SET{2}}}$
the injection-mapping of the set $\SET{1}$ into the set $\SET{2}$
\varfill
$\Injection{\dummy}{\dummy}$\\
%%%%%%%%%%%%%%%%%%%%%%%%%%%%%%%%%%%%%%%%%%%%%%%%%%%%%%%%%%%%%%%%%%%%%%%%%%%%%%%%%%%%%%%%
$\newsymb{\funcprod{\cf_1}{\cf_2}}$
the function-product of the function $\cf_1$ and $\cf_2$
\varfill
$\Cprod{\dummy}{\dummy}$\\
%%%%%%%%%%%%%%%%%%%%%%%%%%%%%%%%%%%%%%%%%%%%%%%%%%%%%%%%%%%%%%%%%%%%%%%%%%%%%%%%%%%%%%%%
$\newsymb{\EqR{\SET{}}}$
the set of all equivalence relations on the set
$\SET{}$
\varfill
$\EqR{\dummy}$\\
%%%%%%%%%%%%%%%%%%%%%%%%%%%%%%%%%%%%%%%%%%%%%%%%%%%%%%%%%%%%%%%%%%%%%%%%%%%%%%%%%%%%%%%%
$\newsymb{\EqClass{\SET{1}}{\SET{2}}}$
quotient-set of
$\SET{1}$
by the equivalence-relation
$\SET{1}$
\varfill
$\EqClass{\dummy}{\dummy}$\\
%%%%%%%%%%%%%%%%%%%%%%%%%%%%%%%%%%%%%%%%%%%%%%%%%%%%%%%%%%%%%%%%%%%%%%%%%%%%%%%%%%%%%%%%
$\newsymb{\pEqclass{\SET{1}}{\SET{2}}}$
equivalence-class of
$\SET{1}$
by the equivalence-relation
$\SET{2}$
\varfill
$\pEqclass{\dummy}{\dummy}$\\
%%%%%%%%%%%%%%%%%%%%%%%%%%%%%%%%%%%%%%%%%%%%%%%%%%%%%%%%%%%%%%%%%%%%%%%%%%%%%%%%%%%%%%%%
$\newsymb{\PEqclass{\SET{1}}{\SET{2}}}$
equivalence-class of
$\SET{2}$
by the equivalence-relation
$\SET{1}$
\varfill
$\PEqclass{\dummy}{\dummy}$\\
%%%%%%%%%%%%%%%%%%%%%%%%%%%%%%%%%%%%%%%%%%%%%%%%%%%%%%%%%%%%%%%%%%%%%%%%%%%%%%%%%%%%%%%%
$\newsymp{\Card{\SET{1}}\cardeq\Card{\SET{2}}}$
$\IF{\SET{1}}{\SET{2}}$
is non-empty.
\varfill
$\Card{\dummy}\cardeq\Card{\dummy}$\\
%%%%%%%%%%%%%%%%%%%%%%%%%%%%%%%%%%%%%%%%%%%%%%%%%%%%%%%%%%%%%%%%%%%%%%%%%%%%%%%%%%%%%%%%
$\newsymb{\CarD{\SET{}}}$
cardinality of
$\SET{}$
\varfill
$\CarD{\dummy}$
%%%%%%%%%%%%%%%%%%%%%%%%%%%%%%%%%%%%%%%%%%%%%%%%%%%%%%%%%%%%%%%%%%%%%%%%%%%%%%%%%%%%%%%%
%%%%%%%%%%%%%%%%%%%%%%%%%%%%%%%%%%%%%%%%%%%%%%%%%%%%%%%%%%%%%%%%%%%%%%%%%%%%%%%%%%%%%%%%
%%%%%%%%%%%%%%%%%%%%%%%%%%%%%%%%%%%%%%%%%%%%%%%%%%%%%%%%%%%%%%%%%%%%%%%%%%%%%%%%%%%%%%%%
%%%%%%%%%%%%%%%%%%%%%%%%%%%%%%%%%%%%%%%%%%%%%%%%%%%%%%%%%%%%%%%%%%%%%%%%%%%%%%%%%%%%%%%%
%%%%%%%%%%%%%%%%%%%%%%%%%%%%%%%%%%%%%%%%%%%%%%%%%%%%%%%%%%%%%%%%%%%%%%%%%%%%%%%%%%%%%%%%
\section*{
Logic
}
%%%%%%%%%%%%%%%%%%%%%%%%%%%%%%%%%%%%%%%%%%%%%%%%%%%%%%%%%%%%%%%%%%%%%%%%%%%%%%%%%%%%%%%%
$\newsymp{\AND{\propos{1}}{\propos{2}}}$
$\propos{1}$
and
$\propos{2}$.
\varfill
$\AND{}{}$\\
%%%%%%%%%%%%%%%%%%%%%%%%%%%%%%%%%%%%%%%%%%%%%%%%%%%%%%%%%%%%%%%%%%%%%%%%%%%%%%%%%%%%%%%%
$\newsymp{\OR{\propos{1}}{\propos{2}}}$
$\propos{1}$
or
$\propos{2}$.
\varfill
$\OR{}{}$\\
%%%%%%%%%%%%%%%%%%%%%%%%%%%%%%%%%%%%%%%%%%%%%%%%%%%%%%%%%%%%%%%%%%%%%%%%%%%%%%%%%%%%%%%%
$\newsymp{{\propos{1}}\then{\propos{2}}}$
if
$\propos{1}$,
then
$\propos{2}$.
\varfill
$\then$\\
%%%%%%%%%%%%%%%%%%%%%%%%%%%%%%%%%%%%%%%%%%%%%%%%%%%%%%%%%%%%%%%%%%%%%%%%%%%%%%%%%%%%%%%%
$\newsymp{{\propos{1}}\thenn{\propos{2}}}$
$\propos{1}$,
if-and-only-if
$\propos{2}$.
\varfill
$\thenn$\\
%%%%%%%%%%%%%%%%%%%%%%%%%%%%%%%%%%%%%%%%%%%%%%%%%%%%%%%%%%%%%%%%%%%%%%%%%%%%%%%%%%%%%%%%
$\newsymp{\negation{\propos{}}}$
$\propos{1}$,
negation of
$\propos{}$.
\varfill
$\negation{}$\\
%%%%%%%%%%%%%%%%%%%%%%%%%%%%%%%%%%%%%%%%%%%%%%%%%%%%%%%%%%%%%%%%%%%%%%%%%%%%%%%%%%%%%%%%
$\newsymp{\Foreach{\SET{1}}{\SET{2}}{\prop{\SET{1}}}}$
for every
$\SET{2}$
in
$\SET{1}$,
$\prop{\SET{1}}$.
\varfill
$\Foreach{\dummy}{\dummy}\dummy$\\
%%%%%%%%%%%%%%%%%%%%%%%%%%%%%%%%%%%%%%%%%%%%%%%%%%%%%%%%%%%%%%%%%%%%%%%%%%%%%%%%%%%%%%%%
$\newsymp{\Exists{\SET{1}}{\SET{2}}{\prop{\SET{1}}}}$
exists
$\SET{2}$
in
$\SET{1}$ such that
$\prop{\SET{1}}$.
\varfill
$\Foreach{\dummy}{\dummy}\dummy$\\
%%%%%%%%%%%%%%%%%%%%%%%%%%%%%%%%%%%%%%%%%%%%%%%%%%%%%%%%%%%%%%%%%%%%%%%%%%%%%%%%%%%%%%%%
%%%%%%%%%%%%%%%%%%%%%%%%%%%%%%%%%%%%%%%%%%%%%%%%%%%%%%%%%%%%%%%%%%%%%%%%%%%%%%%%%%%%%%%%
%%%%%%%%%%%%%%%%%%%%%%%%%%%%%%%%%%%%%%%%%%%%%%%%%%%%%%%%%%%%%%%%%%%%%%%%%%%%%%%%%%%%%%%%
%%%%%%%%%%%%%%%%%%%%%%%%%%%%%%%%%%%%%%%%%%%%%%%%%%%%%%%%%%%%%%%%%%%%%%%%%%%%%%%%%%%%%%%%
%%%%%%%%%%%%%%%%%%%%%%%%%%%%%%%%%%%%%%%%%%%%%%%%%%%%%%%%%%%%%%%%%%%%%%%%%%%%%%%%%%%%%%%%
%%%%%%%%%%%%%%%%%%%%%%%%%%%%%%%%%%%%%%%%%%%%%%%%%%%%%%%%%%%%%%%%%%%%%%%%%%%%%%%%%%%%%%%%
%%%%%%%%%%%%%%%%%%%%%%%%%%%%%%%%%%%%%%%%%%%%%%%%%%%%%%%%%%%%%%%%%%%%%%%%%%%%%%%%%%%%%%%%
\section*{
Group Theory
}
%%%%%%%%%%%%%%%%%%%%%%%%%%%%%%%%%%%%%%%%%%%%%%%%%%%%%%%%%%%%%%%%%%%%%%%%%%%%%%%%%%%%%%%%
$\newsymp{\GHom{\Group{1}}{\Group{2}}}$
the set of all group-homomorphisms from the group $\Group{1}$ targeted to the
group $\Group{2}$
\varfill
$\GHom{\dummy}{\dummy}$\\
%%%%%%%%%%%%%%%%%%%%%%%%%%%%%%%%%%%%%%%%%%%%%%%%%%%%%%%%%%%%%%%%%%%%%%%%%%%%%%%%%%%%%%%%
$\newsymp{\GIsom{\Group{1}}{\Group{2}}}$
the set of all group-isomorphisms from the group $\Group{1}$ targeted to the
group $\Group{2}$
\varfill
$\GIsom{\dummy}{\dummy}$\\
%%%%%%%%%%%%%%%%%%%%%%%%%%%%%%%%%%%%%%%%%%%%%%%%%%%%%%%%%%%%%%%%%%%%%%%%%%%%%%%%%%%%%%%%
$\newsymp{\Subgroups{\Group{}}}$
the set of all subgroups of the group $\Group{}$
\varfill
$\Subgroups{\dummy}$\\
%%%%%%%%%%%%%%%%%%%%%%%%%%%%%%%%%%%%%%%%%%%%%%%%%%%%%%%%%%%%%%%%%%%%%%%%%%%%%%%%%%%%%%%%
$\newsymp{\func{\LCoset{\Group{}}}{\asubgroup{}}}$
the set of all left-cosets of the subgroup $\asubgroup{}$ of the group $\Group{}$
\varfill
$\func{\LCoset{\dummy}}{\dummy}$
%%%%%%%%%%%%%%%%%%%%%%%%%%%%%%%%%%%%%%%%%%%%%%%%%%%%%%%%%%%%%%%%%%%%%%%%%%%%%%%%%%%%%%%%
%%%%%%%%%%%%%%%%%%%%%%%%%%%%%%%%%%%%%%%%%%%%%%%%%%%%%%%%%%%%%%%%%%%%%%%%%%%%%%%%%%%%%%%%
%%%%%%%%%%%%%%%%%%%%%%%%%%%%%%%%%%%%%%%%%%%%%%%%%%%%%%%%%%%%%%%%%%%%%%%%%%%%%%%%%%%%%%%%
%%%%%%%%%%%%%%%%%%%%%%%%%%%%%%%%%%%%%%%%%%%%%%%%%%%%%%%%%%%%%%%%%%%%%%%%%%%%%%%%%%%%%%%%
%%%%%%%%%%%%%%%%%%%%%%%%%%%%%%%%%%%%%%%%%%%%%%%%%%%%%%%%%%%%%%%%%%%%%%%%%%%%%%%%%%%%%%%%
%%%%%%%%%%%%%%%%%%%%%%%%%%%%%%%%%%%%%%%%%%%%%%%%%%%%%%%%%%%%%%%%%%%%%%%%%%%%%%%%%%%%%%%%
%%%%%%%%%%%%%%%%%%%%%%%%%%%%%%%%%%%%%%%%%%%%%%%%%%%%%%%%%%%%%%%%%%%%%%%%%%%%%%%%%%%%%%%%
%%%%%%%%%%%%%%%%%%%%%%%%%%%%%%%%%%%%%%%%%%%%%%%%%%%%%%%%%%%%%%%%%%%%%%%%%%%%%%%%%%%%%%%%
\section*{
Differential Calculus
}
%%%%%%%%%%%%%%%%%%%%%%%%%%%%%%%%%%%%%%%%%%%%%%%%%%%%%%%%%%%%%%%%%%%%%%%%%%%%%%%%%%%%%%%%
$\newsymp{\banachmapdifclass{r}{\NVS{1}}{\NVS{2}}{\U_1}{\U_2}}$
the set of all $r$-times differentiable maps from the Banach-space $\NVS{1}$
to the Banach-space $\NVS{2}$ with domain $\U_1$ and codomain $\U_2$
\varfill
$\banachmapdifclass{\dummy}{\dummy}{\dummy}{\dummy}{\dummy}$\\
%%%%%%%%%%%%%%%%%%%%%%%%%%%%%%%%%%%%%%%%%%%%%%%%%%%%%%%%%%%%%%%%%%%%%%%%%%%%%%%%%%%%%%%%
$\newsymp{\BanachDiff{r}{\NVS{1}}{\NVS{2}}{\U_1}{\U_2}}$
the set of all $r$-times differentiable diffeomorphisms from the Banach-space $\NVS{1}$
to the Banach-space $\NVS{2}$ with domain $\U_1$ and codomain $\U_2$
\varfill
$\BanachDiff{\dummy}{\dummy}{\dummy}{\dummy}{\dummy}$\\
%%%%%%%%%%%%%%%%%%%%%%%%%%%%%%%%%%%%%%%%%%%%%%%%%%%%%%%%%%%%%%%%%%%%%%%%%%%%%%%%%%%%%%%%
$\newsymp{\banachder{\FUNCTION{}}{\NVS{1}}{\NVS{2}}}$
derived map of $\FUNCTION{}$, $\FUNCTION{}$ being a $\difclass{r}$ map
from the Banach-space $\NVS{1}$ to the Banach-space $\NVS{2}$
(with open domain and codomain
in $\NVS{1}$ and $\NVS{2}$, respectively)
\varfill
$\banachder{\dummy}{\dummy}{\dummy}$
%%%%%%%%%%%%%%%%%%%%%%%%%%%%%%%%%%%%%%%%%%%%%%%%%%%%%%%%%%%%%%%%%%%%%%%%%%%%%%%%%%%%%%%%
%%%%%%%%%%%%%%%%%%%%%%%%%%%%%%%%%%%%%%%%%%%%%%%%%%%%%%%%%%%%%%%%%%%%%%%%%%%%%%%%%%%%%%%%
%%%%%%%%%%%%%%%%%%%%%%%%%%%%%%%%%%%%%%%%%%%%%%%%%%%%%%%%%%%%%%%%%%%%%%%%%%%%%%%%%%%%%%%%
%%%%%%%%%%%%%%%%%%%%%%%%%%%%%%%%%%%%%%%%%%%%%%%%%%%%%%%%%%%%%%%%%%%%%%%%%%%%%%%%%%%%%%%%
%%%%%%%%%%%%%%%%%%%%%%%%%%%%%%%%%%%%%%%%%%%%%%%%%%%%%%%%%%%%%%%%%%%%%%%%%%%%%%%%%%%%%%%%
%%%%%%%%%%%%%%%%%%%%%%%%%%%%%%%%%%%%%%%%%%%%%%%%%%%%%%%%%%%%%%%%%%%%%%%%%%%%%%%%%%%%%%%%
%%%%%%%%%%%%%%%%%%%%%%%%%%%%%%%%%%%%%%%%%%%%%%%%%%%%%%%%%%%%%%%%%%%%%%%%%%%%%%%%%%%%%%%%
%%%%%%%%%%%%%%%%%%%%%%%%%%%%%%%%%%%%%%%%%%%%%%%%%%%%%%%%%%%%%%%%%%%%%%%%%%%%%%%%%%%%%%%%
\section*{
Linear Algebra
}
%%%%%%%%%%%%%%%%%%%%%%%%%%%%%%%%%%%%%%%%%%%%%%%%%%%%%%%%%%%%%%%%%%%%%%%%%%%%%%%%%%%%%%%%
$\newsymp{\subvec{\VVS{}}{m}}$
the set of all sets of vectors of all $m$-dimensional vector-subspaces of
the vector-space $\VVS{}$
\varfill
$\subvec{\dummy}{\dummy}$\\
%%%%%%%%%%%%%%%%%%%%%%%%%%%%%%%%%%%%%%%%%%%%%%%%%%%%%%%%%%%%%%%%%%%%%%%%%%%%%%%%%%%%%%%%
$\newsymp{\func{\Vspan{\VVS{}}}{\asubset}}$
the vector-subspace of the vector-space $\VVS{}$ spanned by the subset $\asubset$
of the set of all vectors of $\VVS{}$
\varfill
$\func{\Vspan{\dummy}}{\dummy}$\\
%%%%%%%%%%%%%%%%%%%%%%%%%%%%%%%%%%%%%%%%%%%%%%%%%%%%%%%%%%%%%%%%%%%%%%%%%%%%%%%%%%%%%%%%
$\newsymp{\ovecbasis{\VVS{}}}$
the set of all ordered-bases of the vector-space $\VVS{}$
\varfill
$\ovecbasis{\dummy}$\\
%%%%%%%%%%%%%%%%%%%%%%%%%%%%%%%%%%%%%%%%%%%%%%%%%%%%%%%%%%%%%%%%%%%%%%%%%%%%%%%%%%%%%%%%
$\newsymp{\Lin{\VVS{1}}{\VVS{2}}}$
the set of all linear maps from the vector-space $\VVS{1}$ to the vector-space $\VVS{2}$
\varfill
$\Lin{\dummy}{\dummy}$\\
%%%%%%%%%%%%%%%%%%%%%%%%%%%%%%%%%%%%%%%%%%%%%%%%%%%%%%%%%%%%%%%%%%%%%%%%%%%%%%%%%%%%%%%%
$\newsymp{\Linisom{\VVS{1}}{\VVS{2}}}$
the set of all linear isomorphisms from the vector-space $\VVS{1}$ to the vector-space
$\VVS{2}$
\varfill
$\Linisom{\dummy}{\dummy}$\\
%%%%%%%%%%%%%%%%%%%%%%%%%%%%%%%%%%%%%%%%%%%%%%%%%%%%%%%%%%%%%%%%%%%%%%%%%%%%%%%%%%%%%%%%
$\newsymp{\GL{\VVS{}}{}}$
the set of all linear isomorphisms from the vector-space $\VVS{}$ to itself
\varfill
$\GL{\dummy}{}$\\
%%%%%%%%%%%%%%%%%%%%%%%%%%%%%%%%%%%%%%%%%%%%%%%%%%%%%%%%%%%%%%%%%%%%%%%%%%%%%%%%%%%%%%%%
$\newsymp{\VLin{\VVS{1}}{\VVS{2}}}$
the canonical vector-space of all linear maps from the vector-space
$\VVS{1}$ to the vector-space $\VVS{2}$
\varfill
$\VLin{\dummy}{\dummy}$\\
%%%%%%%%%%%%%%%%%%%%%%%%%%%%%%%%%%%%%%%%%%%%%%%%%%%%%%%%%%%%%%%%%%%%%%%%%%%%%%%%%%%%%%%%
$\newsymp{\NVLin{\NVS{1}}{\NVS{2}}}$
the canonical Banach-space of all linear maps from the finite-dimensional Banach-space
$\NVS{1}$ to the finite-dimensional Banach-space $\VVS{2}$
\varfill
$\NVLin{\dummy}{\dummy}$\\
%%%%%%%%%%%%%%%%%%%%%%%%%%%%%%%%%%%%%%%%%%%%%%%%%%%%%%%%%%%%%%%%%%%%%%%%%%%%%%%%%%%%%%%%
$\newsymp{\Mat{\F}{m}{n}}$
the set of all $m\times n$ matrices over the field $\F$
\varfill
$\Mat{\dummy}{\dummy}{\dummy}$\\
%%%%%%%%%%%%%%%%%%%%%%%%%%%%%%%%%%%%%%%%%%%%%%%%%%%%%%%%%%%%%%%%%%%%%%%%%%%%%%%%%%%%%%%%
$\newsymp{\BMat{\F}{m}{n}}$
the canonical Banach-space of all $m\times n$ matrices over the field $\F$
(endowed with its natural norm)
\varfill
$\BMat{\dummy}{\dummy}{\dummy}$\\
%%%%%%%%%%%%%%%%%%%%%%%%%%%%%%%%%%%%%%%%%%%%%%%%%%%%%%%%%%%%%%%%%%%%%%%%%%%%%%%%%%%%%%%%
$\newsymp{\Det{n}}$
the determinant function on the set of all square $\R$-matrices of degree $n$
\varfill
$\Det{\dummy}$\\
%%%%%%%%%%%%%%%%%%%%%%%%%%%%%%%%%%%%%%%%%%%%%%%%%%%%%%%%%%%%%%%%%%%%%%%%%%%%%%%%%%%%%%%%
$\newsymp{\directsum{\VVS{1}}{\VVS{2}}}$
the direct sum of the vector spaces $\VVS{1}$ and $\VVS{2}$ over the same field
\varfill
$\directsum{\dummy}{\dummy}$
%%%%%%%%%%%%%%%%%%%%%%%%%%%%%%%%%%%%%%%%%%%%%%%%%%%%%%%%%%%%%%%%%%%%%%%%%%%%%%%%%%%%%%%%
%%%%%%%%%%%%%%%%%%%%%%%%%%%%%%%%%%%%%%%%%%%%%%%%%%%%%%%%%%%%%%%%%%%%%%%%%%%%%%%%%%%%%%%%
%%%%%%%%%%%%%%%%%%%%%%%%%%%%%%%%%%%%%%%%%%%%%%%%%%%%%%%%%%%%%%%%%%%%%%%%%%%%%%%%%%%%%%%%
%%%%%%%%%%%%%%%%%%%%%%%%%%%%%%%%%%%%%%%%%%%%%%%%%%%%%%%%%%%%%%%%%%%%%%%%%%%%%%%%%%%%%%%%
%%%%%%%%%%%%%%%%%%%%%%%%%%%%%%%%%%%%%%%%%%%%%%%%%%%%%%%%%%%%%%%%%%%%%%%%%%%%%%%%%%%%%%%%
%%%%%%%%%%%%%%%%%%%%%%%%%%%%%%%%%%%%%%%%%%%%%%%%%%%%%%%%%%%%%%%%%%%%%%%%%%%%%%%%%%%%%%%%
%%%%%%%%%%%%%%%%%%%%%%%%%%%%%%%%%%%%%%%%%%%%%%%%%%%%%%%%%%%%%%%%%%%%%%%%%%%%%%%%%%%%%%%%
%%%%%%%%%%%%%%%%%%%%%%%%%%%%%%%%%%%%%%%%%%%%%%%%%%%%%%%%%%%%%%%%%%%%%%%%%%%%%%%%%%%%%%%%
\section*{Topology}
$\newsymp{\alltopologies{\SET{}}}$
the set of all topologies on the set $\SET{}$
\dotfill
$\alltopologies{\dummy}$\\
%%%%%%%%%%%%%%%%%%%%%%%%%%%%%%%%%%%%%%%%%%%%%%%%%%%%%%%%%%%%%%%%%%%%%%%%%%%%%%%%%%%%%%%%
$\newsymp{\topologyofspace{\Xt}}$
topology of the topological-space $\Xt$
\dotfill
$\topologyofspace{\dummy}$\\
%%%%%%%%%%%%%%%%%%%%%%%%%%%%%%%%%%%%%%%%%%%%%%%%%%%%%%%%%%%%%%%%%%%%%%%%%%%%%%%%%%%%%%%%
$\newsymp{\closedsets{\Xt}}$
the set of all closed sets of the topological-space $\Xt$
\dotfill
$\closedsets{\dummy}$\\
%%%%%%%%%%%%%%%%%%%%%%%%%%%%%%%%%%%%%%%%%%%%%%%%%%%%%%%%%%%%%%%%%%%%%%%%%%%%%%%%%%%%%%%%
$\newsymp{\func{\Cl{\Xt}}{\SET{}}}$
the closure of $\SET{}$ in the topological-space $\Xt$
\dotfill
$\func{\Cl{\dummy}}{\dummy}$\\
%%%%%%%%%%%%%%%%%%%%%%%%%%%%%%%%%%%%%%%%%%%%%%%%%%%%%%%%%%%%%%%%%%%%%%%%%%%%%%%%%%%%%%%%
$\newsymp{\func{\nei{\Xt}}{\asubset}}$
the set of all open sets of the topological-space $\Xt$ including the subset $\asubset$
of $\Xt$ (the set of all neighborhoods of $\asubset$ in $\Xt$)
\dotfill
$\func{\nei{\dummy}}{\dummy}$\\
%%%%%%%%%%%%%%%%%%%%%%%%%%%%%%%%%%%%%%%%%%%%%%%%%%%%%%%%%%%%%%%%%%%%%%%%%%%%%%%%%%%%%%%%
$\newsymp{\func{\IndTop{\Xt}}{\SET{}}}$
the topology on the subset $\SET{}$ of the topological-space $\Xt$ induced (inherited)
from that of $\Xt$
\dotfill
$\func{\IndTop{\dummy}}{\dummy}$\\
%%%%%%%%%%%%%%%%%%%%%%%%%%%%%%%%%%%%%%%%%%%%%%%%%%%%%%%%%%%%%%%%%%%%%%%%%%%%%%%%%%%%%%%%
$\newsymp{\CF{\Xt}{\Yt}}$
the set of all continuous maps from the topological-space $\Xt$ to the topological-space
$\Yt$
\varfill
$\CF{\dummy}{\dummy}$\\
%%%%%%%%%%%%%%%%%%%%%%%%%%%%%%%%%%%%%%%%%%%%%%%%%%%%%%%%%%%%%%%%%%%%%%%%%%%%%%%%%%%%%%%%
$\newsymp{\HOF{\Xt}{\Yt}}$
the set of all homeomorphisms from the topological-space $\Xt$ to the topological-space
$\Yt$
\varfill
$\HOF{\dummy}{\dummy}$\\
%%%%%%%%%%%%%%%%%%%%%%%%%%%%%%%%%%%%%%%%%%%%%%%%%%%%%%%%%%%%%%%%%%%%%%%%%%%%%%%%%%%%%%%%
$\newsymp{\topprod{\Xt}{\Yt}}$
the topological-product of the topological-spaces $\Xt$ and $\Yt$
\varfill
$\topprod{\dummy}{\dummy}$\\
%%%%%%%%%%%%%%%%%%%%%%%%%%%%%%%%%%%%%%%%%%%%%%%%%%%%%%%%%%%%%%%%%%%%%%%%%%%%%%%%%%%%%%%%
$\newsymp{\connecteds{\Xt}}$
the set of all connected sets of the topological-space $\Xt$
\varfill
$\connecteds{\dummy}$\\
%%%%%%%%%%%%%%%%%%%%%%%%%%%%%%%%%%%%%%%%%%%%%%%%%%%%%%%%%%%%%%%%%%%%%%%%%%%%%%%%%%%%%%%%
$\newsymp{\maxcon{\Xt}}$
the set of all connected components of the topological-space $\Xt$
\varfill
$\maxcon{\dummy}$
%%%%%%%%%%%%%%%%%%%%%%%%%%%%%%%%%%%%%%%%%%%%%%%%%%%%%%%%%%%%%%%%%%%%%%%%%%%%%%%%%%%%%%%%
$\newsymp{\TMat{\F}{m}{n}}$
the canonical topological-space of all $m\times n$ matrices over the field $\F$
(with the topology induced by its natural norm)
\varfill
$\TMat{\dummy}{\dummy}{\dummy}$\\
%%%%%%%%%%%%%%%%%%%%%%%%%%%%%%%%%%%%%%%%%%%%%%%%%%%%%%%%%%%%%%%%%%%%%%%%%%%%%%%%%%%%%%%%
$\newsymp{\topR{\R^n}}$
the canonical topological-space of the Euclidean-space $\R^n$
(with the topology induced by its natural norm)
\varfill
$\topR{}$
%%%%%%%%%%%%%%%%%%%%%%%%%%%%%%%%%%%%%%%%%%%%%%%%%%%%%%%%%%%%%%%%%%%%%%%%%%%%%%%%%%%%%%%%
%%%%%%%%%%%%%%%%%%%%%%%%%%%%%%%%%%%%%%%%%%%%%%%%%%%%%%%%%%%%%%%%%%%%%%%%%%%%%%%%%%%%%%%%
%%%%%%%%%%%%%%%%%%%%%%%%%%%%%%%%%%%%%%%%%%%%%%%%%%%%%%%%%%%%%%%%%%%%%%%%%%%%%%%%%%%%%%%%
%%%%%%%%%%%%%%%%%%%%%%%%%%%%%%%%%%%%%%%%%%%%%%%%%%%%%%%%%%%%%%%%%%%%%%%%%%%%%%%%%%%%%%%%
%%%%%%%%%%%%%%%%%%%%%%%%%%%%%%%%%%%%%%%%%%%%%%%%%%%%%%%%%%%%%%%%%%%%%%%%%%%%%%%%%%%%%%%%
%%%%%%%%%%%%%%%%%%%%%%%%%%%%%%%%%%%%%%%%%%%%%%%%%%%%%%%%%%%%%%%%%%%%%%%%%%%%%%%%%%%%%%%%
%%%%%%%%%%%%%%%%%%%%%%%%%%%%%%%%%%%%%%%%%%%%%%%%%%%%%%%%%%%%%%%%%%%%%%%%%%%%%%%%%%%%%%%%
%%%%%%%%%%%%%%%%%%%%%%%%%%%%%%%%%%%%%%%%%%%%%%%%%%%%%%%%%%%%%%%%%%%%%%%%%%%%%%%%%%%%%%%%
\section*{Differential Geometry}
%%%%%%%%%%%%%%%%%%%%%%%%%%%%%%%%%%%%%%%%%%%%%%%%%%%%%%%%%%%%%%%%%%%%%%%%%%%%%%%%%%%%%%%%
$\newsymp{\atlases{r}{\M{}}{\NVS{}}}$
the set of all atlases of differentiablity class $\difclass{r}$
on the set $\M{}$
modeled on the Banach-space $\NVS{}$
\varfill
$\atlases{\dummy}{\dummy}{\dummy}$\\
%%%%%%%%%%%%%%%%%%%%%%%%%%%%%%%%%%%%%%%%%%%%%%%%%%%%%%%%%%%%%%%%%%%%%%%%%%%%%%%%%%%%%%%%
$\newsymp{\maxatlases{r}{\M{}}{\NVS{}}}$
the set of all maximal-atlases of differentiablity class $\difclass{r}$
on the set $\M{}$
modeled on the Banach-space $\NVS{}$
\varfill
$\maxatlases{\dummy}{\dummy}{\dummy}$\\
%%%%%%%%%%%%%%%%%%%%%%%%%%%%%%%%%%%%%%%%%%%%%%%%%%%%%%%%%%%%%%%%%%%%%%%%%%%%%%%%%%%%%%%%
$\newsymp{\func{\maxatlasgen{r}{\M{}}{\NVS{}}}{\atlas{}}}$
the maximal-atlas of differentiablity class $\difclass{r}$
on the set $\M{}$ modeled on
the Banach-space $\NVS{}$, generated by the atlas $\atlas{}$ in
$\atlases{r}{\M{}}{\NVS{}}$
\varfill
$\func{\maxatlasgen{\dummy}{\dummy}{\dummy}}{\dummy}$
%%%%%%%%%%%%%%%%%%%%%%%%%%%%%%%%%%%%%%%%%%%%%%%%%%%%%%%%%%%%%%%%%%%%%%%%%%%%%%%%%%%%%%%%
\section*{Mathematical Environments}
$\newsymp{\blacksquare}$
end of definition
\varfill
$\blacksquare$\\
%%%%%%%%%%%%%%%%%%%%%%%%%%%%%%%%%%%%%%%%%%%%%%%%%%%%%%%%%%%%%%%%%%%%%%%%%%%%%%%%%%%%%%%%
$\newsymp{\square}$
end of theorem, lemma, proposition, or corollary
\varfill
$\square$\\
%%%%%%%%%%%%%%%%%%%%%%%%%%%%%%%%%%%%%%%%%%%%%%%%%%%%%%%%%%%%%%%%%%%%%%%%%%%%%%%%%%%%%%%%
$\newsymp{\Diamond}$
end of the introduction of new fixed objects
\varfill
$\Diamond$
%%%%%%%%%%%%%%%%%%%%%%%%%%%%%%%%%%%%%%%%%%%%%%%%%%%%%%%%%%%%%%%%%%%%%%%%%%%%%%%%%%%%%%%%
%%%%%%%%%%%%%%%%%%%%%%%%%%%%%%%%%%%%%%%%%%%%%%%%%%%%%%%%%%%%%%%%%%%%%%%%%%%%%%%%%%%%%%%%%%%%%%%%%%%%%%%%%%%%%%%%%%%%%%%%%%%%%%%%%%%%%%%%%%%%%%%%%%%%%%%%%%%%%%%%%%%%%%%%%%%%%%%%
%%%%%%%%%%%%%%%%%%%%%%%%%%%%%%%%%%%%%%%%%%%%%%%%%%%%%%%%%%%%%%%%%%%%%%%%%%%%%%%%%%%%%%%%%%%%%%%%%%%%%%%%%%%%%%%%%%%%%%%%%%%%%%%%%%%%%%%%%%%%%%%%%%%%%%%%%%%%%%%%%%%%%%%%%%%%%%%%
%%%%%%%%%%%%%%%%%%%%%%%%%%%%%%%%%%%%%%%%%%%%%%%%%%%%%%%%%%%%%%%%%%%%%%%%%%%%%%%%%%%%%%%%%%%%%%%%%%%%%%%%%%%%%%%%%%%%%%%%%%%%%%%%%%%%%%%%%%%%%%%%%%%%%%%%%%%%%%%%%%%%%%%%%%%%%%%%
%%%%%%%%%%%%%%%%%%%%%%%%%%%%%%%%%%%%%%%%%%%%%%%%%%%%%%%%%%%%%%%%%%%%%%%%%%%%%%%%%%%%%%%%%%%%%%%%%%%%%%%%%%%%%%%%%%%%%%%%%%%%%%%%%%%%%%%%%%%%%%%%%%%%%%%%%%%%%%%%%%%%%%%%%%%%%%%%
\chapteR{Review of Differential Geometry}
\thispagestyle{fancy}
\section{Basic Concepts of smooth manifolds}
\textit{Here, the category of $\difclass{\infty}$ manifolds is decided to be the category that consists of
$\difclass{\infty}$ differentiable-structures without boundary modeled on a Banach-space of any finite non-zero dimension with the
Hausdorff and second-countable underlying topological-space as the objects, and the smooth
mappings between such $\difclass{\infty}$ differentiable-structures as the morphisms.}\\
%%%%%%%%%%%%%%%%%%%%%%%%%%%%%%%%%%%%%%%%
\fixed
$\Man{}=\opair{\M{}}{\maxatlas{}}$ is fixed as an $n$-dimensional and $\difclass{\infty}$ manifold
modeled the Banach-space $\R^n$. So, $\maxatlas{}$ is a $\difclass{\infty}$ maximal-atlas on
$\M{}$ modeled on the Banach-space $\R^{n}$, and the topology on $\M{}$ induced by this maximal-atlas
is Hausdorff and second-countable (i.e. possessing a countable base). The dimension of $\Man{}$
is denoted by $\mandim{\Man{}}$.\\
Also, $\defSet{\Man{i}=\opair{\M{i}}{\maxatlas{i}}}{i\in\Zp}$ is fixed as a collection of manifolds such that
for each positive integer $i$, $\Man{i}$ is an
$m_{i}$-dimensional and $\difclass{\infty}$ manifold modeled on $\R^{m_i}$, where each $m_i$
is a positive integer.
\endfixed
%%%%%%%%%%%%%%%%%%%%%%%%%%%%%%%%%%%%%%%%
\begin{itemize}
\item[\myitem{DG~1.}]
For every point $\point$ of $\Man{}$, any chart $\phi$ of $\Man{}$ whose domain contains
$\point$ is called a $\quotl$chart of $\Man{}$ around $\point$$\quotr$ or
a $\quotl$neighbourhood chart of $\point$ in the manifold $\Man{}$$\quotr$.
For every point $\point$ of $\Man{}$, any chart $\phi$
of $\Man{}$ around $\point$ such that $\func{\phi}{\point}=\zerovec{}$ is called a
$\quotl$chart of $\Man{}$ centered at $\point$$\quotr$, and such a chart exists at every point of $\Man{}$.
%%%%%%%%%%%%%%%%%%%%%%%%%
\item[\myitem{DG~2.}]
For every positive integer $n$,
$\RR^n$ denotes the $n$-dimensional Euclidean-space endowed with its canonical differentiable-structure
constructed upon the Banach-space $\R^n$ that arises
from the trivial atlas consisting merely of the identity map on $\R^{n}$.
So $\identity{\R^n}$ is a chart of the smooth manifold $\RR^n$.
$\RR^n$ is called the $\quotl$canonical differentiable structure of $\R^n$$\quotr$.
%%%%%%%%%%%%%%%%%%%%%%%%%
\item[\myitem{DG~3.}]
The topology on $\M{}$ induced by the maximal-atlas $\maxatlas{}$ is denoted by $\mantop{\Man{}}$,
and the corresponding topological-space is denoted by $\mantops{\Man{}}$, that is,
$\mantops{\Man{}}:=\opair{\M{}}{\mantop{\Man{}}}$. When $\mantop{\Man{}}$ coincides with
an initially given topology on $\M{}$, it is customary to say that $\quotl$the manifold $\Man{}$
is compatible with that topology$\quotr$.
%%%%%%%%%%%%%%%%%%%%%%%%%
\item[\myitem{DG~4.}]
For every positive integer $n$ and every integer $k$ in $\seta{\suc{1}{n}}$, $\projection{n}{k}$
denotes the projection mapping of $\R^{n}$ onto its $k$-th factor. That is, $\projection{n}{k}$
is the element of $\Func{\R^n}{\R}$ such that for every $\mtuple{\x_1}{\x_n}$ in $\R^n$,
$\func{\projection{n}{k}}{\suc{\x_1}{\x_n}}=\x_k$. Furthermore, for every positive integer $n$,
$\Eucbase{n}{}$ denotes the standard ordered-base of $\R^n$ endowed with its canonical linear-structure.
That is, $\Eucbase{n}{}$ is the element of $\Func{\seta{\suc{1}{n}}}{\R^n}$ such that for every $k$
in $\seta{\suc{1}{n}}$, $\func{\projection{n}{j}}{\Eucbase{n}{k}}=\deltaf{k}{j}$.
%%%%%%%%%%%%%%%%%%%%%%%%%
\item[\myitem{DG~5.}]
$\mapdifclass{\infty}{\Man{}}{\Man{1}}$ denotes the set of all mappings
$\cf$ in $\Func{\M{}}{\M{1}}$ such that for every point $\point$ of $\Man{}$, there exists a chart $\phi$ in $\maxatlas{}$
and a chart $\psi$ in $\maxatlas{1}$ such that $\point\in\domain{\phi}$,
%and $\func{\cf}{\point}\in\domain{\psi}$,
$\func{\image{\cf}}{\domain{\phi}}\subseteq\domain{\psi}$, and,
\begin{equation}\label{eqdesmoothnessdefinition}
\cmp{\psi}{\cmp{\cf}{\finv{\phi}}}\in\banachmapdifclass{\infty}{\R^{n}}{\R^{m}}{\funcimage{\phi}}{\funcimage{\psi}}.
\end{equation}
Additionally, given a mapping $\cf$ in $\mapdifclass{\infty}{\Man{}}{\Man{1}}$, for any chart $\phi$ of $\Man{}$
and any chart $\psi$ of $\Man{1}$ such that $\func{\pimage{\cf}}{\domain{\psi}}\cap\domain{\phi}\neq\empty$,
\Ref{eqdesmoothnessdefinition} is satisfied.
Each element of $\mapdifclass{\infty}{\Man{}}{\Man{1}}$ is called a $\quotl$smooth map from $\Man{}$ to $\Man{1}$$\quotr$.
Every smooth map from the manifold $\Man{}$ to the manifold $\Man{1}$ is a continuous map from the underlying topological-space of
$\Man{}$ to the underlying topological-space of $\Man{1}$. That is,
$\mapdifclass{\infty}{\Man{}}{\Man{1}}\subseteq\CF{\mantops{\Man{}}}{\mantops{\Man{1}}}$.
%%%%%%%%%%%%%%%%%%%%%%%%%
\item[\myitem{DG~6.}]
$\Diffeo{\infty}{\Man{}}{\Man{1}}$ denotes the set of all bijective maps $\function{\cf}{\M{}}{\M{1}}$ such that
$\cf$ is a smooth map from $\Man{}$ to $\Man{1}$ and $\finv{\cf}$ is a smooth map from $\Man{1}$ to $\Man{}$. That is,
\begin{equation}
\Diffeo{\infty}{\Man{}}{\Man{1}}:=\defset{\cf}{\IF{\M{}}{\M{1}}}{\cf\in\mapdifclass{\infty}{\Man{}}{\Man{1}},~
{\finv{\cf}}\in\mapdifclass{\infty}{\Man{1}}{\Man{}}}.
\end{equation}
Each element of $\Diffeo{\infty}{\Man{}}{\Man{1}}$ is called an $\quotl$$\infty$-diffeomorphism from $\Man{}$to $\Man{1}$$\quotr$.
It is said that $\quotl$$\Man{}$ is diffeomorphic to $\Man{1}$$\quotr$ iff there exists at least one $\infty$-diffeomorphism
from $\Man{}$ to $\Man{1}$. The existence of diffeomorphism between manifolds clearly induces an equivalence-relation
on a collection of manifolds. Furthermore, evidently
every $\infty$-diffeomorphism from the manifold $\Man{}$ to $\Man{1}$ is a homeomorphism from the underlying topological-space of
$\Man{}$ to the underlying topological-space of $\Man{1}$. That is,
$\Diffeo{\infty}{\Man{}}{\Man{1}}\subseteq\HOF{\mantops{\Man{}}}{\mantops{\Man{1}}}$. So, diffeomorphic manifolds
are homeomorphic, but the converse can not be true in general.\\
The set of all $\infty$-diffeomorphisms from $\Man{}$ to $\Man{}$ is simply denoted by $\Diff{\infty}{\Man{}}$, which
together with the binary operation of function-composition
on them, that is the pair $\opair{\Diff{\infty}{\Man{}}}{\cmp{}{}}$, is obviously a group which is denoted by
$\GDiff{\infty}{\Man{}}$. This group is referred to as the $\quotl$$\infty$-diffeomorphism group of the manifold $\Man{}$$\quotr$, or
the $\quotl$group of $\infty$-automorphisms of $\Man{}$$\quotr$. Each element of the set $\Diff{\infty}{\Man{}}$ is also
referred to as an $\quotl$$\infty$-automorphism of the manifold $\Man{}$$\quotr$.
The composition of any pair of $\infty$-automorphisms of the manifold $\Man{}$
is an $\infty$-automorphism of $\Man{}$.\\
%%%%%%%%%%%%%%%%%%%%%%%%%
Every $\infty$-automorphism of $\Man{}$ brings about a transference of charts of $\Man{}$. Precisely,
given an $\infty$-automorphism $\cf$ of $\Man{}$, for every point $\point$ of $\Man{}$ and every chart $\phi$ of
$\Man{}$ around $\point$, $\cmp{\phi}{\cf}$ is a chart of $\Man{}$ around $\func{\finv{\cf}}{\point}$.
So such transference of charts exhibits a pullback-like behavior when considered pointwise.
Thus for every $\infty$-automorphism $\cf$ of $\Man{}$ and every point $\point$ of $\Man{}$,
it can be defined a mapping $\charttransfer{\Man{}}{\cf}{\point}$ such that,
\begin{align}\label{eqcharttransfer}
&\function{\charttransfer{\Man{}}{\cf}{\point}}{\defset{\psi}{\maxatlas{}}{\point\in\domain{\psi}}}
{\defset{\psi}{\maxatlas{}}{\func{\finv{\cf}}{\point}\in\domain{\psi}}},\cr
&\Foreach{\phi}{\defset{\psi}{\maxatlas{}}{\point\in\domain{\psi}}}
\func{\charttransfer{\Man{}}{\cf}{\point}}{\phi}\eqdef\cmp{\phi}{\cf}.
\end{align}
%%%%%%%%%%%%%%%%%%%%%%%%%
\item[\myitem{DG~7.}]
$\mapdifclass{\infty}{\Man{}}{\RR}$ endowed with its canonical linear-structure is denoted by $\Lmapdifclass{\infty}{\Man{}}{\RR}$
whose addition and scalar-multiplication operations are defined pointwise.
$\rdot$ denotes the usual maltiplication of real-valued functions on $\M{}$, which together with the addition operation
in the linear-structure of $\mapdifclass{\infty}{\Man{}}{\RR}$ induces the structure of a ring on
$\mapdifclass{\infty}{\Man{}}{\RR}$. Also, $\Lmapdifclass{\infty}{\Man{}}{\RR}$
together with the binary operation $\rdot$ is a commutative and associative $\R$-algebra.
Each element of $\mapdifclass{\infty}{\Man{}}{\RR}$ is referred to as a $\quotl$real-valued smooth map on $\Man{}$$\quotr$.
%%%%%%%%%%%%%%%%%%%%%%%%%
\item[\myitem{DG~8.}]
For every $\point$ in $\M{}$, the set of all tangent vectors to $\Man{}$ at $\point$ is denoted by $\tanspace{\point}{\Man{}}$,
and correspondingly the tangent-space of $\Man{}$ at $\point$ is denoted by $\Tanspace{\point}{\Man{}}$ which is a real vector-space,
whose addition and scalar-product operations are defined via any chart around $\point$ and simultaneously independent of the choice of such a chart.
For a given point $\point$ of $\Man{}$, the addition and scalar-product operations of $\Tanspace{\point}{\Man{}}$ are denoted by
$\vsum{\opair{\Man{}}{\point}}$ and $\vsprod{\opair{\Man{}}{\point}}$, respectively. When there is no chance of confusion,
$\vv{1}\vsum{}\vv{2}$ and $\c\vv{}$ can replace $\vv{1}\vsum{\opair{\Man{}}{\point}}\vv{2}$ and $\c\vsprod{\opair{\Man{}}{\point}}\vv{}$, respectively.
For every $\point$ in $\Man{}$ and every chart $\phi$ in $\maxatlas{}$ whose domain contains $\point$, there exists a canonical linear-isomorphism
from $\Tanspace{\point}{\Man{}}$ to $\R^n$ denoted by $\tanspaceiso{\point}{\Man{}}{\phi}$.
The neutral element of addition in the linear-structure of $\Tanspace{\point}{\Man{}}$, that is $\zerovec{\Tanspace{\point}{\Man{}}}$,
is simply denoted by $\zerovec{\point}$.
For every point $\point$ of $\Man{}$, every chart $\phi$ of $\Man{}$, and every $\vv{}$ in $\R^n$,
the element $\func{\finv{\[\tanspaceiso{\point}{\Man{}}{\phi}\]}}{\vv{}}$ of $\tanspace{\point}{\Man{}}$ is called the
$\quotl$vector of the tangent-space of $\Man{}$ at the point $\point$ corresponded to $\vv{}$ with respect to the chart $\phi$$\quotr$.
%%%%%%%%%%%%%%%%%%%%%%%%%
\item[\myitem{DG~9.}]
The set of all tangent vectors to $\Man{}$ is denoted by $\tanbun{\Man{}}$.
The base-point-identifier of $\tanbun{\Man{}}$ is denoted by $\basep{\Man{}}$. That is,
$\basep{\Man{}}$ is defined to be a map from $\tanbun{\Man{}}$ to $\M{}$ such that for every $\point$ in $\M{}$
and every vector $\avec{}$ in $\tanspace{\point}{\Man{}}$, $\func{\basep{\Man{}}}{\avec{}}=\point$.
The tangent-bundle of $\Man{}$ is denoted by $\Tanbun{\Man{}}$, which is a $2n$-dimensional and $\difclass{\infty}$ manifold
built on the set $\tanbun{\Man{}}$ with its canonical differentiable-structure constructed upon the Banach-space $\Cprod{\R^n}{\R^n}$,
obtained from that of $\Man{}$.
Precisely, the set $\defSet{\tanchart{\Man{}}{\phi}}{\phi\in\maxatlas{}}$ is an atlas on $\tanbun{\Man{}}$, where,
\begin{align}\label{eqtangentbundlemaps}
\Foreach{\phi}{\maxatlas{}}
\left\{
\begin{aligned}
&\tanchart{\Man{}}{\phi}\in\Func{\func{\pimage{\basep{\Man{}}}}{\domain{\phi}}}{\Cprod{\funcimage{\phi}}{\R^n}},\cr
&\Foreach{\avec{}}{\func{\pimage{\basep{\Man{}}}}{\domain{\phi}}}
\func{\tanchart{\Man{}}{\phi}}{\avec{}}\eqdef
\opair{\func{\(\cmp{\phi}{\basep{\Man{}}}\)}{\avec{}}}{\func{\(\tanspaceiso{\func{\basep{\Man{}}}{\avec{}}}{\Man{}}{\phi}\)}{\avec{}}}.
\end{aligned}\right.
\end{align}
The $\difclass{\infty}$ maximal-atlas on $\tanbun{\Man{}}$ corresponded to the atlas $\defSet{\tanchart{\Man{}}{\phi}}{\phi\in\maxatlas{}}$
is denoted by $\tanatlas{\Man{}}$ and $\Tanbun{\Man{}}=\opair{\tanbun{\Man{}}}{\tanatlas{\Man{}}}$. For every $\phi$ in $\maxatlas{}$,
$\tanchart{\Man{}}{\phi}$ is referred to as the $\quotl$tangent-bundle chart of $\Man{}$ associated with $\phi$$\quotr$.
%%%%%%%%%%%%%%%%%%%%%%%%%
\item[\myitem{DG~10.}]
The set of all smooth ($\difclass{\infty}$) vector-fields on $\Man{}$ is denoted by $\vecf{\Man{}}{\infty}$.
In other words, $\vecf{\Man{}}{\infty}$ is the set of all maps $\avecf{}$ in $\mapdifclass{\infty}{\Man{}}{\Tanbun{\Man{}}}$
such that $\cmp{\basep{\Man{}}}{\avecf{}}=\identity{\M{}}$.
$\Vecf{\Man{}}{\infty}$ denotes the vector-space obtained from the set $\vecf{\Man{}}{\infty}$ endowed with its canonical linear-structure.
%%%%%%%%%%%%%%%%%%%%%%%%%
\item[\myitem{DG~11.}]
The set of all $\infty$-derivations on $\Man{}$ is denoted by $\Derivation{\Man{}}{\infty}$. In other words,
$\Derivation{\Man{}}{\infty}$ is the set of all linear maps
$\aderivation{}$ in $\Lin{\Lmapdifclass{\infty}{\Man{}}{\RR}}{\Lmapdifclass{\infty}{\Man{}}{\RR}}$ such that
for every pair $\cf$ and $\cg$ of elements of $\mapdifclass{\infty}{\Man{}}{\RR}$,
$\func{\aderivation{}}{\cf\rdot\cg}=\[\func{\aderivation{}}{\cf}\]\rdot\cg+\cf\rdot\[\func{\aderivation{}}{\cg}\]$.
$\LDerivation{\Man{}}{\infty}$ denotes the vector-space obtained from the set $\Derivation{\Man{}}{\infty}$ endowed with its
canonical linear-structure which is inherited from that of $\Lin{\Lmapdifclass{\infty}{\Man{}}{\RR}}{\Lmapdifclass{\infty}{\Man{}}{\RR}}$.
%%%%%%%%%%%%%%%%%%%%%%%%%
\item[\myitem{DG~12.}]
$\derop{\Man{}}{\Man{1}}$ denotes the differential operator on the set $\mapdifclass{\infty}{\Man{}}{\Man{1}}$. That is,
$\derop{\Man{}}{\Man{1}}$ is defined to be the element of
$\Func{\mapdifclass{\infty}{\Man{}}{\Man{1}}}{\mapdifclass{\infty}{\Tanbun{\Man{}}}{\Tanbun{\Man{1}}}}$ such that
for every map $\cf$ in $\mapdifclass{\infty}{\Man{}}{\Man{1}}$, every $\avec{}$ in $\tanbun{\Man{}}$,
and for each chart $\phi$ in $\maxatlas{}$ whose domain contains $\func{\basep{\Man{}}}{\avec{}}$ and
each chart $\psi$ in $\maxatlas{1}$ whose domain contains $\func{\cf}{\func{\basep{\Man{}}}{\avec{}}}$,
\begin{equation}\label{eqdefinitionofdifferentialofamap}
\func{\[\der{\cf}{\Man{}}{\Man{1}}\]}{\avec{}}=
\func{\bigg[\cmp{\finv{\(\tanspaceiso{\func{\cf}{\func{\basep{\Man{}}}{\avec{}}}}{\Man{1}}{\psi}\)}}{\cmp{\bigg(\func{\[\banachder{\(\cmp{\psi}{\cmp{f}{\finv{\phi}}}\)}{\R^n}{\R^m}\]}
{\func{\phi}{\func{\basep{\Man{}}}{\avec{}}}}\bigg)}{\tanspaceiso{\func{\basep{\Man{}}}{\avec{}}}{\Man{}}{\phi}}}\bigg]}{\avec{}}.
\end{equation}
%%%%%%%%%%%%%%%%%
It can be easily seen that for a smooth map $\cf$ in $\mapdifclass{\infty}{\Man{}}{\Man{1}}$,
$\der{\cf}{\Man{}}{\Man{1}}$ maps $\tanspace{\point}{\Man{}}$ into $\tanspace{\func{\cf}{\point}}{\Man{1}}$
for every point $\point$ of $\Man{}$. In other words,
\begin{equation}\label{eqtangentmapbasepoint}
\Foreach{\avec{}}{\tanbun{\Man{}}}
\func{\basep{\Man{}}}{\func{\[\der{\cf}{\Man{}}{\Man{1}}\]}{\avec{}}}=
\func{\[\cmp{\cf}{\basep{\Man{}}}\]}{\avec{}}.
\end{equation}
Furthermore, for every $\cf$ in $\mapdifclass{\infty}{\Man{}}{\Man{1}}$ and each point $\point$ of $\Man{}$,
the restriction of $\der{\cf}{\Man{}}{\Man{1}}$ to $\tanspace{\point}{\Man{}}$ is a linear map from
$\Tanspace{\point}{\Man{}}$ to $\Tanspace{\func{\cf}{\point}}{\Man{1}}$, that is,
\begin{equation}
\Foreach{\point}{\M{}}
\func{\res{\der{\cf}{\Man{}}{\Man{1}}}}{\tanspace{\point}{\Man{}}}\in\Lin{\Tanspace{\point}{\Man{}}}{\Tanspace{\func{\cf}{\point}}{\Man{1}}}.
\end{equation}
%%%%%%%%%%%%%%%%%
The differential operator has local behavior. This means, for any $\cf$ and $\cg$ in $\mapdifclass{\infty}{\Man{}}{\RR}$
that coincide in an open set $\U$ of $\mantops{\Man{}}$ (an element of $\mantop{\Man{}}$),
$\der{\cf}{\Man{}}{\Man{1}}$ and $\der{\cg}{\Man{}}{\Man{1}}$ coincide in $\func{\pimage{\basep{\Man{}}}}{\U}$.\\
The chain rule of differentiation (in the category of smooth manifolds) asserts that for every $\cf$ in
$\mapdifclass{\infty}{\Man{}}{\Man{1}}$ and every $\cg$ in $\mapdifclass{\infty}{\Man{1}}{\Man{2}}$,
\begin{equation}\label{eqchainrule}
\der{\(\cmp{\cg}{\cf}\)}{\Man{}}{\Man{2}}=\cmp{\(\der{\cg}{\Man{1}}{\Man{2}}\)}{\(\der{\cf}{\Man{}}{\Man{1}}\)}.
\end{equation}
So since for every $\cf$ in $\Diffeo{\infty}{\Man{}}{\Man{1}}$, $\cmp{\finv{\cf}}{\cf}=\identity{\M{}}$, and
$\der{\identity{\M{}}}{\Man{}}{\Man{}}=\identity{\tanbun{\Man{}}}$, clearly,
\begin{equation}\label{eqdiffeomorphismdifrule}
\Foreach{\cf}{\Diffeo{\infty}{\Man{}}{\Man{1}}}
\finv{\(\der{\cf}{\Man{}}{\Man{1}}\)}=\der{\finv{\cf}}{\Man{1}}{\Man{}}.
\end{equation}
\\$\der{\cf}{\Man{}}{\Man{1}}$ can simply be denoted by $\derr{\cf}{}{}$ when there is no ambiguity about the underlying
source and target manifolds $\Man{}$ and $\Man{1}$.
%%%%%%%%%%%%%%%%%%%%%%%%%
\item[\myitem{DG~13.}]
$\Rderop{\Man{}}$ denotes the derivative operator on $\mapdifclass{\infty}{\Man{}}{\RR}$. That is,
$\Rderop{\Man{}}$ is defined to be the element of
$\Func{\mapdifclass{\infty}{\Man{}}{\RR}}{\mapdifclass{\infty}{\Tanbun{\Man{}}}{\RR}}$ such that
for every smooth map $\cf$ in $\mapdifclass{\infty}{\Man{}}{\RR}$, every $\avec{}$ in $\tanbun{\Man{}}$,
\begin{align}\label{eqdefderivativeoperator}
\func{\[\Rder{\cf}{\Man{}}\]}{\avec{}}:=&\func{\tanspaceiso{\func{\cf}{\func{\basep{\Man{}}}{\avec{}}}}{\RR}
{\identity{\R}}}{\func{\[\der{\cf}{\Man{}}{\RR}\]}{\avec{}}}\cr
=&\func{\[\cmp{\bigg(\func{\[\banachder{\(\cmp{f}{\finv{\phi}}\)}{\R^n}{\R}\]}
{\func{\phi}{\func{\basep{\Man{}}}{\avec{}}}}\bigg)}{\tanspaceiso{\func{\basep{\Man{}}}{\avec{}}}{\Man{}}{\phi}}\]}{\avec{}},
\end{align}
where $\phi$ can be any of charts in $\maxatlas{}$ whose domain contains $\func{\basep{\Man{}}}{\point}$.
For every $\cf$ in $\mapdifclass{\infty}{\Man{}}{\RR}$, $\Rderop{\Man{}}$
%$\Rder{\cf}{\Man{}}$ is a map from $\tanbun{\Man{}}$ to $\R$ which
can be regarded as a replacement for the differential of $\cf$.
%For every $\cf$ in $\mapdifclass{\infty}{\Man{}}{\RR}$, $\Rder{\cf}{\Man{}}\in\mapdifclass{\infty}{\Tanbun{\Man{}}}{\RR}$.
Equipping $\mapdifclass{\infty}{\Man{}}{\RR}$ and $\mapdifclass{\infty}{\Tanbun{\Man{}}}{\RR}$ with
their natural linear-structure, $\Rderop{\Man{}}$ is a linear map as a direct result of a similar property for
the derivative operator of real-valued smooth maps on Banach-spaces. That is,
$\Rderop{\Man{}}\in\Lin{\mapdifclass{\infty}{\Man{}}{\RR}}{\mapdifclass{\infty}{\Tanbun{\Man{}}}{\RR}}$.
As another crucial property of $\Rderop{\Man{}}$, for every $\cf$ and $\cg$ in $\mapdifclass{\infty}{\Man{}}{\RR}$,
and every $\avec{}$ in $\tanbun{\Man{}}$,
\begin{equation}\label{Jacobipropertyofderivative}
\func{\[\func{\Rderop{\Man{}}}{\cf\rdot\cg}\]}{\avec{}}=\[\func{\(\cmp{\cg}{\basep{\Man{}}}\)}{\avec{}}\]\[\func{\Rder{\cf}{\Man{}}}{\avec{}}\]+
\[\func{\(\cmp{\cf}{\basep{\Man{}}}\)}{\avec{}}\]\[\func{\Rder{\cg}{\Man{}}}{\avec{}}\].
\end{equation}
As a trivial consequence of the general case,
for every $\cf$ in $\mapdifclass{\infty}{\Man{}}{\RR}$ and each point $\point$ of $\Man{}$,
the restriction of $\Rder{\cf}{\Man{}}$ to $\tanspace{\point}{\Man{}}$ is a linear map from
$\Tanspace{\point}{\Man{}}$ to $\R$, that is,
$\func{\res{\Rder{\cf}{\Man{}}}}{\tanspace{\point}{\Man{}}}\in\Lin{\Tanspace{\point}{\Man{}}}{\R}$.
%%%%%%%%%%%%%%%%%%%%%%%%%%%
\\According to \Ref{eqtangentmapbasepoint}, \Ref{eqchainrule}, and \Ref{eqdefderivativeoperator},
it is evident that,
\begin{equation}\label{eqspecialchainrule1}
\Foreach{\cf}{\mapdifclass{\infty}{\Man{1}}{\RR}}
\Foreach{\cg}{\mapdifclass{\infty}{\Man{}}{\Man{1}}}
\Rder{\(\cmp{\cf}{\cg}\)}{\Man{}}=\cmp{\(\Rder{\cf}{\Man{1}}\)}{\(\der{\cg}{\Man{}}{\Man{1}}\)}.
\end{equation}
This is called the $\quotl$special chain-rule of differentiation$\quotr$, here.
%%%%%%%%%%%%%%%%%%%%%%%%%%%
$\Rder{\cf}{\Man{}}$ can simply be denoted by $\Rder{\cf}{}$ when there is no ambiguity about the underlying
manifold $\Man{}$, and it can be simply referred to as the $\quotl$derivative of $\cf$$\quotr$. $\func{\Rder{\cf}{\Man{}}}{\avec{}}$
can simply be referred to as the $\quotl$derivative of $\cf$ in the direction $\avec{}$$\quotr$.
%%%%%%%%%%%%%%%%%%%%%%%%%
\item[\myitem{DG~14.}]
The set $\atlas{}:=\defSet{\funcprod{\phi}{\psi}}{\phi\in\maxatlas{},~\psi\in\maxatlas{1}}$
is an element of $\atlases{\infty}{\Cprod{\M{}}{\M{1}}}{\R^{n+m}}$, that is, a $\difclass{\infty}$ atlas on $\Cprod{\M{}}{\M{1}}$
constructed upon $\R^{n+m}$. The topology on $\Cprod{\M{}}{\M{1}}$ induced by the $\difclass{\infty}$ maximal-atlas on
$\Cprod{\M{}}{\M{1}}$ constructed upon $\R^{n+m}$ that is generated by the atlas $\atlas{}$
(that is, $\func{\maxatlasgen{\infty}{\Cprod{\M{}}{\M{1}}}{\R^{n+m}}}{\atlas{}}$),
is the topological-product of the underlying topological-spaces of the manifolds $\Man{}$ and $\Man{1}$, that is,
$\topprod{\mantops{\Man{}}}{\mantops{\Man{1}}}$. Since $\mantops{\Man{}}$ and $\mantops{\Man{1}}$
are Hausdorff and second-countable spaces, so is the topological space $\topprod{\mantops{\Man{}}}{\mantops{\Man{1}}}$.
Consequently, the set $\Cprod{\M{}}{\M{1}}$ endowed with the $\difclass{\infty}$ maximal-atlas
$\func{\maxatlasgen{\infty}{\Cprod{\M{}}{\M{1}}}{\R^{n+m}}}{\atlas{}}$ is a differentiable-structure
whose underlying topological-space is Hausdorff and second-countable. Thus,
$\opair{\Cprod{\M{}}{\M{1}}}{\func{\maxatlasgen{\infty}{\Cprod{\M{}}{\M{1}}}{\R^{n+m}}}{\atlas{}}}$ is a manifold.
$\manprod{\Man{}}{\Man{1}}$ denotes this manifold and is referred to as
the $\quotl$manifold-product of the manifolds $\Man{}$ and $\Man{1}$$\quotr$.
%%%%%%%%%%%%%%%%%%%%%%%%%
\item[\myitem{DG~15.}]
For every $\cf$ in $\mapdifclass{\infty}{\Man{}}{\Man{2}}$ and every $\cg$ in $\mapdifclass{\infty}{\Man{1}}{\Man{3}}$,
the product mapping $\funcprod{\cf}{\cg}$ is a smooth map from the manifold-product of $\Man{}$ and $\Man{1}$ to the
manifold-product of $\Man{2}$ and $\Man{3}$, that is,
$\(\funcprod{\cf}{\cg}\)\in\mapdifclass{\infty}{\manprod{\Man{}}{\Man{1}}}{\manprod{\Man{2}}{\Man{3}}}$.
%%%%%%%%%%%%%%%%%%%%%%%%%
\item[\myitem{DG~16.}]
Every smooth map $\cf$ from $\Man{}$ to $\Man{1}$ whose differential becomes a a linear monomorphism from
$\tanspace{\point}{\Man{}}$ to $\tanspace{\func{\cf}{\point}}{\Man{1}}$ when restricted to $\tanspace{\point}{\Man{}}$
for every point $\point$ of $\Man{}$,
is referred to as a $\quotl$$\infty$-immersion from the manifold $\Man{}$ to the manifold $\Man{1}$$\quotr$.
The set of all $\infty$-immersions from $\Man{}$ to $\Man{1}$ is denoted by $\immersion{\infty}{\Man{}}{\Man{1}}$.\\
Every $\infty$-immersion from $\Man{}$ to $\Man{1}$ which is additionally a topological embedding of
$\mantops{\Man{}}$ (the underlying topological-space of $\Man{}$) into $\mantops{\Man{1}}$ (the underlying
topological-space of $\Man{1}$), is referred to as an $\quotl$$\infty$-embedding of $\Man{}$ into $\Man{1}$$\quotr$.
The set of all $\infty$-embedding from $\Man{}$ to $\Man{1}$ is denoted by $\embedding{\infty}{\Man{}}{\Man{1}}$.\\
Every smooth map $\cf$ from $\Man{}$ to $\Man{1}$ whose differential becomes a a linear epimorphism from
$\tanspace{\point}{\Man{}}$ and $\tanspace{\func{\cf}{\point}}{\Man{1}}$ when restricted to $\tanspace{\point}{\Man{}}$
for every point $\point$ of $\Man{}$,
is referred to as a $\quotl$$\infty$-submersion from the manifold $\Man{}$ to the manifold $\Man{1}$$\quotr$.
The set of all $\infty$-submersions from $\Man{}$ to $\Man{1}$ is denoted by $\submersion{\infty}{\Man{}}{\Man{1}}$.
%%%%%%%%%%%%%%%%%%%%%%%%%
\item[\myitem{DG~17.}]
$\Man{1}$ is said to be an $\quotl$($m_{1}$-dimensional) $\infty$-immersed submanifold of $\Man{}$$\quotr$ or
simply an $\quotl$immersed submanifold of $\Man{}$$\quotr$ iff $\M{1}\subseteq\M{}$ and
the injection of $\M{1}$ into $\M{}$ is an $\infty$-immersion from $\Man{1}$ to $\Man{}$, that is
$\Injection{\M{1}}{\M{}}\in\immersion{\infty}{\Man{1}}{\Man{}}$.\\
$\Man{1}$ is said to be an $\quotl$($m_{1}$-dimensional) $\infty$-embedded submanifold of $\Man{}$$\quotr$ or
simply an $\quotl$embedded submanifold of $\Man{}$$\quotr$ iff $\M{1}\subseteq\M{}$ and
the injection of $\M{1}$ into $\M{}$ is an $\infty$-embedding from $\Man{1}$ to $\Man{}$, that is
$\Injection{\M{1}}{\M{}}\in\embedding{\infty}{\Man{1}}{\Man{}}$. It is clear that every
embedded submanifold of $\Man{}$ is necessarily an immersed submanifold of $\Man{}$.\\
given a subset $\SET{}$ of $\M{}$, there is at most one corresponded differentiable structure modeled on $\R^{q}$ for some $q\in\Zp$
that makes it a $\difclass{\infty}$ manifold which is an embedded submanifold of $\Man{}$;
If this unique differentiable structure exists, it is inherited from the maximal atlas of $\Man{}$ in a
canonical way. The set of all subsets of $\M{}$ accepting a differentiable structure modeled on a Euclidean-space
that makes it an embedded submanifold of $\Man{}$ is denoted by $\Emsubman{\Man{}}$.
For every $\SET{}$ in $\Emsubman{\Man{}}$, the unique $\difclass{\infty}$
manifold with the set of points $\SET{}$ that is an embedded submanifold of $\Man{}$ is denoted by $\emsubman{\Man{}}{\SET{}}$ or $\subman{\Man{}}{\SET{}}$.
Each element of $\Emsubman{\Man{}}$ is called a $\quotl$embedded set of the manifold $\Man{}$$\quotr$.
Every open set of $\Man{}$ is an embedded set of $\Man{}$, that is $\mantop{\Man{}}\subseteq\Emsubman{\Man{}}$.
Furthermore, every embedded set of $\Man{}$ is a locally-closed set of the underlying topological-space of $\Man{1}$,
which means every embedded set of $\Man{}$ is open in its closure (with respect to its topology inherited from the topology of $\Man{}$).\\
If $\Man{2}$ is an immersed submanifold of $\Man{1}$ and $\Man{1}$ is an immersed submanifold of $\Man{}$, then
$\Man{2}$ is an immersed submanifold of $\Man{}$. Additionally, given subsets $\SET{1}$ and $\SET{2}$ of $\M{}$
such that $\SET{2}\subseteq\SET{1}$, if $\SET{1}\in\Emsubman{\Man{}}$ and $\SET{2}\in\Emsubman{\emsubman{\Man{}}{\SET{1}}}$,
then $\SET{2}\in\Emsubman{\Man{}}$.
%%%%%%%%%%%%%%%%%%%%%%%%%
\item[\myitem{DG~18.}]
For every smooth map $\cf$ from $\Man{}$ to $\Man{1}$, and every element $\SET{}$ of $\Emsubman{\Man{}}$ (every embedded set of $\Man{}$),
the domain-restriction of $\cf$ to $\SET{}$ is also a smooth map from the embedded submanifold $\emsubman{\Man{}}{\SET{}}$ of
$\Man{}$ to $\Man{1}$.\\
For every map $\cf$ from $\Man{}$ to $\Man{1}$, and every element $\SET{}$ of $\Emsubman{\Man{1}}$ (every embedded set of $\Man{1}$)
that includes the image of $\cf$,
$\cf$ is a smooth map from $\Man{}$ to $\Man{1}$ if and only if
the codomain-restriction of $\cf$ to $\SET{}$ is a smooth map from the $\Man{}$ to the embedded submanifold
$\emsubman{\Man{1}}{\SET{}}$ of $\Man{1}$.\\
The embedded sets of the manifold $\Man{}$ is in a one-to-one correspondence with the embedded sets of
any manifold diffeomorphic to $\Man{}$. Actually,
Every $\infty$-diffeomorphism $\cf$ from $\Man{}$ to $\Man{1}$ maps an embedded set of $\Man{}$ to an embedded set of
$\Man{1}$. Therefeore the restriction (that is simultaneously domain-restriction and codomain-restriction) of
a smooth map $\cf$ from $\Man{}$ to $\Man{1}$ to an embedded set $\SET{}$ of $\Man{}$ is again a smooth map
from the embedded submanifold $\emsubman{\Man{}}{\SET{}}$ of $\Man{}$ to the embedded submanifold
$\emsubman{\Man{1}}{\func{\image{\cf}}{\SET{}}}$ of $\Man{1}$.
%%%%%%%%%%%%%%%%%%%%%%%%%
\item[\myitem{DG~19.}]
Given an $\infty$-immersion $\cf$ from $\Man{}$ to $\Man{1}$, for every point $\point$ of $\Man{}$,
there exists an open set $\U$ of the manifold $\Man{}$ such that $\func{\image{\cf}}{\U}$ is an embedded
set of the manifold $\Man{1}$, and the restriction of $\cf$ to $\U$ is an $\infty$-diffeomorphism from
the embedded submanifold $\emsubman{\Man{}}{\U}$ to the embedded submanifold $\emsubman{\Man{1}}{\func{\image{\cf}}{\U}}$
of $\Man{1}$.
%%%%%%%%%%%%%%%%%%%%%%%%%
\item[\myitem{DG~20.}]
Given a positive integer $n$, and a non-empty open set $\U$ of $\topR{\R^n}$, trivially $\U$ is an embedded
set of the canonical smooth manifold of $\R^n$, and the manifold $\emsubman{\RR^n}{\U}$ is simply denoted by $\Ropenman{\U}{n}$,
and is called the $\quotl$canonical differentiable structure of $\U$ inherited from $\RR^n$$\quotr$.
It is also a trivial fact that $\Ropenman{\U}{n}$ is a $\difclass{\infty}$ differentiable structure modeled on $\R^n$, the singleton
$\seta{\identity{\U}}$ being an atlas and consequently $\identity{\U}$ a chart of it.
%%%%%%%%%%%%%%%%%%%%%%%%%
\item[\myitem{DG~21.}]
The notion of partial derivative of differentiable mappings from a product of Banach spaces to a Banach space
as introduced in \cite{Cartan}, can be generalized to the case of smooth mappings from a product of manifolds
to a manifold. Here, only the simplest case is considered, where the domain of a smooth mapping is the product of
a pair of manifolds.\\
$\point_1$ and $\point_2$ are taken as points of the manifolds $\Man{1}$ and $\Man{2}$, respectively.
The mappings $\function{\leftparinj{\M{1}}{\M{2}}{\point_2}}{\M{1}}{\Cprod{\M{1}}{\M{2}}}$ and
$\function{\rightparinj{\M{1}}{\M{2}}{\point_2}}{\M{2}}{\Cprod{\M{1}}{\M{2}}}$ are defined as,
\begin{align}
&\Foreach{\x}{\M{1}}
\func{\leftparinj{\M{1}}{\M{2}}{\point_2}}{\x}\eqdef\opair{\x}{\point_2},\\
&\Foreach{\x}{\M{2}}
\func{\rightparinj{\M{1}}{\M{2}}{\point_1}}{\x}\eqdef\opair{\point_1}{\x}.
\end{align}
It is a trivial fact that $\leftparinj{\M{1}}{\M{2}}{\point_2}$ is a smooth map from $\Man{1}$ to
$\manprod{\Man{1}}{\Man{2}}$, and $\rightparinj{\M{1}}{\M{2}}{\point_1}$ a smooth map
from $\Man{2}$ to $\manprod{\Man{1}}{\Man{2}}$. The mapping $\prodmantan{\Man{1}}{\Man{2}}{\point_1}{\point_2}$
is defined as,
\begin{align}
&\prodmantan{\Man{1}}{\Man{2}}{\point_1}{\point_2}\indef\Func{\Cprod{\tanspace{\point_1}{\Man{1}}}{\tanspace{\point_2}{\Man{2}}}}
{\tanspace{\opair{\point_1}{\point_2}}{\manprod{\Man{1}}{\Man{2}}}},\cr
&\begin{aligned}
&\Foreach{\opair{\vv{1}}{\vv{2}}}{\Cprod{\tanspace{\point_1}{\Man{1}}}{\tanspace{\point_2}{\Man{2}}}}\cr
&\func{\prodmantan{\Man{1}}{\Man{2}}{\point_1}{\point_2}}{\binary{\vv{1}}{\vv{2}}}\eqdef
\func{\[\der{\leftparinj{\M{1}}{\M{2}}{\point_2}}{\Man{1}}{\manprod{\Man{1}}{\Man{2}}}\]}{\vv{1}}+
\func{\[\der{\rightparinj{\M{1}}{\M{2}}{\point_1}}{\Man{2}}{\manprod{\Man{1}}{\Man{2}}}\]}{\vv{2}}.
\end{aligned}\cr
&{}
\end{align}
Actually, $\prodmantan{\Man{1}}{\Man{2}}{\point_1}{\point_2}$ is a bijection and hence the tangent-space of $\manprod{\Man{1}}{\Man{2}}$
at the point $\opair{\point_1}{\point_2}$ can be identified with $\Cprod{\tanspace{\point_1}{\Man{1}}}{\tanspace{\point_2}{\Man{2}}}$
via $\prodmantan{\Man{1}}{\Man{2}}{\point_1}{\point_2}$.
\\$\cf$ is taken as a smooth map from the manifold-product of $\Man{1}$ and $\Man{2}$ to the manifold $\Man{}$, that is an element of
$\mapdifclass{\infty}{\manprod{\Man{1}}{\Man{2}}}{\Man{}}$.
\begin{align}\label{EQpartialsLemma}
&\Foreach{\opair{\vv{1}}{\vv{2}}}{\Cprod{\tanspace{\point_1}{\Man{1}}}{\tanspace{\point_2}{\Man{2}}}}\cr
&\begin{aligned}
\func{\[\der{\cf}{\manprod{\Man{1}}{\Man{2}}}{\Man{}}\]}{\func{\prodmantan{\Man{1}}{\Man{2}}{\point_1}{\point_2}}{\binary{\vv{1}}{\vv{2}}}}=
&\hskip0.5\baselineskip\func{\[\der{\(\cmp{\cf}{\leftparinj{\M{1}}{\M{2}}{\point_2}}\)}{\Man{1}}{\Man{}}\]}{\vv{1}}\cr
&+\func{\[\der{\(\cmp{\cf}{\rightparinj{\M{1}}{\M{2}}{\point_1}}\)}{\Man{2}}{\Man{}}\]}{\vv{2}},
\end{aligned}
\end{align}
where the addition takes place in the tangent-space of $\manprod{\Man{1}}{\Man{2}}$ at $\func{\cf}{\binary{\point_1}{\point_2}}$.
\end{itemize}
\chapteR{Multi-linear Algebra}
\thispagestyle{fancy}
%%%%%%%%%%%%%
\noindent
\textit{Here the existence of multiplicative identity is assumed as an axiom within the definition of the algebraic structure of a ring.}\\
\textit{All occurrences of modules are taken to be left modules.}
%%%%%%%%%%%%%%%%%%%%%%%%%%%%%%%%%%%%%%%%%%%%%%%%%%%%%%%%%%%%%%%%%%%%%%%%%%%%%%%%%%%%%%%%%%%%%%%%%%%%%%%%%%%%%%%%%%%%%%%%%%%%%%%%%%%%%%%%%%%%%
\section{Tensors}
\subsection{Tensors on Modules}
\fixed
$\CRing{}$ is fixed as a commutative ring, and $\CModule{}=\quadruple{\Cmodule{}}{+}{\times}{\CRing{}}$ as a
$\CRing{}$-module.
\endfixed
%%%%%%%%%%%%%%%%%%%%%%%%%%%%%%%%%%%%%%%%%%%%%%%%%%%%%%%%%%%%%%%%%%%%%%%%%%%%%%%%%%%%%%%%%%%%%%%%%%%%%%%%%%%%%%%%%%%%%%%%%%%%%%%%%%%%%%%%%%%%%
\definition
Let $r$ and $s$ be non-negative integers such that $r+s>0$, and let
$\mlinmap{}$ be a function from $\Cprod{\multiprod{\CModule{}}{r}}{\multiprod{\Vdual{\CModule{}}}{s}}$ to $\CRing{}$.
$\mlinmap{}$ is called a $\quotl$$\opair{r}{s}$-tensor on $\CModule{}$$\quotr$ or a
$\quotl$multi-linear form (of type $\opair{r}{s}$) on $\CModule{}$$\quotr$ iff $\mlinmap{}$ is linear in each factor, that is,
\begin{align}
&\Foreach{j}{\seta{\suc{1}{r}}}
\Foreach{\mtuple{\vv{1}}{\binary{\vv{r}}{\ww{}}}}{\multiprod{\CModule{}}{r+1}}
\Foreach{\mtuple{\uu{1}}{\uu{s}}}{\multiprod{\Vdual{\CModule{}}}{s}}
\Foreach{c}{\CRing{}}\cr
&\hskip\baselineskip\begin{aligned}
&~\func{\mlinmap{}}{\suc{\vv{1}}{\vv{j}+c\ww{}}\binary{\suc{}{\vv{r}}}{\suc{\uu{1}}{\uu{s}}}}\cr
=&~\func{\mlinmap{}}{\binary{\suc{\vv{1}}{\vv{r}}}{\suc{\uu{1}}{\uu{s}}}}+
c\func{\mlinmap{}}{\binary{\suc{\vv{1}}{\binary{\vv{j-1}}{\ww{}}}}{\binary{\suc{\vv{j+1}}{\vv{r}}}{\suc{\uu{1}}{\uu{s}}}}},
\end{aligned}
\end{align}
and
\begin{align}
&\Foreach{j}{\seta{\suc{1}{s}}}
\Foreach{\mtuple{\vv{1}}{\vv{r}}}{\multiprod{\CModule{}}{r}}
\Foreach{\mtuple{\uu{1}}{\binary{\uu{s}}{\ww{}}}}{\multiprod{\Vdual{\CModule{}}}{s+1}}
\Foreach{c}{\CRing{}}\cr
&\hskip\baselineskip\begin{aligned}
&~\func{\mlinmap{}}{\binary{\suc{\vv{1}}{\vv{r}}}{\suc{\uu{1}}{}\suc{\uu{j}+c\ww{}}{\uu{s}}}}\cr
=&~\func{\mlinmap{}}{\binary{\suc{\vv{1}}{\vv{r}}}{\suc{\uu{1}}{\uu{s}}}}+
c\func{\mlinmap{}}{\binary{\suc{\vv{1}}{\vv{r}}}{\binary{\suc{\uu{1}}{\binary{\uu{j-1}}{\ww{}}}}{\suc{\uu{j+1}}{\uu{s}}}}}.
\end{aligned}
\end{align}
The set of all $\opair{r}{s}$-tensors on $\CModule{}$ is denoted by $\Tensors{r}{s}{\CModule{}}$.
Additionally, by convention, $\Tensors{0}{0}{\CModule{}}:=\CRing{}$.\\
\caution
Any $\opair{r}{s}$-tensor on $\CModule{}$ is also called a $\quotl$$r$-covariant and $s$-contravariant tensor on $\CModule{}$$\quotr$.
Any element of $\Tensors{r}{0}{\CModule{}}$ is called a $\quotl$covariant tensor of rank $r$ on $\CModule{}$$\quotr$,
and any element of $\Tensors{0}{s}{\CModule{}}$ is called a $\quotl$contravariant tensor of rank $s$ on $\CModule{}$$\quotr$.
\endef
%%%%%%%%%%%%%%%%%%%%%%%%%%%%%%%%%%%%%%%%%%%%%%%%%%%%%%%%%%%%%%%%%%%%%%%%%%%%%%%%%%%%%%%%%%%%%%%%%%%%%%%%%%%%%%%%%%%%%%%%%%%%%%%%%%%%%%%%%%%%%
\remark
\textit{The notion of a multi-linear form in the class of vector-spaces is adopted as that of the class of modules,
considering that any vector-space over a field is inherently a module over a commutative ring.}
\endremark
%%%%%%%%%%%%%%%%%%%%%%%%%%%%%%%%%%%%%%%%%%%%%%%%%%%%%%%%%%%%%%%%%%%%%%%%%%%%%%%%%%%%%%%%%%%%%%%%%%%%%%%%%%%%%%%%%%%%%%%%%%%%%%%%%%%%%%%%%%%%%
\theorem
Let $r$ and $s$ be non-negative integers. Endowing $\Tensors{r}{s}{\CModule{}}$ with the operations of addition and
scalar multiplication defined as
\begin{align}
&\Foreach{\mtuple{\vv{1}}{\vv{r}}}{{\CModule{}}^{\times r}}
\Foreach{\mtuple{\uu{1}}{\uu{s}}}{{\Vdual{\CModule{}}}^{\times s}}\cr
&\begin{cases}
\Foreach{\opair{\alpha}{\beta}}{\Cprod{\Tensors{r}{s}{\CModule{}}}{\Tensors{r}{s}{\CModule{}}}}
\begin{aligned}
&~\func{\[\alpha+\beta\]}{\binary{\suc{\vv{1}}{\vv{r}}}{\suc{\uu{1}}{\uu{s}}}}\cr
\eqdef&~\func{\alpha}{\binary{\suc{\vv{1}}{\vv{r}}}{\suc{\uu{1}}{\uu{s}}}}+
\func{\beta}{\binary{\suc{\vv{1}}{\vv{r}}}{\suc{\uu{1}}{\uu{s}}}},
\end{aligned}\cr
\Foreach{\opair{c}{\alpha}}{\Cprod{\CRing{}}{\Tensors{r}{s}{\CModule{}}}}
\func{\(c\alpha\)}{\binary{\suc{\vv{1}}{\vv{r}}}{\suc{\uu{1}}{\uu{s}}}}\eqdef
c\func{\alpha}{\binary{\suc{\vv{1}}{\vv{r}}}{\suc{\uu{1}}{\uu{s}}}},
\end{cases}\cr
&{}
\end{align}
when $r+s>0$, and defined as the same addition and multiplication operations of the ring $\CRing{}$ when
$r=s=0$,
turns $\Tensors{r}{s}{\CModule{}}$ into a $\CRing{}$-module.
\proof
It is trivial to check the axioms of the module structure.
\endthm
%%%%%%%%%%%%%%%%%%%%%%%%%%%%%%%%%%%%%%%%%%%%%%%%%%%%%%%%%%%%%%%%%%%%%%%%%%%%%%%%%%%%%%%%%%%%%%%%%%%%%%%%%%%%%%%%%%%%%%%%%%%%%%%%%%%%%%%%%%%%%
\definition
Let $r$ and $s$ be non-negative integers. The $\CRing{}$-module obtained by endowing
$\Tensors{r}{s}{\CModule{}}$ with the natural operations of addition and scalar multiplication
defined in the theorem above, will be denoted by $\MTensors{r}{s}{\CModule{}}$.
\endef
%%%%%%%%%%%%%%%%%%%%%%%%%%%%%%%%%%%%%%%%%%%%%%%%%%%%%%%%%%%%%%%%%%%%%%%%%%%%%%%%%%%%%%%%%%%%%%%%%%%%%%%%%%%%%%%%%%%%%%%%%%%%%%%%%%%%%%%%%%%%%
\definition
The direct sum of all $\MTensors{r}{s}{\CModule{}}$-s when $r$ and $s$ range over all non-negative integers simoltaneously,
will simply be denoted by $\DTensors{\CModule{}}$. That is,
\begin{equation}
\DTensors{\CModule{}}:=\Dsum{\defSet{\MTensors{r}{s}{\CModule{}}}{r,s=0,1,\ldots}}.
\end{equation}
$\DTensors{\CModule{}}$ is referred to as the $\quotl$space of tensors on $\CModule{}$$\quotr$.
For every pair $r$ and $s$ of non-negative integers, there is a natural injection of $\Tensors{r}{s}{\CModule{}}$ into
$\DTensors{\CModule{}}$ via which every $\opair{r}{s}$-tensor on $\CModule{}$ can be identified with a unique element of
$\DTensors{\CModule{}}$.
\endef
%%%%%%%%%%%%%%%%%%%%%%%%%%%%%%%%%%%%%%%%%%%%%%%%%%%%%%%%%%%%%%%%%%%%%%%%%%%%%%%%%%%%%%%%%%%%%%%%%%%%%%%%%%%%%%%%%%%%%%%%%%%%%%%%%%%%%%%%%%%%%
\definition
The binary operation $\tensor{\CModule{}}$ on $\DTensors{\CModule{}}$ is defined to be the unique bilinear map
$\function{\tensor{\CModule{}}}{\Cprod{\DTensors{\CModule{}}}{\DTensors{\CModule{}}}}{\DTensors{\CModule{}}}$
such that for every non-negative integers $r$, $s$, $p$, $q$, if $r+s>0$ and $p+q>0$,
\begin{align}
&\Foreach{\opair{\alpha}{\beta}}{\Cprod{\Tensors{r}{s}{\CModule{}}}{\Tensors{p}{q}{\CModule{}}}}
\Foreach{\suc{\vv{1}}{\vv{r+p}}}{{\Cmodule{}}^{\times{r+p}}}
\Foreach{\suc{\uu{1}}{\uu{s+q}}}{{\Vdual{\CModule{}}}^{\times{s+q}}}\cr
&\hskip\baselineskip\begin{aligned}
&~\func{\(\alpha\tensor{\CModule{}}\beta\)}{\binary{\suc{\vv{1}}{\vv{r+p}}}{\suc{\uu{1}}{\uu{s+q}}}}\cr
\eqdef&~\func{\alpha}{\binary{\suc{\vv{1}}{\vv{r}}}{\suc{\uu{1}}{\uu{s}}}}
\func{\beta}{\binary{\suc{\vv{r+1}}{\vv{r+p}}}{\suc{\uu{s+1}}{\uu{s+q}}}},
\end{aligned}
\end{align}
and if $r+s=0$,
\begin{equation}
\Foreach{\opair{\alpha}{\beta}}{\Cprod{\Tensors{r}{s}{\CModule{}}}{\Tensors{p}{q}{\CModule{}}}}
\alpha\tensor{\CModule{}}\beta\eqdef\alpha\beta,
\end{equation}
and if $p+q=0$,
\begin{equation}
\Foreach{\opair{\alpha}{\beta}}{\Cprod{\Tensors{r}{s}{\CModule{}}}{\Tensors{p}{q}{\CModule{}}}}
\alpha\tensor{\CModule{}}\beta\eqdef\beta\alpha.
\end{equation}
$\tensor{\CModule{}}$ is called the $\quotl$tensor operation of the $\CRing{}$-module $\CModule{}$$\quotr$.
For any non-negative integers $r$ and $s$, each element of $\Tensors{r}{s}{\CModule{}}$ is also called a
$\quotl$(type) $\opair{r}{s}$ simple tensor (on $\CModule{}$)$\quotr$. Also for every $\alpha$ and $\beta$
in $\DTensors{\CModule{}}$, $\alpha\tensor{\CModule{}}\beta$ is called the $\quotl$tensor product of the tensors $\alpha$
and $\beta$$\quotr$.
When there is no ambiguity about the underlying module, $\tensor{\CModule{}}$ can simply be denoted by $\tensor{}$.\\
\caution
It is trivial that for any $\alpha\in\Tensors{r}{s}{\CModule{}}\subseteq\DTensors{\CModule{}}$ and any
$\beta\in\Tensors{p}{q}{\CModule{}}\subseteq\DTensors{\CModule{}}$,
$\alpha\tensor{\CModule{}}\beta\in\Tensors{r+p}{q+s}{\CModule{}}\subseteq\DTensors{\CModule{}}$ and $\tensor{\CModule{}}$ acts bilinearly
on $\Cprod{\Tensors{r}{s}{\CModule{}}}{\Tensors{p}{q}{\CModule{}}}$. We presumed the triviality of these facts prior to the definition.
\endef
%%%%%%%%%%%%%%%%%%%%%%%%%%%%%%%%%%%%%%%%%%%%%%%%%%%%%%%%%%%%%%%%%%%%%%%%%%%%%%%%%%%%%%%%%%%%%%%%%%%%%%%%%%%%%%%%%%%%%%%%%%%%%%%%%%%%%%%%%%%%%
\theorem\label{thmtensoroperationofmoduleisassociative}
$\tensor{\CModule{}}$ is an associative binary operation.
\proof
Let $\alpha$, $\beta$, and $\gamma$ be arbitrary elements of $\DTensors{\CModule{}}$.\\
Suppose that $\alpha$, $\beta$
and $\gamma$ are simple tensors. Then there exists non-negative integers $l$, $m$, $p$, $q$, $r$, and $s$ such that
$\alpha\in\Tensors{l}{m}{\CModule{}}$, $\beta\in\Tensors{p}{q}{\CModule{}}$, and
$\gamma\in\Tensors{r}{s}{\CModule{}}$. Clearly,
for every $\mtuple{\vv{1}}{\vv{l+p+r}}\in{\CModule{}}^{\times l+p+r}$ and every
$\mtuple{\uu{1}}{\uu{m+q+s}}\in{\Vdual{\CModule{}}}^{\times m+q+s}$,
\begin{align}
&\hskip0.6\baselineskip\func{\[\alpha\tensor{}\(\beta\tensor{}\gamma\)\]}{\binary{\suc{\vv{1}}{\vv{l+p+r}}}{\suc{\uu{1}}{\uu{m+q+s}}}}\cr
&=~\func{\alpha}{\binary{\suc{\vv{1}}{\vv{l}}}{\suc{\uu{1}}{\uu{m}}}}
\[\func{\beta\tensor{}\gamma}{\binary{\suc{\vv{l+1}}{\vv{l+p+r}}}{\suc{\uu{m+1}}{\uu{m+q+s}}}}\]\cr
&=~\func{\alpha}{\binary{\suc{\vv{1}}{\vv{l}}}{\suc{\uu{1}}{\uu{m}}}}\cr
&\hskip0.6\baselineskip\[\func{\beta}{\binary{\suc{\vv{l+1}}{\vv{l+p}}}{\suc{\uu{m+1}}{\uu{m+q}}}}
\func{\gamma}{\binary{\suc{\vv{l+p+1}}{\vv{l+p+r}}}{\suc{\uu{m+q+1}}{\uu{m+q+s}}}}\],\cr
&{}
\end{align}
and similarly,
\begin{align}
&\hskip0.6\baselineskip\func{\[\(\alpha\tensor{}\beta\)\tensor{}\gamma\]}{\binary{\suc{\vv{1}}{\vv{l+p+r}}}{\suc{\uu{1}}{\uu{m+q+s}}}}\cr
&=~\[\func{\alpha}{\binary{\suc{\vv{1}}{\vv{l}}}{\suc{\uu{1}}{\uu{m}}}}
\func{\beta}{\binary{\suc{\vv{l+1}}{\vv{l+p}}}{\suc{\uu{m+1}}{\uu{m+q}}}}\]\cr
&\hskip0.6\baselineskip\func{\gamma}{\binary{\suc{\vv{l+p+1}}{\vv{l+p+r}}}{\suc{\uu{m+q+1}}{\uu{m+q+s}}}}.\cr
&{}
\end{align}
Hence, considering the associativity of the product operation of the ring $\CRing{}$, it becomes evident that,
\begin{equation}
\alpha\tensor{}\(\beta\tensor{}\gamma\)=\(\alpha\tensor{}\beta\)\tensor{}\gamma.
\end{equation}
%\func{\[\(\alpha\tensor{}\beta\)\tensor{}\gamma\]}{\binary{\suc{\vv{1}}{\vv{l+p+r}}}{\suc{\uu{1}}{\uu{m+q+s}}}}
So the tensor operation is associative when restricted to simple tensors.\\
Generally there exists simple tensors of different types $\suc{\alpha_1}{\alpha_m}$,
simple tensors of different types $\suc{\beta_1}{\beta_n}$, and
simple tensors of different types $\suc{\gamma_1}{\gamma_p}$
in $\DTensors{\CModule{}}$ such that,
\begin{align}
\alpha&=\alpha_1+\cdots+\alpha_m\cr
\beta&=\beta_1+\cdots+\beta_n\cr
\gamma&=\gamma_1+\cdots+\gamma_p.
\end{align}
Hence according to the definition of tensor operation, considering the bilinearity of the tensor operation and the commutativity of
addition in the space of tensors, and the associativity of the tensor operation when restricted to simple tensors,
\begin{align}
\alpha\tensor{}\(\beta\tensor{}\gamma\)&=
\(\sum_{i=1}^{m}\alpha_i\)\tensor{}\[\(\sum_{j=1}^{n}\beta_j\)\tensor{}\(\sum_{k=1}^{p}\gamma_k\)\]\cr
&=\(\sum_{i=1}^{m}\alpha_i\)\tensor{}\(\sum_{j=1}^{n}\sum_{k=1}^{p}\beta_j\tensor{}\gamma_k\)\cr
&=\sum_{i=1}^{m}\sum_{j=1}^{n}\sum_{k=1}^{p}\alpha_i\tensor{}\(\beta_j\tensor{}\gamma_k\)\cr
&=\sum_{i=1}^{m}\sum_{j=1}^{n}\sum_{k=1}^{p}\(\alpha_i\tensor{}\beta_j\)\tensor{}\gamma_k\cr
&=\(\alpha\tensor{}\beta\)\tensor{}\gamma.
\end{align}
\endthm
%%%%%%%%%%%%%%%%%%%%%%%%%%%%%%%%%%%%%%%%%%%%%%%%%%%%%%%%%%%%%%%%%%%%%%%%%%%%%%%%%%%%%%%%%%%%%%%%%%%%%%%%%%%%%%%%%%%%%%%%%%%%%%%%%%%%%%%%%%%%%
\definition
$\opair{\DTensors{\CModule{}}}{\tensor{\CModule{}}}$ is called the $\quotl$tensor algebra of $\CModule{}$$\quotr$.
\endef
%%%%%%%%%%%%%%%%%%%%%%%%%%%%%%%%%%%%%%%%%%%%%%%%%%%%%%%%%%%%%%%%%%%%%%%%%%%%%%%%%%%%%%%%%%%%%%%%%%%%%%%%%%%%%%%%%%%%%%%%%%%%%%%%%%%%%%%%%%%%%
%%%%%%%%%%%%%%%%%%%%%%%%%%%%%%%%%%%%%%%%%%%%%%%%%%%%%%%%%%%%%%%%%%%%%%%%%%%%%%%%%%%%%%%%%%%%%%%%%%%%%%%%%%%%%%%%%%%%%%%%%%%%%%%%%%%%%%%%%%%%%
\subsection{Tensors on Finite-Dimensional Vector-Spaces}
\fixed
$\algfield{}$ is fixed as a field of characteristic zero, and
$\vecs{}$, $\vecs{1}$, and $\vecs{2}$ as $\algfield{}$-vector-spaces with finite dimensions $n$, $n_1$, and $n_2$ respectively.
$\vsbase{}$ is fixed as an ordered basis of $\vecs{}$, and $\Bdual{\vsbase{}}{}$ is taken to be its dual basis, which is an ordered basis
for $\Vdual{\vecs{}}$. $\Bddual{\vsbase{}}{}$ denotes the dual basis of $\Bdual{\vsbase{}}{}$, which can be
identified with $\vsbase{}$ via the natural linear isomorphism
$\map{\function{\iota}{\vecs{}}{\Vddual{\vecs{}}}}{\ww{}}{\[\map{\func{\iota}{\ww{}}}{\alpha}{\func{\alpha}{\ww{}}}\]}$.
\endfixed
%%%%%%%%%%%%%%%%%%%%%%%%%%%%%%%%%%%%%%%%%%%%%%%%%%%%%%%%%%%%%%%%%%%%%%%%%%%%%%%%%%%%%%%%%%%%%%%%%%%%%%%%%%%%%%%%%%%%%%%%%%%%%%%%%%%%%%%%%%%%%
\corollary
Let $r$ and $s$ be non-negative integers. $\MTensors{r}{s}{\vecs{}}$ is a finite-dimensional vector-space
with dimension $n^{r+s}$. This vector-space is called the $\quotl$vector-space of (all) $\opair{r}{s}$ tensors
on $\vecs{}$.
\endcor
%%%%%%%%%%%%%%%%%%%%%%%%%%%%%%%%%%%%%%%%%%%%%%%%%%%%%%%%%%%%%%%%%%%%%%%%%%%%%%%%%%%%%%%%%%%%%%%%%%%%%%%%%%%%%%%%%%%%%%%%%%%%%%%%%%%%%%%%%%%%%
\definition
For every non-negative integers $r$ and $s$ such that $r+s>0$, the function
$\function{\tensorvsbase{r}{s}{}}{\seta{\suc{1}{n^{r+s}}}}{\Tensors{r}{s}{\vecs{}}}$
is defined as,
\begin{equation}
\Foreach{j}{n^{r+s}}
\tensorvsbase{r}{s}{j}\eqdef\Bdual{\vsbase{}}{j_1}\tensor{}\cdots\tensor{}\Bdual{\vsbase{}}{j_r}
\tensor{}\vsbase{j_{r+1}}\tensor{}\cdots\tensor{}\vsbase{j_{r+s}},
\end{equation}
where the string $\mtuple{j_1}{j_{r+s}}$ is the unique representation of $j$ relative to the radix $n$, that is,
\begin{align}
&j=j_1\(n^{r+s-1}\)+j_2\(n^{r+s-2}\)+\cdots+j_{r+s},\cr
&\Foreach{i}{\seta{\suc{1}{r+s}}}0\leq j_i\leq n.
\end{align}
\endef
%%%%%%%%%%%%%%%%%%%%%%%%%%%%%%%%%%%%%%%%%%%%%%%%%%%%%%%%%%%%%%%%%%%%%%%%%%%%%%%%%%%%%%%%%%%%%%%%%%%%%%%%%%%%%%%%%%%%%%%%%%%%%%%%%%%%%%%%%%%%%
\theorem
Let $r$ and $s$ be non-nagative integers such that $r+s>0$.
$\tensorvsbase{r}{s}{}$ is an ordered basis for the $\algfield{}$-vector-space $\MTensors{r}{s}{\vecs{}}$,
and every $\alpha\in\Tensors{r}{s}{\vecs{}}$ has the unique representation with respect to this basis as,
\begin{equation}
\alpha=\sum_{j=1}^{n+s}\func{\alpha}{\binary{\suc{\vsbase{j_1}}{\vsbase{j_r}}}{\suc{\Bdual{\vsbase{}}{j_{r+1}}}{\Bdual{\vsbase{}}{j_{r+s}}}}}
\tensorvsbase{r}{s}{j},
\end{equation}
where the string $\mtuple{j_1}{j_{r+s}}$ is the representation of $j$ relative to the radix $n$.
\proof
For every $\alpha\in\Tensors{r}{s}{\vecs{}}$,
\begin{align}
&\begin{aligned}
\Foreach{\mtuple{\vv{1}}{\vv{r+s}}}{\Cprod{\multiprod{\vecs{}}{r}}{\multiprod{\Vdual{\vecs{}}}{s}}}
\end{aligned}\cr
&\begin{aligned}
&~\func{\alpha}{\suc{\vv{1}}{\vv{r+s}}}\cr
=&~\func{\alpha}{\binary
{{\suc{\sum_{j_1=1}^{n}\[\func{\Bdual{\vsbase{}}{j_1}}{\vv{1}}\]\vsbase{j_1}}{\sum_{j_r=1}^{n}\[\func{\Bdual{\vsbase{}}{j_r}}{\vv{r}}\]\vsbase{j_r}}}}
{{\suc{\sum_{j_{r+1}=1}^{n}\[\func{\vsbase{j_{r+1}}}{\vv{r+1}}\]\Bdual{\vsbase{}}{j_{r+1}}}{\sum_{j_{r+s}=1}^{n}\[\func{\vsbase{j_{r+s}}}{\vv{r+s}}\]\Bdual{\vsbase{}}{j_{r+s}}}}}}\cr
=&~\sum_{j_1=1}^{n}\cdots\sum_{j_{r+s}=1}^{n}\[\func{\alpha}{\binary{\suc{\vsbase{j_1}}{\vsbase{j_r}}}{\suc{\Bdual{\vsbase{}}{j_{r+1}}}{\Bdual{\vsbase{}}{j_{r+s}}}}}\]
\[\func{\Bdual{\vsbase{}}{j_1}}{\vv{1}}\cdots{\func{\Bdual{\vsbase{}}{j_r}}{\vv{r}}}\func{\vsbase{j_{r+1}}}{\vv{r+1}}\cdots\func{\vsbase{j_{r+s}}}{\vv{r+s}}\]\cr
=&~\sum_{j_1=1}^{n}\cdots\sum_{j_{r+s}=1}^{n}\[\func{\alpha}{\binary{\suc{\vsbase{j_1}}{\vsbase{j_r}}}{\suc{\Bdual{\vsbase{}}{j_{r+1}}}{\Bdual{\vsbase{}}{j_{r+s}}}}}\]
\[\func{\Bdual{\vsbase{}}{j_1}\tensor{}\cdots\tensor{}\Bdual{\vsbase{}}{j_r}
\tensor{}\vsbase{j_{r+1}}\tensor{}\cdots\tensor{}\vsbase{j_{r+s}}}{\suc{\vv{1}}{\vv{r+s}}}\]\cr
=&~\func{\[\sum_{j=1}^{n^{r+s}}\func{\alpha}{\binary{\suc{\vsbase{j_1}}{\vsbase{j_r}}}{\suc{\Bdual{\vsbase{}}{j_{r+1}}}{\Bdual{\vsbase{}}{j_{r+s}}}}}
\tensorvsbase{r}{s}{j}\]}{\suc{\vv{1}}{\vv{r+s}}},
\end{aligned}\cr
&{}
\end{align}
and hence,
\begin{equation}
\alpha=\sum_{j=1}^{n^{r+s}}\func{\alpha}{\binary{\suc{\vsbase{j_1}}{\vsbase{j_r}}}{\suc{\Bdual{\vsbase{}}{j_{r+1}}}{\Bdual{\vsbase{}}{j_{r+s}}}}}
\tensorvsbase{r}{s}{j}.
\end{equation}
So it becomes clear that $\defSet{\tensorvsbase{r}{s}{j}}{j=\suc{1}{n^{r+s}}}$ spans the vector-space
$\MTensors{r}{s}{\vecs{}}$.\\
Now let $\suc{c_1}{\c_{n^{r+s}}}$ be elements of $\algfield{}$, and suppose that,
\begin{equation}
\sum_{j=1}^{n^{r+s}}c_{j}\tensorvsbase{r}{s}{j}=\zerovec{},
\end{equation}
where $\zerovec{}$ denotes the neutral element of the vector-space $\MTensors{r}{s}{\vecs{}}$. Then clearly,
\begin{align}
\Foreach{k}{\seta{\suc{1}{n^{r+s}}}}
0&=\func{\[\sum_{j=1}^{n^{r+s}}c_{j}\tensorvsbase{r}{s}{j}\]}{\vsbase{k_1}\tensor{}\cdots\tensor{}\vsbase{k_r}\tensor{}
\Bdual{\vsbase{}}{k_{r+1}}\tensor{}\cdots\tensor{}\Bdual{\vsbase{}}{k_{r+s}}}\cr
&=c_k\tensorvsbase{r}{s}{k},
\end{align}
and therefore considering that each $\tensorvsbase{r}{s}{k}$ ($k$ ranging over all possible indices) is a non-zero element of
$\Tensors{r}{s}{\vecs{}}$,
\begin{equation}
\Foreach{k}{\seta{\suc{1}{n^{r+s}}}}
c_k=0.
\end{equation}
Therefore it is inferred that $\defSet{\tensorvsbase{r}{s}{j}}{j=\suc{1}{n^{r+s}}}$ is a linearly-independent subset of
$\MTensors{r}{s}{\vecs{}}$.
\endthm
%%%%%%%%%%%%%%%%%%%%%%%%%%%%%%%%%%%%%%%%%%%%%%%%%%%%%%%%%%%%%%%%%%%%%%%%%%%%%%%%%%%%%%%%%%%%%%%%%%%%%%%%%%%%%%%%%%%%%%%%%%%%%%%%%%%%%%%%%%%%%
\corollary
The dimension of the $\algfield{}$-vector-space $\MTensors{r}{s}{\vecs{}}$ is $n^{r+s}$.
\endcor
%%%%%%%%%%%%%%%%%%%%%%%%%%%%%%%%%%%%%%%%%%%%%%%%%%%%%%%%%%%%%%%%%%%%%%%%%%%%%%%%%%%%%%%%%%%%%%%%%%%%%%%%%%%%%%%%%%%%%%%%%%%%%%%%%%%%%%%%%%%%%
%%%%%%%%%%%%%%%%%%%%%%%%%%%%%%%%%%%%%%%%%%%%%%%%%%%%%%%%%%%%%%%%%%%%%%%%%%%%%%%%%%%%%%%%%%%%%%%%%%%%%%%%%%%%%%%%%%%%%%%%%%%%%%%%%%%%%%%%%%%%%
%%%%%%%%%%%%%%%%%%%%%%%%%%%%%%%%%%%%%%%%%%%%%%%%%%%%%%%%%%%%%%%%%%%%%%%%%%%%%%%%%%%%%%%%%%%%%%%%%%%%%%%%%%%%%%%%%%%%%%%%%%%%%%%%%%%%%%%%%%%%%
\section{Pullbacks of Linear Maps}
%%%%%%%%%%%%%%%%%%%%%%%%%%%%%%%%%%%%%%%%
\fixed
In this section the dimensional restriction of the fixed $\algfield{}$-vector-spaces
$\vecs{}$, $\vecs{1}$, and $\vecs{2}$ is relaxed and they can be infinite-dimensional.
\endfixed
%%%%%%%%%%%%%%%%%%%%%%%%%%%%%%%%%%%%%%%%%%%%%%%%%%%%%%%%%%%%%%%%%%%%%%%%%%%%%%%%%%%%%%%%%%%%%%%%%%%%%%%%%%%%%%%%%%%%%%%%%%%%%%%%%%%%%%%%%%%%%
\definition
Let $L$ be an element of $\Lin{\vecs{}}{\vecs{1}}$ (a linear map from $\vecs{}$ to $\vecs{1}$). The mapping
$\function{\dualpb{L}}{\Vdual{\vecs{1}}}{\Vdual{\vecs{}}}$ is defined as,
\begin{equation}
\Foreach{\alpha}{\Vdual{\vecs{1}}}
\Foreach{\vv{}}{\vecs{}}
\func{\[\func{\dualpb{L}}{\alpha}\]}{\vv{}}\eqdef\func{\alpha}{\func{L}{\vv{}}}.
\end{equation}
$\dualpb{L}$ is referred to as the $\quotl$dual-pullback of $L$$\quotr$.
\endef
%%%%%%%%%%%%%%%%%%%%%%%%%%%%%%%%%%%%%%%%%%%%%%%%%%%%%%%%%%%%%%%%%%%%%%%%%%%%%%%%%%%%%%%%%%%%%%%%%%%%%%%%%%%%%%%%%%%%%%%%%%%%%%%%%%%%%%%%%%%%%
\proposition
Let $L\in\Lin{\vecs{}}{\vecs{1}}$ and $T\in\Lin{\vecs{1}}{\vecs{2}}$.
\begin{itemize}
\item
\begin{equation}
\dualpb{\(\cmp{T}{L}\)}=\cmp{\dualpb{L}}{\dualpb{T}}.
\end{equation}
\item
If $L\in\Linisom{\vecs{}}{\vecs{1}}$, then $\dualpb{L}\in\Linisom{\Vdual{\vecs{1}}}{\Vdual{\vecs{}}}$, and
\begin{equation}
\finv{\(\dualpb{L}\)}=\dualpb{\(\finv{L}\)}.
\end{equation}
\end{itemize}
\proof
It is a direct consequence of the definition of the dual-pullback.
\endpro
%%%%%%%%%%%%%%%%%%%%%%%%%%%%%%%%%%%%%%%%%%%%%%%%%%%%%%%%%%%%%%%%%%%%%%%%%%%%%%%%%%%%%%%%%%%%%%%%%%%%%%%%%%%%%%%%%%%%%%%%%%%%%%%%%%%%%%%%%%%%%
\definition
Let $r$ and $s$ be non-negative integers. Suppose that $\vecs{}$ and $\vecs{1}$ are linearly isomorphic and
let $L$ be an element of $\Linisom{\vecs{}}{\vecs{1}}$ (a linear isomorphism from $\vecs{}$ to $\vecs{1}$).
The mapping $\function{\Vpullback{L}{r}{s}}{\Tensors{r}{s}{\vecs{1}}}{\Tensors{r}{s}{\vecs{}}}$ is defined as,
\begin{itemize}
\item
if $r+s>0$,
\begin{align}
&\Foreach{\beta}{\Tensors{r}{s}{\vecs{1}}}
\Foreach{\mtuple{\vv{1}}{\vv{r+s}}}{\Cprod{\multiprod{\vecs{}}{r}}{\multiprod{\Vdual{\vecs{}}}{s}}}\cr
&\func{\[\func{\Vpullback{L}{r}{s}}{\beta}\]}{\suc{\vv{1}}{\vv{r+s}}}\eqdef
\func{\beta}{\binary{\suc{\func{L}{\vv{1}}}{\func{L}{\vv{r}}}}
{\suc{\func{\dualpb{\(\finv{L}\)}}{\vv{r+1}}}{\func{\dualpb{\(\finv{L}\)}}{\vv{r+s}}}}},\cr
&{}
\end{align}
\item
if $r=s=0$,
\begin{equation}
\Vpullback{L}{r}{s}\eqdef\identity{\algfield{}}.
\end{equation}
\end{itemize}
$\Vpullback{L}{r}{s}$ is called the $\quotl$$\opair{r}{s}$-pullback of the linear isomorphism $L\in\Linisom{\vecs{}}{\vecs{1}}$$\quotr$.
\endef
%%%%%%%%%%%%%%%%%%%%%%%%%%%%%%%%%%%%%%%%%%%%%%%%%%%%%%%%%%%%%%%%%%%%%%%%%%%%%%%%%%%%%%%%%%%%%%%%%%%%%%%%%%%%%%%%%%%%%%%%%%%%%%%%%%%%%%%%%%%%%
\theorem\label{thmrspullbackofcompositionoflinearisomorphisms}
Let $r$ and $s$ be non-negative integers.
Suppose that $\vecs{}$, $\vecs{1}$, and $\vecs{2}$ are linearly isomorphic and
let $L$ be an element of $\Linisom{\vecs{}}{\vecs{1}}$ and
$T$ an element of $\Linisom{\vecs{1}}{\vecs{2}}$.
\begin{equation}
\Vpullback{\(\identity{\vecs{}}\)}{r}{s}=\identity{\Tensors{r}{s}{\vecs{}}},
\end{equation}
and
\begin{equation}
\Vpullback{\(\cmp{T}{L}\)}{r}{s}=\cmp{\Vpullback{L}{r}{s}}{\Vpullback{T}{r}{s}}.
\end{equation}
\proof
It is an immediate result of the definition of $\opair{r}{s}$-pullback of linear isomorphisms.
\endthm
%%%%%%%%%%%%%%%%%%%%%%%%%%%%%%%%%%%%%%%%%%%%%%%%%%%%%%%%%%%%%%%%%%%%%%%%%%%%%%%%%%%%%%%%%%%%%%%%%%%%%%%%%%%%%%%%%%%%%%%%%%%%%%%%%%%%%%%%%%%%%
\theorem\label{thminverseofpullbackofrstensor}
Let $r$ and $s$ be non-negative integers. Suppose that $\vecs{}$ and $\vecs{1}$ are linearly isomorphic and
let $L$ be an element of $\Linisom{\vecs{}}{\vecs{1}}$ (a linear isomorphism from $\vecs{}$ to $\vecs{1}$).
The mapping $\function{\Vpullback{L}{r}{s}}{\Tensors{r}{s}{\vecs{1}}}{\Tensors{r}{s}{\vecs{}}}$ is a linear isomorphism
from $\MTensors{r}{s}{\vecs{1}}$ to $\MTensors{r}{s}{\vecs{}}$. That is, $\Vpullback{L}{r}{s}$ is a bijection and,
\begin{equation}
\Foreach{\opair{\beta_1}{\beta_2}}{\Cprod{\Tensors{r}{s}{\vecs{1}}}{\Tensors{r}{s}{\vecs{1}}}}
\Foreach{c}{\algfield{}}
\func{\Vpullback{L}{r}{s}}{c\beta_1+\beta_2}=c\func{\Vpullback{L}{r}{s}}{\beta_1}+
\func{\Vpullback{L}{r}{s}}{\beta_2}.
\end{equation}
Furthermore,
\begin{equation}
\finv{\(\Vpullback{L}{r}{s}\)}=\Vpullback{\(\finv{L}\)}{r}{s}.
\end{equation}
\proof
$\Vpullback{L}{0}{0}=\identity{\algfield{}}$, and hence $\Vpullback{L}{0}{0}$
is a linear isomorphism from the $\algfield{}$-vector-space $\algfield{}$ to itself, that is
a linear isomorphism from $\Tensors{0}{0}{\vecs{1}}$ to $\Tensors{0}{0}{\vecs{}}$. So the assertion is true when $r=s=0$.\\
Now suppose that $r+s>0$.
The linearity of $\Vpullback{L}{r}{s}$ follows directly from the definition of $\opair{r}{s}$-pullback of a linear isomorphism
and the definition of the linear structure of the space of $\opair{r}{s}$ tensors on a vector-space.\\
The final assertion is a consequence of \refthm{thmrspullbackofcompositionoflinearisomorphisms}. Let us give
a direct proof, which actually mimics the idea for proving \refthm{thmrspullbackofcompositionoflinearisomorphisms}.
Let $L$ be a linear isomorphism from $\vecs{}$ to $\vecs{1}$. Clearly
$\function{\cmp{\Vpullback{L}{r}{s}}{\Vpullback{\(\finv{L}\)}{r}{s}}}{\Tensors{r}{s}{\vecs{}}}{\Tensors{r}{s}{\vecs{}}}$, and
\begin{align}
&\begin{aligned}
\Foreach{\alpha}{\Tensors{r}{s}{\vecs{}}}
\Foreach{\mtuple{\vv{1}}{\vv{r+s}}}{\Cprod{\multiprod{\vecs{}}{r}}{\multiprod{\Vdual{\vecs{}}}{s}}}
\end{aligned}\cr
&\begin{aligned}
&~\func{\[\func{\(\cmp{\Vpullback{L}{r}{s}}{\Vpullback{\(\finv{L}\)}{r}{s}}\)}{\alpha}\]}{\suc{\vv{1}}{\vv{r+s}}}\cr
=&~
%%%%
\func{\[\func{\Vpullback{\(\finv{L}\)}{r}{s}}{\alpha}\]}{\binary{\suc{\func{L}{\vv{1}}}{\func{L}{\vv{r}}}}
{\suc{\func{\dualpb{\(\finv{L}\)}}{\vv{r+1}}}{\func{\dualpb{\(\finv{L}\)}}{\vv{r+s}}}}}\cr
=&~\func{\alpha}{\binary{\suc{\func{\finv{L}}{\func{L}{\vv{1}}}}{\func{\finv{L}}{\func{L}{\vv{r}}}}}
{\suc{\func{\dualpb{L}}{\func{\dualpb{\(\finv{L}\)}}{\vv{r+1}}}}{\func{\dualpb{L}}{\func{\dualpb{\(\finv{L}\)}}{\vv{r+s}}}}}}\cr
=&~\func{\alpha}{\suc{\vv{1}}{\vv{r+s}}}.
\end{aligned}\cr
&{}
\end{align}
Therefore,
\begin{equation}
\cmp{\Vpullback{L}{r}{s}}{\Vpullback{\(\finv{L}\)}{r}{s}}=
\identity{\Tensors{r}{s}{\vecs{}}}.
\end{equation}
Similarly, it can be seen that,
\begin{equation}
\cmp{\Vpullback{\(\finv{L}\)}{r}{s}}{\Vpullback{L}{r}{s}}=
\identity{\Tensors{r}{s}{\vecs{1}}}.
\end{equation}
Thus, $\Vpullback{L}{r}{s}$ is a bijection and
$\finv{\(\Vpullback{L}{r}{s}\)}=\Vpullback{\(\finv{L}\)}{r}{s}$.
\endthm
%%%%%%%%%%%%%%%%%%%%%%%%%%%%%%%%%%%%%%%%%%%%%%%%%%%%%%%%%%%%%%%%%%%%%%%%%%%%%%%%%%%%%%%%%%%%%%%%%%%%%%%%%%%%%%%%%%%%%%%%%%%%%%%%%%%%%%%%%%%%%
\definition
Let $r$ be a non-negative integer.
Let $L$ be an element of $\Lin{\vecs{}}{\vecs{1}}$ (a linear map from $\vecs{}$ to $\vecs{1}$).
The mapping $\function{\Vpullbackcov{L}{r}}{\Tensors{r}{0}{\vecs{1}}}{\Tensors{r}{0}{\vecs{}}}$ is defined as,
\begin{itemize}
\item
if $r>0$,
\begin{align}
&\Foreach{\beta}{\Tensors{r}{0}{\vecs{1}}}
\Foreach{\mtuple{\vv{1}}{\vv{r}}}{\multiprod{\vecs{}}{r}}\cr
&\func{\[\func{\Vpullbackcov{L}{r}}{\beta}\]}{\suc{\vv{1}}{\vv{r}}}\eqdef
\func{\beta}{\suc{\func{L}{\vv{1}}}{\func{L}{\vv{r}}}},
\end{align}
\item
if $r=s=0$,
\begin{equation}
\Vpullbackcov{L}{r}\eqdef\identity{\algfield{}}.
\end{equation}
\end{itemize}
$\Vpullbackcov{L}{r}$ is called the $\quotl$$r$-covariant--pullback of the linear map $L\in\Lin{\vecs{}}{\vecs{1}}$$\quotr$.
\endef
%%%%%%%%%%%%%%%%%%%%%%%%%%%%%%%%%%%%%%%%%%%%%%%%%%%%%%%%%%%%%%%%%%%%%%%%%%%%%%%%%%%%%%%%%%%%%%%%%%%%%%%%%%%%%%%%%%%%%%%%%%%%%%%%%%%%%%%%%%%%%
\theorem\label{thmrpullbackofcompositionoflinearmaps}
Let $r$ be a non-negative integer.
Let $L$ be an element of $\Lin{\vecs{}}{\vecs{1}}$ and
$T$ an element of $\Lin{\vecs{1}}{\vecs{2}}$.
\begin{equation}
\Vpullbackcov{\(\identity{\vecs{}}\)}{r}=\identity{\Tensors{r}{0}{\vecs{}}},
\end{equation}
and
\begin{equation}
\Vpullbackcov{\(\cmp{T}{L}\)}{r}=\cmp{\Vpullbackcov{L}{r}}{\Vpullbackcov{T}{r}}.
\end{equation}
\proof
It is an immediate result of the definition of $r$-pullback of linear maps.
\endthm
%%%%%%%%%%%%%%%%%%%%%%%%%%%%%%%%%%%%%%%%%%%%%%%%%%%%%%%%%%%%%%%%%%%%%%%%%%%%%%%%%%%%%%%%%%%%%%%%%%%%%%%%%%%%%%%%%%%%%%%%%%%%%%%%%%%%%%%%%%%%%
%%%%%%%%%%%%%%%%%%%%%%%%%%%%%%%%%%%%%%%%%%%%%%%%%%%%%%%%%%%%%%%%%%%%%%%%%%%%%%%%%%%%%%%%%%%%%%%%%%%%%%%%%%%%%%%%%%%%%%%%%%%%%%%%%%%%%%%%%%%%%
\theorem\label{thmpullbackofcovtensorislinear}
Let $r$ be a non-negative integer.
Let $L$ be an element of $\Lin{\vecs{}}{\vecs{1}}$ (a linear map from $\vecs{}$ to $\vecs{1}$).
The mapping $\function{\Vpullbackcov{L}{r}}{\Tensors{r}{0}{\vecs{1}}}{\Tensors{r}{0}{\vecs{}}}$ is a linear map
from $\MTensors{r}{0}{\vecs{1}}$ to $\MTensors{r}{0}{\vecs{}}$. That is,
\begin{equation}
\Foreach{\opair{\beta_1}{\beta_2}}{\Cprod{\Tensors{r}{0}{\vecs{1}}}{\Tensors{r}{0}{\vecs{1}}}}
\Foreach{c}{\algfield{}}
\func{\Vpullbackcov{L}{r}}{c\beta_1+\beta_2}=c\func{\Vpullbackcov{L}{r}}{\beta_1}+
\func{\Vpullbackcov{L}{r}}{\beta_2}.
\end{equation}
\proof
The linearity of $\Vpullbackcov{L}{r}$ is an immediate consequence of the definition of $r$-covariant--pullback of a linear map
and the definition of the linear structure of the space of $r$-covariant tensors on a vector-space.
\endthm
%%%%%%%%%%%%%%%%%%%%%%%%%%%%%%%%%%%%%%%%%%%%%%%%%%%%%%%%%%%%%%%%%%%%%%%%%%%%%%%%%%%%%%%%%%%%%%%%%%%%%%%%%%%%%%%%%%%%%%%%%%%%%%%%%%%%%%%%%%%%%
\definition
Suppose that $\vecs{}$ and $\vecs{1}$ are linearly isomorphic and
Let $L$ be an element of $\Linisom{\vecs{}}{\vecs{1}}$ (a linear isomorphism from $\vecs{}$ to $\vecs{1}$).
The mapping $\function{\VPullback{L}}{\DTensors{\vecs{1}}}{\DTensors{\vecs{}}}$ is defined term-wise as,
\begin{equation}
\Foreach{\beta}{\DTensors{\vecs{1}}}
{\[\func{\VPullback{L}}{\beta}\]}_{\opair{r}{s}}:=\func{\Vpullback{L}{r}{s}}{\beta_{\opair{r}{s}}}.
\end{equation}
\endef
%%%%%%%%%%%%%%%%%%%%%%%%%%%%%%%%%%%%%%%%%%%%%%%%%%%%%%%%%%%%%%%%%%%%%%%%%%%%%%%%%%%%%%%%%%%%%%%%%%%%%%%%%%%%%%%%%%%%%%%%%%%%%%%%%%%%%%%%%%%%%
\corollary
Suppose that $\vecs{}$ and $\vecs{1}$ are linearly isomorphic and
let $L$ be an element of $\Linisom{\vecs{}}{\vecs{1}}$ (a linear isomorphism from $\vecs{}$ to $\vecs{1}$).
$\VPullback{L}$ is a linear isomorphism from $\DTensors{\vecs{1}}$ to $\DTensors{\vecs{}}$.
\endcor
%%%%%%%%%%%%%%%%%%%%%%%%%%%%%%%%%%%%%%%%%%%%%%%%%%%%%%%%%%%%%%%%%%%%%%%%%%%%%%%%%%%%%%%%%%%%%%%%%%%%%%%%%%%%%%%%%%%%%%%%%%%%%%%%%%%%%%%%%%%%%
%%%%%%%%%%%%%%%%%%%%%%%%%%%%%%%%%%%%%%%%%%%%%%%%%%%%%%%%%%%%%%%%%%%%%%%%%%%%%%%%%%%%%%%%%%%%%%%%%%%%%%%%%%%%%%%%%%%%%%%%%%%%%%%%%%%%%%%%%%%%%
%%%%%%%%%%%%%%%%%%%%%%%%%%%%%%%%%%%%%%%%%%%%%%%%%%%%%%%%%%%%%%%%%%%%%%%%%%%%%%%%%%%%%%%%%%%%%%%%%%%%%%%%%%%%%%%%%%%%%%%%%%%%%%%%%%%%%%%%%%%%%
\section{Tensor Fields on Banach-Spaces}
%%%%%%%%%%%%%%%%%%%%%%%%%%%%%%%%%%%%%%%%
\fixed
In this section, each of the $\algfield{}$-vector-spaces $\vecs{}$, $\vecs{1}$, and $\vecs{2}$ is assumed
to be equipped with a norm that turns it into a Banach-space. Any of these Banach-spaces
will be denoted by the same notation as that of its underlying vector-space. $\OO{}$, $\OO{1}$, and $\OO{2}$
are fixed as non-empty open sets of $\vecs{}$, $\vecs{1}$, and $\vecs{2}$, respectively. Notice that,
the considered vector-spaces are decided to be finite-dimensional in this section.
\endfixed
\remark
\textit{
In this section, the dimensional limitation imposed on vector-spaces is a vital constraint,
because the norm defined here for the whole space of $\opair{r}{s}$ tensors (for some non-negative integers $r$ and $s$)
may not be well-defined in general when the dimension of the considered vector-space is infinite. In the most general case,
this norm is well-defined for a sub-space of all tensors which is the space of continuous tensors. However,
when the dimension of the underlying vector-space is finite, this sub-space coincides with the whole space of tensors
and thus the definition here makes sense. \cite[Chapter 1]{Cartan} contains a thorough study of this problem.}
%%%%%%%%%%%%%%%%%%%%%%%%%%%%%%%%%%%%%%%%%%%%%%%%%%%%%%%%%%%%%%%%%%%%%%%%%%%%%%%%%%%%%%%%%%%%%%%%%%%%%%%%%%%%%%%%%%%%%%%%%%%%%%%%%%%%%%%%%%%%%
\definition
The function $\function{\norm{}{\Vdual{\vecs{}}}}{\Vdual{\vecs{}}}{\R}$ is defined as,
\begin{align}
\Foreach{\alpha}{\Vdual{\vecs{}}}
\norm{\alpha}{\Vdual{\vecs{}}}\eqdef
\sup{\defSet{\abs{\func{\alpha}{\vv{}}}}{\norm{\vv{}}{}\leq 1}}.
\end{align}
\endef
%%%%%%%%%%%%%%%%%%%%%%%%%%%%%%%%%%%%%%%%%%%%%%%%%%%%%%%%%%%%%%%%%%%%%%%%%%%%%%%%%%%%%%%%%%%%%%%%%%%%%%%%%%%%%%%%%%%%%%%%%%%%%%%%%%%%%%%%%%%%%
\proposition
$\norm{}{\Vdual{\vecs{}}}$ is a norm on $\Vdual{\vecs{}}$. This norm is referred to as the
$\quotl$dual norm of $\vecs{}$$\quotr$. The Banach-space obtained by endowing $\Vdual{\vecs{}}$ with
this norm will also be denoted by $\Vdual{\vecs{}}$.
\endpro
%%%%%%%%%%%%%%%%%%%%%%%%%%%%%%%%%%%%%%%%%%%%%%%%%%%%%%%%%%%%%%%%%%%%%%%%%%%%%%%%%%%%%%%%%%%%%%%%%%%%%%%%%%%%%%%%%%%%%%%%%%%%%%%%%%%%%%%%%%%%%
\definition
The function $\function{\norm{}{\MTensors{r}{s}{\vecs{}}}}{\Tensors{r}{s}{\vecs{}}}{\R}$ is defined as,
\begin{align}
&\Foreach{\alpha}{\Tensors{r}{s}{\vecs{}}}\cr
&\begin{aligned}
&~\norm{\alpha}{\MTensors{r}{s}{\vecs{}}}\cr
\eqdef&~\sup{\defSet{\abs{\func{\alpha}{\suc{\vv{1}}{\vv{r+s}}}}}{\bigg(\suc{\norm{\vv{1}}{}\leq 1}{\norm{\vv{r}}{}\leq 1},~
\suc{\norm{\vv{r+1}}{\Vdual{\vecs{}}}\leq 1}{\norm{\vv{r+s}}{\Vdual{\vecs{}}}\leq 1}\bigg)}}.
\end{aligned}\cr
&{}
\end{align}
\endef
%%%%%%%%%%%%%%%%%%%%%%%%%%%%%%%%%%%%%%%%%%%%%%%%%%%%%%%%%%%%%%%%%%%%%%%%%%%%%%%%%%%%%%%%%%%%%%%%%%%%%%%%%%%%%%%%%%%%%%%%%%%%%%%%%%%%%%%%%%%%%
\proposition
$\norm{}{\MTensors{r}{s}{\vecs{}}}$ is a norm on $\MTensors{r}{s}{\vecs{}}$. This norm is
referred to as the $\quotl$$\opair{r}{s}$ tensor norm of $\vecs{}$$\quotr$. The Banach-space obtained by endowing
$\MTensors{r}{s}{\vecs{}}$ with this norm will also be denoted by $\MTensors{r}{s}{\vecs{}}$.
\endpro
%%%%%%%%%%%%%%%%%%%%%%%%%%%%%%%%%%%%%%%%%%%%%%%%%%%%%%%%%%%%%%%%%%%%%%%%%%%%%%%%%%%%%%%%%%%%%%%%%%%%%%%%%%%%%%%%%%%%%%%%%%%%%%%%%%%%%%%%%%%%%
\remark
\textit{All differentiability properties and the values of derivatives of any mapping between Banach-spaces
coincide for any two choices of equivalent norms on both source and target Banach-spaces.
Thus, knowing that all norms of any finite-dimensional Banach-space are equivalent norms that turn it to a Banach-space,
all differentiability properties and the values of derivatives of any mapping between finite-dimensional
Banach-spaces must be independent of the choice of a norm for both the source and target Banach-spaces.
So, specifying a particular norm for a finite-dimensional Banach-space is just a choice of a representative
among all norms, when differentiability is the main objective of the study. Moreover, since any finite-dimensional
vector-space possesses a norm definitely, the study of differentiability of mappings between
finite-dimensional vector-spaces automatically makes sense without choosing any specific norm on them.
This is in particular true for the spaces of tensors on any finite-dimensional vector-space, but invoking a particular norm
is certainly necessary for calculations and investigations about differentiability. The norm chosen here for the spaces of
tensors on a finite-dimensional vector-space is the standard norm that is defined in the most general setting.
The proofs of these assertions are all available in \cite{Cartan}.}
%%%%%%%%%%%%%%%%%%%%%%%%%%%%%%%%%%%%%%%%%%%%%%%%%%%%%%%%%%%%%%%%%%%%%%%%%%%%%%%%%%%%%%%%%%%%%%%%%%%%%%%%%%%%%%%%%%%%%%%%%%%%%%%%%%%%%%%%%%%%%
\definition\label{defspaceofrstensorfieldsonopensetsofvectorspaces}
Let $r$ and $s$ be non-negative integers.
Any smooth (infinitely differentiable) map from $\OO{}\subseteq\vecs{}$ to $\MTensors{r}{s}{\vecs{}}$ is called a
$\quotl$(type) $\opair{r}{s}$ tensor field on the open subset $\OO{}$ of $\vecs{}$$\quotr$.
The set of all $\opair{r}{s}$ tensor fields on $\OO{}\subseteq\vecs{}$ will be denoted by
$\TFB{r}{s}{\vecs{}}{\OO{}}$. That is,
\begin{equation}
\TFB{r}{s}{\vecs{}}{\OO{}}:=
\banachmapdifclass{\infty}{\vecs{}}{\MTensors{r}{s}{\vecs{}}}{\OO{}}{\Tensors{r}{s}{\vecs{}}}
\end{equation}
Each element of $\TFB{r}{0}{\vecs{}}{\OO{}}$ is called alternatively a
$\quotl$covariant tensor field of rank $r$ on $\OO{}\subseteq{\vecs{}}$$\quotr$.
\endef
%%%%%%%%%%%%%%%%%%%%%%%%%%%%%%%%%%%%%%%%%%%%%%%%%%%%%%%%%%%%%%%%%%%%%%%%%%%%%%%%%%%%%%%%%%%%%%%%%%%%%%%%%%%%%%%%%%%%%%%%%%%%%%%%%%%%%%%%%%%%%
\proposition
The set $\TFB{r}{s}{\vecs{}}{\OO{}}$ endowed with the addition and scalar multiplication
operations defined pointwise as
\begin{align}
&\Foreach{\opair{\atf{1}}{\atf{2}}}{\Cprod{\TFB{r}{s}{\vecs{}}{\OO{}}}{\TFB{r}{s}{\vecs{}}{\OO{}}}}
\Foreach{\point}{\OO{}}
\func{\(\atf{1}+\atf{2}\)}{\point}\eqdef\func{\atf{1}}{\point}+\func{\atf{2}}{\point}\cr
&\Foreach{\opair{c}{\atf{}}}{\Cprod{\algfield{}}{\TFB{r}{s}{\vecs{}}{\OO{}}}}
\Foreach{\point}{\OO{}}
\func{\(c\atf{}\)}{\point}\eqdef c\[\func{\atf{}}{\point}\],
\end{align}
is a $\algfield{}$-vector-space. This vector-space will be denoted by $\VTFB{r}{s}{\vecs{}}{\OO{}}$.
\endpro
%%%%%%%%%%%%%%%%%%%%%%%%%%%%%%%%%%%%%%%%%%%%%%%%%%%%%%%%%%%%%%%%%%%%%%%%%%%%%%%%%%%%%%%%%%%%%%%%%%%%%%%%%%%%%%%%%%%%%%%%%%%%%%%%%%%%%%%%%%%%%
%%%%%%%%%%%%%%%%%%%%%%%%%%%%%%%%%%%%%%%%%%%%
\definition
The direct sum of all $\VTFB{r}{s}{\vecs{}}{\OO{}}$-s when $r$ and $s$ range over all non-negative integers simultaneously,
will simply be denoted by $\DVTFB{\vecs{}}{\OO{}}$. That is,
\begin{equation}
\DVTFB{\vecs{}}{\OO{}}:=\Dsum{\defSet{\VTFB{r}{s}{\vecs{}}{\OO{}}}{r,s=0,1,\ldots}}.
\end{equation}
$\DVTFB{\vbundle{}}{\OO{}}$ is referred to as the $\quotl$space of tensor fields on $\OO{}\subseteq\vecs{}$$\quotr$.
For every pair $r$ and $s$ of non-negative integers, there is a natural injection of $\TFB{r}{s}{\vecs{}}{\OO{}}$ into
$\DVTFB{\vecs{}}{\OO{}}$ via which every $\opair{r}{s}$ tensor field on $\OO{}\subseteq\vecs{}$ can be
identified with a unique element of $\DVTFB{\vecs{}}{\OO{}}$.
\endef
%%%%%%%%%%%%%%%%%%%%%%%%%%%%%%%%%%%%%%%%%%%%%%%%%%%%%%%%%%%%%%%%%%%%%%%%%%%%%%%%%%%%%%%%%%%%%%%%%%%%%%%%%%%%%%%%%%%%%%%%%%%%%%%%
\definition
The binary operation $\tftensor{\OO{}}$ on $\DVTFB{\vecs{}}{\OO{}}$ is defined to be the unique bilinear map
$\function{\tftensor{\vecs{}}{\OO{}}}{\Cprod{\DVTFB{\vecs{}}{\OO{}}}{\DVTFB{\vecs{}}{\OO{}}}}{\DVTFB{\vecs{}}{\OO{}}}$
such that for every non-negative integers $r$, $s$, $p$, $q$,
%if $r+s>0$ and $p+q>0$,
%%$\func{\image{\tensor{\CModule{}}}}{}\subseteq
\begin{align}
&\Foreach{\opair{\atf{}}{\atff{}}}{\Cprod{\TFB{r}{s}{\vecs{}}{\OO{}}}{\TFB{p}{q}{\vecs{}}{\OO{}}}}
\Foreach{\point}{\OO{}}
\func{\(\atf{}\tftensor{\OO{}}\atff{}\)}{\point}\eqdef
\func{\atf{}}{\point}\tensor{}\func{\atff{}}{\point},
%%%%%%%%%%%%%%%%
\end{align}
where the tensor product $\func{\atf{}}{\point}\tensor{}\func{\atff{}}{\point}$ takes place
in the space of tensors on $\vecs{}$, that is $\DTensors{\vecs{}}$.
%%%%%%%%%
$\tftensor{\OO{}}$ is called the $\quotl$tensor field operation of the open set $\OO{}\subseteq\vecs{}$$\quotr$.
For any non-negative integers $r$ and $s$, each element of $\TFB{r}{s}{\vecs{}}{\OO{}}$ is also called a
$\quotl$(type) $\opair{r}{s}$ simple tensor field on $\OO{}$$\quotr$. Also for every $\atf{}$ and $\atff{}$
in $\DVTFB{\vecs{}}{\OO{}}$, $\atf{}\tftensor{\OO{}}\atff{}$ is called the $\quotl$tensor product of the tensor fields $\atf{}$
and $\atff{}$$\quotr$.
When there is no ambiguity about the underlying open set $\OO{}$ of $\vecs{}$, $\tftensor{\OO{}}$ can
simply be denoted by $\tftensor{}$.\\
\caution
It is easy to verify that for any $\atf{}\in\TFB{r}{s}{\vecs{}}{\OO{}}\subseteq\DVTFB{\vecs{}}{\OO{}}$ and any
$\atff{}\in\TFB{p}{q}{\vecs{}}{\OO{}}\subseteq\DVTFB{\vecs{}}{\OO{}}$, we have
$\atf{}\tftensor{\OO{}}\atff{}\in\TFB{r+p}{q+s}{\vecs{}}{\OO{}}\subseteq\DVTFB{\vecs{}}{\OO{}}$ and $\tftensor{\OO{}}$
acts bilinearly on\\
$\Cprod{\TFB{r}{s}{\vecs{}}{\OO{}}}{\TFB{p}{q}{\vecs{}}{\OO{}}}$. We presumed the triviality of
these facts prior to the definition.
\endef
%%%%%%%%%%%%%%%%%%%%%%%%%%%%%%%%%%%%%%%%%%%%%%%%%%%%%%%%%%%%%%%%%%%%%%%%%%%%%%%%%%%%%%%%%%%%%%%%%%%%%%%%%%%%%%%%%%%%%%%%%%%%%%%%
\proposition
$\tftensor{\OO{}}$ is an associative and bilinear binary operation on $\DVTFB{\vecs{}}{\OO{}}$.
So the pair $\opair{\DVTFB{\vecs{}}{\OO{}}}{\tftensor{\OO{}}}$ is an algebra.
\proof
The bilinearity of $\tftensor{\OO{}}$ lies in the definition of it.\\
Now let $\atf{1}\in\TFB{r}{s}{\vecs{}}{\OO{}}$, $\atf{2}\in\TFB{p}{1}{\vecs{}}{\OO{}}$, and
$\atf{3}\in\TFB{l}{m}{\vecs{}}{\OO{}}$.
According to \refthm{thmtensoroperationofmoduleisassociative},
\begin{align}
\Foreach{\point}{\OO{}}
\func{\[\atf{1}\tftensor{}\(\atf{2}\tftensor{}\atf{3}\)\]}{\point}&=
\func{\atf{1}}{\point}\tensor{}\(\func{\atf{2}}{\point}\tensor{}\func{\atf{3}}{\point}\)\cr
&=\(\func{\atf{1}}{\point}\tensor{}\func{\atf{2}}{\point}\)\tensor{}\func{\atf{3}}{\point}\cr
&=\func{\[\(\atf{1}\tftensor{}\atf{2}\)\tftensor{}\atf{3}\]}{\point},
\end{align}
and thus,
\begin{equation}
\atf{1}\tftensor{}\(\atf{2}\tftensor{}\atf{3}\)=
\(\atf{1}\tftensor{}\atf{2}\)\tftensor{}\atf{3}.
\end{equation}
So $\tftensor{\OO{}}$ is associative when restricted to simple tensor fields. The verification of associativity
for general tensor fields in $\DVTFB{\vecs{}}{\OO{}}$ can be obtained in a way completely similar to what is used in
the proof of \refthm{thmtensoroperationofmoduleisassociative}.
\endpro
%%%%%%%%%%%%%%%%%%%%%%%%%%%%%%%%%%%%%%%%%%%%%%%%%%%%%%%%%%%%%%%%%%%%%%%%%%%%%%%%%%%%%%%%%%%%%%%%%%%%%%%%%%%%%%%%%%%%%%%%%%%%%%%%
\definition
$\opair{\DVTFB{\vecs{}}{\OO{}}}{\tftensor{\OO{}}}$ is called the
$\quotl$tensor algebra of tensor fields on the open set $\OO{}$ of $\vecs{}$$\quotr$.
\endef
%%%%%%%%%%%%%%%%%%%%%%%%%%%%%%%%%%%%%%%%%%%%%%%%%%%%%%%%%%%%%%%%%%%%%%%%%%%%%%%%%%%%%%%%%%%%%%%%%%%%%%%%%%%%%%%%%%%%%%%%%%%%%%%%
\section{Pullbacks of Smooth Maps}
%%%%%%%%%%%%%%%%%%%%%%%%%%%%%%%%%%%%%
\proposition
Let $r$ and $s$ be non-negative integers, and let $f$ be a smooth diffeomorphism from $\OO{1}\subseteq\vecs{1}$
to $\OO{2}\subseteq\vecs{2}$ (that is, an element of
$\BanachDiff{\infty}{\vecs{1}}{\vecs{2}}{\OO{1}}{\OO{2}}$).
The image of the assignment $\TFB{r}{s}{\vecs{2}}{\OO{2}}\ni\atf{}\mapsto\bar{\atf{}}\in
\Func{\OO{1}}{\Tensors{r}{s}{\vecs{1}}}$ such that
\begin{align}
\Foreach{\point}{\OO{1}}
\func{\bar{\atf{}}}{\point}\eqdef
\func{\Vpullback{\[\func{\banachder{f}{\vecs{1}}{\vecs{2}}}{\point}\]}{r}{s}}{\func{\atf{}}{\func{f}{\point}}},
\end{align}
is included in $\TFB{r}{s}{\vecs{1}}{\OO{1}}$. In other words, for every $\atf{}\in\VTFB{r}{s}{\vecs{2}}{\OO{2}}$,
$\bar{\atf{}}$ lies in $\VTFB{r}{s}{\vecs{1}}{\OO{1}}$.
$\banachder{f}{\vecs{1}}{\vecs{2}}$ denotes the derived map of $f$ (with respect to the background
Banach-spaces $\vecs{1}$ and $\vecs{2}$), and it can be denoted simply by $\banachderr{f}$ when the
ambient Banach-spaces are clearly recognized.\\
\caution
Note that this assignment is well-defined since the value of the derived map of a diffeomorphism
at every point is a linear isomorphism.
\proof
It is left to the reader as an exercise.
\endpro
%%%%%%%%%%%%%%%%%%%%%%%%%%%%%%%%%%%%%%%%%%%%%%%%%%%%%%%%%%%%%%%%%%%%%%%%%%%%%%%%%%%%%%%%%%%%%%%%%%%%%%%%%%%%%%%%%%%%%%%%%%%%%%%%
\definition\label{defrspullbackofdiffeomorphism}
Let $r$ and $s$ be non-negative integers, and let $f$ be a smooth diffeomorphism from $\OO{1}\subseteq\vecs{1}$
to $\OO{2}\subseteq\vecs{2}$ (that is, an element of
$\BanachDiff{\infty}{\vecs{1}}{\vecs{2}}{\OO{1}}{\OO{2}}$). The mapping
$\function{\TFpullback{f}{r}{s}}{\TFB{r}{s}{\vecs{2}}{\OO{2}}}{\TFB{r}{s}{\vecs{1}}{\OO{1}}}$ is defined as
\begin{align}
&\Foreach{\atf{}}{\TFB{r}{s}{\vecs{2}}{\OO{2}}}
\Foreach{\point}{\OO{1}}\cr
&\func{\[\func{\TFpullback{f}{r}{s}}{\atf{}}\]}{\point}\eqdef
\func{\Vpullback{\[\func{\banachder{f}{\vecs{1}}{\vecs{2}}}{\point}\]}{r}{s}}{\func{\atf{}}{\func{f}{\point}}},
\end{align}
$\TFpullback{f}{r}{s}$ is referred to as the $\quotl$$\opair{r}{s}$-pullback of the diffeomorphism
$f\in\BanachDiff{\infty}{\vecs{1}}{\vecs{2}}{\OO{1}}{\OO{2}}$$\quotr$.
\endef
%%%%%%%%%%%%%%%%%%%%%%%%%%%%%%%%%%%%%%%%%%%%%%%%%%%%%%%%%%%%%%%%%%%%%%%%%%%%%%%%%%%%%%%%%%%%%%%%%%%%%%%%%%%%%%%%%%%%%%%%%%%%%%%%
\proposition
Let $r$ and $s$ be non-negative integers, and let $f$ be a smooth diffeomorphism from $\OO{1}\subseteq\vecs{1}$
to $\OO{2}\subseteq\vecs{2}$ (that is, an element of
$\BanachDiff{\infty}{\vecs{1}}{\vecs{2}}{\OO{1}}{\OO{2}}$).
\begin{align}
&\Foreach{\atf{}}{\TFB{r}{s}{\vecs{2}}{\OO{2}}}
\Foreach{\point}{\OO{1}}
\Foreach{\mtuple{\vv{1}}{\vv{r+s}}}{\Cprod{\multiprod{\vecs{}}{r}}{\multiprod{\Vdual{\vecs{}}}{s}}}\cr
&\begin{aligned}
&\hskip0.6\baselineskip\func{\(\func{\[\func{\TFpullback{f}{r}{s}}{\atf{}}\]}{\point}\)}{\suc{\vv{1}}{\vv{r+s}}}\cr
&=
\func{\[\func{\atf{}}{\func{f}{\point}}\]}{\binary{\suc{\func{\[\func{\banachderr{f}}{\point}\]}{\vv{1}}}{\func{\[\func{\banachderr{f}}{\point}\]}{\vv{r}}}}
{\suc{\func{\dualpb{\[\finv{\(\func{\banachderr{f}}{\point}\)}\]}}{\vv{r+1}}}{\func{\dualpb{\[\finv{\(\func{\banachderr{f}}{\point}\)}\]}}{\vv{r+s}}}}}.
\end{aligned}\cr
&{}
\end{align}
\proof
It is trivial.
\endpro
%%%%%%%%%%%%%%%%%%%%%%%%%%%%%%%%%%%%%%%%%%%%%%%%%%%%%%%%%%%%%%%%%%%%%%%%%%%%%%%%%%%%%%%%%%%%%%%%%%%%%%%%%%%%%%%%%%%%%%%%%%%%%%%%
\theorem\label{thmrspullbackofcompositionofdiffeomorphisms}
Let $r$ and $s$ be non-negative integers.
Let $f$ be an element of $\BanachDiff{\infty}{\vecs{}}{\vecs{1}}{\OO{}}{\OO{1}}$ and
$g$ an element of $\BanachDiff{\infty}{\vecs{1}}{\vecs{2}}{\OO{1}}{\OO{2}}$.
\begin{equation}
\TFpullback{\(\identity{\OO{}}\)}{r}{s}=\identity{\TFB{r}{s}{\vecs{}}{\OO{}}},
\end{equation}
and
\begin{equation}
\TFpullback{\(\cmp{g}{f}\)}{r}{s}=\cmp{\TFpullback{f}{r}{s}}{\TFpullback{g}{r}{s}}.
\end{equation}
\proof
The first assertion is trivial.\\
According to \refthm{thmrspullbackofcompositionoflinearisomorphisms},
\refdef{defrspullbackofdiffeomorphism}, and the chain rule of differentiation,
\begin{align}
&\Foreach{\atf{}}{\TFB{r}{s}{\vecs{2}}{\OO{2}}}
\Foreach{\point}{\OO{}}\cr
&\begin{aligned}
\func{\[\func{\TFpullback{\(\cmp{g}{f}\)}{r}{s}}{\atf{}}\]}{\point}&=
\func{\Vpullback{\[\func{\banachderr{\(\cmp{g}{f}\)}}{\point}\]}{r}{s}}{\func{\atf{}}{\func{g}{\func{f}{\point}}}}\cr
&=\func{\Vpullback{\bigg(\cmp{\[\func{\banachderr{g}}{\func{f}{\point}}\]}{\[\func{\banachderr{f}}{\point}\]}\bigg)}{r}{s}}
{\func{\atf{}}{\func{g}{\func{f}{\point}}}}\cr
&=\func{\bigg(\cmp{\Vpullback{\[\func{\banachderr{f}}{\point}\]}{r}{s}}{\Vpullback{\[\func{\banachderr{g}}{\func{f}{\point}}\]}{r}{s}}\bigg)}
{\func{\atf{}}{\func{g}{\func{f}{\point}}}}\cr
&=\func{\Vpullback{\[\func{\banachderr{f}}{\point}\]}{r}{s}}{\func{\[\func{\TFpullback{g}{r}{s}}{\atf{}}\]}{\func{f}{\point}}}\cr
&=\func{\[\func{\TFpullback{f}{r}{s}}{\func{\TFpullback{g}{r}{s}}{\atf{}}}\]}{\point}\cr
&=\func{\[\func{\(\cmp{\TFpullback{f}{r}{s}}{\TFpullback{g}{r}{s}}\)}{\atf{}}\]}{\point},
\end{aligned}\cr
&{}
\end{align}
and thus,
\begin{equation}
\TFpullback{\(\cmp{g}{f}\)}{r}{s}=\cmp{\TFpullback{f}{r}{s}}{\TFpullback{g}{r}{s}}.
\end{equation}
\endthm
%%%%%%%%%%%%%%%%%%%%%%%%%%%%%%%%%%%%%%%%%%%%%%%%%%%%%%%%%%%%%%%%%%%%%%%%%%%%%%%%%%%%%%%%%%%%%%%%%%%%%%%%%%%%%%%%%%%%%%%%%%%%%%%%
\theorem\label{thmpullbackofrstensorfieldisislinear}
Let $r$ and $s$ be non-negative integers, and
let $f$ be a smooth diffeomorphism from $\OO{1}\subseteq\vecs{1}$ to $\OO{2}\subseteq\vecs{2}$
(that is, an element of $\BanachDiff{\infty}{\vecs{1}}{\vecs{2}}{\OO{1}}{\OO{2}}$).
%%%%%%%%
The mapping $\TFpullback{f}{r}{s}$ is a linear isomorphism from
$\VTFB{r}{s}{\vecs{2}}{\OO{2}}$ to $\VTFB{r}{s}{\vecs{1}}{\OO{1}}$. That is $\TFpullback{f}{r}{s}$ is a bijection
and,
%%%%%%%%
\begin{align}
&\Foreach{\opair{\atf{1}}{\atf{2}}}{\Cprod{\TFB{r}{s}{\vecs{2}}{\OO{2}}}{\TFB{r}{s}{\vecs{2}}{\OO{2}}}}
\Foreach{c}{\algfield{}}\cr
&\func{\TFpullback{f}{r}{s}}{c\atf{1}+\atf{2}}=c\func{\TFpullback{f}{r}{s}}{\atf{1}}+
\func{\TFpullback{f}{r}{s}}{\atf{2}}.
\end{align}
Furthermore,
\begin{equation}
\finv{\(\TFpullback{f}{r}{s}\)}=\TFpullback{\(\finv{f}\)}{r}{s}.
\end{equation}
\proof
The linearity of $\TFpullback{f}{r}{s}$ follows directly from 
\refthm{thminverseofpullbackofrstensor}, \refdef{defspaceofrstensorfieldsonopensetsofvectorspaces},
and \refdef{defrspullbackofdiffeomorphism}. It is also easy to verify that $\TFpullback{\(\finv{f}\)}{r}{s}$
is both a left and right inverse for $\TFpullback{f}{r}{s}$, considering that
$\func{\banachderr{\finv{f}}}{\func{f}{\point}}$ is the inverse of $\func{\banachderr{f}}{\point}$
for every $\point\in\OO{1}$.
\endthm
%%%%%%%%%%%%%%%%%%%%%%%%%%%%%%%%%%%%%%%%%%%%%%%%%%%%%%%%%%%%%%%%%%%%%%%%%%%%%%%%%%%%%%%%%%%%%%%%%%%%%%%%%%%%%%%%%%%%%%%%%%%%%%%%
\proposition
Let $r$ be a non-negative integer, and $f$ an element of
$\banachmapdifclass{\infty}{\vecs{1}}{\vecs{2}}{\OO{1}}{\OO{2}}$ (that is an infinitely differentiable map
from $\OO{1}\subseteq\vecs{1}$ to $\OO{2}\subseteq\vecs{2}$).
The image of the assignment $\TFB{r}{0}{\vecs{2}}{\OO{2}}\ni\atf{}\mapsto\bar{\atf{}}\in
\Func{\OO{1}}{\Tensors{r}{0}{\vecs{1}}}$ such that
\begin{align}
\Foreach{\point}{\OO{1}}
\func{\bar{\atf{}}}{\point}\eqdef
\func{\Vpullbackcov{\[\func{\banachder{f}{\vecs{1}}{\vecs{2}}}{\point}\]}{r}}{\func{\atf{}}{\func{f}{\point}}},
\end{align}
is included in $\TFB{r}{0}{\vecs{1}}{\OO{1}}$. In other words, for every $\atf{}\in\VTFB{r}{0}{\vecs{2}}{\OO{2}}$,
$\bar{\atf{}}$ lies in $\VTFB{r}{0}{\vecs{1}}{\OO{1}}$.
\proof
It is left to the reader as an exercise.
\endpro
%%%%%%%%%%%%%%%%%%%%%%%%%%%%%%%%%%%%%%%%%%%%%%%%%%%%%%%%%%%%%%%%%%%%%%%%%%%%%%%%%%%%%%%%%%%%%%%%%%%%%%%%%%%%%%%%%%%%%%%%%%%%%%%%
\definition\label{defcovariantpullbackofsmoothmap}
Let $r$ be a non-negative integer, and $f$ an element of
$\banachmapdifclass{\infty}{\vecs{1}}{\vecs{2}}{\OO{1}}{\OO{2}}$ (that is a smooth map from $\OO{1}\subseteq\vecs{1}$
to $\OO{2}\subseteq\vecs{2}$). The mapping
$\function{\TFpullbackcov{f}{r}}{\TFB{r}{0}{\vecs{2}}{\OO{2}}}{\TFB{r}{0}{\vecs{1}}{\OO{1}}}$ is defined as
\begin{align}
&\Foreach{\atf{}}{\TFB{r}{0}{\vecs{2}}{\OO{2}}}
\Foreach{\point}{\OO{1}}\cr
&\func{\[\func{\TFpullbackcov{f}{r}}{\atf{}}\]}{\point}\eqdef
\func{\Vpullbackcov{\[\func{\banachder{f}{\vecs{1}}{\vecs{2}}}{\point}\]}{r}}{\func{\atf{}}{\func{f}{\point}}}.
\end{align}
$\TFpullbackcov{f}{r}$ is referred to as the $\quotl$$r$-pullback of the smooth map
$f\in\banachmapdifclass{\infty}{\vecs{1}}{\vecs{2}}{\OO{1}}{\OO{2}}$$\quotr$.
\endef
%%%%%%%%%%%%%%%%%%%%%%%%%%%%%%%%%%%%%%%%%%%%%%%%%%%%%%%%%%%%%%%%%%%%%%%%%%%%%%%%%%%%%%%%%%%%%%%%%%%%%%%%%%%%%%%%%%%%%%%%%%%%%%%%
\proposition
Let $r$ be a non-negative integer, and $f$ an element of
$\banachmapdifclass{\infty}{\vecs{1}}{\vecs{2}}{\OO{1}}{\OO{2}}$ (that is a smooth map from $\OO{1}\subseteq\vecs{1}$
to $\OO{2}\subseteq\vecs{2}$).
\begin{align}
&\Foreach{\atf{}}{\TFB{r}{0}{\vecs{2}}{\OO{2}}}
\Foreach{\point}{\OO{1}}
\Foreach{\mtuple{\vv{1}}{\vv{r}}}{\multiprod{\vecs{}}{r}}\cr
&\begin{aligned}
&\hskip0.6\baselineskip\func{\(\func{\[\func{\TFpullbackcov{f}{r}}{\atf{}}\]}{\point}\)}{\suc{\vv{1}}{\vv{r}}}\cr
&=
\func{\[\func{\atf{}}{\func{f}{\point}}\]}{\suc{\func{\[\func{\banachderr{f}}{\point}\]}{\vv{1}}}{\func{\[\func{\banachderr{f}}{\point}\]}{\vv{r}}}}.
\end{aligned}\cr
&{}
\end{align}
\proof
It is trivial.
\endpro
%%%%%%%%%%%%%%%%%%%%%%%%%%%%%%%%%%%%%%%%%%%%%%%%%%%%%%%%%%%%%%%%%%%%%%%%%%%%%%%%%%%%%%%%%%%%%%%%%%%%%%%%%%%%%%%%%%%%%%%%%%%%%%%%
\theorem\label{thmrpullbackofcompositionofsmoothmaps}
Let $r$ be a non-negative integer.
Let $f$ be an element of $\banachmapdifclass{\infty}{\vecs{}}{\vecs{1}}{\OO{}}{\OO{1}}$ and
$g$ an element of $\banachmapdifclass{\infty}{\vecs{1}}{\vecs{2}}{\OO{1}}{\OO{2}}$.
\begin{equation}
\TFpullbackcov{\(\identity{\OO{}}\)}{r}=\identity{\TFB{r}{0}{\vecs{}}{\OO{}}},
\end{equation}
and
\begin{equation}
\TFpullbackcov{\(\cmp{g}{f}\)}{r}=\cmp{\TFpullbackcov{f}{r}}{\TFpullbackcov{g}{r}}.
\end{equation}
\proof
The proof is similar to that of \refthm{thmrspullbackofcompositionofdiffeomorphisms}, as a result of
\refthm{thmrpullbackofcompositionoflinearmaps},
\refdef{defcovariantpullbackofsmoothmap}, and the chain rule of differentiation.
\endthm
%%%%%%%%%%%%%%%%%%%%%%%%%%%%%%%%%%%%%%%%%%%%%%%%%%%%%%%%%%%%%%%%%%%%%%%%%%%%%%%%%%%%%%%%%%%%%%%%%%%%%%%%%%%%%%%%%%%%%%%%%%%%%%%%%%%%%%%%%%%%%
\theorem\label{thmpullbackofcovtensorfieldisislinear}
Let $r$ be a non-negative integer.
Let $f$ be an element of $\banachmapdifclass{\infty}{\vecs{1}}{\vecs{2}}{\OO{1}}{\OO{2}}$
(that is an infinitely differentiable map from $\OO{1}\subseteq\vecs{1}$ to $\OO{2}\subseteq\vecs{2}$).
$\TFpullbackcov{f}{r}$ is a linear map
from $\VTFB{r}{0}{\vecs{2}}{\OO{2}}$ to $\VTFB{r}{0}{\vecs{1}}{\OO{1}}$. That is,
\begin{equation}
\Foreach{\opair{\atf{1}}{\atf{2}}}{\Cprod{\TFB{r}{0}{\vecs{2}}{\OO{2}}}{\TFB{r}{0}{\vecs{2}}{\OO{2}}}}
\Foreach{c}{\algfield{}}
\func{\TFpullbackcov{f}{r}}{c\atf{1}+\atf{2}}=c\func{\TFpullbackcov{f}{r}}{\atf{1}}+
\func{\TFpullbackcov{f}{r}}{\atf{2}}.
\end{equation}
\proof
The linearity of $\TFpullback{f}{r}{s}$ follows directly from 
\refthm{thmpullbackofcovtensorislinear}, \refdef{defcovariantpullbackofsmoothmap},
and \refdef{defrspullbackofdiffeomorphism}.
\endthm
%%%%%%%%%%%%%%%%%%%%%%%%%%%%%%%%%%%%%%%%%%%%%%%%%%%%%%%%%%%%%%%%%%%%%%%%%%%%%%%%%%%%%%%%%%%%%%%%%%%%%%%%%%%%%%%%%%%%%%%%%%%%%%%%
\definition
Let $f$ be an element of $\BanachDiff{\infty}{\vecs{1}}{\vecs{2}}{\OO{1}}{\OO{2}}$.
The mapping $\function{\TFPullback{f}}{\DVTFB{\vecs{2}}{\OO{2}}}{\DVTFB{\vecs{1}}{\OO{1}}}$ is defined term-wise as,
\begin{equation}
\Foreach{\atf{}}{\DVTFB{\vecs{2}}{\OO{2}}}
{\[\func{\TFPullback{f}}{\atf{}}\]}_{\opair{r}{s}}:=\func{\TFpullback{f}{r}{s}}{\atf{\opair{r}{s}}}.
\end{equation}
\endef
%%%%%%%%%%%%%%%%%%%%%%%%%%%%%%%%%%%%%%%%%%%%%%%%%%%%%%%%%%%%%%%%%%%%%%%%%%%%%%%%%%%%%%%%%%%%%%%%%%%%%%%%%%%%%%%%%%%%%%%%%%%%%%%%%%%%%%%%%%%%%
\corollary
Let $f$ be an element of $\BanachDiff{\infty}{\vecs{1}}{\vecs{2}}{\OO{1}}{\OO{2}}$.
$\TFPullback{f}$ is a linear isomorphism from $\DVTFB{\vecs{2}}{\OO{2}}$ to $\DVTFB{\vecs{1}}{\OO{1}}$.
\endcor
\chapteR{Smooth Fiber Bundles}
\thispagestyle{fancy}
%%%%%%%%%%%%%%%%%%%%%%%%%%%%%%%%%%%%%%%%%%%%%%%%%%%%%%%%%%%%%%%%%%%%%%%%%%%%%%%%%%%%%%%%%%%%%%%%%%%%%%%%%%%%%%%%%%%%%%%%%%%%%%%%
\section{Basic Structure of Smooth Fiber Bundles}
%%%%%%%%%%%%%%%%%%%%%%%%%%%%%%%%%%%%%%%%%%%%%%%%%%%%%%%%%%%%%%%%%%%%%%%%%%%%%%%%%%%%%%%%%%%%%%%%%%%%%%%%%%%%%%%%%%%%%%%%%%%%%%%%
\definition
Each
$\fbtotal{}=\opair{\Tot{}}{\maxatlas{\Tot{}}}$,
$\fbbase{}=\opair{\B{}}{\maxatlas{\B{}}}$,
$\fbfiber{}=\opair{\Fib{}}{\maxatlas{\Fib{}}}$ is taken as a $\difclass{\infty}$ manifold, and
$\fbprojection{}$ as an element of $\mapdifclass{\infty}{\fbtotal{}}{\fbbase{}}$, that is
a smooth map from $\fbtotal{}$ to $\fbbase{}$.
%%%%%%%%
The quadruple $\quadruple{\fbtotal{}}{\fbprojection{}}{\fbbase{}}{\fbfiber{}}$ is referred to as a $\quotl$smooth fiber bundle$\quotr$
iff the following property is satisfied.
\begin{itemize}
\item[\myitem{FB~1.}]
For every point $\point$ of $\fbbase{}$ there exists an open set $\U$ of $\Man{}$ containing $\point$, and
a diffeomorphism $\phi$ from $\subman{\fbtotal{}}{\func{\pimage{\fbprojection{}}}{\U}}$ to $\manprod{\subman{\fbbase{}}{\U}}{\fbfiber{}}$,
that is an element $\phi$ of
$\Diffeo{\infty}{\subman{\fbtotal{}}{\func{\pimage{\fbprojection{}}}{\U}}}{\manprod{\subman{\fbbase{}}{\U}}{\fbfiber{}}}$,
such that,
\begin{equation}
\Foreach{\point}{\U}
\func{\(\cmp{\proj{\U}{\Fib{}}{1}}{\phi}\)}{\point}=\func{\fbprojection{}}{\point},
\end{equation}
which means the following diagram commutes.
\begin{center}
\vskip0.5\baselineskip
\hskip-2\baselineskip
\begin{tikzcd}[row sep=6em, column sep=6em]
& \func{\pimage{\fbprojection{}}}{\U}
\arrow{r}{\phi}
\arrow[swap]{d}{\func{\res{\fbprojection{}}}{\func{\pimage{\fbprojection{}}}{\U}}}
& \Cprod{\U}{\Fib{}}
\arrow{dl}{\proj{\U}{\Fib{}}{1}} \\
& \U
\end{tikzcd}
\end{center}
\end{itemize}
Note that since $\U$ is an open set of $\fbbase{}$ and $\fbprojection{}$ is continuous,
$\func{\pimage{\fbprojection{}}}{\U}$ is an open set of $\fbtotal{}$, and therefore trivially
$\U$ and $\func{\pimage{\fbprojection{}}}{\U}$ are embedded submanifolds of $\fbbase{}$ and
$\fbtotal{}$, respectively. Thus, $\U$ and $\func{\pimage{\fbprojection{}}}{\U}$ inherit
a canonical differentiable structure from $\fbbase{}$ and $\fbtotal{}$, denoted here by
$\subman{\fbbase{}}{\U}$ and $\subman{\fbtotal{}}{\func{\pimage{\fbprojection{}}}{\U}}$, respectively.\\
When $\quadruple{\fbtotal{}}{\fbprojection{}}{\fbbase{}}{\fbfiber{}}$ is a smooth fiber bundle,
$\fbtotal{}$, $\fbprojection{}$, $\fbbase{}$, and $\fbfiber{}$ are referred to as the
$\quotl$total space$\quotr$, the $\quotl$projection$\quotr$, the $\quotl$base space$\quotr$,
and the $\quotl$(typical) fiber$\quotr$ of the smooth fiber bundle
$\quadruple{\fbtotal{}}{\fbprojection{}}{\fbbase{}}{\fbfiber{}}$, respectively.
For every point $\point$ of $\fbbase{}$, $\func{\pimage{\fbprojection{}}}{\seta{\point}}$
is called the $\quotl$fiber of the fiber bundle $\quadruple{\fbtotal{}}{\fbprojection{}}{\fbbase{}}{\fbfiber{}}$
over $\point$$\quotr$.\\
For convenience, the restriction of $\fbprojection{}$ to $\func{\pimage{\fbprojection{}}}{\U}$, that is
$\func{\res{\fbprojection{}}}{\func{\pimage{\fbprojection{}}}{\U}}$,
can be alternatively denoted by $\reS{\fbprojection{}}{\func{\pimage{\fbprojection{}}}{\U}}$.
\endef
%%%%%%%%%%%%%%%%%%%%%%%%%%%%%%%%%%%%%%%%%%%%%%%%%%%%%%%%%%%%%%%%%%%%%%%%%%%%%%%%%%%%%%%%%%%%%%%%%%%%%%%%%%%%%%%%%%%%%%%%%%%%%%%%
\fixed
$\fbundle{}=\quadruple{\fbtotal{}}{\fbprojection{}}{\fbbase{}}{\fbfiber{}}$ is fixed as a smooth fiber bundle, where
$\fbtotal{}=\opair{\Tot{}}{\maxatlas{\Tot{}}}$,
$\fbbase{}=\opair{\B{}}{\maxatlas{\B{}}}$,
$\fbfiber{}=\opair{\Fib{}}{\maxatlas{\Fib{}}}$ are $\difclass{\infty}$ manifolds
modeled on the Banach-spaces $\R^{n_{\Tot{}}}$, $\R^{n_{\B{}}}$, and $\R^{n_{\Fib{}}}$, respectively.\\
Also, $\defSet{\fbundle{i}=\quadruple{\fbtotal{i}}{\fbprojection{i}}{\fbbase{i}}{\fbfiber{i}}}{i\in\Zp}$ is fixed as a
collection of smooth fiber bundles such that
for each positive integer $i$, $\fbtotal{i}=\opair{\Tot{i}}{\maxatlas{\Tot{i}}}$,
$\fbbase{i}=\opair{\B{i}}{\maxatlas{\B{i}}}$,
$\fbfiber{i}=\opair{\Fib{i}}{\maxatlas{\Fib{i}}}$ are $\difclass{\infty}$ manifolds
modeled on the Banach-spaces $\R^{n_{\Tot{i}}}$, $\R^{n_{\B{i}}}$, and $\R^{n_{\Fib{i}}}$, respectively.\\
\endfixed
%%%%%%%%%%%%%%%%%%%%%%%%%%%%%%%%%%%%%%%%%%%%%%%%%%%%%%%%%%%%%%%%%%%%%%%%%%%%%%%%%%%%%%%%%%%%%%%%%%%%%%%%%%%%%%%%%%%%%%%%%%%%%%%%
\definition
%\begin{itemize}
%\item[$\cdot$]
$\phi$ is taken as a mapping from $\func{\pimage{\fbprojection{}}}{\U}$ to $\Cprod{\U}{\Fib{}}$, for some non-empty open set
$\U$ of $\fbbase{}$.
$\phi$ is referred to as a $\quotl$local trivialization (or FB-chart) of the fiber bundle $\fbundle{}$$\quotr$ iff
$\phi$ is a diffeomorphism from $\subman{\fbtotal{}}{\func{\pimage{\fbprojection{}}}{\U}}$ to
$\manprod{\subman{\fbbase{}}{\U}}{\fbfiber{}}$ such that $\cmp{\proj{\U}{\Fib{}}{1}}{\phi}=
\reS{\fbprojection{}}{\func{\pimage{\fbprojection{}}}{\U}}$.\\
The set of all local trivializations of $\fbundle{}$ is denoted by $\fbatlas{\fbundle{}}$. That is,
\begin{align}
&~\fbatlas{\fbundle{}}\cr
:=&~\defSet{\function{\phi}{\func{\pimage{\fbprojection{}}}{\U}}{\Cprod{\U}{\Fib{}}}}{\[\U\in\mantop{\fbbase{}},~
\phi\in\Diffeo{\infty}{\subman{\fbtotal{}}{\func{\pimage{\fbprojection{}}}{\U}}}{\manprod{\subman{\fbbase{}}{\U}}{\fbfiber{}}},~
\cmp{\proj{\U}{\Fib{}}{1}}{\phi}=
\reS{\fbprojection{}}{\func{\pimage{\fbprojection{}}}{\U}}\]}.\cr
&{}
\end{align}
$\fbatlas{\fbundle{}}$ is referred to as the $\quotl$FB-atlas of the fiber bundle $\fbundle{}$$\quotr$.\\
Any local trivialization $\phi$ of $\fbundle{}$ can be alternatively denoted by $\opair{\U}{\phi}$, where
$\U$ is the unique open set of $\fbbase{}$ such that $\domain{\phi}=\func{\pimage{\fbprojection{}}}{\U}$.
Any subset $\maxatlas{}$ of $\fbatlas{\fbundle{}}$ such that $\defSet{\U}{\opair{\U}{\phi}\in\maxatlas{}}$
is an open cover of $\fbbase{}$ is called an $\quotl$atlas of the smooth fiber bundle $\fbundle{}$$\quotr$.
%\end{itemize}
\endef
%%%%%%%%%%%%%%%%%%%%%%%%%%%%%%%%%%%%%%%%%%%%%%%%%%%%%%%%%%%%%%%%%%%%%%%%%%%%%%%%%%%%%%%%%%%%%%%%%%%%%%%%%%%%%%%%%%%%%%%%%%%%%%%%
\corollary
$\defSet{\domain{\phi}}{\phi\in\fbatlas{\fbundle{}}}$ is an open cover of $\fbbase{}$.
\endcor
%%%%%%%%%%%%%%%%%%%%%%%%%%%%%%%%%%%%%%%%%%%%%%%%%%%%%%%%%%%%%%%%%%%%%%%%%%%%%%%%%%%%%%%%%%%%%%%%%%%%%%%%%%%%%%%%%%%%%%%%%%%%%%%%
\lemma\label{lemlocaltrivializationcontraction}
Let $\opair{\U}{\phi}$ be a local trivialization of $\fbundle{}$, and $\V$ a subset of $\U$. The restriction of
$\phi$ to $\func{\pimage{\fbprojection{}}}{\V}$ is again a local trivialization of $\fbundle{}$.
\proof
Let $\p{\phi}$ denote the restriction of $\phi$ to $\func{\pimage{\fbprojection{}}}{\V}$. Then clearly
$\function{\p{\phi}}{\func{\pimage{\fbprojection{}}}{\V}}{\Cprod{\V}{\Fib{}}}$, and considering that
$\func{\pimage{\fbprojection{}}}{\V}$ and $\Cprod{\V}{\Fib{}}$ are open sets of
$\subman{\fbtotal{}}{\func{\pimage{\fbprojection{}}}{\U}}$ and $\manprod{\subman{\fbbase{}}{\U}}{\fbfiber{}}$ respectively,
it is immediately inferred that $\p{\phi}$
must be a diffeomorphism from $\subman{\fbtotal{}}{\func{\pimage{\fbprojection{}}}{\V}}$ to
$\manprod{\subman{\fbbase{}}{\V}}{\fbfiber{}}$ as the restriction of the diffeomorphism $\phi$
to an open set.
\endlem
%%%%%%%%%%%%%%%%%%%%%%%%%%%%%%%%%%%%%%%%%%%%%%%%%%%%%%%%%%%%%%%%%%%%%%%%%%%%%%%%%%%%%%%%%%%%%%%%%%%%%%%%%%%%%%%%%%%%%%%%%%%%%%%%
\theorem
Each $\Man{1}$ and $\Man{2}$ is taken as a $\difclass{\infty}$ manifold.
$\quadruple{\manprod{\Man{1}}{\Man{2}}}{\proj{\Man{1}}{\Man{2}}{1}}{\Man{1}}{\Man{2}}$ is a smooth fiber bundle.\\
\caution
This type of smooth fiber bundles are called $\quotl$trivial smooth fiber bundles$\quotr$.
\proof
It is trivial
\endthm
%%%%%%%%%%%%%%%%%%%%%%%%%%%%%%%%%%%%%%%%%%%%%%%%%%%%%%%%%%%%%%%%%%%%%%%%%%%%%%%%%%%%%%%%%%%%%%%%%%%%%%%%%%%%%%%%%%%%%%%%%%%%%%%%
\theorem
$\fbprojection{}$ is a submersion from $\fbtotal{}$ to $\fbbase{}$.
\proof
Let $\point$ be an arbitrary point pf $\fbbase{}$, and $\phi$ a local trivialization of $\fbundle{}$ such that
$\point\in\U:=\domain{\phi}$. Since $\proj{\U}{\Fib{}}{1}$ and $\phi$ are trivially submersions, so is
$\cmp{\proj{\U}{\Fib{}}{1}}{\phi}$, and thus considering that $\reS{\fbprojection{}}{\func{\pimage{\fbprojection{}}}{\U}}=
\cmp{\proj{\U}{\Fib{}}{1}}{\phi}$, it becomes evident that $\reS{\fbprojection{}}{\func{\pimage{\fbprojection{}}}{\U}}$
is a submersion from $\subman{\fbtotal{}}{\func{\pimage{\fbprojection{}}}{\U}}$
to $\subman{\fbbase{}}{\U}$. Thus considering that $\U$ and $\func{\pimage{\fbprojection{}}}{\U}$
are open submanifolds of $\fbbase{}$ and $\fbtotal{}$ respectively, $\fbprojection{}$ must be a submersion at every point of
$\func{\pimage{\fbprojection{}}}{\U}$, and in particular at $\point$.\\
Therefore $\fbprojection{}$ is a submersion at every points of $\fbtotal{}$, and hence a submersion from
$\fbtotal{}$ to $\fbbase{}$.
\endthm
%%%%%%%%%%%%%%%%%%%%%%%%%%%%%%%%%%%%%%%%%%%%%%%%%%%%%%%%%%%%%%%%%%%%%%%%%%%%%%%%%%%%%%%%%%%%%%%%%%%%%%%%%%%%%%%%%%%%%%%%%%%%%%%%
\theorem
For every point $\point$ of $\fbbase{}$, the fiber of $\fbundle{}$ over $\point$ is an embedded set of $\fbtotal{}$. That is,
\begin{align}
\Foreach{\point}{\fbbase{}}
%\begin{cases}
\func{\pimage{\fbprojection{}}}{\seta{\point}}\in\Emsubman{\fbtotal{}},
%\Diffeo{\infty}{\subman{\fbtotal{}}{\func{\pimage{\fbprojection{}}}{\seta{\point}}}}{\fbfiber{}}\neq\empty.
%\end{cases}
\end{align}
\proof
Let $\point$ be an arbitrary point pf $\fbbase{}$, and $\phi$ a local trivialization of $\fbundle{}$ such that
$\point\in\U:=\domain{\phi}$. Since $\phi$ is a bijection from $\func{\pimage{\fbprojection{}}}{\U}$ to
$\Cprod{\U}{\Fib{}}$, and $\reS{\fbprojection{}}{\func{\pimage{\fbprojection{}}}{\U}}=
\cmp{\proj{\U}{\Fib{}}{1}}{\phi}$, it is obvious that $\func{\image{\phi}}{\func{\pimage{\fbprojection{}}}{\seta{\point}}}=
\Cprod{\seta{\point}}{\Fib{}}$. It is also known that $\Cprod{\seta{\point}}{\Fib{}}$ is an embedded set of
$\manprod{\subman{\fbbase{}}{\U}}{\fbfiber{}}$, and hence considering that $\phi$ is a diffeomorphism from
$\subman{\fbtotal{}}{\func{\pimage{\fbprojection{}}}{\U}}$
to $\manprod{\subman{\fbbase{}}{\U}}{\fbfiber{}}$, $\func{\pimage{\fbprojection{}}}{\seta{\point}}=
\func{\pimage{\phi}}{\Cprod{\seta{\point}}{\Fib{}}}$ must be
an embedded set of $\subman{\fbtotal{}}{\func{\pimage{\fbprojection{}}}{\U}}$. Therefore considering that
$\func{\pimage{\fbprojection{}}}{\U}$ is itself an embedded set of the manifold $\fbtotal{}$ (because it is an open set of
$\fbtotal{}$), it becomes evident that $\func{\pimage{\fbprojection{}}}{\seta{\point}}$ is an embedded set of $\fbtotal{}$.
So $\func{\pimage{\fbprojection{}}}{\seta{\point}}$ inherits canonically a differentiable structure from that of $\fbtotal{}$,
endowed with which it becomes a manifold denoted by $\subman{\fbtotal{}}{\func{\pimage{\fbprojection{}}}{\seta{\point}}}$.\\
\endthm
%%%%%%%%%%%%%%%%%%%%%%%%%%%%%%%%%%%%%%%%%%%%%%%%%%%%%%%%%%%%%%%%%%%%%%%%%%%%%%%%%%%%%%%%%%%%%%%%%%%%%%%%%%%%%%%%%%%%%%%%%%%%%%%%
\definition
For every local trivialization $\opair{\U}{\phi}$ of the smooth fiber bundle $\fbundle{}$, let
\begin{equation}
\plt{\fbundle{}}{\phi}:=\cmp{\proj{\U}{\Fib{}}{2}}{\phi},
\end{equation}
which is referred to as the $\quotl$principal part of the local trivialization $\phi$ of $\fbundle{}$$\quotr$.\\
When there is no confusion about the underlying smooth fiber bundle, for every $\point\in\U$, the restriction of
$\plt{\fbundle{}}{\phi}$ to $\func{\pimage{\fbprojection{}}}{\seta{\point}}$, that is
$\func{\res{\plt{\fbundle{}}{\phi}}}{\func{\pimage{\fbprojection{}}}{\seta{\point}}}$, simply can be denoted by
$\reS{\plt{\fbundle{}}{\phi}}{\func{\pimage{\fbprojection{}}}{\seta{\point}}}$ or $\pltfib{\phi}{\point}$.
\endef
%%%%%%%%%%%%%%%%%%%%%%%%%%%%%%%%%%%%%%%%%%%%%%%%%%%%%%%%%%%%%%%%%%%%%%%%%%%%%%%%%%%%%%%%%%%%%%%%%%%%%%%%%%%%%%%%%%%%%%%%%%%%%%%%
\theorem
For every local trivialization $\opair{\U}{\phi}$ of $\fbundle{}$, $\plt{\fbundle{}}{\phi}$ is a smooth map from
$\subman{\fbtotal{}}{\func{\pimage{\fbprojection{}}}{\U}}$ to $\fbfiber{}$, and
$\phi=\opair{\reS{\fbprojection{}}{\func{\pimage{\fbprojection{}}}{\U}}}{\plt{\fbundle{}}{\phi}}$.
\begin{align}
\Foreach{\opair{\U}{\phi}}{\fbatlas{\fbundle{}}}
\begin{cases}
\plt{\fbundle{}}{\phi}\in\mapdifclass{\infty}{\subman{\fbtotal{}}{\func{\pimage{\fbprojection{}}}{\U}}}{\fbfiber{}},\cr
\Foreach{\x}{\func{\pimage{\fbprojection{}}}{\U}}\func{\phi}{\x}=\opair{\func{\fbprojection{}}{\x}}{\func{\plt{\fbundle{}}{\phi}}{\x}}.                    	
\end{cases}
\end{align}
\proof
It is trivial according to the definition of a local trivialization of a smooth fiber bundle and its principal part.
\endthm
%%%%%%%%%%%%%%%%%%%%%%%%%%%%%%%%%%%%%%%%%%%%%%%%%%%%%%%%%%%%%%%%%%%%%%%%%%%%%%%%%%%%%%%%%%%%%%%%%%%%%%%%%%%%%%%%%%%%%%%%%%%%%%%%
\theorem
Let $\opair{\U}{\phi}$ be a local trivialization of $\fbundle{}$. For every $\point\in\U$,
the restriction of $\plt{\fbundle{}}{\phi}$ to $\func{\pimage{\fbprojection{}}}{\seta{\point}}$ is a
diffeomorphism from $\subman{\fbtotal{}}{\func{\pimage{\fbprojection{}}}{\seta{\point}}}$ to $\fbfiber{}$.
That is,
\begin{equation}
\reS{\plt{\fbundle{}}{\phi}}{\func{\pimage{\fbprojection{}}}{\seta{\point}}}\in
\Diffeo{\infty}{\subman{\fbtotal{}}{\func{\pimage{\fbprojection{}}}{\seta{\point}}}}{\fbfiber{}},
\end{equation}
where $\reS{\plt{\fbundle{}}{\phi}}{\func{\pimage{\fbprojection{}}}{\seta{\point}}}:=
\func{\res{\plt{\fbundle{}}{\phi}}}{\func{\pimage{\fbprojection{}}}{\seta{\point}}}$.
\proof
Let $\point$ be an arbitrary element of $\U$.
Since $\subman{\fbtotal{}}{\func{\pimage{\fbprojection{}}}{\seta{\point}}}$ is an embedded set of
$\subman{\fbtotal{}}{\func{\pimage{\fbprojection{}}}{\U}}$, and $\plt{\fbundle{}}{\phi}$
is a smooth map from $\subman{\fbtotal{}}{\func{\pimage{\fbprojection{}}}{\U}}$ to $\fbfiber{}$,
clearly $\reS{\plt{\fbundle{}}{\phi}}{\func{\pimage{\fbprojection{}}}{\seta{\point}}}$
must be a smooth map from $\subman{\fbtotal{}}{\func{\pimage{\fbprojection{}}}{\seta{\point}}}$
to $\fbfiber{}$. In addition, considering that $\phi$ is a bijection from $\func{\pimage{\fbprojection{}}}{\U}$
to $\Cprod{\U}{\Fib{}}$ and $\cmp{\proj{\U}{\Fib{}}{1}}{\phi}=\reS{\fbprojection{}}{\func{\pimage{\fbprojection{}}}{\U}}$,
it is clear that $\cmp{\proj{\U}{\Fib{}}{2}}{\phi}$ forms a one-to-one corepondence between
$\func{\pimage{\fbprojection{}}}{\seta{\point}}$ and $\Fib{}$. That is,
$\reS{\plt{\fbundle{}}{\phi}}{\func{\pimage{\fbprojection{}}}{\seta{\point}}}$ is a bijection from
$\func{\pimage{\fbprojection{}}}{\seta{\point}}$ to $\Fib{}$. Moreover, define the map
$\function{\psi_{\point}}{\Cprod{\seta{\point}}{\Fib{}}}{\func{\pimage{\fbprojection{}}}{\seta{\point}}}$ as the restriction of $\finv{\phi}$
to $\Cprod{\seta{\point}}{\Fib{}}$, and the map $\function{\xi_{\point}}{\Fib{}}{\Cprod{\seta{\point}}{\Fib{}}}$ as
$\map{\xi_{\point}}{\x}{\opair{\point}{\x}}$. Clearly $\Cprod{\seta{\point}}{\Fib{}}$ is an embedded submanifold of
$\manprod{\subman{\fbbase{}}{\U}}{\fbfiber{}}$ and $\xi_{\point}$ is a smooth map from $\fbfiber{}$ to
$\Cprod{\seta{\point}}{\Fib{}}$. Moreover, considering that $\finv{\phi}$ is a smooth map from
$\manprod{\subman{\fbbase{}}{\U}}{\fbfiber{}}$ to $\subman{\fbtotal{}}{\func{\pimage{\fbprojection{}}}{\U}}$,
$\psi_{\point}$ must be a smooth map from $\Cprod{\seta{\point}}{\Fib{}}$ to
$\subman{\fbtotal{}}{\func{\pimage{\fbprojection{}}}{\seta{\point}}}$. Therefore according to the
composition rule of smooth maps, $\cmp{\psi_{\point}}{\xi_{\point}}$ is a smooth map from $\fbfiber{}$ to
$\subman{\fbtotal{}}{\func{\pimage{\fbprojection{}}}{\seta{\point}}}$. Furthermore, it can be easily verified that
$\finv{\(\reS{\plt{\fbundle{}}{\phi}}{\func{\pimage{\fbprojection{}}}{\seta{\point}}}\)}=\cmp{\psi_{\point}}{\xi_{\point}}$.
Therefore, it is inferred that $\reS{\plt{\fbundle{}}{\phi}}{\func{\pimage{\fbprojection{}}}{\seta{\point}}}$ is a
diffeomorphism from $\subman{\fbtotal{}}{\func{\pimage{\fbprojection{}}}{\seta{\point}}}$ to $\fbfiber{}$.
\endthm
%%%%%%%%%%%%%%%%%%%%%%%%%%%%%%%%%%%%%%%%%%%%%%%%%%%%%%%%%%%%%%%%%%%%%%%%%%%%%%%%%%%%%%%%%%%%%%%%%%%%%%%%%%%%%%%%%%%%%%%%%%%%%%%%
\corollary
Let $\opair{\U}{\phi}$ be a local trivialization of $\fbundle{}$.
\begin{equation}
\Foreach{\x}{\func{\pimage{\fbprojection{}}}{\U}}
\func{\phi}{\x}=\opair{\func{\fbprojection{}}{\x}}{\func{\plt{\fbundle{}}{\phi}}{\x}},
\end{equation}
and in particular,
\begin{equation}
\Foreach{\point}{\U}
\Foreach{\x}{\func{\pimage{\fbprojection{}}}{\seta{\point}}}
\func{\phi}{\x}=\opair{\point}{\func{\pltfib{\phi}{\point}}{\x}}.
\end{equation}
Also, for every point $\point\in\U$, $\subman{\fbtotal{}}{\func{\pimage{\fbprojection{}}}{\seta{\point}}}$
is diffeomorphic to $\fbfiber{}$, $\pltfib{\phi}{\point}$ being an instance for such a diffeomorphism.
\endcor
%%%%%%%%%%%%%%%%%%%%%%%%%%%%%%%%%%%%%%%%%%%%%%%%%%%%%%%%%%%%%%%%%%%%%%%%%%%%%%%%%%%%%%%%%%%%%%%%%%%%%%%%%%%%%%%%%%%%%%%%%%%%%%%%
\definition
Let $\opair{\U}{\phi}$ and $\opair{\V}{\psi}$ be a pair of local trivializations of the smooth fiber bundle $\fbundle{}$
such that $\U\cap\V\neq\empty$.
The mapping $\function{\transition{\fbundle{}}{\phi}{\psi}}{\U\cap\V}{\Diff{\infty}{\fbfiber{}}}$ is defined as,
\begin{equation}
\Foreach{\point}{\U\cap\V}
\func{\transition{\fbundle{}}{\phi}{\psi}}{\point}\eqdef
\cmp{\(\reS{\plt{\fbundle{}}{\phi}}{\func{\pimage{\fbprojection{}}}{\seta{\point}}}\)}
{\finv{\(\reS{\plt{\fbundle{}}{\psi}}{\func{\pimage{\fbprojection{}}}{\seta{\point}}}\)}}=
\cmp{\pltfib{\phi}{\point}}{\finv{\(\pltfib{\psi}{\point}\)}}.
\end{equation}
$\transition{\fbundle{}}{\phi}{\psi}$ is referred to as the $\quotl$transition map of the local trivializations
$\phi$ and $\psi$ of the smooth fiber bundle $\fbundle{}$$\quotr$. When the underlying smooth fiber bundle is clearly
recognized, $\transition{\fbundle{}}{\phi}{\psi}$ can simply be denoted by
$\transition{}{\phi}{\psi}$.
\endef
%%%%%%%%%%%%%%%%%%%%%%%%%%%%%%%%%%%%%%%%%%%%%%%%%%%%%%%%%%%%%%%%%%%%%%%%%%%%%%%%%%%%%%%%%%%%%%%%%%%%%%%%%%%%%%%%%%%%%%%%%%%%%%%%
\theorem
Let $\opair{\U}{\phi}$ and $\opair{\V}{\psi}$ be a pair of local trivializations of the smooth fiber bundle $\fbundle{}$
such that $\U\cap\V\neq\empty$. $\cmp{\phi}{\finv{\psi}}$ is a mapping from $\Cprod{\(\U\cap\V\)}{\Fib{}}$ to
$\Cprod{\(\U\cap\V\)}{\Fib{}}$, and,
\begin{equation}
\Foreach{\opair{\point}{z}}{\Cprod{\(\U\cap\V\)}{\Fib{}}}
\func{\(\cmp{\phi}{\finv{\psi}}\)}{\binary{\point}{z}}=
\opair{\point}{\func{\[\func{\transition{\fbundle{}}{\phi}{\psi}}{\point}\]}{z}}.
\end{equation}
\proof
It is trivial.
\endthm
%%%%%%%%%%%%%%%%%%%%%%%%%%%%%%%%%%%%%%%%%%%%%%%%%%%%%%%%%%%%%%%%%%%%%%%%%%%%%%%%%%%%%%%%%%%%%%%%%%%%%%%%%%%%%%%%%%%%%%%%%%%%%%%%
\theorem
Let $\opair{\U}{\phi}$, $\opair{\V}{\psi}$, and $\opair{\WW{}}{\eta}$ be local trivializations of
the smooth fiber bundle $\fbundle{}$ such that $\U\cap\V\cap\WW{}\neq\empty$.
\begin{align}
&\Foreach{\point}{\U}\func{\transition{\fbundle{}}{\phi}{\phi}}{\point}=\identity{\Fib{}},\\
&\Foreach{\point}{\U\cap\V}\func{\transition{\fbundle{}}{\phi}{\psi}}{\point}=
\finv{\[\func{\transition{\fbundle{}}{\psi}{\phi}}{\point}\]},\\
&\Foreach{\point}{\U\cap\V\cap\WW{}}\cmp{\cmp{\func{\transition{\fbundle{}}{\phi}{\psi}}{\point}}
{\func{\transition{\fbundle{}}{\psi}{\eta}}{\point}}}{\func{\transition{\fbundle{}}{\eta}{\phi}}{\point}}=\identity{\Fib{}}.
\end{align}
\proof
It is trivial.
\endthm
%%%%%%%%%%%%%%%%%%%%%%%%%%%%%%%%%%%%%%%%%%%%%%%%%%%%%%%%%%%%%%%%%%%%%%%%%%%%%%%%%%%%%%%%%%%%%%%%%%%%%%%%%%%%%%%%%%%%%%%%%%%%%%%%
\definition
Let $\function{\fbsec{}}{\B{}}{\Tot{}}$ be a smooth map from $\fbbase{}$ to $\fbtotal{}$.
$\fbsec{}$ is referred to as a $\quotl$(global) section of the smooth fiber bundle $\fbundle{}$$\quotr$ iff
$\cmp{\fbprojection{}}{\fbsec{}}=\identity{\B{}}$.\\
The set of all global sections of $\fbundle{}$ is denoted by $\fbsections{\fbundle{}}$. That is,
\begin{equation}
\fbsections{\fbundle{}}:=\defset{\fbsec{}}{\mapdifclass{\infty}{\fbbase{}}{\fbtotal{}}}
{\cmp{\fbprojection{}}{\fbsec{}}=\identity{\B{}}}.
\end{equation}
Now let $\U$ be a non-empty open set of $\fbbase{}$, and $\function{\fbsec{\U}}{\U}{\Tot{}}$ a smooth map from
$\subman{\fbbase{}}{\U}$ to $\fbtotal{}$. $\fbsec{\U}$ is referred to as a $\quotl$local section of the smooth
fiber bundle $\fbundle{}$ over $\U$$\quotr$ iff $\cmp{\fbprojection{}}{\fbsec{\U}}=\Injection{\U}{\B{}}$,
which means $\fbsec{\U}$ maps every point $\point$ of $\U$ into the fiber of $\fbundle{}$ over $\point$
(that is $\func{\pimage{\fbprojection{}}}{\seta{\point}}$).\\
The set of all local sections of $\fbundle{}$ over $\U$ is denoted by $\fbsectionsl{\fbundle{}}{\U}$.
\begin{equation}
\fbsectionsl{\fbundle{}}{\U}:=\defset{\fbsec{}}{\mapdifclass{\infty}{\subman{\fbbase{}}{\U}}{\fbtotal{}}}
{\cmp{\fbprojection{}}{\fbsec{}}=\Injection{\U}{\B{}}}.
\end{equation}
\endef
%%%%%%%%%%%%%%%%%%%%%%%%%%%%%%%%%%%%%%%%%%%%%%%%%%%%%%%%%%%%%%%%%%%%%%%%%%%%%%%%%%%%%%%%%%%%%%%%%%%%%%%%%%%%%%%%%%%%%%%%%%%%%%%%
%%%%%%%%%%%%%%%%%%%%%%%%%%%%%%%%%%%%%%%%%%%%%%%%%%%%%%%%%%%%%%%%%%%%%%%%%%%%%%%%%%%%%%%%%%%%%%%%%%%%%%%%%%%%%%%%%%%%%%%%%%%%%%%%
%%%%%%%%%%%%%%%%%%%%%%%%%%%%%%%%%%%%%%%%%%%%%%%%%%%%%%%%%%%%%%%%%%%%%%%%%%%%%%%%%%%%%%%%%%%%%%%%%%%%%%%%%%%%%%%%%%%%%%%%%%%%%%%%
\section{Smooth Fiber Bundle Morphisms}
\definition
Let $f$ be a smooth map from $\fbtotal{1}$ to $\fbtotal{2}$. $f$ is referred to as a $\quotl$smooth fiber bundle morphism
from the smooth fiber bundle $\fbundle{1}$ to the smooth fiber bundle $\fbundle{2}$$\quotr$ iff there exists a smooth map
$g$ from $\fbbase{1}$ to $\fbbase{2}$ such that $\cmp{\fbprojection{2}}{f}=\cmp{g}{\fbprojection{1}}$, that is the following
diagram commutes.
\begin{center}
\vskip0.5\baselineskip
\hskip-2\baselineskip
\begin{tikzcd}[row sep=6em, column sep=6em]
& \Tot{1}
\arrow{r}{f}
\arrow[swap]{d}{\fbprojection{1}}
& \Tot{2}
\arrow{d}{\fbprojection{2}} \\
& \B{1}
\arrow[swap]{r}{g}
& \B{2}
\end{tikzcd}
\end{center}
When $f$ is a morphism from $\fbundle{1}$ to $\fbundle{2}$, it can be easily verified that there exists only one
smooth map $g$ from $\fbbase{1}$ to $\fbbase{2}$ such that $\cmp{\fbprojection{2}}{f}=\cmp{g}{\fbprojection{1}}$,
and in addition $f$ maps the fiber of $\fbundle{1}$ over $\point$ into the fiber of $\fbundle{2}$ over $\func{g}{\point}$,
for every $\point\in\B{1}$. So, when $f$ is a morphism from $\fbundle{1}$ to $\fbundle{2}$, and
the source and and target smooth fiber bundles are clearly identified, we will simply denote by $\fbmorb{f}$
the unique smooth map from $\fbbase{1}$ to $\fbbase{2}$ such that
$\cmp{\fbprojection{2}}{f}=\cmp{\fbmorb{f}}{\fbprojection{1}}$; in this case, $f$ is called a $\quotl$morphism
from $\fbundle{1}$ to $\fbundle{2}$ along $\fbmorb{f}$$\quotr$.\\
The set of all morphisms from $\fbundle{1}$ to $\fbundle{2}$ is denoted by $\fbmorphisms{\fbundle{1}}{\fbundle{2}}$.
That is,
\begin{equation}
\fbmorphisms{\fbundle{1}}{\fbundle{2}}:=\defset{f}{\mapdifclass{\infty}{\fbtotal{1}}{\fbtotal{2}}}
{\(\Exists{g}{\mapdifclass{\infty}{\fbbase{1}}{\fbbase{2}}}\cmp{\fbprojection{2}}{f}=\cmp{g}{\fbprojection{1}}\)}.
\end{equation}
\endef
%%%%%%%%%%%%%%%%%%%%%%%%%%%%%%%%%%%%%%%%%%%%%%%%%%%%%%%%%%%%%%%%%%%%%%%%%%%%%%%%%%%%%%%%%%%%%%%%%%%%%%%%%%%%%%%%%%%%%%%%%%%%%%%%
\corollary\label{corfbmorphismmapsfiberintofiber}
Let $f$ be a morphism from $\fbundle{1}$ to $\fbundle{2}$. For every point $\point$ of $\fbbase{1}$,
$f$ maps the fiber of $\fbundle{1}$ over $\point$ into the fiber of $\fbundle{2}$ over $\func{\fbmorb{f}}{\point}$.
That is,
\begin{equation}
\Foreach{\point}{\B{1}}
\func{\image{f}}{\func{\pimage{\fbprojection{1}}}{\seta{\point}}}\subseteq
\func{\pimage{\fbprojection{2}}}{\seta{\func{\fbmorb{f}}{\point}}}.
\end{equation}
\proof
It is an immediate consequence of the definition of morphisms between smooth fiber bundles, and the fact that
the projection of a vector bundle is surjective.
\endcor
%%%%%%%%%%%%%%%%%%%%%%%%%%%%%%%%%%%%%%%%%%%%%%%%%%%%%%%%%%%%%%%%%%%%%%%%%%%%%%%%%%%%%%%%%%%%%%%%%%%%%%%%%%%%%%%%%%%%%%%%%%%%%%%%
\proposition
$\identity{\Tot{}}$ is a morphism from $\fbundle{}$ to itself, that is
$\identity{\Tot{}}\in\fbmorphisms{\fbundle{}}{\fbundle{}}$. Furtheremore,
$\fbmorb{\identity{\Tot{}}}=\identity{\B{}}$.
\proof
It is trivial.
\endpro
%%%%%%%%%%%%%%%%%%%%%%%%%%%%%%%%%%%%%%%%%%%%%%%%%%%%%%%%%%%%%%%%%%%%%%%%%%%%%%%%%%%%%%%%%%%%%%%%%%%%%%%%%%%%%%%%%%%%%%%%%%%%%%%%
\proposition\label{procompositionoffbmorphisms}
Let $f_1$ be a morphism from $\fbundle{1}$ to $\fbundle{2}$, and $f_2$ be a morphism from $\fbundle{2}$ to $\fbundle{3}$.
$\cmp{f_2}{f_1}$ is a morphism from $\fbundle{1}$ to $\fbundle{3}$ along $\cmp{\fbmorb{f_2}}{\fbmorb{f_1}}$.
\proof
According to the definition of the morphisms between smooth fiber bundles, we have the following commutative diagram,
which along with the consideration that the composition of smooth maps is again smooth, immediately implies the assertion.
\begin{center}
\vskip0.5\baselineskip
\hskip-2\baselineskip
\begin{tikzcd}[row sep=6em, column sep=6em]
& \Tot{1}
\arrow{r}{f_1}
\arrow[swap]{d}{\fbprojection{1}}
& \Tot{2}
\arrow{r}{f_2}
\arrow{d}{\fbprojection{2}}
& \Tot{3}
\arrow{d}{\fbprojection{3}} \\
& \B{1}
\arrow[swap]{r}{\fbmorb{f_1}}
& \B{2}
\arrow[swap]{r}{\fbmorb{f_2}}
& \B{3}
\end{tikzcd}
\end{center}
\endpro
%%%%%%%%%%%%%%%%%%%%%%%%%%%%%%%%%%%%%%%%%%%%%%%%%%%%%%%%%%%%%%%%%%%%%%%%%%%%%%%%%%%%%%%%%%%%%%%%%%%%%%%%%%%%%%%%%%%%%%%%%%%%%%%%
\remark
\textit{According to these propositions, it is inferred that any collection of smooth fiber bundles along with the set of
all possible morphisms between them forms a category,
having the ordinary composition of morphisms as the composition rule of its arrows.}
\endremark
%%%%%%%%%%%%%%%%%%%%%%%%%%%%%%%%%%%%%%%%%%%%%%%%%%%%%%%%%%%%%%%%%%%%%%%%%%%%%%%%%%%%%%%%%%%%%%%%%%%%%%%%%%%%%%%%%%%%%%%%%%%%%%%%
\definition
Let $f$ be a morphism from $\fbundle{1}$ to $\fbundle{2}$. $f$ is referred to as an $\quotl$(smooth fiber bundle)
isomorphism from $\fbundle{1}$ to $\fbundle{2}$$\quotr$ iff there exists a morphism $g$ from $\fbundle{2}$ to
$\fbundle{1}$ such that $\cmp{f}{g}=\identity{\B{2}}$, and $\cmp{g}{f}=\identity{\B{1}}$.\\
The set of all isomorphisms from $\fbundle{1}$ to $\fbundle{2}$ is denoted by $\fbisomorphisms{\fbundle{1}}{\fbundle{2}}$.
\endef
%%%%%%%%%%%%%%%%%%%%%%%%%%%%%%%%%%%%%%%%%%%%%%%%%%%%%%%%%%%%%%%%%%%%%%%%%%%%%%%%%%%%%%%%%%%%%%%%%%%%%%%%%%%%%%%%%%%%%%%%%%%%%%%%
\theorem\label{thmfbisomorphisms}
Let $f$ be a morphism from $\fbundle{1}$ to $\fbundle{2}$. $f$ is an isomorphism from $\fbundle{1}$ to $\fbundle{2}$
if and only if $f$ is a diffeomorphism from $\fbtotal{1}$ to $\fbtotal{2}$ and $\fbmorb{f}$ is a diffeomorphism from
$\fbbase{1}$ to $\fbbase{2}$. That is,
\begin{equation}
\fbisomorphisms{\fbundle{1}}{\fbundle{2}}=
\defset{f}{\fbmorphisms{\fbundle{1}}{\fbundle{2}}}{\(f\in\Diffeo{\infty}{\fbtotal{1}}{\fbtotal{2}},~
\fbmorb{f}\in\Diffeo{\infty}{\fbbase{1}}{\fbbase{2}}\)}.
\end{equation}
\proof
It is trivial.
\endthm
\chapteR{Smooth Vector Bundles}
\thispagestyle{fancy}
%%%%%%%%%%%%%%%%%%%%%%%%%%%%%%%%%%%%%%%%%%%%%%%%%%%%%%%%%%%%%%%%%%%%%%%%%%%%%%%%%%%%%%%%%%%%%%%%%%%%%%%%%%%%%%%%%%%%%%%%%%%%%%%%
\textit{$\algfield{}$ refers to one of the fields $\R$ or $\C$.
Given a $\algfield{}$-vector-space $\vbfiber{}$ with the finite dimension $n$, $\vecsmanifold{\vbfiber{}}$ stands for the
canonical differentiable structure of $\vbfiber{}$ endowed with which it becomes a $\difclass{\infty}$
manifold modeled on $\R^l$. When $\algfield{}=\R$, $l$ coincides with $n$, but when $\algfield{}=\C$,
we have $l=2n$.}
%%%%%%%%%%%%%%%%%%%%%%%%%%%%%%%%%%%%%%%%%%%%%%%%%%%%%%%%%%%%%%%%%%%%%%%%%%%%%%%%%%%%%%%%%%%%%%%%%%%%%%%%%%%%%%%%%%%%%%%%%%%%%%%%
\section{Basic Structur of Smooth Vector Bundles}
%%%%%%%%%%%%%%%%%%%%%%%%%%%%%%%%%%%%%%%%%%%%%%%%%%%%%%%%%%%%%%%%%%%%%%%%%%%%%%%%%%%%%%%%%%%%%%%%%%%%%%%%%%%%%%%%%%%%%%%%%%%%%%%%
\definition
Let each $\vbtotal{}=\opair{\vTot{}}{\maxatlas{\vTot{}}}$ and
$\vbbase{}=\opair{\vB{}}{\maxatlas{\vB{}}}$ be a $\difclass{\infty}$ manifold, and
$\vbprojection{}$ an element of $\mapdifclass{\infty}{\vbtotal{}}{\vbbase{}}$, that is
a smooth map from $\vbtotal{}$ to $\vbbase{}$. Also let $\vbfiber{}$ be a finite-dimensional $\algfield{}$-vector-space,
and $\vbatlas{}$ a collection of mappings $\function{\phi}{\func{\pimage{\vbprojection{}}}{\U}}{\Cprod{\U}{\vbfiber{}}}$
with $\U$ a non-empty open set of $\vbbase{}$.
%%%%%%%%
The quintuple $\quintuple{\vbtotal{}}{\vbprojection{}}{\vbbase{}}{\vbfiber{}}{\vbatlas{}}$ is referred to as a
$\quotl$smooth vector bundle$\quotr$ iff the following propertirs are satisfied.
\begin{itemize}
\item[\myitem{VB~1.}]
$\quadruple{\vbtotal{}}{\vbprojection{}}{\vbbase{}}{\vecsmanifold{\vbfiber{}}}$ is a smooth fiber bundle.
\item[\myitem{VB~2.}]
$\vbatlas{}$ is a maximal subset of the FB-atlas of
$\quadruple{\vbtotal{}}{\vbprojection{}}{\vbbase{}}{\vecsmanifold{\vbfiber{}}}$
such that the set of all domains of mappings of $\vbatlas{}$ forms an open covering of $\vbbase{}$
(that is a maximal atlas of the smooth fiber bundle
$\quadruple{\vbtotal{}}{\vbprojection{}}{\vbbase{}}{\vecsmanifold{\vbfiber{}}}$), and
for every pair $\opair{\U}{\phi}$ and $\opair{\V}{\psi}$ of local trivializations contained in $\vbatlas{}$
such that $\U\cap\V\neq\empty$, the value of the transition map of $\phi$ and $\psi$ (in the sense of the smooth fiber bundle
$\quadruple{\vbtotal{}}{\vbprojection{}}{\vbbase{}}{\vecsmanifold{\vbfiber{}}}$) at any point $\point\in\U\cap\V$
is a linear-isomorphism from the vector-space $\vbfiber{}$ to itself. That is,
\begin{equation}
\Foreach{\opair{\opair{\U}{\phi}}{\opair{\V}{\psi}}}
{\Cprod{\vbatlas{}}{\vbatlas{}}}
\bigg(\Foreach{\point}{\U\cap\V}
\func{\transition{\fbundle{}}{\phi}{\psi}}{\point}\in\GL{\vbfiber{}}{}\bigg),
\end{equation}
where $\fbundle{}:=\quadruple{\vbtotal{}}{\vbprojection{}}{\vbbase{}}{\vecsmanifold{\vbfiber{}}}$.
\end{itemize}
When $\vbundle{}:=\quintuple{\vbtotal{}}{\vbprojection{}}{\vbbase{}}{\vbfiber{}}{\vbatlas{}}$ is a smooth vector bundle:
\begin{itemize}
\item[$\cdot$]
$\vbtotal{}$, $\vbprojection{}$, $\vbbase{}$, $\vbfiber{}$, and $\vbatlas{}$ are called the
$\quotl$total space of the smooth vector bundle $\vbundle{}$$\quotr$,
the $\quotl$projection of $\vbundle{}$$\quotr$, the $\quotl$base space of $\vbundle{}$$\quotr$,
the $\quotl$(typical) fiber of $\vbundle{}$$\quotr$, and the
$\quotl$VB-atlas (or maximal-atlas) of the smooth vector bundle $\vbundle{}$$\quotr$, respectively.
\item[$\cdot$]
its $\quotl$rank$\quotr$ is defined to be the dimension of the $\algfield{}$-vector-space $\vbfiber{}$.
\item[$\cdot$]
it is also called a $\quotl$real (resp. complex) vector bundle$\quotr$ if $\algfield{}=\R$ (resp.
$\algfield{}=\C$).
\item[$\cdot$]
the underlying smooth fiber bundle structure of $\vbundle{}$, that is
$\quadruple{\vbtotal{}}{\vbprojection{}}{\vbbase{}}{\vecsmanifold{\vbfiber{}}}$, will simply be denoted by
$\Fvbundle{\vbundle{}}$.
\item[$\cdot$]
each element of $\vbatlas{}$ is called a $\quotl$local trivialization (or a VB-chart) of the smooth vector bundle $\vbundle{}$$\quotr$.
\item[$\cdot$]
for every point $\point$ of $\fbbase{}$, $\func{\pimage{\vbprojection{}}}{\seta{\point}}$ is called the
$\quotl$fiber of the smooth vector bundle $\vbundle{}$ over $\point$$\quotr$.
\item[$\cdot$]
any subset of $\vbatlas{}$ possessing necessarily all features stated in $[$\myitem{VB~2}$]$ except for being maximal,
is called an $\quotl$atlas of $\vbundle{}$$\quotr$.
\end{itemize}
\endef
%%%%%%%%%%%%%%%%%%%%%%%%%%%%%%%%%%%%%%%%%%%%%%%%%%%%%%%%%%%%%%%%%%%%%%%%%%%%%%%%%%%%%%%%%%%%%%%%%%%%%%%%%%%%%%%%%%%%%%%%%%%%%%%%
\fixed
$\vbundle{}=\quintuple{\vbtotal{}}{\vbprojection{}}{\vbbase{}}{\vbfiber{}}{\vbatlas{}}$ is fixed as a smooth vector bundle
of rank $d$,
where $\vbtotal{}=\opair{\vTot{}}{\maxatlas{\vTot{}}}$ and
$\vbbase{}=\opair{\vB{}}{\maxatlas{\vB{}}}$ are $\difclass{\infty}$ manifolds
modeled on the Banach-spaces $\R^{n_{\vTot{}}}$ and $\R^{n_{\vB{}}}$, respectively.\\
Also, $\defSet{\vbundle{i}=\quintuple{\vbtotal{i}}{\vbprojection{i}}{\vbbase{i}}{\vbfiber{i}}{\vbatlas{i}}}{i\in\Z}$
is fixed as a collection of smooth vector bundles such that
for each positive integer $i$, $\vbtotal{i}=\opair{\vTot{i}}{\maxatlas{\vTot{i}}}$,
$\vbbase{i}=\opair{\vB{i}}{\maxatlas{\vB{i}}}$ are $\difclass{\infty}$ manifolds
modeled on the Banach-spaces $\R^{n_{\vTot{i}}}$ and $\R^{n_{\vB{i}}}$, respectively. The rank of $\vbundle{i}$
is taken to be $r_i$.\\
\endfixed
%%%%%%%%%%%%%%%%%%%%%%%%%%%%%%%%%%%%%%%%%%%%%%%%%%%%%%%%%%%%%%%%%%%%%%%%%%%%%%%%%%%%%%%%%%%%%%%%%%%%%%%%%%%%%%%%%%%%%%%%%%%%%%%%
\remark
\textit{When $\algfield{}=\R$, the rank of the smooth vector bundle
$\vbundle{}=\quintuple{\vbtotal{}}{\vbprojection{}}{\vbbase{}}{\vbfiber{}}{\vbatlas{}}$ coincides with the rank of
the smooth fiber bundle $\quadruple{\vbtotal{}}{\vbprojection{}}{\vbbase{}}{\vecsmanifold{\vbfiber{}}}$.
When $\algfield{}=\C$, the rank of the smooth fiber bundle
$\quadruple{\vbtotal{}}{\vbprojection{}}{\vbbase{}}{\vecsmanifold{\vbfiber{}}}$
is twice the rank of $\vbundle{}$.}
\endremark
%%%%%%%%%%%%%%%%%%%%%%%%%%%%%%%%%%%%%%%%%%%%%%%%%%%%%%%%%%%%%%%%%%%%%%%%%%%%%%%%%%%%%%%%%%%%%%%%%%%%%%%%%%%%%%%%%%%%%%%%%%%%%%%%
\definition
\begin{itemize}
\item[$\cdot$]
For every local trivialization $\opair{\U}{\phi}$ of the smooth vector bundle $\vbundle{}$ ($\phi\in\vbatlas{}$), let
$\plt{\vbundle{}}{\phi}$ be the same as $\plt{\Fvbundle{\vbundle{}}}{\phi}$ (in the sense of smooth fiber bundles),
which is referred to as the $\quotl$principal part of the local trivialization $\phi$ of $\vbundle{}$$\quotr$.\\
When there is no confusion about the underlying smooth vector bundle, for every $\point\in\U$, the restriction of
$\plt{\vbundle{}}{\phi}$ to $\func{\pimage{\vbprojection{}}}{\seta{\point}}$, that is
$\func{\res{\plt{\vbundle{}}{\phi}}}{\func{\pimage{\vbprojection{}}}{\seta{\point}}}$, simply can be denoted by
$\reS{\plt{\vbundle{}}{\phi}}{\func{\pimage{\vbprojection{}}}{\seta{\point}}}$ or $\pltfib{\phi}{\point}$.
\item[$\cdot$]
Let $\opair{\U}{\phi}$ and $\opair{\V}{\psi}$ be a pair of local trivializations of the smooth vector bundle $\vbundle{}$
such that $\U\cap\V\neq\empty$.
The mapping $\function{\transition{\vbundle{}}{\phi}{\psi}}{\U\cap\V}{\GL{\vbfiber{}}{}}$ is defined to be the same as
$\transition{\Fvbundle{\vbundle{}}}{\phi}{\psi}$ (in the sense of smooth fiber bundles).
$\transition{\vbundle{}}{\phi}{\psi}$ is referred to as the $\quotl$transition map of the local trivializations
$\phi$ and $\psi$ of the smooth vector bundle $\vbundle{}$$\quotr$. When the underlying smooth vector bundle is clearly
recognized, $\transition{\vbundle{}}{\phi}{\psi}$ can simply be denoted by $\transition{}{\phi}{\psi}$.
\end{itemize}
\endef
%%%%%%%%%%%%%%%%%%%%%%%%%%%%%%%%%%%%%%%%%%%%%%%%%%%%%%%%%%%%%%%%%%%%%%%%%%%%%%%%%%%%%%%%%%%%%%%%%%%%%%%%%%%%%%%%%%%%%%%%%%%%%%%%
\theorem\label{thmvbfiberslinearstrucure}
Let $\point$ be a point of $\vbbase{}$, and let $\opair{\U}{\phi}$ and $\opair{\V}{\psi}$
be a pair of local trivializations of $\vbundle{}$ such that $\point\in\U\cap\V$.
\begin{align}
&\Foreach{\opair{a}{b}}{\Cprod{\algfield{}}{\algfield{}}}
\Foreach{\opair{u}{v}}{\Cprod{\func{\pimage{\vbprojection{}}}{\seta{\point}}}{\func{\pimage{\vbprojection{}}}{\seta{\point}}}}\cr
&\hskip\baselineskip
\begin{aligned}
&\hskip0.5\baselineskip
\func{\finv{\(\reS{\plt{\vbundle{}}{\phi}}{\func{\pimage{\vbprojection{}}}{\seta{\point}}}\)}}{
a\[\func{\reS{\plt{\vbundle{}}{\phi}}{\func{\pimage{\vbprojection{}}}{\seta{\point}}}}{u}\]+
b\[\func{\reS{\plt{\vbundle{}}{\phi}}{\func{\pimage{\vbprojection{}}}{\seta{\point}}}}{v}\]
}\cr
&=\func{\finv{\(\reS{\plt{\vbundle{}}{\psi}}{\func{\pimage{\vbprojection{}}}{\seta{\point}}}\)}}{
a\[\func{\reS{\plt{\vbundle{}}{\psi}}{\func{\pimage{\vbprojection{}}}{\seta{\point}}}}{u}\]+
\[b\func{\reS{\plt{\vbundle{}}{\psi}}{\func{\pimage{\vbprojection{}}}{\seta{\point}}}}{v}\]
},
\end{aligned}
\end{align}
or denoted simply,
\begin{align}
&\Foreach{\opair{a}{b}}{\Cprod{\algfield{}}{\algfield{}}}
\Foreach{\opair{u}{v}}{\Cprod{\func{\pimage{\vbprojection{}}}{\seta{\point}}}{\func{\pimage{\vbprojection{}}}{\seta{\point}}}}\cr
&\hskip\baselineskip
\func{\finv{\(\pltfib{\phi}{\point}\)}}{
a\[\func{\pltfib{\phi}{\point}}{u}\]+
b\[\func{\pltfib{\phi}{\point}}{v}\]
}=
\func{\finv{\(\pltfib{\psi}{\point}\)}}{
a\[\func{\pltfib{\psi}{\point}}{u}\]+
b\[\func{\pltfib{\psi}{\point}}{v}\]
}.
\end{align}
\proof
Since $\cmp{\pltfib{\phi}{\point}}{\finv{\(\pltfib{\psi}{\point}\)}}=\transition{\vbundle{}}{\phi}{\psi}\in\GL{\vbfiber{}}{}$,
clearly for every $a$ and $b$ in $\algfield{}$ and every $u$ and $v$ in $\func{\pimage{\vbprojection{}}}{\seta{\point}}$,
\begin{align}
&~\func{\finv{\(\pltfib{\psi}{\point}\)}}{
a\[\func{\pltfib{\psi}{\point}}{u}\]+
b\[\func{\pltfib{\psi}{\point}}{v}\]}\cr
=&~
\func{\[\cmp{\cmp{\finv{\(\pltfib{\phi}{\point}\)}}{\pltfib{\phi}{\point}}}{\finv{\(\pltfib{\psi}{\point}\)}}\]}{
a\[\func{\pltfib{\psi}{\point}}{u}\]+
b\[\func{\pltfib{\psi}{\point}}{v}\]}\cr
=&~\func{\finv{\(\pltfib{\phi}{\point}\)}}
{a\[\func{\cmp{\cmp{\pltfib{\phi}{\point}}{\finv{\(\pltfib{\psi}{\point}\)}}}{\pltfib{\psi}{\point}}}{u}\]+
b\[\func{\cmp{\cmp{\pltfib{\phi}{\point}}{\finv{\(\pltfib{\psi}{\point}\)}}}{\pltfib{\psi}{\point}}}{v}\]}\cr
=&~\func{\finv{\(\pltfib{\phi}{\point}\)}}{
a\[\func{\pltfib{\phi}{\point}}{u}\]+
b\[\func{\pltfib{\phi}{\point}}{v}\]}.
\end{align}
\endthm
%%%%%%%%%%%%%%%%%%%%%%%%%%%%%%%%%%%%%%%%%%%%%%%%%%%%%%%%%%%%%%%%%%%%%%%%%%%%%%%%%%%%%%%%%%%%%%%%%%%%%%%%%%%%%%%%%%%%%%%%%%%%%%%%
\definition\label{defvbfiberslinearstrucure}
Let $\point$ be a point of $\vbbase{}$. $\fibervecs{\vbundle{}}{\point}$
is defined to be the $\algfield{}$-vector-space formed by endowing $\func{\pimage{\vbprojection{}}}{\seta{\point}}$
with the linear structure:
\begin{align}
&\Foreach{\opair{a}{b}}{\Cprod{\algfield{}}{\algfield{}}}
\Foreach{\opair{u}{v}}{\Cprod{\func{\pimage{\vbprojection{}}}{\seta{\point}}}{\func{\pimage{\vbprojection{}}}{\seta{\point}}}}\cr
&\hskip\baselineskip
au+bv\eqdef\func{\finv{\(\pltfib{\phi}{\point}\)}}{
a\[\func{\pltfib{\phi}{\point}}{u}\]+
b\[\func{\pltfib{\phi}{\point}}{v}\]},
\end{align}
where $\phi$ is any element of $\vbatlas{}$ (local trivialization of $\vbundle{}$) such that $\point\in\domain{\phi}$.
$\fibervecs{\vbundle{}}{\point}$ is referred to as the $\quotl$(canonical) fiber space of $\vbundle{}$ over $\point$$\quotr$.
When the underlying smooth vector bundle is clearly identified, $\fibervecs{\vbundle{}}{\point}$
can simply be denoted by $\fibervecs{}{\point}$.\\
\caution
This linear structure is well-defined in the shadow of \refthm{thmvbfiberslinearstrucure}.
\endef
%%%%%%%%%%%%%%%%%%%%%%%%%%%%%%%%%%%%%%%%%%%%%%%%%%%%%%%%%%%%%%%%%%%%%%%%%%%%%%%%%%%%%%%%%%%%%%%%%%%%%%%%%%%%%%%%%%%%%%%%%%%%%%%%
\theorem
Let $\opair{\U}{\phi}$ be an element of $\vbatlas{}$ (a local trivialization of $\vbundle{}$).
For every $\point\in\U$, $\reS{\plt{\vbundle{}}{\phi}}{\func{\pimage{\vbprojection{}}}{\seta{\point}}}$ is
a linear-isomorphism from $\fibervecs{\vbundle{}}{\point}$ to $\vbfiber{}$.
\proof
It is already known that $\pltfib{\phi}{\point}$ is a bijection from $\func{\pimage{\vbprojection{}}}{\seta{\point}}$
to $\vbfiber{}$, actually a diffeomorphism from $\subman{\vbtotal{}}{\func{\pimage{\vbprojection{}}}{\seta{\point}}}$
to $\vecsmanifold{\vbfiber{}}$. In addition, the linearity of $\pltfib{\phi}{\point}$ is obvious according to
\refdef{defvbfiberslinearstrucure}.
\endthm
%%%%%%%%%%%%%%%%%%%%%%%%%%%%%%%%%%%%%%%%%%%%%%%%%%%%%%%%%%%%%%%%%%%%%%%%%%%%%%%%%%%%%%%%%%%%%%%%%%%%%%%%%%%%%%%%%%%%%%%%%%%%%%%%
\corollary
For every point $\point$ of $\vbbase{}$, the fiber space of $\vbundle{}$ over $\point$, that is
$\fibervecs{\vbundle{}}{\point}$, is linearly isomorphic to $\vbfiber{}$.
\endcor
%%%%%%%%%%%%%%%%%%%%%%%%%%%%%%%%%%%%%%%%%%%%%%%%%%%%%%%%%%%%%%%%%%%%%%%%%%%%%%%%%%%%%%%%%%%%%%%%%%%%%%%%%%%%%%%%%%%%%%%%%%%%%%%%
\theorem
For every point $\point$ of $\vbbase{}$, there exists an open set $\U$ of $\vbbase{}$ containing $\point$,
a chart $\phi$ of the manifold $\vbbase{}$, and a local trivialization $\psi$ of the smooth vector vundle
$\vbundle{}$ such that $\func{\image{\vbprojection{}}}{\domain{\psi}}=\domain{\phi}=\U$. In other words, for
every point $\point$ of $\vbbase{}$, there exists an open set $\U$ of $\vbbase{}$ containing $\point$,
a chart $\opair{\U}{\phi}$ of $\vbbase{}$, and a local trivialization $\opair{\U}{\psi}$ of $\vbundle{}$.
\proof
Let $\point$ be a point of $\vbbase{}$. Trivially there exists a chart $\opair{\U}{\phi}$ of $\vbbase{}$
and a local trivialization $\opair{\V}{\psi}$ of $\vbundle{}$ such that $\point\in\U\cap\V$.
Let $\p{\phi}$ denote the restriction of $\phi$ to $\U\cap\V$ and $\p{\psi}$ denote the restriction of
$\psi$ to $\func{\pimage{\vbprojection{}}}{\U\cap\V}$. Then clearly $\opair{\U\cap\V}{\p{\phi}}$ is a chart of $\vbbase{}$
and $\opair{\U\cap\V}{\p{\psi}}$ is a local trivialization of $\vbundle{}$ as a result of
\reflem{lemlocaltrivializationcontraction}.
\endthm
%%%%%%%%%%%%%%%%%%%%%%%%%%%%%%%%%%%%%%%%%%%%%%%%%%%%%%%%%%%%%%%%%%%%%%%%%%%%%%%%%%%%%%%%%%%%%%%%%%%%%%%%%%%%%%%%%%%%%%%%%%%%%%%%
%%%%%%%%%%%%%%%%%%%%%%%%%%%%%%%%%%%%%%%%%%%%%%%%%%%%%%%%%%%%%%%%%%%%%%%%%%%%%%%%%%%%%%%%%%%%%%%%%%%%%%%%%%%%%%%%%%%%%%%%%%%%%%%%
%%%%%%%%%%%%%%%%%%%%%%%%%%%%%%%%%%%%%%%%%%%%%%%%%%%%%%%%%%%%%%%%%%%%%%%%%%%%%%%%%%%%%%%%%%%%%%%%%%%%%%%%%%%%%%%%%%%%%%%%%%%%%%%%
%%%%%%%%%%%%%%%%%%%%%%%%%%%%%%%%%%%%%%%%%%%%%%%%%%%%%%%%%%%%%%%%%%%%%%%%%%%%%%%%%%%%%%%%%%%%%%%%%%%%%%%%%%%%%%%%%%%%%%%%%%%%%%%%
\section{Sections of a Smooth Vector Bundle}
%%%%%%%%%%%%%%%%%%%%%%%%%%%%%%%%%%%%%%%%%%%%%%%%%%%%%%%%%%%%%%%%%%%%%%%%%%%%%%%%%%%%%%%%%%%%%%%%%%%%%%%%%%%%%%%%%%%%%%%%%%%%%%%%
\definition
Let $\function{\vbsec{}}{\vB{}}{\vTot{}}$ be a smooth map from $\vbbase{}$ to $\vbtotal{}$.
$\vbsec{}$ is referred to as a $\quotl$(global) section of the smooth vector bundle $\vbundle{}$$\quotr$ iff
it is a global section of the underlying smooth fiber bundle of $\vbundle{}$,
that is $\cmp{\vbprojection{}}{\vbsec{}}=\identity{\vB{}}$.\\
The set of all global sections of $\vbundle{}$ is denoted by $\vbsections{\vbundle{}}$. That is,
\begin{equation}
\vbsections{\vbundle{}}:=\defset{\vbsec{}}{\mapdifclass{\infty}{\vbbase{}}{\vbtotal{}}}
{\cmp{\vbprojection{}}{\vbsec{}}=\identity{\vB{}}}.
\end{equation}
Now let $\U$ be a non-empty open set of $\vbbase{}$, and $\function{\vbsec{\U}}{\U}{\vTot{}}$ a smooth map from
$\subman{\vbbase{}}{\U}$ to $\vbtotal{}$. $\vbsec{\U}$ is referred to as a $\quotl$local section of the smooth
vector bundle $\vbundle{}$ over $\U$$\quotr$ iff it is a local section of the underlying smooth fiber bundle of $\vbundle{}$,
that is $\cmp{\fbprojection{}}{\fbsec{\U}}=\Injection{\U}{\B{}}$.
%which means $\fbsec{\U}$ maps every point $\point$ of $\U$ into the fiber of $\fbundle{}$ over $\point$
%(that is $\func{\pimage{\fbprojection{}}}{\seta{\point}}$).
\\
The set of all local sections of $\vbundle{}$ over $\U$ is denoted by $\vbsectionsl{\vbundle{}}{\U}$.
\begin{equation}
\vbsectionsl{\vbundle{}}{\U}:=\defset{\vbsec{}}{\mapdifclass{\infty}{\subman{\vbbase{}}{\U}}{\vbtotal{}}}
{\cmp{\vbprojection{}}{\vbsec{}}=\Injection{\U}{\vB{}}}.
\end{equation}
\endef
%%%%%%%%%%%%%%%%%%%%%%%%%%%%%%%%%%%%%%%%%%%%%%%%%%%%%%%%%%%%%%%%%%%%%%%%%%%%%%%%%%%%%%%%%%%%%%%%%%%%%%%%%%%%%%%%%%%%%%%%%%%%%%%%
\definition
The operations
$\function{\vbsadd{\vbundle{}}}{\Cprod{\vbsections{\vbundle{}}}{\vbsections{\vbundle{}}}}{\vbsections{\vbundle{}}}$,
$\function{\vbsprod{\vbundle{}}}{\Cprod{\smoothring{\infty}{\vbbase{}}}{\vbsections{\vbundle{}}}}{\vbsections{\vbundle{}}}$, and
$\function{\vbsscalprod{\vbundle{}}}{\Cprod{\algfield{}}{\vbsections{\vbundle{}}}}{\vbsections{\vbundle{}}}$
are defined as,
\begin{align}
&\Foreach{\opair{\vbsec{1}}{\vbsec{2}}}{\Cprod{\vbsections{\vbundle{}}}{\vbsections{\vbundle{}}}}
\Foreach{\point}{\vB{}}
\func{\[\vbsec{1}\vbsadd{\vbundle{}}\vbsec{2}\]}{\point}\eqdef\func{\vbsec{1}}{\point}+\func{\vbsec{2}}{\point},\cr
%%%%%%%%%%%
&\Foreach{\opair{f}{\vbsec{}}}{\Cprod{\smoothring{\infty}{\vbbase{}}}{\vbsections{\vbundle{}}}}
\Foreach{\point}{\vB{}}
\func{\[f\vbsprod{\vbundle{}}\vbsec{}\]}{\point}\eqdef\func{f}{\point}\func{\vbsec{}}{\point},\cr
%%%%%%%%%%%
&\Foreach{\opair{c}{\vbsec{}}}{\Cprod{\algfield{}}{\vbsections{\vbundle{}}}}
\Foreach{\point}{\vB{}}
\func{\[c\vbsscalprod{\vbundle{}}\vbsec{}\]}{\point}\eqdef c\[\func{\vbsec{}}{\point}\].
\end{align}
Note that the addition of $\func{\vbsec{1}}{\point}$ and $\func{\vbsec{2}}{\point}$
and the product of $\func{f}{\point}$ and $\func{\vbsec{}}{\point}$ take place in $\fibervecs{\vbundle{}}{\point}$,
that is the fiber space of $\vbundle{}$ over $\point$. When the context is clear enough, $\vbsadd{\vbundle{}}$
can be denoted by $+$, $f\vbsprod{\vbundle{}}\vbsec{}$ by $f\vbsec{}$, and
$c\vbsscalprod{\vbundle{}}\vbsec{}$ by $c\vbsec{}$ simply.\\
\caution
$\smoothring{\infty}{\vbbase{}}$ refers to the commutative ring of all real-valued smooth maps on $\vbbase{}$.
\endef
%%%%%%%%%%%%%%%%%%%%%%%%%%%%%%%%%%%%%%%%%%%%%%%%%%%%%%%%%%%%%%%%%%%%%%%%%%%%%%%%%%%%%%%%%%%%%%%%%%%%%%%%%%%%%%%%%%%%%%%%%%%%%%%%
\theorem
The quadruple $\quadruple{\vbsections{\vbundle{}}}{\vbsadd{\vbundle{}}}{\vbsprod{\vbundle{}}}{\smoothring{\infty}{\vbbase{}}}$
has the structure of a module over a commutative ring. We will make reference to this module by the same notation
$\vbsections{\vbundle{}}$. The quadruple
$\quadruple{\vbsections{\vbundle{}}}{\vbsadd{\vbundle{}}}{\vbsscalprod{\vbundle{}}}{\algfield{}}$
has the structure of a $\algfield{}$-vector-space.
\proof
It is straightforward to check the axioms for module and vector-space structures.
\endthm
%%%%%%%%%%%%%%%%%%%%%%%%%%%%%%%%%%%%%%%%%%%%%%%%%%%%%%%%%%%%%%%%%%%%%%%%%%%%%%%%%%%%%%%%%%%%%%%%%%%%%%%%%%%%%%%%%%%%%%%%%%%%%%%%
\definition
For every point $\point$ of $\vbbase{}$, an ordered basis of the vector-space $\fibervecs{\vbundle{}}{\point}$
is called a $\quotl$frame of $\vbundle{}$ at $\point$$\quotr$.\\
Let $\U$ be an open set of $\vbbase{}$. Any sequence $\suc{\vbsec{1}}{\vbsec{d}}$ of local sections of $\vbundle{}$
over $\U$ such that for every $\point\in\U$ the sequence $\suc{\func{\vbsec{1}}{\point}}{\func{\vbsec{d}}{\point}}$
is an ordered basis of $\fibervecs{\vbundle{}}{\point}$ (that is a frame of $\vbundle{}$ at $\point$), is referred to as a
$\quotl$local frame field of $\vbundle{}$ over $\U$$\quotr$. The set of all local frame fields of
$\vbundle{}$ over $\U$ is denoted by $\lframe{\vbundle{}}{\U}$. Any element of $\lframe{\vbundle{}}{\vB{}}$ is alternatively
called a $\quotl$global frame field of the smooth vector bundle $\vbundle{}$.
\endef
%%%%%%%%%%%%%%%%%%%%%%%%%%%%%%%%%%%%%%%%%%%%%%%%%%%%%%%%%%%%%%%%%%%%%%%%%%%%%%%%%%%%%%%%%%%%%%%%%%%%%%%%%%%%%%%%%%%%%%%%%%%%%%%%
\theorem
Let $\opair{\U}{\phi}$ be a local trivialization of the smooth vector bundle $\vbundle{}$, and let
$\suc{\vsbase{1}}{\vsbase{d}}$ be an ordered basis of fiber space $\vbfiber{}$.\\
The sequence of maps
$\suc{\function{\vbsec{1}}{\U}{\vTot{}}}{\function{\vbsec{d}}{\U}{\vTot{}}}$ defined as
\begin{equation}
\Foreach{i}{\seta{\suc{1}{d}}}
\Foreach{\point}{\U}
\func{\vbsec{i}}{\point}\eqdef\func{\finv{\(\pltfib{\phi}{\point}\)}}{\vsbase{i}}
=\func{\finv{\phi}}{\binary{\point}{\vsbase{i}}},
\end{equation}
is a local frame field of $\vbundle{}$ over $\U$.
\proof
Clearly each $\vsbase{i}$ is smooth, and thus a local section of $\vbundle{}$ over $\U$.
Moreover, according to the fact that each $\finv{\(\pltfib{\phi}{\point}\)}$
is a linear isomorphism from $\fibervecs{\vbundle{}}{\point}$ to $\vbfiber{}$, it is evident that for every $\point\in\U$,
the sequence $\suc{\func{\vbsec{1}}{\point}}{\func{\vbsec{d}}{\point}}$ is an ordered base of $\fibervecs{\vbundle{}}{\point}$.
\endthm
%%%%%%%%%%%%%%%%%%%%%%%%%%%%%%%%%%%%%%%%%%%%%%%%%%%%%%%%%%%%%%%%%%%%%%%%%%%%%%%%%%%%%%%%%%%%%%%%%%%%%%%%%%%%%%%%%%%%%%%%%%%%%%%%
%%%%%%%%%%%%%%%%%%%%%%%%%%%%%%%%%%%%%%%%%%%%%%%%%%%%%%%%%%%%%%%%%%%%%%%%%%%%%%%%%%%%%%%%%%%%%%%%%%%%%%%%%%%%%%%%%%%%%%%%%%%%%%%%
%%%%%%%%%%%%%%%%%%%%%%%%%%%%%%%%%%%%%%%%%%%%%%%%%%%%%%%%%%%%%%%%%%%%%%%%%%%%%%%%%%%%%%%%%%%%%%%%%%%%%%%%%%%%%%%%%%%%%%%%%%%%%%%%
%%%%%%%%%%%%%%%%%%%%%%%%%%%%%%%%%%%%%%%%%%%%%%%%%%%%%%%%%%%%%%%%%%%%%%%%%%%%%%%%%%%%%%%%%%%%%%%%%%%%%%%%%%%%%%%%%%%%%%%%%%%%%%%%
%%%%%%%%%%%%%%%%%%%%%%%%%%%%%%%%%%%%%%%%%%%%%%%%%%%%%%%%%%%%%%%%%%%%%%%%%%%%%%%%%%%%%%%%%%%%%%%%%%%%%%%%%%%%%%%%%%%%%%%%%%%%%%%%
\section{Tensor Fields on a Smooth Vector Bundle}
%%%%%%%%%%%%%%%%%%%%%%%%%%%%%%%%
\definition
Let $r$ and $s$ be non-negative integers, and $\U$ a subset of $\vB{}$.
\begin{equation}
\vbTensors{r}{s}{\vbundle{}}{\U}:=
\Union{\point}{\U}{\Tensors{r}{s}{\fibervecs{\vbundle{}}{\point}}}.
\end{equation}
Specifically, $\vbTensors{r}{s}{\vbundle{}}{\vB{}}$ can alternatively be denoted by $\VBTensors{r}{s}{\vbundle{}}$.
\endef
%%%%%%%%%%%%%%%%%%%%%%%%%%%%%%%%%%%%%%%%%%%%%%%%%%%%%%%%%%%%%%%%%%%%%%%%%%%%%%%%%%%%%%%%%%%%%%%%%%%%%%%%%%%%%%%%%%%%%%%%%%%%%%%%
\definition
\begin{equation}
\vbchartlocalt{\vbundle{}}:=
\defSet{\opair{\phi}{\psi}}{\[\phi\in\maxatlas{\vB{}},~
\psi\in\vbatlas{},~\domain{\phi}=\func{\image{\vbprojection{}}}{\domain{\psi}}\]}.
\end{equation}
An element $\opair{\phi}{\psi}$ of $\vbchartlocalt{\vbundle{}}$ can alternatively be denoted by
$\btriple{\U}{\phi}{\psi}$ where $\U:=\domain{\phi}=\func{\image{\vbprojection{}}}{\domain{\psi}}$.
\endef
%%%%%%%%%%%%%%%%%%%%%%%%%%%%%%%%%%%%%%%%%%%%%%%%%%%%%%%%%%%%%%%%%%%%%%%%%%%%%%%%%%%%%%%%%%%%%%%%%%%%%%%%%%%%%%%%%%%%%%%%%%%%%%%%
\lemma
The set $\defsets{\U}{\vB{}}{\[\Exists{\opair{\phi}{\psi}}{\vbchartlocalt{\vbundle{}}}\U=\domain{\phi}=
\func{\image{\vbprojection{}}}{\domain{\psi}}\]}$ is an open covering of $\vbbase{}$.
\proof
It is an immediate consequence of \reflem{lemlocaltrivializationcontraction}.
\endlem
%%%%%%%%%%%%%%%%%%%%%%%%%%%%%%%%%%%%%%%%%%%%%%%%%%%%%%%%%%%%%%%%%%%%%%%%%%%%%%%%%%%%%%%%%%%%%%%%%%%%%%%%%%%%%%%%%%%%%%%%%%%%%%%%
\definition
Let $r$ and $s$ be non-negative integers. We associate to each element $\btriple{\U}{\phi}{\psi}$ of
$\vbchartlocalt{\vbundle{}}$ the mapping
$\function{\vbtensorchart{r}{s}{\phi}{\psi}}{\vbTensors{r}{s}{\vbundle{}}{\U}}
{\Cprod{\func{\phi}{\U}}{\Tensors{r}{s}{\vbfiber{}}}}$ defined as,
\begin{align}
&\begin{aligned}
\Foreach{\point}{\U}
\Foreach{\alpha}{\Tensors{r}{s}{\fibervecs{\vbundle{}}{\point}}}
\end{aligned}\cr
&\begin{aligned}
\func{\vbtensorchart{r}{s}{\phi}{\psi}}{\alpha}\eqdef
\opair{\func{\phi}{\point}}{\func{\[\Vpullback{\(\finv{\(\pltfib{\psi}{\point}\)}\)}{r}{s}\]}{\alpha}}.
\end{aligned}
\end{align}
Furthermore, we associate to each element $\opair{\U}{\psi}$ of
$\vbatlas{}$ the mapping
$\function{\vbtensorlocalt{r}{s}{\psi}}{\vbTensors{r}{s}{\vbundle{}}{\U}}
{\Cprod{\U}{\Tensors{r}{s}{\vbfiber{}}}}$ defined as,
\begin{align}
&\begin{aligned}
\Foreach{\point}{\U}
\Foreach{\alpha}{\Tensors{r}{s}{\fibervecs{\vbundle{}}{\point}}}
\end{aligned}\cr
&\begin{aligned}
\func{\vbtensorlocalt{r}{s}{\psi}}{\alpha}\eqdef
\opair{\point}{\func{\[\Vpullback{\(\finv{\(\pltfib{\psi}{\point}\)}\)}{r}{s}\]}{\alpha}}.
\end{aligned}
\end{align}
\endef
%%%%%%%%%%%%%%%%%%%%%%%%%%%%%%%%%%%%%%%%%%%%%%%%%%%%%%%%%%%%%%%%%%%%%%%%%%%%%%%%%%%%%%%%%%%%%%%%%%%%%%%%%%%%%%%%%%%%%%%%%%%%%%%%
\textcolor{Blue}{\lemma
Let $r$ and $s$ be non-negative integers.
The set $\defSet{\vbtensorchart{r}{s}{\phi}{\psi}}{\opair{\phi}{\psi}\in\vbchartlocalt{\vbundle{}}}$
is a $\difclass{\infty}$ atlas on $\VBTensors{r}{s}{\vbundle{}}$ modeled on
the Banach-space $\Cprod{\R^{n_{\vB{}}}}{\MTensors{r}{s}{\vbfiber{}}}$.}
\endlem
%%%%%%%%%%%%%%%%%%%%%%%%%%%%%%%%%%%%%%%%%%%%%%%%%%%%%%%%%%%%%%%%%%%%%%%%%%%%%%%%%%%%%%%%%%%%%%%%%%%%%%%%%%%%%%%%%%%%%%%%%%%%%%%%
\definition
Let $r$ and $s$ be non-negative integers.
The $\difclass{\infty}$ maximal-atlas on $\VBTensors{r}{s}{\vbundle{}}$ modeled on
the Banach-space $\Cprod{\R^{n_{\vB{}}}}{\MTensors{r}{s}{\vbfiber{}}}$ generated by the atlas
$\defSet{\vbtensorchart{r}{s}{\phi}{\psi}}{\opair{\phi}{\psi}\in\vbchartlocalt{\vbundle{}}}$ will be denoted by
$\vbtensormaxatlas{\vbundle{}}$. That is,
\begin{equation}
\vbtensormaxatlas{\vbundle{}}:=
\maxatlasgen{\infty}{\VBTensors{r}{s}{\vbundle{}}}{\Cprod{\R^{n_{\vB{}}}}{\Tensors{r}{s}{\vbfiber{}}}}{
\defSet{\vbtensorchart{r}{s}{\phi}{\psi}}{\opair{\phi}{\psi}\in\vbchartlocalt{\vbundle{}}}}.
\end{equation}
\endef
%%%%%%%%%%%%%%%%%%%%%%%%%%%%%%%%%%%%%%%%%%%%%%%%%%%%%%%%%%%%%%%%%%%%%%%%%%%%%%%%%%%%%%%%%%%%%%%%%%%%%%%%%%%%%%%%%%%%%%%%%%%%%%%%
\textcolor{Blue}{\theorem
Let $r$ and $s$ be non-negative integers.\\
The differentiable structure $\opair{\VBTensors{r}{s}{\vbundle{}}}{\vbtensormaxatlas{\vbundle{}}}$
is a $\difclass{\infty}$ manifold, which means the topology induced by the maximal atlas $\vbtensormaxatlas{\vbundle{}}$
on $\VBTensors{r}{s}{\vbundle{}}$ is Hausdorff and second-countable.}
\endthm
%%%%%%%%%%%%%%%%%%%%%%%%%%%%%%%%%%%%%%%%%%%%%%%%%%%%%%%%%%%%%%%%%%%%%%%%%%%%%%%%%%%%%%%%%%%%%%%%%%%%%%%%%%%%%%%%%%%%%%%%%%%%%%%%
\definition
Let $r$ and $s$ be non-negative integers.
The manifold $\opair{\VBTensors{r}{s}{\vbundle{}}}{\vbtensormaxatlas{\vbundle{}}}$ is called the
$\quotl$$\opair{r}{s}$-tensor-bundle of the smooth vector bundle $\vbundle{}$$\quotr$,
which will be denoted by $\VBTensorsMan{r}{s}{\vbundle{}}$.
\endef
%%%%%%%%%%%%%%%%%%%%%%%%%%%%%%%%%%%%%%%%%%%%%%%%%%%%%%%%%%%%%%%%%%%%%%%%%%%%%%%%%%%%%%%%%%%%%%%%%%%%%%%%%%%%%%%%%%%%%%%%%%%%%%%%
\definition
Let $r$ and $s$ be non-negative integers.
The mapping $\function{\vbtensorprojection{\vbundle{}}{r}{s}}{\VBTensors{r}{s}{\vbundle{}}}{\vB{}}$ is
defined as,
\begin{equation}
\Foreach{\point}{\vB{}}
\Foreach{\alpha}{\Tensors{r}{s}{\fibervecs{\vbundle{}}{\point}}}
\func{\vbtensorprojection{\vbundle{}}{r}{s}}{\alpha}\eqdef\point.
\end{equation}
\endef
%%%%%%%%%%%%%%%%%%%%%%%%%%%%%%%%%%%%%%%%%%%%%%%%%%%%%%%%%%%%%%%%%%%%%%%%%%%%%%%%%%%%%%%%%%%%%%%%%%%%%%%%%%%%%%%%%%%%%%%%%%%%%%%%
\lemma
Let $r$ and $s$ be non-negative integers.
The quadruple
$\quadruple{\VBTensorsMan{r}{s}{\vbundle{}}}{\vbtensorprojection{\vbundle{}}{r}{s}}{\vbbase{}}
{\vecsmanifold{\MTensors{r}{s}{\vbfiber{}}}}$
is a smooth fiber bundle.
\proof
Let $\point$ be an arbitrary point of $\vbbase{}$. Let $\btriple{\U}{\phi}{\psi}$ be an element of
$\vbchartlocalt{\vbundle{}}$ that contains $\point$. Clearly $\funcprod{\finv{\phi}}{\identity{\Tensors{r}{s}{\vbfiber{}}}}$
is a diffeomorphism from $\Cprod{\func{\phi}{\U}}{\MTensors{r}{s}{\vbfiber{}}}$ (with its canonical differentiable structure
inherited from that of $\vecsmanifold{\Cprod{\R^{n_{\vB{}}}}{\MTensors{r}{s}{\vbfiber{}}}}$) to
$\Cprod{\subman{\vbbase{}}{\U}}{\vecsmanifold{\MTensors{r}{s}{\vbfiber{}}}}$. In addition, it is also trivial that
$\vbtensorchart{r}{s}{\phi}{\psi}$ is a diffeomorphism from
$\subman{\VBTensorsMan{r}{s}{\vbundle{}}}{\vbTensors{r}{s}{\vbundle{}}{\U}}$ to
$\Cprod{\func{\phi}{\U}}{\MTensors{r}{s}{\vbfiber{}}}$, because it is a chart of the manifold
$\VBTensorsMan{r}{s}{\vbundle{}}$.
Thus, considering that the composition of a pair of diffeomorphisms as again a diffeomorphism,
$\vbtensorlocalt{r}{s}{\psi}=
\cmp{\(\funcprod{\finv{\phi}}{\identity{\Tensors{r}{s}{\vbfiber{}}}}\)}{\vbtensorchart{r}{s}{\phi}{\psi}}$
is a diffeomorphism from $\subman{\VBTensorsMan{r}{s}{\vbundle{}}}{\vbTensors{r}{s}{\vbundle{}}{\U}}$ to
$\Cprod{\subman{\vbbase{}}{\U}}{\vecsmanifold{\MTensors{r}{s}{\vbfiber{}}}}$.
Moreover, it is evident that
$\func{\pimage{\(\vbtensorprojection{\vbundle{}}{r}{s}\)}}{\U}=
\vbTensors{r}{s}{\vbundle{}}{\U}$, and the following diagram is commutative.
\begin{center}
\vskip0.5\baselineskip
\hskip-2\baselineskip
\begin{tikzcd}[row sep=6em, column sep=6em]
& \vbTensors{r}{s}{\vbundle{}}{\U}
\arrow{r}{\vbtensorlocalt{r}{s}{\psi}}
\arrow[swap]{d}{\vbtensorprojection{\vbundle{}}{r}{s}}
& \Cprod{\U}{\Tensors{r}{s}{\vbfiber{}}}
\arrow{dl}{\proj{\U}{\Tensors{r}{s}{\vbfiber{}}}{1}} \\
& \U
\end{tikzcd}
\end{center}
\endlem
%%%%%%%%%%%%%%%%%%%%%%%%%%%%%%%%%%%%%%%%%%%%%%%%%%%%%%%%%%%%%%%%%%%%%%%%%%%%%%%%%%%%%%%%%%%%%%%%%%%%%%%%%%%%%%%%%%%%%%%%%%%%%%%%
\lemma
Let $r$ and $s$ be non-negative integers.
The set $\defSet{\vbtensorlocalt{r}{s}{\psi}}{\psi\in\vbatlas{}}$ is an atlas of the smooth fiber bundle
$\quadruple{\VBTensorsMan{r}{s}{\vbundle{}}}{\vbtensorprojection{\vbundle{}}{r}{s}}{\vbbase{}}
{\vecsmanifold{\MTensors{r}{s}{\vbfiber{}}}}$.
\proof
Let $\opair{\U}{\psi}$ be an arbitrary local trivialization of the smooth vector bundle
$\vbundle{}$.
Clearly, $\func{\pimage{\(\vbtensorprojection{\vbundle{}}{r}{s}\)}}{\U}=
\vbTensors{r}{s}{\vbundle{}}{\U}$ and the diagram in the proof of the previous lemma is again commutative. There just remains
to show that $\vbtensorlocalt{r}{s}{\psi}$ is a diffeomorphism from
$\subman{\VBTensorsMan{r}{s}{\vbundle{}}}{\vbTensors{r}{s}{\vbundle{}}{\U}}$ to
$\Cprod{\subman{\vbbase{}}{\U}}{\vecsmanifold{\MTensors{r}{s}{\vbfiber{}}}}$. This is true because according
to the previous lemma,there exists an open covering of $\U$ that the restriction of
$\vbtensorlocalt{r}{s}{\psi}$ to each of the elements of that covering is a diffeomorphism.\\
Thus, for every local trivialization $\opair{\U}{\psi}$ of $\vbundle{}$, $\opair{\U}{\vbtensorlocalt{r}{s}{\psi}}$ is
a local trivialization of the  smooth fiber bundle
$\quadruple{\VBTensorsMan{r}{s}{\vbundle{}}}{\vbtensorprojection{\vbundle{}}{r}{s}}{\vbbase{}}
{\vecsmanifold{\MTensors{r}{s}{\vbfiber{}}}}$. Therefore,
$\defSet{\vbtensorlocalt{r}{s}{\psi}}{\psi\in\vbatlas{}}$ is an atlas of the smooth fiber bundle
$\quadruple{\VBTensorsMan{r}{s}{\vbundle{}}}{\vbtensorprojection{\vbundle{}}{r}{s}}{\vbbase{}}
{\vecsmanifold{\MTensors{r}{s}{\vbfiber{}}}}$.
\endlem
%%%%%%%%%%%%%%%%%%%%%%%%%%%%%%%%%%%%%%%%%%%%%%%%%%%%%%%%%%%%%%%%%%%%%%%%%%%%%%%%%%%%%%%%%%%%%%%%%%%%%%%%%%%%%%%%%%%%%%%%%%%%%%%%
\proposition\label{provbtensorlocaltplt}
Let $r$ and $s$ be non-negative integers.
Let $\opair{\U}{\psi}$ be a local trivialization of the smooth vector bundle $\vbundle{}$,
and $\point$ a point of $\U$. The principal part
of the local trivialization $\vbtensorlocalt{r}{s}{\psi}$ (of the smooth fiber bundle
$\quadruple{\VBTensorsMan{r}{s}{\vbundle{}}}{\vbtensorprojection{\vbundle{}}{r}{s}}{\vbbase{}}
{\vecsmanifold{\MTensors{r}{s}{\vbfiber{}}}}$) when restricted to
$\func{\finv{\(\vbtensorprojection{\vbundle{}}{r}{s}\)}}{\seta{\point}}=
\Tensors{r}{s}{\fibervecs{\vbundle{}}{\point}}$ coincides with
$\Vpullback{\(\finv{\(\pltfib{\psi}{\point}\)}\)}{r}{s}$. That is,
\begin{equation}
\func{\res{\plt{}{\vbtensorlocalt{r}{s}{\psi}}}}{\func{\finv{\(\vbtensorprojection{\vbundle{}}{r}{s}\)}}{\seta{\point}}}=
\pltfib{\vbtensorlocalt{r}{s}{\psi}}{\point}=
\Vpullback{\(\finv{\(\pltfib{\psi}{\point}\)}\)}{r}{s}.
\end{equation}
\proof
It is trivial according to the definition of $\vbtensorlocalt{r}{s}{\psi}$.
\endpro
%%%%%%%%%%%%%%%%%%%%%%%%%%%%%%%%%%%%%%%%%%%%%%%%%%%%%%%%%%%%%%%%%%%%%%%%%%%%%%%%%%%%%%%%%%%%%%%%%%%%%%%%%%%%%%%%%%%%%%%%%%%%%%%%
\proposition
Let $r$ and $s$ be non-negative integers.
For every pair $\opair{\U}{\psi}$ and $\opair{\V}{\eta}$ of local trivialization of the smooth vector bundle $\vbundle{}$,
$\function{\cmp{\vbtensorlocalt{r}{s}{\psi}}{\finv{\(\vbtensorlocalt{r}{s}{\eta}\)}}}
{\Cprod{\(\U\cap\V\)}{\Tensors{r}{s}{\vbfiber{}}}}{\Cprod{\(\U\cap\V\)}{\Tensors{r}{s}{\vbfiber{}}}}$, and
the value of the transition map of the local trivializations $\vbtensorlocalt{r}{s}{\psi}$ and $\vbtensorlocalt{r}{s}{\eta}$
of the fiber bundle $\quadruple{\VBTensorsMan{r}{s}{\vbundle{}}}{\vbtensorprojection{\vbundle{}}{r}{s}}{\vbbase{}}
{\vecsmanifold{\MTensors{r}{s}{\vbfiber{}}}}$ at every point $\point$ of its domain $\U\cap\V$ is a linear isomorphism
from $\MTensors{r}{s}{\vbfiber{}}$ to itself. That is,
\begin{align}
\Foreach{\point}{\U\cap\V}
\func{\[\transition{}{\vbtensorlocalt{r}{s}{\psi}}{\vbtensorlocalt{r}{s}{\eta}}\]}{\point}\in
\GL{\MTensors{r}{s}{\vbfiber{}}}{}.
\end{align}
\proof
Let $\opair{\U}{\phi}$ and $\opair{\V}{\psi}$ be an arbitrary pair of local trivializations of $\vbundle{}$,
and $\point$ an arbitrary element of $\U\cap\V$.
According to the definition of the transition map of a pair of local trivializations of a smooth fiber bundle,
\refpro{provbtensorlocaltplt}, and \refthm{thminverseofpullbackofrstensor},
\begin{align}
\func{\[\transition{}{\vbtensorlocalt{r}{s}{\psi}}{\vbtensorlocalt{r}{s}{\eta}}\]}{\point}&=
\cmp{\pltfib{\vbtensorlocalt{r}{s}{\psi}}{\point}}{\finv{\(\pltfib{\vbtensorlocalt{r}{s}{\eta}}{\point}\)}}\cr
&=\cmp{\Vpullback{\(\finv{\(\pltfib{\psi}{\point}\)}\)}{r}{s}}
{\finv{\[\Vpullback{\(\finv{\(\pltfib{\eta}{\point}\)}\)}{r}{s}\]}}\cr
&=\cmp{\Vpullback{\(\finv{\(\pltfib{\psi}{\point}\)}\)}{r}{s}}
{\Vpullback{\(\pltfib{\eta}{\point}\)}{r}{s}}.
\end{align}
Since $\pltfib{\eta}{\point}\in\Linisom{\fibervecs{\vbundle{}}{\point}}{\vbfiber{}}$,
clearly
$\Vpullback{\(\pltfib{\eta}{\point}\)}{r}{s}\in
\Linisom{\MTensors{r}{s}{\vbfiber{}}}{\MTensors{r}{s}{\fibervecs{\vbundle{}}{\point}}}$.
Similarlt, it is evident that
$\Vpullback{\(\finv{\(\pltfib{\psi}{\point}\)}\)}{r}{s}\in
\Linisom{\MTensors{r}{s}{\fibervecs{\vbundle{}}{\point}}}{\MTensors{r}{s}{\vbfiber{}}}$. Therefore,
\begin{equation}
\func{\[\transition{}{\vbtensorlocalt{r}{s}{\psi}}{\vbtensorlocalt{r}{s}{\eta}}\]}{\point}=
\cmp{\Vpullback{\(\finv{\(\pltfib{\psi}{\point}\)}\)}{r}{s}}
{\Vpullback{\(\pltfib{\eta}{\point}\)}{r}{s}}\in\GL{\MTensors{r}{s}{\vbfiber{}}}{}.
\end{equation}
\endlem
%%%%%%%%%%%%%%%%%%%%%%%%%%%%%%%%%%%%%%%%%%%%%%%%%%%%%%%%%%%%%%%%%%%%%%%%%%%%%%%%%%%%%%%%%%%%%%%%%%%%%%%%%%%%%%%%%%%%%%%%%%%%%%%%
\definition
Let $r$ and $s$ be non-negative integers. We will denote by $\vbtensoratlas{r}{s}{\vbundle{}}$
the maximal atlas of the smooth fiber bundle
$\quadruple{\VBTensorsMan{r}{s}{\vbundle{}}}{\vbtensorprojection{\vbundle{}}{r}{s}}{\vbbase{}}
{\vecsmanifold{\MTensors{r}{s}{\vbfiber{}}}}$
including $\defSet{\vbtensorlocalt{r}{s}{\psi}}{\psi\in\vbatlas{}}$, endowed with which, the fiber bundle
$\quadruple{\VBTensorsMan{r}{s}{\vbundle{}}}{\vbtensorprojection{\vbundle{}}{r}{s}}{\vbbase{}}
{\vecsmanifold{\MTensors{r}{s}{\vbfiber{}}}}$ becomes a smooth vector bundle.
\endef
%%%%%%%%%%%%%%%%%%%%%%%%%%%%%%%%%%%%%%%%%%%%%%%%%%%%%%%%%%%%%%%%%%%%%%%%%%%%%%%%%%%%%%%%%%%%%%%%%%%%%%%%%%%%%%%%%%%%%%%%%%%%%%%%
\corollary
Let $r$ and $s$ be non-negative integers.\\
The quintuple
$\quintuple{\VBTensorsMan{r}{s}{\vbundle{}}}{\vbtensorprojection{\vbundle{}}{r}{s}}{\vbbase{}}
{\vecsmanifold{\MTensors{r}{s}{\vbfiber{}}}}{\vbtensoratlas{r}{s}{\vbundle{}}}$ is a smooth vector bundle.
\endcor
%%%%%%%%%%%%%%%%%%%%%%%%%%%%%%%%%%%%%%%%%%%%%%%%%%%%%%%%%%%%%%%%%%%%%%%%%%%%%%%%%%%%%%%%%%%%%%%%%%%%%%%%%%%%%%%%%%%%%%%%%%%%%%%%
\definition
Let $r$ and $s$ be non-negative integers.
We will denote by $\vbtensorbundle{r}{s}{\vbundle{}}$ the smooth vector bundle
$\quintuple{\VBTensorsMan{r}{s}{\vbundle{}}}{\vbtensorprojection{\vbundle{}}{r}{s}}{\vbbase{}}
{\vecsmanifold{\MTensors{r}{s}{\vbfiber{}}}}{\vbtensoratlas{r}{s}{\vbundle{}}}$, which is referred to as the
$\quotl$$\opair{r}{s}$ tensor bundle of the smooth vector bundle $\vbundle{}$$\quotr$.
\endef
%%%%%%%%%%%%%%%%%%%%%%%%%%%%%%%%%%%%%%%%%%%%%%%%%%%%%%%%%%%%%%%%%%%%%%%%%%%%%%%%%%%%%%%%%%%%%%%%%%%%%%%%%%%%%%%%%%%%%%%%%%%%%%%%
\corollary
Let $r$ and $s$ be non-negative integers.
For every point $\point$ of $\vbbase{}$, the fiber space of $\vbtensorbundle{r}{s}{\vbundle{}}$ over $\point$ equals
the vector-space of all $\opair{r}{s}$ tensors on the fiber space of $\vbundle{}$ over $\point$. That is,
\begin{align}
\Foreach{\point}{\vB{}}
\fibervecs{\vbtensorbundle{r}{s}{\vbundle{}}}{\point}=
\MTensors{r}{s}{\fibervecs{\vbundle{}}{\point}}.
\end{align}
\endcor
%%%%%%%%%%%%%%%%%%%%%%%%%%%%%%%%%%%%%%%%%%%%%%%%%%%%%%%%%%%%%%%%%%%%%%%%%%%%%%%%%%%%%%%%%%%%%%%%%%%%%%%%%%%%%%%%%%%%%%%%%%%%%%%%
\definition
Let $r$ and $s$ be non-negative integers. Any section of the $\opair{r}{s}$ tensor bundle of $\vbundle{}$
is called a $\quotl$$\opair{r}{s}$ tensor field of the smooth vector bundle $\vbundle{}$$\quotr$.
The set of all $\opair{r}{s}$ tensor fields of $\vbundle{}$ will be denoted by $\TF{r}{s}{\vbundle{}}$. That is,
\begin{equation}
\TF{r}{s}{\vbundle{}}:=\vbsections{\vbtensorbundle{r}{s}{\vbundle{}}}.
\end{equation}
The canonical $\smoothring{\infty}{\vbbase{}}$-module strucure associated with $\TF{r}{s}{\vbundle{}}$
will be denoted by\\
$\MTF{r}{s}{\vbundle{}}$, and the canonical linear structure associated with
$\TF{r}{s}{\vbundle{}}$ will be denoted by $\VTF{r}{s}{\vbundle{}}$.
\endef
%%%%%%%%%%%%%%%%%%%%%%%%%%%%%%%%%%%%%%%%%%%%%%%%%%%%%%%%%%%%%%%%%%%%%%%%%%%%%%%%%%%%%%%%%%%%%%%%%%%%%%%%%%%%%%%%%%%%%%%%%%%%%%%%
%%%%%%%%%%%%%%%%%%%%%%%%%%%%%%%%%%%%%%%%%%%%%%%%%%%%%%%%%%%%%%%%%%%%%%%%%%%%%%%%%%%%%%%%%%%%%%%%%%%%%%%%%%%%%%%%%%%%%%%%%%%%%%%%
%%%%%%%%%%%%%%%%%%%%%%%%%%%%%%%%%%%%%%%%%%%%%%%%%%%%%%%%%%%%%%%%%%%%%%%%%%%%%%%%%%%%%%%%%%%%%%%%%%%%%%%%%%%%%%%%%%%%%%%%%%%%%%%%
\subsection{Tensor Algebra of Tensor Fields on Vector Bundles}
%%%%%%%%%%%%%%%%%%%%%%%%%%%%%%%%%%%%%%%%%%%%
\definition
The direct sum of all $\VTF{r}{s}{\vbundle{}}$-s when $r$ and $s$ range over all non-negative integers simultaneously,
will simply be denoted by $\DVTF{\vbundle{}}$. That is,
\begin{equation}
\DVTF{\vbundle{}}:=\Dsum{\defSet{\VTF{r}{s}{\vbundle{}}}{r,s=0,1,\ldots}}.
\end{equation}
$\DVTF{\vbundle{}}$ is referred to as the $\quotl$space of tensor fields on the smooth vector bundle $\vbundle{}$$\quotr$.
For every pair $r$ and $s$ of non-negative integers, there is a natural injection of $\TF{r}{s}{\vbundle{}}$ into
$\DVTF{\CModule{}}$ via which every $\opair{r}{s}$ tensor field on $\vbundle{}$ can be identified with a unique element of
$\DVTF{\vbundle{}}$.
\endef
%%%%%%%%%%%%%%%%%%%%%%%%%%%%%%%%%%%%%%%%%%%%%%%%%%%%%%%%%%%%%%%%%%%%%%%%%%%%%%%%%%%%%%%%%%%%%%%%%%%%%%%%%%%%%%%%%%%%%%%%%%%%%%%%
\definition
The binary operation $\vbtensor{\vbundle{}}$ on $\DVTF{\vbundle{}}$ is defined to be the unique bilinear map
$\function{\vbtensor{\vbundle{}}}{\Cprod{\DVTF{\vbundle{}}}{\DVTF{\vbundle{}}}}{\DVTF{\vbundle{}}}$
such that for every non-negative integers $r$, $s$, $p$, $q$,
%if $r+s>0$ and $p+q>0$,
%%$\func{\image{\tensor{\CModule{}}}}{}\subseteq
\begin{align}
&\Foreach{\opair{\atf{}}{\atff{}}}{\Cprod{\TF{r}{s}{\vbundle{}}}{\TF{p}{q}{\vbundle{}}}}
\Foreach{\point}{\vB{}}
\func{\(\atf{}\vbtensor{\vbundle{}}\atff{}\)}{\point}\eqdef
\func{\atf{}}{\point}\tensor{}\func{\atff{}}{\point},
\end{align}
where the tensor product $\func{\atf{}}{\point}\tensor{}\func{\atff{}}{\point}$ takes place
in the  space of tensors on the fiber space of $\vbundle{}$ over $\point$, that is $\DTensors{\fibervecs{\vbundle{}}{\point}}$.
%%%%%%%%%
$\vbtensor{\vbundle{}}$ is called the $\quotl$tensor (or tensor field) operation of the smooth vector bundle $\vbundle{}$$\quotr$.
For any non-negative integers $r$ and $s$, each element of $\TF{r}{s}{\vbundle{}}$ is also called a
$\quotl$(type) $\opair{r}{s}$ simple tensor field on $\vbundle{}$$\quotr$. Also for every $\atf{}$ and $\atff{}$
in $\DVTF{\vbundle{}}$, $\atf{}\vbtensor{\vbundle{}}\atff{}$ is called the $\quotl$tensor product of the tensor fields $\atf{}$
and $\atff{}$$\quotr$.
When there is no ambiguity about the underlying smooth vector bundle, $\vbtensor{\vbundle{}}$ can
simply be denoted by $\vbtensor{}$.\\
\caution
It is easy to verify that for any $\atf{}\in\TF{r}{s}{\vbundle{}}\subseteq\DVTF{\vbundle{}}$ and any
$\atff{}\in\TF{p}{q}{\vbundle{}}\subseteq\DVTF{\vbundle{}}$, we have
$\atf{}\vbtensor{\vbundle{}}\atff{}\in\TF{r+p}{q+s}{\vbundle{}}\subseteq\DVTF{\vbundle{}}$ and $\vbtensor{\vbundle{}}$
acts bilinearly on $\Cprod{\TF{r}{s}{\vbundle{}}}{\TF{p}{q}{\vbundle{}}}$. We presumed the triviality of
these facts prior to the definition.
\endef
%%%%%%%%%%%%%%%%%%%%%%%%%%%%%%%%%%%%%%%%%%%%%%%%%%%%%%%%%%%%%%%%%%%%%%%%%%%%%%%%%%%%%%%%%%%%%%%%%%%%%%%%%%%%%%%%%%%%%%%%%%%%%%%%
\proposition
$\vbtensor{\vbundle{}}$ is an associative and bilinear binary operation on $\DVTF{\vbundle{}}$.
So the pair $\opair{\DVTF{\vbundle{}}}{\vbtensor{\vbundle{}}}$ is an algebra.
\proof
The bilinearity of $\vbtensor{\vbundle{}}$ lies in the definition of it.\\
Now let $\atf{1}\in\TF{r}{s}{\vbundle{}}$, $\atf{2}\in\TF{p}{q}{\vbundle{}}$, and $\atf{3}\in\TF{l}{m}{\vbundle{}}$.
According to \refthm{thmtensoroperationofmoduleisassociative},
\begin{align}
\Foreach{\point}{\vB{}}
\func{\[\atf{1}\vbtensor{}\(\atf{2}\vbtensor{}\atf{3}\)\]}{\point}&=
\func{\atf{1}}{\point}\tensor{}\(\func{\atf{2}}{\point}\tensor{}\func{\atf{3}}{\point}\)\cr
&=\(\func{\atf{1}}{\point}\tensor{}\func{\atf{2}}{\point}\)\tensor{}\func{\atf{3}}{\point}\cr
&=\func{\[\(\atf{1}\vbtensor{}\atf{2}\)\vbtensor{}\atf{3}\]}{\point},
\end{align}
and thus,
\begin{equation}
\atf{1}\vbtensor{}\(\atf{2}\vbtensor{}\atf{3}\)=
\(\atf{1}\vbtensor{}\atf{2}\)\vbtensor{}\atf{3}.
\end{equation}
So $\vbtensor{\vbundle{}}$ is associative when restricted to simple tensor fields. The verification of associativity
for general tensor fields in $\DVTF{\vbundle{}}$ can be achieved in a manner completely similar to what is used in
the proof of \refthm{thmtensoroperationofmoduleisassociative}.
\endpro
%%%%%%%%%%%%%%%%%%%%%%%%%%%%%%%%%%%%%%%%%%%%%%%%%%%%%%%%%%%%%%%%%%%%%%%%%%%%%%%%%%%%%%%%%%%%%%%%%%%%%%%%%%%%%%%%%%%%%%%%%%%%%%%%
\definition
$\opair{\DVTF{\vbundle{}}}{\vbtensor{\vbundle{}}}$ is called the
$\quotl$tensor algebra of tensor fields on the smooth vector bundle $\vbundle{}$$\quotr$.
\endef
%%%%%%%%%%%%%%%%%%%%%%%%%%%%%%%%%%%%%%%%%%%%%%%%%%%%%%%%%%%%%%%%%%%%%%%%%%%%%%%%%%%%%%%%%%%%%%%%%%%%%%%%%%%%%%%%%%%%%%%%%%%%%%%%
%%%%%%%%%%%%%%%%%%%%%%%%%%%%%%%%%%%%%%%%%%%%%%%%%%%%%%%%%%%%%%%%%%%%%%%%%%%%%%%%%%%%%%%%%%%%%%%%%%%%%%%%%%%%%%%%%%%%%%%%%%%%%%%%
%%%%%%%%%%%%%%%%%%%%%%%%%%%%%%%%%%%%%%%%%%%%%%%%%%%%%%%%%%%%%%%%%%%%%%%%%%%%%%%%%%%%%%%%%%%%%%%%%%%%%%%%%%%%%%%%%%%%%%%%%%%%%%%%
%%%%%%%%%%%%%%%%%%%%%%%%%%%%%%%%%%%%%%%%%%%%%%%%%%%%%%%%%%%%%%%%%%%%%%%%%%%%%%%%%%%%%%%%%%%%%%%%%%%%%%%%%%%%%%%%%%%%%%%%%%%%%%%%
\section{Smooth Vector Bundle Morphisms}
%%%%%%%%%%%%%%%%%%%%%%%%%%%%%%%%%%%%%%%%%%
\definition\label{defvbmorphism}
Let $f$ be an element of $\fbmorphisms{\Fvbundle{\vbundle{1}}}{\Fvbundle{\vbundle{2}}}$, that is
a smooth fiber bundle morphism from the underlying smooth fiber bundle structure of $\vbundle{1}$
to the underlying smooth fiber bundle structure of $\vbundle{2}$. $f$ is referred to as a
$\quotl$smooth vector bundle morphism from the smooth vector bundle $\vbundle{1}$
to the smooth vector bundle $\vbundle{2}$$\quotr$ iff for every point $\point$ of the base space of $\vbundle{1}$,
the restriction of $f$ to the fiber space of $\vbundle{1}$ over $\point$ is a linear map from
$\fibervecs{\vbundle{1}}{\point}$ to $\fibervecs{\vbundle{2}}{\func{\fbmorb{f}}{\point}}$,
where $\fbmorb{f}$ denotes the unique smooth map from $\vbbase{1}$ to $\vbbase{2}$ such that
$\cmp{\vbprojection{2}}{f}=\cmp{\fbmorb{f}}{\vbprojection{1}}$.\\
The set of all smooth vector bundle morphisms from $\vbundle{1}$ to $\vbundle{2}$ will be denoted by
$\vbmorphisms{\vbundle{1}}{\vbundle{2}}$, that is,
\begin{align}
&~\vbmorphisms{\vbundle{1}}{\vbundle{2}}\cr
:=&~\defset{f}{\fbmorphisms{\Fvbundle{\vbundle{1}}}{\Fvbundle{\vbundle{2}}}}
{\[\Foreach{\point}{\vbbase{1}}\func{\res{f}}{\func{\pimage{\vbprojection{1}}}{\seta{\point}}}\in
\Lin{\fibervecs{\vbundle{1}}{\point}}{\fibervecs{\vbundle{2}}{\func{\fbmorb{f}}{\point}}}\]}.\cr
&{}
\end{align}
\endef
%%%%%%%%%%%%%%%%%%%%%%%%%%%%%%%%%%%%%%%%%%%%%%%%%%%%%%%%%%%%%%%%%%%%%%%%%%%%%%%%%%%%%%%%%%%%%%%%%%%%%%%%%%%%%%%%%%%%%%%%%%%%%%%%
\proposition
$\identity{\vTot{}}$ is a morphism from $\vbundle{}$ to itself, that is
$\identity{\vTot{}}\in\vbmorphisms{\vbundle{}}{\vbundle{}}$. Furtheremore,
$\fbmorb{\identity{\vTot{}}}=\identity{\vB{}}$.
\proof
It is trivial.
\endpro
%%%%%%%%%%%%%%%%%%%%%%%%%%%%%%%%%%%%%%%%%%%%%%%%%%%%%%%%%%%%%%%%%%%%%%%%%%%%%%%%%%%%%%%%%%%%%%%%%%%%%%%%%%%%%%%%%%%%%%%%%%%%%%%%
\proposition
Let $f_1$ be a morphism from $\vbundle{1}$ to $\vbundle{2}$, and $f_2$ be a morphism from $\vbundle{2}$ to $\vbundle{3}$.
$\cmp{f_2}{f_1}$ is a morphism from $\vbundle{1}$ to $\vbundle{3}$ along $\cmp{\fbmorb{f_2}}{\fbmorb{f_1}}$.
\proof
According to \refpro{procompositionoffbmorphisms}, $\cmp{f_2}{f_1}$ is a smooth fiber bundle morphism from
$\Fvbundle{\vbundle{1}}$ to $\Fvbundle{\vbundle{3}}$, and $\fbmorb{\cmp{f_2}{f_1}}=\cmp{\fbmorb{f_2}}{\fbmorb{f_1}}$.\\
Additionally, considering the definition of smooth vector bundle morphisms,
$\reS{f_1}{\func{\pimage{\vbprojection{1}}}{\seta{\point}}}\in
\Lin{\fibervecs{\vbundle{1}}{\point}}{\fibervecs{\vbundle{2}}{\func{\fbmorb{f_1}}{\point}}}$ and
$\reS{f_2}{\func{\pimage{\vbprojection{2}}}{\seta{\func{\fbmorb{f_1}}{\point}}}}\in
\Lin{\fibervecs{\vbundle{2}}{\func{\fbmorb{f_1}}{\point}}}{\fibervecs{\vbundle{3}}
{\func{\(\cmp{\fbmorb{f_2}}{\fbmorb{f_1}}\)}{\point}}}$ for every point $\point$ of $\vbbase{1}$.
Therefore, considering that composition of a pair of linear maps is a linear map, and further
considering \refcor{corfbmorphismmapsfiberintofiber} (that is,
$\func{\image{f_1}}{\func{\pimage{\fbprojection{1}}}{\seta{\point}}}\subseteq
\func{\pimage{\fbprojection{2}}}{\seta{\func{\fbmorb{f_1}}{\point}}}$ for every $\point\in\vB{1}$),
it becomes evident that for every point $\point$ of $\vbbase{1}$,
\begin{align}
\reS{\(\cmp{f_2}{f_1}\)}{\func{\pimage{\vbprojection{1}}}{\seta{\point}}}=
\cmp{\(\reS{f_2}{\func{\pimage{\vbprojection{2}}}{\seta{\func{\fbmorb{f_1}}{\point}}}}\)}
{\(\reS{f_1}{\func{\pimage{\vbprojection{1}}}{\seta{\point}}}\)}
\in\Lin{\fibervecs{\vbundle{1}}{\point}}{\fibervecs{\vbundle{3}}
{\func{\(\fbmorb{\cmp{f_2}{f_1}}\)}{\point}}}.
\end{align}
Therefore, $\cmp{f_2}{f_1}\in\vbmorphisms{\vbundle{1}}{\vbundle{3}}$.
\endpro
%%%%%%%%%%%%%%%%%%%%%%%%%%%%%%%%%%%%%%%%%%%%%%%%%%%%%%%%%%%%%%%%%%%%%%%%%%%%%%%%%%%%%%%%%%%%%%%%%%%%%%%%%%%%%%%%%%%%%%%%%%%%%%%%
\remark
\textit{According to these propositions, it is inferred that any collection of smooth vector bundles along with the set of
all possible morphisms between them forms a category,
having the ordinary composition of morphisms as the composition rule of its arrows.
Such a category goes under the name of a $\quotl$category of smooth vector bundles$\quotr$.}
\endremark
%%%%%%%%%%%%%%%%%%%%%%%%%%%%%%%%%%%%%%%%%%%%%%%%%%%%%%%%%%%%%%%%%%%%%%%%%%%%%%%%%%%%%%%%%%%%%%%%%%%%%%%%%%%%%%%%%%%%%%%%%%%%%%%%
\definition\label{defvbisomorphism}
Let $f$ be a morphism from $\vbundle{1}$ to $\vbundle{2}$. $f$ is referred to as an $\quotl$(smooth vector bundle)
isomorphism from $\vbundle{1}$ to $\vbundle{2}$$\quotr$ iff there exists a morphism $g$ from $\vbundle{2}$ to
$\vbundle{1}$ such that $\cmp{f}{g}=\identity{\B{2}}$, and $\cmp{g}{f}=\identity{\B{1}}$.\\
The set of all isomorphisms from $\vbundle{1}$ to $\vbundle{2}$ will be denoted by $\vbisomorphisms{\fbundle{1}}{\fbundle{2}}$.\\
When $\vbisomorphisms{\fbundle{1}}{\fbundle{2}}\neq\empty$, the smooth vector bundles $\vbundle{1}$ and $\vbundle{2}$
are said to be $\quotl$isomorphic (in the sense of smooth vector bundles)$\quotr$.
\endef
%%%%%%%%%%%%%%%%%%%%%%%%%%%%%%%%%%%%%%%%%%%%%%%%%%%%%%%%%%%%%%%%%%%%%%%%%%%%%%%%%%%%%%%%%%%%%%%%%%%%%%%%%%%%%%%%%%%%%%%%%%%%%%%%
\theorem
Let $f$ be a morphism from $\vbundle{1}$ to $\vbundle{2}$. $f$ is an isomorphism from $\vbundle{1}$ to $\vbundle{2}$
if and only if $f$ is a diffeomorphism from $\vbtotal{1}$ to $\vbtotal{2}$ and $\fbmorb{f}$ is a diffeomorphism from
$\vbbase{1}$ to $\vbbase{2}$. That is,
\begin{equation}
\vbisomorphisms{\vbundle{1}}{\vbundle{2}}=
\defset{f}{\vbmorphisms{\vbundle{1}}{\vbundle{2}}}{\(f\in\Diffeo{\infty}{\vbtotal{1}}{\vbtotal{2}},~
\fbmorb{f}\in\Diffeo{\infty}{\vbbase{1}}{\vbbase{2}}\)}.
\end{equation}
\proof
It is trivial according to \refthm{thmfbisomorphisms}, \refdef{defvbmorphism}, and \refdef{defvbisomorphism}.
\endthm
%%%%%%%%%%%%%%%%%%%%%%%%%%%%%%%%%%%%%%%%%%%%%%%%%%%%%%%%%%%%%%%%%%%%%%%%%%%%%%%%%%%%%%%%%%%%%%%%%%%%%%%%%%%%%%%%%%%%%%%%%%%%%%%%
\corollary
Let $f$ be an isomorphism from $\vbundle{1}$ to $\vbundle{2}$. For every point $\point$ of $\vbbase{1}$,
the restriction of $f$ to $\func{\pimage{\vbprojection{1}}}{\seta{\point}}$ is a linear isomorphism from
$\fibervecs{\vbundle{1}}{\point}$ to $\fibervecs{\vbundle{2}}{\func{f}{\point}}$. That is,
\begin{equation}
\Foreach{\point}{\vB{1}}
\func{\res{f}}{\func{\pimage{\vbprojection{1}}}{\seta{\point}}}\in
\Linisom{\fibervecs{\vbundle{1}}{\point}}{\fibervecs{\vbundle{2}}{\func{f}{\point}}}.
\end{equation}
\proof
It is an immediate consequence of \refdef{defvbisomorphism}.
\endcor
%%%%%%%%%%%%%%%%%%%%%%%%%%%%%%%%%%%%%%%%%%%%%%%%%%%%%%%%%%%%%%%%%%%%%%%%%%%%%%%%%%%%%%%%%%%%%%%%%%%%%%%%%%%%%%%%%%%%%%%%%%%%%%%%
%%%%%%%%%%%%%%%%%%%%%%%%%%%%%%%%%%%%%%%%%%%%%%%%%%%%%%%%%%%%%%%%%%%%%%%%%%%%%%%%%%%%%%%%%%%%%%%%%%%%%%%%%%%%%%%%%%%%%%%%%%%%%%%%
%%%%%%%%%%%%%%%%%%%%%%%%%%%%%%%%%%%%%%%%%%%%%%%%%%%%%%%%%%%%%%%%%%%%%%%%%%%%%%%%%%%%%%%%%%%%%%%%%%%%%%%%%%%%%%%%%%%%%%%%%%%%%%%%
\section{Pullbacks of Vector Bundle Morphisms}
%%%%%%%%%%%%%%%%%%%%%%%
\proposition
Let $r$ and $s$ be non-negative integers. Suppose that $\vbundle{1}$ and $\vbundle{2}$ are isomorphic, and
let $f$ be an element of $\vbisomorphisms{\vbundle{1}}{\vbundle{2}}$ (a smooth vector bundle isomorphism from
$\vbundle{1}$ to $\vbundle{2}$). The image of the assignment
$\TF{r}{s}{\vbundle{2}}\ni\atf{}\mapsto\bar{\atf{}}\in\Func{\vB{1}}{\VBTensors{r}{s}{\vbundle{}}}$ such that
\begin{align}
\Foreach{\atf{}}{\TF{r}{s}{\vbundle{2}}}
\Foreach{\point}{\vB{1}}
\func{\bar{\atf{}}}{\point}\eqdef
\func{\Vpullback{\(\reS{f}{\func{\pimage{\vbprojection{1}}}{\seta{\point}}}\)}{r}{s}}{\func{\atf{}}{\func{f}{\point}}},
\end{align}
is included in $\TF{r}{s}{\vbundle{1}}$.
In other words, for every $\atf{}\in\TF{r}{s}{\vbundle{2}}$,
$\bar{\atf{}}$ lies in $\TF{r}{s}{\vbundle{1}}$.
\proof
It is left to the reader as an exercise.
\endpro
%%%%%%%%%%%%%%%%%%%%%%%%%%%%%%%%%%%%%%%%%%%%%%%%%%%%%%%%%%%%%%%%%%%%%%%%%%%%%%%%%%%%%%%%%%%%%%%%%%%%%%%%%%%%%%%%%%%%%%%%%%%%%%%%
\definition\label{defrsVBpullbackofVBisomorphism}
Let $r$ and $s$ be non-negative integers. Suppose that $\vbundle{1}$ and $\vbundle{2}$ are isomorphic, and
let $f$ be an element of
$\vbisomorphisms{\vbundle{1}}{\vbundle{2}}$ (a smooth vector bundle isomorphism from
$\vbundle{1}$ to $\vbundle{2}$). The mapping $\function{\VBpullback{f}{r}{s}}{\TF{r}{s}{\vbundle{2}}}{\TF{r}{s}{\vbundle{1}}}$
is defined as
\begin{align}
\Foreach{\atf{}}{\TF{r}{s}{\vbundle{2}}}
\Foreach{\point}{\vB{1}}
\func{\[\func{\VBpullback{f}{r}{s}}{\atf{}}\]}{\point}\eqdef
\func{\Vpullback{\(\reS{f}{\func{\pimage{\vbprojection{1}}}{\seta{\point}}}\)}{r}{s}}{\func{\atf{}}{\func{f}{\point}}}.
\end{align}
$\VBpullback{f}{r}{s}$ is referred to as the $\quotl$$\opair{r}{s}$-pullback of the (smooth vector bundle) isomorphism
$f\in\vbisomorphisms{\vbundle{1}}{\vbundle{2}}$$\quotr$.
\endef
%%%%%%%%%%%%%%%%%%%%%%%%%%%%%%%%%%%%%%%%%%%%%%%%%%%%%%%%%%%%%%%%%%%%%%%%%%%%%%%%%%%%%%%%%%%%%%%%%%%%%%%%%%%%%%%%%%%%%%%%%%%%%%%%
\proposition
Let $r$ and $s$ be non-negative integers, and $f$ an element of
$\vbisomorphisms{\vbundle{1}}{\vbundle{2}}$ (a smooth vector bundle isomorphism from
$\vbundle{1}$ to $\vbundle{2}$).
\begin{align}
&\Foreach{\atf{}}{\TF{r}{s}{\vbundle{2}}}\cr
&\Foreach{\point}{\vB{1}}
\Foreach{\mtuple{\vv{1}}{\vv{r+s}}}{\Cprod{\multiprod{\fibervecs{\vbundle{2}}{\point}}{r}}{\multiprod{\Vdual{\fibervecs{\vbundle{2}}{\point}}}{s}}}\cr
&\begin{aligned}
&~\func{\(\func{\[\func{\VBpullback{f}{r}{s}}{\atf{}}\]}{\point}\)}{\suc{\vv{1}}{\vv{r+s}}}\cr
=&~
\func{\[\func{\atf{}}{\func{f}{\point}}\]}
{\begin{aligned}
%\binary{
&\suc{\func{\(\reS{f}{\func{\pimage{\vbprojection{1}}}{\seta{\point}}}\)}{\vv{1}}}
{\func{\(\reS{f}{\func{\pimage{\vbprojection{1}}}{\seta{\point}}}\)}{\vv{r}}},\cr
%}
%%
%{
&\suc{\func{\dualpb{\[\finv{\(\reS{f}{\func{\pimage{\vbprojection{1}}}{\seta{\point}}}\)}\]}}{\vv{r+1}}}
{\func{\dualpb{\[\finv{\(\reS{f}{\func{\pimage{\vbprojection{1}}}{\seta{\point}}}\)}\]}}{\vv{r+s}}}
%}
\end{aligned}
}.
\end{aligned}\cr
&{}
\end{align}
\proof
It is trivial.
\endpro
%%%%%%%%%%%%%%%%%%%%%%%%%%%%%%%%%%%%%%%%%%%%%%%%%%%%%%%%%%%%%%%%%%%%%%%%%%%%%%%%%%%%%%%%%%%%%%%%%%%%%%%%%%%%%%%%%%%%%%%%%%%%%%%%
\theorem\label{thmrsVBpullbackofcompositionofdiffeomorphisms}
Let $r$ and $s$ be non-negative integers.
Let $f$ be an element of $\vbisomorphisms{\vbundle{}}{\vbundle{1}}$ and
$g$ an element of $\vbisomorphisms{\vbundle{1}}{\vbundle{2}}$.
\begin{equation}
\VBpullback{\(\identity{\vTot{}}\)}{r}{s}=\identity{\TF{r}{s}{\vbundle{}}},
\end{equation}
and
\begin{equation}
\VBpullback{\(\cmp{g}{f}\)}{r}{s}=\cmp{\VBpullback{f}{r}{s}}{\VBpullback{g}{r}{s}}.
\end{equation}
\proof
According to \refthm{thmrspullbackofcompositionoflinearisomorphisms} and
\refdef{defrsVBpullbackofVBisomorphism}, it is trivial.
\endthm
%%%%%%%%%%%%%%%%%%%%%%%%%%%%%%%%%%%%%%%%%%%%%%%%%%%%%%%%%%%%%%%%%%%%%%%%%%%%%%%%%%%%%%%%%%%%%%%%%%%%%%%%%%%%%%%%%%%%%%%%%%%%%%%%
\theorem\label{thmpullbackofrsVBtensorfieldisislinear}
Let $r$ and $s$ be non-negative integers, and
let $f$ be an element of $\vbisomorphisms{\vbundle{1}}{\vbundle{2}}$.
%%%%%%%%
The mapping $\VBpullback{f}{r}{s}$ is a linear isomorphism from
$\TF{r}{s}{\vbundle{2}}$ to $\TF{r}{s}{\vbundle{1}}$. That is $\VBpullback{f}{r}{s}$ is a bijection
and,
%%%%%%%%
\begin{equation}
\Foreach{\opair{\atf{1}}{\atf{2}}}{\Cprod{\TF{r}{s}{\vbundle{2}}}{\TF{r}{s}{\vbundle{2}}}}
\Foreach{c}{\algfield{}}
\func{\VBpullback{f}{r}{s}}{c\atf{1}+\atf{2}}=c\func{\VBpullback{f}{r}{s}}{\atf{1}}+
\func{\VBpullback{f}{r}{s}}{\atf{2}}.
\end{equation}
Furthermore,
\begin{equation}
\finv{\(\VBpullback{f}{r}{s}\)}=\VBpullback{\(\finv{f}\)}{r}{s}.
\end{equation}
\proof
The linearity of $\VBpullback{f}{r}{s}$ follows directly from 
\refthm{thminverseofpullbackofrstensor}, \refdef{defrsVBpullbackofVBisomorphism},
and the canonical linear structures of $\VTF{r}{s}{\vbundle{1}}$ and $\VTF{r}{s}{\vbundle{2}}$.
It is also easy to verify that $\TFpullback{\(\finv{f}\)}{r}{s}$
is both a left and right inverse for $\TFpullback{f}{r}{s}$, considering that
$\reS{\finv{f}}{\func{\pimage{\vbprojection{2}}}{\seta{\func{f}{\point}}}}$ is the inverse of
$\reS{f}{\func{\pimage{\vbprojection{1}}}{\seta{\point}}}$
for every $\point\in\vB{1}$.
\endthm
%%%%%%%%%%%%%%%%%%%%%%%%%%%%%%%%%%%%%%%%%%%%%%%%%%%%%%%%%%%%%%%%%%%%%%%%%%%%%%%%%%%%%%%%%%%%%%%%%%%%%%%%%%%%%%%%%%%%%%%%%%%%%%%%
%%%%%%%%%%%%%%%%%%%%%%%%%%%%%%%%%%%%%%%%%%%%%%%%%%%%%%%%%%%%%%%%%%%%%%%%%%%%%%%%%%%%%%%%%%%%%%%%%%%%%%%%%%%%%%%%%%%%%%%%%%%%%%%%
\proposition
Let $r$ be a non-negative integer, and let $f$ be an element of $\vbmorphisms{\vbundle{1}}{\vbundle{2}}$
(a smooth vector bundle morphism from $\vbundle{1}$ to $\vbundle{2}$). The image of the assignment
$\TF{r}{0}{\vbundle{2}}\ni\atf{}\mapsto\bar{\atf{}}\in\Func{\vB{1}}{\VBTensors{r}{0}{\vbundle{}}}$ such that
\begin{align}
\Foreach{\atf{}}{\TF{r}{0}{\vbundle{2}}}
\Foreach{\point}{\vB{1}}
%\func{\[\func{\VBpullback{f}{r}{s}}{\atf{}}\]}{\point}\eqdef
\func{\bar{\atf{}}}{\point}\eqdef
\func{\Vpullbackcov{\(\reS{f}{\func{\pimage{\vbprojection{1}}}{\seta{\point}}}\)}{r}}{\func{\atf{}}{\func{f}{\point}}},
\end{align}
is included in $\TF{r}{0}{\vbundle{1}}$. In other words, for every $\atf{}\in\TF{r}{0}{\vbundle{2}}$,
$\bar{\atf{}}$ lies in $\TF{r}{0}{\vbundle{1}}$.
\proof
It is left to the reader as an exercise.
\endpro
%%%%%%%%%%%%%%%%%%%%%%%%%%%%%%%%%%%%%%%%%%%%%%%%%%%%%%%%%%%%%%%%%%%%%%%%%%%%%%%%%%%%%%%%%%%%%%%%%%%%%%%%%%%%%%%%%%%%%%%%%%%%%%%%
\definition\label{defcovariantpullbackofVBmorphism}
Let $r$ be a non-negative integer, and let $f$ be an element of $\vbmorphisms{\vbundle{1}}{\vbundle{2}}$
(a smooth vector bundle morphism from $\vbundle{1}$ to $\vbundle{2}$). The mapping
$\function{\VBpullbackcov{f}{r}}{\TF{r}{0}{\vbundle{2}}}{\TF{r}{0}{\vbundle{1}}}$ is defined as
\begin{align}
&\Foreach{\atf{}}{\TF{r}{0}{\vbundle{2}}}
\Foreach{\point}{\vB{1}}\cr
&\func{\[\func{\VBpullbackcov{f}{r}}{\atf{}}\]}{\point}\eqdef
\func{\Vpullbackcov{\(\reS{f}{\func{\pimage{\vbprojection{1}}}{\seta{\point}}}\)}{r}}{\func{\atf{}}{\func{f}{\point}}}.
\end{align}
$\VBpullbackcov{f}{r}$ is referred to as the $\quotl$$r$-pullback of the smooth vector bundle morphism
$f\in\vbmorphisms{\vbundle{1}}{\vbundle{2}}$$\quotr$.
\endef
%%%%%%%%%%%%%%%%%%%%%%%%%%%%%%%%%%%%%%%%%%%%%%%%%%%%%%%%%%%%%%%%%%%%%%%%%%%%%%%%%%%%%%%%%%%%%%%%%%%%%%%%%%%%%%%%%%%%%%%%%%%%%%%%
\proposition
Let $r$ be a non-negative integer, and let $f$ be an element of $\vbmorphisms{\vbundle{1}}{\vbundle{2}}$
(a smooth vector bundle morphism from $\vbundle{1}$ to $\vbundle{2}$).
\begin{align}
&\Foreach{\atf{}}{\TF{r}{0}{\vbundle{2}}}
\Foreach{\point}{\vB{1}}
\Foreach{\mtuple{\vv{1}}{\vv{r}}}{\multiprod{\fibervecs{\vbundle{2}}{\point}}{r}}\cr
&\func{\(\func{\[\func{\VBpullbackcov{f}{r}}{\atf{}}\]}{\point}\)}{\suc{\vv{1}}{\vv{r}}}
=
\func{\[\func{\atf{}}{\func{f}{\point}}\]}
{\suc{\func{\(\reS{f}{\func{\pimage{\vbprojection{1}}}{\seta{\point}}}\)}{\vv{1}}}
{\func{\(\reS{f}{\func{\pimage{\vbprojection{1}}}{\seta{\point}}}\)}{\vv{r}}}}.
\end{align}
\proof
It is trivial.
\endpro
%%%%%%%%%%%%%%%%%%%%%%%%%%%%%%%%%%%%%%%%%%%%%%%%%%%%%%%%%%%%%%%%%%%%%%%%%%%%%%%%%%%%%%%%%%%%%%%%%%%%%%%%%%%%%%%%%%%%%%%%%%%%%%%%
\theorem\label{thmcovpullbackofcompositionofVBmorphisms}
Let $r$ be a non-negative integer.
Let $f$ be an element of $\vbmorphisms{\vbundle{}}{\vbundle{1}}$ and
$g$ an element of $\vbmorphisms{\vbundle{1}}{\vbundle{2}}$.
\begin{equation}
\VBpullbackcov{\(\identity{\vTot{}}\)}{r}=\identity{\TF{r}{0}{\vbundle{}}},
\end{equation}
and
\begin{equation}
\VBpullbackcov{\(\cmp{g}{f}\)}{r}=\cmp{\VBpullbackcov{f}{r}}{\VBpullbackcov{g}{r}}.
\end{equation}
\proof
According to \refthm{thmrpullbackofcompositionoflinearmaps} and
\refdef{defcovariantpullbackofVBmorphism}, it is trivial.
\endthm
%%%%%%%%%%%%%%%%%%%%%%%%%%%%%%%%%%%%%%%%%%%%%%%%%%%%%%%%%%%%%%%%%%%%%%%%%%%%%%%%%%%%%%%%%%%%%%%%%%%%%%%%%%%%%%%%%%%%%%%%%%%%%%%%%%%%%%%%%%%%%
\theorem\label{thmpullbackofVBcovtensorfieldisislinear}
Let $r$ be a non-negative integer,
and let $f$ be an element of $\vbmorphisms{\vbundle{1}}{\vbundle{2}}$
(a smooth vector bundle morphism from $\vbundle{1}$ to $\vbundle{2}$).
$\VBpullbackcov{f}{r}$ is a linear map
from $\VTF{r}{0}{\vbundle{2}}$ to $\VTF{r}{0}{\vbundle{1}}$. That is,
\begin{equation}
\Foreach{\opair{\atf{1}}{\atf{2}}}{\Cprod{\TF{r}{s}{\vbundle{2}}}{\TF{r}{0}{\vbundle{2}}}}
\Foreach{c}{\algfield{}}
\func{\VBpullbackcov{f}{r}}{c\atf{1}+\atf{2}}=c\func{\VBpullbackcov{f}{r}}{\atf{1}}+
\func{\VBpullbackcov{f}{r}}{\atf{2}}.
\end{equation}
\proof
The linearity of $\VBpullbackcov{f}{r}$ follows directly from 
\refthm{thmpullbackofcovtensorislinear}, \refdef{defcovariantpullbackofVBmorphism},
and the canonical linear structures of $\VTF{r}{0}{\vbundle{1}}$ and $\VTF{r}{0}{\vbundle{2}}$.
\endthm
%%%%%%%%%%%%%%%%%%%%%%%%%%%%%%%%%%%%%%%%%%%%%%%%%%%%%%%%%%%%%%%%%%%%%%%%%%%%%%%%%%%%%%%%%%%%%%%%%%%%%%%%%%%%%%%%%%%%%%%%%%%%%%%%
\definition
Let $f$ be an element of $\vbisomorphisms{\vbundle{1}}{\vbundle{2}}$
(a smooth vector bundle isomorphism from $\vbundle{1}$ to $\vbundle{2}$).
The mapping $\function{\VBPullback{f}}{\DVTF{\vbundle{2}}}{\DVTF{\vbundle{1}}}$ is defined term-wise as,
\begin{equation}
\Foreach{\atf{}}{\DVTF{\vbundle{2}}}
{\[\func{\VBPullback{f}}{\atf{}}\]}_{\opair{r}{s}}:=\func{\VBpullback{f}{r}{s}}{\atf{\opair{r}{s}}}.
\end{equation}
\endef
%%%%%%%%%%%%%%%%%%%%%%%%%%%%%%%%%%%%%%%%%%%%%%%%%%%%%%%%%%%%%%%%%%%%%%%%%%%%%%%%%%%%%%%%%%%%%%%%%%%%%%%%%%%%%%%%%%%%%%%%%%%%%%%%
\corollary
Let $f$ be an element of $\vbisomorphisms{\vbundle{1}}{\vbundle{2}}$.
$\VBPullback{f}$ is a linear isomorphism from $\DVTF{\vbundle{2}}$ to $\DVTF{\vbundle{1}}$.
\endcor
%%%%%%%%%%%%%%%%%%%%%%%%%%%%%%%%%%%%%%%%%%%%%%%%%%%%%%%%%%%%%%%%%%%%%%%%%%%%%%%%%%%%%%%%%%%%%%%%%%%%%%%%%%%%%%%%%%%%%%%%%%%%%%%%
%%%%%%%%%%%%%%%%%%%%%%%%%%%%%%%%%%%%%%%%%%%%%%%%%%%%%%%%%%%%%%%%%%%%%%%%%%%%%%%%%%%%%%%%%%%%%%%%%%%%%%%%%%%%%%%%%%%%%%%%%%%%%%%%
%%%%%%%%%%%%%%%%%%%%%%%%%%%%%%%%%%%%%%%%%%%%%%%%%%%%%%%%%%%%%%%%%%%%%%%%%%%%%%%%%%%%%%%%%%%%%%%%%%%%%%%%%%%%%%%%%%%%%%%%%%%%%%%%
%%%%%%%%%%%%%%%%%%%%%%%%%%%%%%%%%%%%%%%%%%%%%%%%%%%%%%%%%%%%%%%%%%%%%%%%%%%%%%%%%%%%%%%%%%%%%%%%%%%%%%%%%%%%%%%%%%%%%%%%%%%%%%%%
%%%%%%%%%%%%%%%%%%%%%%%%%%%%%%%%%%%%%%%%%%%%%%%%%%%%%%%%%%%%%%%%%%%%%%%%%%%%%%%%%%%%%%%%%%%%%%%%%%%%%%%%%%%%%%%%%%%%%%%%%%%%%%%%
%%%%%%%%%%%%%%%%%%%%%%%%%%%%%%%%%%%%%%%%%%%%%%%%%%%%%%%%%%%%%%%%%%%%%%%%%%%%%%%%%%%%%%%%%%%%%%%%%%%%%%%%%%%%%%%%%%%%%%%%%%%%%%%%
%%%%%%%%%%%%%%%%%%%%%%%%%%%%%%%%%%%%%%%%%%%%%%%%%%%%%%%%%%%%%%%%%%%%%%%%%%%%%%%%%%%%%%%%%%%%%%%%%%%%%%%%%%%%%%%%%%%%%%%%%%%%%%%%
%%%%%%%%%%%%%%%%%%%%%%%%%%%%%%%%%%%%%%%%%%%%%%%%%%%%%%%%%%%%%%%%%%%%%%%%%%%%%%%%%%%%%%%%%%%%%%%%%%%%%%%%%%%%%%%%%%%%%%%%%%%%%%%%
%%%%%%%%%%%%%%%%%%%%%%%%%%%%%%%%%%%%%%%%%%%%%%%%%%%%%%%%%%%%%%%%%%%%%%%%%%%%%%%%%%%%%%%%%%%%%%%%%%%%%%%%%%%%%%%%%%%%%%%%%%%%%%%%
%%%%%%%%%%%%%%%%%%%%%%%%%%%%%%%%%%%%%%%%%%%%%%%%%%%%%%%%%%%%%%%%%%%%%%%%%%%%%%%%%%%%%%%%%%%%%%%%%%%%%%%%%%%%%%%%%%%%%%%%%%%%%%%%
%%%%%%%%%%%%%%%%%%%%%%%%%%%%%%%%%%%%%%%%%%%%%%%%%%%%%%%%%%%%%%%%%%%%%%%%%%%%%%%%%%%%%%%%%%%%%%%%%%%%%%%%%%%%%%%%%%%%%%%%%%%%%%%%
%%%%%%%%%%%%%%%%%%%%%%%%%%%%%%%%%%%%%%%%%%%%%%%%%%%%%%%%%%%%%%%%%%%%%%%%%%%%%%%%%%%%%%%%%%%%%%%%%%%%%%%%%%%%%%%%%%%%%%%%%%%%%%%%
\section{Smooth Vector Sub-bundle}
%%%%%%%%%%%%%%%%%%%%%%%%%%%%%%%%%%%%%%
\textcolor{Blue}{\theorem
Let $\p{\vbtotal{}}$ be a regular submanifold of $\vbtotal{}$, and let $\p{\vbfiber{}}$ be a
vector sub-space of $\vbfiber{}$. Suppose that for every point $\point$ of $\vbbase{}$ there exists a local trivialization
$\opair{\U}{\phi}$ of the smooth vector bundle $\vbundle{}$ such that
\begin{equation}
\func{\image{\phi}}{\func{\finv{\vbprojection{}}}{\U}\cap\p{\vbtotal{}}}=\Cprod{\U}{\p{\vbfiber{}}}.
\end{equation}
Then the quadruple $\quadruple{\p{\vbtotal{}}}{\Res{\vbprojection{}}{\p{\vbtotal{}}}}
{\vbbase{}}{\vecsmanifold{\p{\vbfiber{}}}}$ is a smooth fiber bundle, having the set\\
$\defSet{\Res{\phi}{\func{\finv{\vbprojection{}}}{\U}\cap\p{\vbtotal{}}}}
{\opair{\U}{\phi}\in\vbatlas{}}$ as a smooth fiber bundle atlas. Furthermore, for every local trivializations
$\opair{\U}{\phi}$ and $\opair{\V}{\psi}$ of $\vbundle{}$, the transition map of
$\Res{\phi}{\func{\pimage{\vbprojection{}}}{\U}\cap\p{\vbtotal{}}}$ and
$\Res{\psi}{\func{\pimage{\vbprojection{}}}{\V}\cap\p{\vbtotal{}}}$ is an element of $\GL{\p{\vbfiber{}}}{}$ for every
$\point\in\U\cap\V$. That is,
\begin{align}
&\Foreach{\opair{\opair{\U}{\phi}}{\opair{\V}{\psi}}}{\Cprod{\vbatlas{}}{\vbatlas{}}}
\Foreach{\point}{\U\cap\V}\cr
&\begin{aligned}
\func{\(\transition{}{\Res{\phi}{\func{\pimage{\vbprojection{}}}{\U}\cap\p{\vbtotal{}}}}
{\Res{\psi}{\func{\pimage{\vbprojection{}}}{\V}\cap\p{\vbtotal{}}}}\)}{\point}\in
\GL{\p{\vbfiber{}}}{}.
\end{aligned}
\end{align}}
\endthm
%%%%%%%%%%%%%%%%%%%%%%%%%%%%%%%%%%%%%%%%%%%%%%%%%%%%%%%%%%%%%%%%%%%%%%%%%%%%%%%%%%%%%%%%%%%%%%%%%%%%%%%%%%%%%%%%%%%%%%%%%%%%%%%%
\definition
In continuation of the previous theorem,
the maximal subset of $\vbatlas{}$ including
$\defSet{\Res{\phi}{\func{\finv{\vbprojection{}}}{\U}\cap\p{\vbtotal{}}}}
{\opair{\U}{\phi}\in\vbatlas{}}$, endowed with which the smooth fiber bundle
$\quadruple{\p{\vbtotal{}}}{\Res{\vbprojection{}}{\p{\vbtotal{}}}}
{\vbbase{}}{\vecsmanifold{\p{\vbfiber{}}}}$ becomes a
smooth vector bundle, will be denoted by
$\subvbatlas{\vbundle{}}{\p{\vbtotal{}}}{\p{\vbfiber{}}}$. Moreover, the smooth vector bundle
$\quintuple{\p{\vbtotal{}}}{\Res{\vbprojection{}}{\p{\vbtotal{}}}}
{\vbbase{}}{\p{\vbfiber{}}}{\subvbatlas{\vbundle{}}{\p{\vbtotal{}}}{\p{\vbfiber{}}}}$
will be called a $\quotl$smooth vector sub-bundle of the vector bundle $\vbundle{}$$\quotr$.
\endef
%%%%%%%%%%%%%%%%%%%%%%%%%%%%%%%%%%%%%%%%%%%%%%%%%%%%%%%%%%%%%%%%%%%%%%%%%%%%%%%%%%%%%%%%%%%%%%%%%%%%%%%%%%%%%%%%%%%%%%%%%%%%%%%%
\textcolor{Blue}{\theorem
Let $\WW{\point}$ be a vector subspace of $\fibervecs{\vbundle{}}{\point}$ for every point $\point$ of $\vbbase{}$,
and let $\displaystyle\WW{}:=\Union{\point}{\vB{}}{\WW{\point}}$.
$\WW{}$ is a regular submanifold of $\vbtotal{}$ and the quintuple
$\quintuple{\subman{\vbtotal{}}{\WW{}}}{\reS{\vbprojection{}}{\WW{}}}{\vbbase{}}{\p{\vbfiber{}}}
{\subvbatlas{\vbundle{}}{\p{\vbtotal{}}}{\p{\vbfiber{}}}}$ is a smooth vector sub-bundle of $\vbundle{}$ with rank $l$
for some vector subspace $\p{\vbfiber{}}$ of $\vbfiber{}$
if and only if for every $\point\in\vB{}$ there exists a neighborhood $\U$ of $\point$ and a system of
$l$ local sections $\suc{\vbsec{1}}{\vbsec{l}}\in\vbsectionsl{\vbundle{}}{\U}$ such that for each $q\in\U$ the set
$\seta{\suc{\func{\vbsec{1}}{q}}{\func{\vbsec{l}}{q}}}$ is base of the vector space $\WW{\point}$.
}
\endthm
%%%%%%%%%%%%%%%%%%%%%%%%%%%%%%%%%%%%%%%%%%%%%%%%%%%%%%%%%%%%%%%%%%%%%%%%%%%%%%%%%%%%%%%%%%%%%%%%%%%%%%%%%%%%%%%%%%%%%%%%%%%%%%%%
%%%%%%%%%%%%%%%%%%%%%%%%%%%%%%%%%%%%%%%%%%%%%%%%%%%%%%%%%%%%%%%%%%%%%%%%%%%%%%%%%%%%%%%%%%%%%%%%%%%%%%%%%%%%%%%%%%%%%%%%%%%%%%%%
%%%%%%%%%%%%%%%%%%%%%%%%%%%%%%%%%%%%%%%%%%%%%%%%%%%%%%%%%%%%%%%%%%%%%%%%%%%%%%%%%%%%%%%%%%%%%%%%%%%%%%%%%%%%%%%%%%%%%%%%%%%%%%%%
%%%%%%%%%%%%%%%%%%%%%%%%%%%%%%%%%%%%%%%%%%%%%%%%%%%%%%%%%%%%%%%%%%%%%%%%%%%%%%%%%%%%%%%%%%%%%%%%%%%%%%%%%%%%%%%%%%%%%%%%%%%%%%%%
%%%%%%%%%%%%%%%%%%%%%%%%%%%%%%%%%%%%%%%%%%%%%%%%%%%%%%%%%%%%%%%%%%%%%%%%%%%%%%%%%%%%%%%%%%%%%%%%%%%%%%%%%%%%%%%%%%%%%%%%%%%%%%%%
\section{Base-Restriction of a Smooth Vector Bundle}
\fixed
$\M{}$ is fixed as a regular submanifold of $\vbbase{}$.
\endfixed
%%%%%%%%%%%%%%%%%%%%%%%%%%%%%%%%%%%%%%%%%%%%%%%%%%%%%%%%%%%%%%%%%%%%%%%%%%%%%%%%%%%%%%%%%%%%%%%%%%%%%%%%%%%%%%%%%%%%%%%%%%%%%%%%
\textcolor{Blue}{\theorem
\begin{itemize}
\item
$\func{\pimage{\vbprojection{}}}{\M{}}$ is a regular submanifold of $\vbtotal{}$.
\item
$\quadruple{\subman{\vbtotal{}}{\func{\pimage{\vbprojection{}}}{\M{}}}}
{\Res{\vbprojection{}}{\func{\pimage{\vbprojection{}}}{\M{}}}}
{\subman{\vbbase{}}{\M{}}}{\vbfiber{}}$ is a smooth fiber bundle, and the set\\
$\defSet{\Res{\phi}{\func{\pimage{\vbprojection{}}}{\M{}\cap\U}}}{\opair{\U}{\phi}\in\vbatlas{}}$ is an atlas of this
fiber bundle.
\item
For every pair $\opair{\U}{\phi}$ and $\opair{\V}{\psi}$ of local trivializations of the smooth vector bundle
$\vbundle{}$,
\begin{equation}
\Foreach{\point}{\func{\pimage{\vbprojection{}}}{\M{}\cap\U\cap\V}}
\func{\transition{}{\Res{\phi}{\func{\pimage{\vbprojection{}}}{\M{}\cap\U}}}
{\Res{\psi}{\func{\pimage{\vbprojection{}}}{\M{}\cap\V}}}}{\point}\in\GL{\vbfiber{}}{}.
\end{equation}
\end{itemize}
}
\endthm
%%%%%%%%%%%%%%%%%%%%%%%%%%%%%%%%%%%%%%%%%%%%%%%%%%%%%%%%%%%%%%%%%%%%%%%%%%%%%%%%%%%%%%%%%%%%%%%%%%%%%%%%%%%%%%%%%%%%%%%%%%%%%%%%
\definition
We will denote by $\vbResatlas{\vbundle{}}{\M{}}$
the maximal atlas of the smooth fiber bundle\\
$\quadruple{\subman{\vbtotal{}}{\func{\pimage{\vbprojection{}}}{\M{}}}}
{\Res{\vbprojection{}}{\func{\pimage{\vbprojection{}}}{\M{}}}}
{\subman{\vbbase{}}{\M{}}}{\vbfiber{}}$
including
$\defSet{\Res{\phi}{\func{\pimage{\vbprojection{}}}{\M{}\cap\U}}}{\opair{\U}{\phi}\in\vbatlas{}}$,
endowed with which, this fiber bundle becomes a smooth vector bundle.
\endef
%%%%%%%%%%%%%%%%%%%%%%%%%%%%%%%%%%%%%%%%%%%%%%%%%%%%%%%%%%%%%%%%%%%%%%%%%%%%%%%%%%%%%%%%%%%%%%%%%%%%%%%%%%%%%%%%%%%%%%%%%%%%%%%%
\corollary
The quintuple
$\quintuple{\subman{\vbtotal{}}{\func{\pimage{\vbprojection{}}}{\M{}}}}
{\Res{\vbprojection{}}{\func{\pimage{\vbprojection{}}}{\M{}}}}
{\subman{\vbbase{}}{\M{}}}{\vbfiber{}}
{\vbResatlas{\vbundle{1}}{\M{}}}$ is a smooth vector bundle.
\endcor
%%%%%%%%%%%%%%%%%%%%%%%%%%%%%%%%%%%%%%%%%%%%%%%%%%%%%%%%%%%%%%%%%%%%%%%%%%%%%%%%%%%%%%%%%%%%%%%%%%%%%%%%%%%%%%%%%%%%%%%%%%%%%%%%
\definition
We will denote by $\vbResbundle{\vbundle{}}{\M{}}$ the smooth vector bundle\\
$\quintuple{\subman{\vbtotal{}}{\func{\pimage{\vbprojection{}}}{\M{}}}}
{\Res{\vbprojection{}}{\func{\pimage{\vbprojection{}}}{\M{}}}}
{\subman{\vbbase{}}{\M{}}}{\vbfiber{}}
{\vbResatlas{\vbundle{1}}{\M{}}}$, which will be referred to as the $\quotl$restriction of the smooth vector bundle
$\vbundle{}$ to $\M{}$$\quotr$.
\endef
%%%%%%%%%%%%%%%%%%%%%%%%%%%%%%%%%%%%%%%%%%%%%%%%%%%%%%%%%%%%%%%%%%%%%%%%%%%%%%%%%%%%%%%%%%%%%%%%%%%%%%%%%%%%%%%%%%%%%%%%%%%%%%%%
\corollary
For every point $\point\in\M{}$,
the fiber space of $\vbResbundle{\vbundle{}}{\M{}}$ over $\point$ is the same as the fiber space of
$\vbundle{}$ at $\point$. That is,
\begin{align}
\Foreach{\point}{\M{}}
\fibervecs{\vbResbundle{\vbundle{}}{\M{}}}{\point}=
\fibervecs{\vbundle{}}{\point}.
\end{align}
\endcor
\section{Operations on Smooth Vector Bundles}
\subsection{Hom Bundle}
\fixed
In this section,
$\vbundle{1}=\quintuple{\vbtotal{1}}{\vbprojection{1}}{\vbbase{}}{\vbfiber{}}{\vbatlas{1}}$ and
$\vbundle{2}=\quintuple{\vbtotal{2}}{\vbprojection{2}}{\vbbase{}}{\vbfiber{}}{\vbatlas{2}}$
are fixed as a pair of smooth vector bundles
of rank $d$,
where $\vbtotal{i}=\opair{\vTot{i}}{\maxatlas{\vTot{i}}}$ for $i=1,~2$, and
$\vbbase{}=\opair{\vB{}}{\maxatlas{\vB{}}}$ are $\difclass{\infty}$ manifolds
modeled on the Banach-spaces $\R^{n_{\vTot{i}}}$ ($i=1,~2$) and $\R^{n_{\vB{}}}$, respectively.
\endfixed
%%%%%%%%%%%%%%%%%%%%%%%%%%%%%%%%%%%%%%%%%%%%%%%%%%%%%%%%%%%%%%
\definition
Let $\U$ be a subset of $\vB{}$.
\begin{equation}
\vbHom{\vbundle{1}}{\vbundle{2}}{\U}:=
\Union{\point}{\U}{\Lin{\fibervecs{\vbundle{1}}{\point}}{\fibervecs{\vbundle{2}}{\point}}}.
\end{equation}
Specifically, $\vbHom{\vbundle{1}}{\vbundle{2}}{\vB{}}$ can alternatively be denoted by $\VBHom{\vbundle{1}}{\vbundle{2}}$.
\endef
%%%%%%%%%%%%%%%%%%%%%%%%%%%%%%%%%%%%%%%%%%%%%%%%%%%%%%%%%%%%%%%%%%%%%%%%%%%%%%%%%%%%%%%%%%%%%%%%%%%%%%%%%%%%%%%%%%%%%%%%%%%%%%%%
%%%%%%%%%%%%%%%%%%%%%%%%%%%%%%%%%%%%%%%%%%%%%%%%%%%%%%%%%%%%%%%%%%%%%%%%%%%%%%%%%%%%%%%%%%%%%%%%%%%%%%%%%%%%%%%%%%%%%%%%%%%%%%%%
\definition
\begin{align}
&~\vbchartlocaltt{\vbundle{1}}{\vbundle{2}}\cr
:=&~\defSet{\triple{\phi}{\psi_1}{\psi_2}}{\[\phi\in\maxatlas{\vB{}},~
\psi_1\in\vbatlas{1},~\psi_2\in\vbatlas{2},~\domain{\phi}=\func{\image{\vbprojection{1}}}{\domain{\psi_1}}=
\func{\image{\vbprojection{2}}}{\domain{\psi_2}}\]}.\cr
&{}
\end{align}
An element $\triple{\phi}{\psi_1}{\psi_2}$ of $\vbchartlocaltt{\vbundle{1}}{\vbundle{2}}$ can alternatively be denoted by
$\btriple{\U}{\phi}{\binary{\psi_1}{\psi_2}}$ where $\U:=\domain{\phi}=\func{\image{\vbprojection{1}}}{\domain{\psi_1}}=
\func{\image{\vbprojection{2}}}{\domain{\psi_2}}$.
\endef
%%%%%%%%%%%%%%%%%%%%%%%%%%%%%%%%%%%%%%%%%%%%%%%%%%%%%%%%%%%%%%%%%%%%%%%%%%%%%%%%%%%%%%%%%%%%%%%%%%%%%%%%%%%%%%%%%%%%%%%%%%%%%%%%
\lemma
$\defsets{\U}{\vB{}}{\[\Exists{\triple{\phi}{\psi_1}{\psi_2}}{\vbchartlocaltt{\vbundle{1}}{\vbundle{2}}}\U=\domain{\phi}=
\func{\image{\vbprojection{i}}}{\domain{\psi_i}}~(i=1,~2)\]}$ is an open covering of $\vbbase{}$.
\proof
It is an immediate consequence of \reflem{lemlocaltrivializationcontraction}.
\endlem
%%%%%%%%%%%%%%%%%%%%%%%%%%%%%%%%%%%%%%%%%%%%%%%%%%%%%%%%%%%%%%%%%%%%%%%%%%%%%%%%%%%%%%%%%%%%%%%%%%%%%%%%%%%%%%%%%%%%%%%%%%%%%%%%
\definition\label{defHombundlelocaltrivializations}
We associate to each element $\btriple{\U}{\phi}{\binary{\psi_1}{\psi_2}}$ of
$\vbchartlocaltt{\vbundle{1}}{\vbundle{2}}$ the mapping
$\function{\vbHomchart{\phi}{\psi_1}{\psi_2}}{\vbHom{\vbundle{1}}{\vbundle{2}}{\U}}
{\Cprod{\func{\phi}{\U}}{\Endo{\vbfiber{}}}}$ defined as,
\begin{align}
&\begin{aligned}
\Foreach{\point}{\U}
\Foreach{\alpha}{\Lin{\fibervecs{\vbundle{1}}{\point}}{\fibervecs{\vbundle{2}}{\point}}}
\end{aligned}\cr
&\begin{aligned}
\func{\vbHomchart{\phi}{\psi_1}{\psi_2}}{\alpha}\eqdef
\opair{\func{\phi}{\point}}{\cmp{\pltfib{\psi_2}{\point}}{\cmp{\alpha}{\finv{\(\pltfib{\psi_1}{\point}\)}}}}.
\end{aligned}
\end{align}
Furthermore, we associate to each element $\btriple{\U}{\phi}{\binary{\psi_1}{\psi_2}}$ of
$\vbchartlocaltt{\vbundle{1}}{\vbundle{2}}$\\
the mapping
$\function{\vbHomlocalt{\phi}{\psi_1}{\psi_2}}{\vbHom{\vbundle{1}}{\vbundle{2}}{\U}}
{\Cprod{\U}{\Endo{\vbfiber{}}}}$ defined as,
\begin{align}
&\begin{aligned}
\Foreach{\point}{\U}
\Foreach{\alpha}{\Lin{\fibervecs{\vbundle{1}}{\point}}{\fibervecs{\vbundle{2}}{\point}}}
\end{aligned}\cr
&\begin{aligned}
\func{\vbHomlocalt{\phi}{\psi_1}{\psi_2}}{\alpha}\eqdef
\opair{\point}{\cmp{\pltfib{\psi_2}{\point}}{\cmp{\alpha}{\finv{\(\pltfib{\psi_1}{\point}\)}}}}.
\end{aligned}
\end{align}
\endef
%%%%%%%%%%%%%%%%%%%%%%%%%%%%%%%%%%%%%%%%%%%%%%%%%%%%%%%%%%%%%%%%%%%%%%%%%%%%%%%%%%%%%%%%%%%%%%%%%%%%%%%%%%%%%%%%%%%%%%%%%%%%%%%%
\textcolor{Blue}{\lemma
The set $\defSet{\vbHomchart{\phi}{\psi_1}{\psi_2}}{\triple{\phi}{\psi_1}{\psi_2}\in\vbchartlocaltt{\vbundle{1}}{\vbundle{2}}}$
is a $\difclass{\infty}$ atlas on $\VBHom{\vbundle{1}}{\vbundle{2}}$ modeled on
the Banach-space $\Cprod{\R^{n_{\vB{}}}}{\Endo{\vbfiber{}}}$.}
\endlem
%%%%%%%%%%%%%%%%%%%%%%%%%%%%%%%%%%%%%%%%%%%%%%%%%%%%%%%%%%%%%%%%%%%%%%%%%%%%%%%%%%%%%%%%%%%%%%%%%%%%%%%%%%%%%%%%%%%%%%%%%%%%%%%%
\definition
The $\difclass{\infty}$ maximal-atlas on $\VBHom{\vbundle{1}}{\vbundle{2}}$ modeled on
the Banach-space $\Cprod{\R^{n_{\vB{}}}}{\Endo{\vbfiber{}}}$ generated by the atlas
$\defSet{\vbHomchart{\phi}{\psi_1}{\psi_2}}{\triple{\phi}{\psi_1}{\psi_2}\in\vbchartlocaltt{\vbundle{1}}{\vbundle{2}}}$
will be denoted by
$\vbHommaxatlas{\vbundle{1}}{\vbundle{2}}$. That is,
\begin{equation}
\vbHommaxatlas{\vbundle{1}}{\vbundle{2}}:=
\maxatlasgen{\infty}{\VBHom{\vbundle{1}}{\vbundle{2}}}{\Cprod{\R^{n_{\vB{}}}}{\Endo{\vbfiber{}}}}{
\defSet{\vbHomchart{\phi}{\psi_1}{\psi_2}}{\triple{\phi}{\psi_1}{\psi_2}\in\vbchartlocaltt{\vbundle{1}}{\vbundle{2}}}}.
\end{equation}
\endef
%%%%%%%%%%%%%%%%%%%%%%%%%%%%%%%%%%%%%%%%%%%%%%%%%%%%%%%%%%%%%%%%%%%%%%%%%%%%%%%%%%%%%%%%%%%%%%%%%%%%%%%%%%%%%%%%%%%%%%%%%%%%%%%%
\textcolor{Blue}{\theorem
The differentiable structure $\opair{\VBHom{\vbundle{1}}{\vbundle{2}}}{\vbHommaxatlas{\vbundle{1}}{\vbundle{2}}}$
is a $\difclass{\infty}$ manifold, which means the topology induced by the maximal atlas
$\vbHommaxatlas{\vbundle{1}}{\vbundle{2}}$
on $\VBHom{\vbundle{1}}{\vbundle{2}}$ is Hausdorff and second-countable.}
\endthm
%%%%%%%%%%%%%%%%%%%%%%%%%%%%%%%%%%%%%%%%%%%%%%%%%%%%%%%%%%%%%%%%%%%%%%%%%%%%%%%%%%%%%%%%%%%%%%%%%%%%%%%%%%%%%%%%%%%%%%%%%%%%%%%%
\definition
The manifold $\opair{\VBHom{\vbundle{1}}{\vbundle{2}}}{\vbHommaxatlas{\vbundle{1}}{\vbundle{2}}}$
will be denoted by $\vbHomMan{\vbundle{1}}{\vbundle{2}}$.
\endef
%%%%%%%%%%%%%%%%%%%%%%%%%%%%%%%%%%%%%%%%%%%%%%%%%%%%%%%%%%%%%%%%%%%%%%%%%%%%%%%%%%%%%%%%%%%%%%%%%%%%%%%%%%%%%%%%%%%%%%%%%%%%%%%%
\definition
The mapping $\function{\vbHomprojection{\vbundle{1}}{\vbundle{2}}}{\VBHom{\vbundle{1}}{\vbundle{2}}}{\vB{}}$ is
defined as,
\begin{equation}
\Foreach{\point}{\vB{}}
\Foreach{\alpha}{\Lin{\fibervecs{\vbundle{1}}{\point}}{\fibervecs{\vbundle{2}}{\point}}}
\func{\vbHomprojection{\vbundle{1}}{\vbundle{2}}}{\alpha}\eqdef\point.
\end{equation}
\endef
%%%%%%%%%%%%%%%%%%%%%%%%%%%%%%%%%%%%%%%%%%%%%%%%%%%%%%%%%%%%%%%%%%%%%%%%%%%%%%%%%%%%%%%%%%%%%%%%%%%%%%%%%%%%%%%%%%%%%%%%%%%%%%%%
\lemma\label{lemHombundleisafiberbundle}
The quadruple
$\quadruple{\vbHomMan{\vbundle{1}}{\vbundle{2}}}{\vbHomprojection{\vbundle{1}}{\vbundle{2}}}{\vbbase{}}
{\Endo{\vbfiber{}}}$
is a smooth fiber bundle.
\proof
Let $\point$ be an arbitrary point of $\vbbase{}$. Let $\btriple{\U}{\phi}{\binary{\psi_1}{\psi_2}}$ be an element of
$\vbchartlocaltt{\vbundle{1}}{\vbundle{2}}$ such that $\point\in\U$.
Clearly $\funcprod{\finv{\phi}}{\identity{\VBHom{\vbundle{1}}{\vbundle{2}}}}$
is a diffeomorphism from $\Cprod{\func{\phi}{\U}}{\Endo{\vbfiber{}}}$
(with its canonical differentiable structure
inherited from that of $\vecsmanifold{\Cprod{\R^{n_{\vB{}}}}{\VBHom{\vbundle{1}}{\vbundle{2}}}}$) to
$\Cprod{\subman{\vbbase{}}{\U}}{\vecsmanifold{\VBHom{\vbundle{1}}{\vbundle{2}}}}$. In addition, it is also trivial that
$\vbHomchart{\phi}{\psi_1}{\psi_2}$ is a diffeomorphism from
$\subman{\vbHomMan{\vbundle{1}}{\vbundle{2}}}{\vbHom{\vbundle{1}}{\vbundle{2}}{\U}}$ to
$\Cprod{\func{\phi}{\U}}{\Endo{\vbfiber{}}}$, because it is a chart of the manifold
$\vbHomMan{\vbundle{1}}{\vbundle{2}}$.
Thus, considering that the composition of a pair of diffeomorphisms as again a diffeomorphism,
$\vbHomlocalt{\phi}{\psi_1}{\psi_2}=
\cmp{\(\funcprod{\finv{\phi}}{\identity{\VBHom{\vbundle{1}}{\vbundle{2}}}}\)}{\vbHomchart{\phi}{\psi_1}{\psi_2}}$
is a diffeomorphism from $\subman{\vbHomMan{\vbundle{1}}{\vbundle{2}}}{\vbHom{\vbundle{1}}{\vbundle{2}}{\U}}$ to
$\Cprod{\subman{\vbbase{}}{\U}}{\vecsmanifold{\VBHom{\vbundle{1}}{\vbundle{2}}}}$.
Moreover, it is evident that
$\func{\pimage{\(\vbHomprojection{\vbundle{1}}{\vbundle{2}}\)}}{\U}=
\vbHom{\vbundle{1}}{\vbundle{2}}{\U}$, and the following diagram is commutative.
\begin{center}
\vskip0.5\baselineskip
\hskip-2\baselineskip
\begin{tikzcd}[row sep=6em, column sep=6em]
& \vbHom{\vbundle{1}}{\vbundle{2}}{\U}
\arrow{r}{\vbHomlocalt{\phi}{\psi_1}{\psi_2}}
\arrow[swap]{d}{\vbHomprojection{\vbundle{1}}{\vbundle{2}}}
& \Cprod{\U}{\Endo{\vbfiber{}}}
\arrow{dl}{\proj{\U}{\Endo{\vbfiber{}}}{1}} \\
& \U
\end{tikzcd}
\end{center}
\endlem
%%%%%%%%%%%%%%%%%%%%%%%%%%%%%%%%%%%%%%%%%%%%%%%%%%%%%%%%%%%%%%%%%%%%%%%%%%%%%%%%%%%%%%%%%%%%%%%%%%%%%%%%%%%%%%%%%%%%%%%%%%%%%%%%
\corollary
The set $\defSet{\vbHomlocalt{\phi}{\psi_1}{\psi_2}}{\triple{\phi}{\psi_1}{\psi_2}\in\vbchartlocaltt{\vbundle{1}}{\vbundle{2}}}$
is an atlas of the smooth fiber bundle
$\quadruple{\vbHomMan{\vbundle{1}}{\vbundle{2}}}{\vbHomprojection{\vbundle{1}}{\vbundle{2}}}{\vbbase{}}
{\Endo{\vbfiber{}}}$.
\endcor
%%%%%%%%%%%%%%%%%%%%%%%%%%%%%%%%%%%%%%%%%%%%%%%%%%%%%%%%%%%%%%%%%%%%%%%%%%%%%%%%%%%%%%%%%%%%%%%%%%%%%%%%%%%%%%%%%%%%%%%%%%%%%%%%
\proposition
Let $\btriple{\U}{\phi}{\binary{\psi_1}{\psi_2}}$ and $\btriple{\V}{\p{\phi}}{\binary{\eta_1}{\eta_2}}$ be a pair of elements of
$\vbchartlocaltt{\vbundle{1}}{\vbundle{2}}$.
%such that there exists a pair $\opair{\U}{\phi}$ and $\opair{\V}{\p{\phi}}$ of charts of the manifold $\vbbase{}$,
$\function{\cmp{\vbHomlocalt{\phi}{\psi_1}{\psi_2}}{\finv{\(\vbHomlocalt{\p{\phi}}{\eta_1}{\eta_2}\)}}}
{\Cprod{\(\U\cap\V\)}{\Endo{\vbfiber{}}}}{\Cprod{\(\U\cap\V\)}{\Endo{\vbfiber{}}}}$, and
the value of the transition map of the local trivializations $\vbHomlocalt{\phi}{\psi_1}{\psi_2}$ and
$\vbHomlocalt{\p{\phi}}{\eta_1}{\eta_2}$
of the fiber bundle $\quadruple{\vbHomMan{\vbundle{1}}{\vbundle{2}}}{\vbHomprojection{\vbundle{1}}{\vbundle{2}}}{\vbbase{}}
{\vecsmanifold{\Endo{\vbfiber{}}}}$ at every point $\point$ of its domain $\U\cap\V$ is a linear isomorphism
from $\Endo{\vbfiber{}}$ to itself. That is,
\begin{align}
\Foreach{\point}{\U\cap\V}
\func{\[\transition{}{\vbHomlocalt{\phi}{\psi_1}{\psi_2}}{\vbHomlocalt{\p{\phi}}{\eta_1}{\eta_2}}\]}{\point}\in
\GL{\Endo{\vbfiber{}}}{}.
\end{align}
\proof
Let $\point$ be an arbitrary element of $\U\cap\V$.
According to the definition of the transition map of a pair of local trivializations of a smooth fiber bundle,
and \refdef{defHombundlelocaltrivializations},
\begin{align}
\Foreach{\alpha}{\Endo{\vbfiber{}}}
\func{\(\func{\[\transition{}{\vbHomlocalt{\phi}{\psi_1}{\psi_2}}{\vbHomlocalt{\p{\phi}}{\eta_1}{\eta_2}}\]}{\point}\)}{\alpha}&=
\func{\[\cmp{\pltfib{\vbHomlocalt{\phi}{\psi_1}{\psi_2}}{\point}}{\finv{\(\pltfib{\vbHomlocalt{\p{\phi}}{\eta}{\eta}}{\point}\)}}\]}{\alpha}\cr
&=\cmp{\pltfib{\psi_2}{\point}}{\cmp{
\cmp{\finv{\(\pltfib{\eta_2}{\point}\)}}{\cmp{\alpha}{\pltfib{\eta_1}{\point}}}
}{\finv{\(\pltfib{\psi_1}{\point}\)}}},
\end{align}
which obviously implies that
$\func{\[\transition{}{\vbHomlocalt{\phi}{\psi_1}{\psi_2}}{\vbHomlocalt{\p{\phi}}{\eta_1}{\eta_2}}\]}{\point}\in
\GL{\Endo{\vbfiber{}}}{}$, because
$\cmp{\pltfib{\psi_2}{\point}}{\finv{\(\pltfib{\eta_2}{\point}\)}}$ and
$\cmp{\pltfib{\eta_1}{\point}}{\finv{\(\pltfib{\psi_1}{\point}\)}}$ are linear automorphisms of $\vbfiber{}$.
\endlem
%%%%%%%%%%%%%%%%%%%%%%%%%%%%%%%%%%%%%%%%%%%%%%%%%%%%%%%%%%%%%%%%%%%%%%%%%%%%%%%%%%%%%%%%%%%%%%%%%%%%%%%%%%%%%%%%%%%%%%%%%%%%%%%%
\definition
We will denote by $\vbHomatlas{\vbundle{1}}{\vbundle{2}}$
the maximal atlas of the smooth fiber bundle\\
$\quadruple{\vbHomMan{\vbundle{1}}{\vbundle{2}}}{\vbHomprojection{\vbundle{1}}{\vbundle{2}}}{\vbbase{}}
{\vecsmanifold{\Endo{\vbfiber{}}}}$
including $\defSet{\vbHomlocalt{\phi}{\psi_1}{\psi_2}}{\triple{\phi}{\psi_1}{\psi_2}\in\vbchartlocaltt{\vbundle{1}}{\vbundle{2}}}$,
endowed with which, the fiber bundle
$\quadruple{\vbHomMan{\vbundle{1}}{\vbundle{2}}}{\vbHomprojection{\vbundle{1}}{\vbundle{2}}}{\vbbase{}}
{\vecsmanifold{\Endo{\vbfiber{}}}}$ becomes a smooth vector bundle.
\endef
%%%%%%%%%%%%%%%%%%%%%%%%%%%%%%%%%%%%%%%%%%%%%%%%%%%%%%%%%%%%%%%%%%%%%%%%%%%%%%%%%%%%%%%%%%%%%%%%%%%%%%%%%%%%%%%%%%%%%%%%%%%%%%%%
\corollary
The quintuple
$\quintuple{\vbHomMan{\vbundle{1}}{\vbundle{2}}}{\vbHomprojection{\vbundle{1}}{\vbundle{2}}}{\vbbase{}}
{\vecsmanifold{\Endo{\vbfiber{}}}}{\vbHomatlas{\vbundle{1}}{\vbundle{2}}}$ is a smooth vector bundle.
\endcor
%%%%%%%%%%%%%%%%%%%%%%%%%%%%%%%%%%%%%%%%%%%%%%%%%%%%%%%%%%%%%%%%%%%%%%%%%%%%%%%%%%%%%%%%%%%%%%%%%%%%%%%%%%%%%%%%%%%%%%%%%%%%%%%%
\definition
We will denote by $\vbHombundle{\vbundle{1}}{\vbundle{2}}$ the smooth vector bundle\\
$\quintuple{\vbHomMan{\vbundle{1}}{\vbundle{2}}}{\vbHomprojection{\vbundle{1}}{\vbundle{2}}}{\vbbase{}}
{\vecsmanifold{\Endo{\vbfiber{}}}}{\vbHomatlas{\vbundle{1}}{\vbundle{2}}}$, which is referred to as the
$\quotl$Hom bundle of the smooth vector bundles $\vbundle{1}$ and $\vbundle{2}$$\quotr$.
\endef
%%%%%%%%%%%%%%%%%%%%%%%%%%%%%%%%%%%%%%%%%%%%%%%%%%%%%%%%%%%%%%%%%%%%%%%%%%%%%%%%%%%%%%%%%%%%%%%%%%%%%%%%%%%%%%%%%%%%%%%%%%%%%%%%
\corollary
For every point $\point$ of $\vbbase{}$, the fiber space of $\vbHombundle{\vbundle{1}}{\vbundle{2}}$ over $\point$ equals
the vector-space of all linear maps from the fiber space of $\vbundle{1}$ over $\point$ to the
the fiber space of $\vbundle{2}$ over $\point$. That is,
\begin{align}
\Foreach{\point}{\vB{}}
\fibervecs{\vbHombundle{\vbundle{1}}{\vbundle{2}}}{\point}=
\Lin{\fibervecs{\vbundle{1}}{\point}}{\fibervecs{\vbundle{2}}{\point}}.
\end{align}
\endcor
%%%%%%%%%%%%%%%%%%%%%%%%%%%%%%%%%%%%%%%%%%%%%%%%%%%%%%%%%%%%%%%%%%%%%%%%%%%%%%%%%%%%%%%%%%%%%%%%%%%%%%%%%%%%%%%%%%%%%%%%%%%%%%%%
%%%%%%%%%%%%%%%%%%%%%%%%%%%%%%%%%%%%%%%%%%%%%%%%%%%%%%%%%%%%%%%%%%%%%%%%%%%%%%%%%%%%%%%%%%%%%%%%%%%%%%%%%%%%%%%%%%%%%%%%%%%%%%%%
%%%%%%%%%%%%%%%%%%%%%%%%%%%%%%%%%%%%%%%%%%%%%%%%%%%%%%%%%%%%%%%%%%%%%%%%%%%%%%%%%%%%%%%%%%%%%%%%%%%%%%%%%%%%%%%%%%%%%%%%%%%%%%%%
\subsection{Dual Bundle}
%%%%%%%%%%%%%%%%%%%%%%%%%
\definition
The $\opair{1}{0}$ tensor bundle of the smooth vector bundle $\vbundle{}$, that is
$\vbtensorbundle{1}{0}{\vbundle{}}=
\quintuple{\VBTensorsMan{1}{0}{\vbundle{}}}{\vbtensorprojection{\vbundle{}}{1}{0}}{\vbbase{}}
{\vecsmanifold{\MTensors{1}{0}{\vbfiber{}}}}{\vbtensoratlas{1}{0}{\vbundle{}}}$,  is referred to as the
$\quotl$dual bundle of the smooth vector bundle $\vbundle{}$$\quotr$ and will be denoted by
$\vbDualbundle{\vbundle{}}$, alternatively.
\endef
%%%%%%%%%%%%%%%%%%%%%%%%%%%%%%%%%%%%%%%%%%%%%%%%%%%%%%%%%%%%%%%%%%%%%%%%%%%%%%%%%%%%%%%%%%%%%%%%%%%%%%%%%%%%%%%%%%%%%%%%%%%%%%%%
\definition
Let $\U$ be a non-empty open set of $\vbbase{}$, and let $\function{\vbsec{}}{\seta{\suc{1}{d}}}{\vbsectionsl{\vbundle{}}{\U}}$
be a local frame field of $\vbundle{}$ over $\U$, that is an element of $\lframe{\vbundle{}}{\U}$. The sequence
$\function{\lframedual{\vbsec{}}{}}{\seta{\suc{1}{d}}}{\Func{\U}{\VBTensorsMan{1}{0}{\vbundle{}}}}$ is definded as the unique
such sequence that satisfies
\begin{align}
\Foreach{\opair{i}{j}}{\Cprod{\seta{\suc{1}{d}}}{\seta{\suc{1}{d}}}}
\Foreach{\point}{\U}
\func{\[\func{\lframedual{\vbsec{}}{i}}{\point}\]}{\func{\vbsec{j}}{\point}}\eqdef
\deltaf{i}{j}.
\end{align}
In other words, at every point $\point$ of $\U$, the sequence $\suc{\func{\lframedual{\vbsec{}}{1}}{\point}}
{\func{\lframedual{\vbsec{}}{d}}{\point}}$ is the dual of the ordered basis
$\suc{\func{\vbsec{1}}{\point}}{\func{\vbsec{d}}{\point}}$ of $\fibervecs{\vbundle{}}{\point}$, by definition.
$\lframedual{\vbsec{}}{}$ is referred to as the $\quotl$dual of the local frame field $\vbsec{}$ of $\vbundle{}$
over $\U$$\quotr$.
\endef
%%%%%%%%%%%%%%%%%%%%%%%%%%%%%%%%%%%%%%%%%%%%%%%%%%%%%%%%%%%%%%%%%%%%%%%%%%%%%%%%%%%%%%%%%%%%%%%%%%%%%%%%%%%%%%%%%%%%%%%%%%%%%%%%
\textcolor{Blue}{\theorem
Let $\U$ be a non-empty open set of $\vbbase{}$, and let $\function{\vbsec{}}{\seta{\suc{1}{d}}}{\vbsectionsl{\vbundle{}}{\U}}$
be a local frame field of $\vbundle{}$ over $\U$, that is an element of $\lframe{\vbundle{}}{\U}$.
For every $i$ in $\seta{\suc{1}{d}}$, $\lframedual{\vbsec{}}{i}$ is a $\opair{1}{0}$ tensor field of the smooth vector bundle
$\vbResbundle{\vbundle{}}{\U}$, or equivalently a local section of the dual bundle (that is, the $\opair{1}{0}$ tensor bundle)
of $\vbundle{}$ over $\U$.
That is,
\begin{equation}
\Foreach{i}{\seta{\suc{1}{d}}}
\lframedual{\vbsec{}}{i}\in\vbsectionsl{\vbDualbundle{\vbundle{}}}{\U}.
\end{equation}
}
\endthm
%%%%%%%%%%%%%%%%%%%%%%%%%%%%%%%%%%%%%%%%%%%%%%%%%%%%%%%%%%%%%%%%%%%%%%%%%%%%%%%%%%%%%%%%%%%%%%%%%%%%%%%%%%%%%%%%%%%%%%%%%%%%%%%%
\textcolor{Blue}{\theorem
Let $\atf{}$ be an element of $\TF{r}{s}{\vbundle{}}$, that is an $\opair{r}{s}$ tensor field of $\vbundle{}$.
Let $\U$ be a non-empty open set of $\vbbase{}$, and let $\function{\vbsec{}}{\seta{\suc{1}{r}}}{\vbsectionsl{\vbundle{}}{\U}}$
be a local frame field of $\vbundle{}$ over $\U$.
\begin{align}
\reS{\atf{}}{\U}=\sum_{\suc{j_1}{j_r}=1}^{d}\sum_{\suc{k_1}{k_s}=1}^{d}
\atf{\suc{j_1}{j_r}}^{\suc{k_1}{k_s}}
\(\lframedual{\vbsec{}}{j_1}\vbtensor{}\lframedual{\vbsec{}}{j_r}\vbtensor{}
\vbsec{k_1}\vbtensor{}\cdots\vbtensor{}\vbsec{k_s}\),
\end{align}
where for every $\suc{j_1}{j_r},\cdots\suc{k_1}{k_s}$ in $\seta{\suc{1}{d}}$,
the mapping $\function{\atf{\suc{j_1}{j_r}}^{\suc{k_1}{k_s}}}{\U}{\algfield{}}$ is defined as,
\begin{align}
\Foreach{\point}{\U}
\func{\atf{\suc{j_1}{j_r}}^{\suc{k_1}{k_s}}}{\point}\eqdef
\func{\[\func{\atf{}}{\point}\]}{\binary{\suc{\func{\vbsec{j_1}}{\point}}{\func{\vbsec{j_r}}{\point}}}
{\suc{\func{\lframedual{\vbsec{}}{k_1}}{\point}}{\func{\lframedual{\vbsec{}}{k_s}}{\point}}}},
\end{align}
which is a smooth map from $\subman{\vbbase{}}{\U}$ to $\algfield{}$, that is an element of
$\mapdifclass{\infty}{\subman{\vbbase{}}{\U}}{\vecsmanifold{\algfield{}}}$.\\
\caution
Notice that the tensor products here, take place in the tensor algebra of the vector bundle
$\vbResbundle{\vbundle{}}{\U}$. So $\vbtensor{}$ stands for $\vbtensor{\vbResbundle{\vbundle{}}{\U}}$ here.\\
\caution
Such a local representation of a tensor field is called a
$\quotl$local expression of the tensor field $\atf{}$ in terms of the local frame field of $\vbsec{}$ of $\vbundle{}$$\quotr$.
}
\endthm
%%%%%%%%%%%%%%%%%%%%%%%%%%%%%%%%%%%%%%%%%%%%%%%%%%%%%%%%%%%%%%%%%%%%%%%%%%%%%%%%%%%%%%%%%%%%%%%%%%%%%%%%%%%%%%%%%%%%%%%%%%%%%%%%
%%%%%%%%%%%%%%%%%%%%%%%%%%%%%%%%%%%%%%%%%%%%%%%%%%%%%%%%%%%%%%%%%%%%%%%%%%%%%%%%%%%%%%%%%%%%%%%%%%%%%%%%%%%%%%%%%%%%%%%%%%%%%%%%
%%%%%%%%%%%%%%%%%%%%%%%%%%%%%%%%%%%%%%%%%%%%%%%%%%%%%%%%%%%%%%%%%%%%%%%%%%%%%%%%%%%%%%%%%%%%%%%%%%%%%%%%%%%%%%%%%%%%%%%%%%%%%%%%
\subsection{Whitney-Sum of Smooth Vector Bundles}
%%%%%%%%%%%%%%%%%%%%%%%%
\fixed
In this section,
$\vbundle{1}=\quintuple{\vbtotal{1}}{\vbprojection{1}}{\vbbase{1}}{\vbfiber{}}{\vbatlas{1}}$ and
$\vbundle{2}=\quintuple{\vbtotal{2}}{\vbprojection{2}}{\vbbase{2}}{\vbfiber{}}{\vbatlas{2}}$
are fixed as a pair of smooth vector bundles
of rank $d_i$,
where $\vbtotal{i}=\opair{\vTot{i}}{\maxatlas{\vTot{i}}}$ for $i=1,~2$, and
$\vbbase{}=\opair{\vB{}}{\maxatlas{\vB{}}}$ are $\difclass{\infty}$ manifolds
modeled on the Banach-spaces $\R^{n_{\vTot{i}}}$ ($i=1,~2$) and $\R^{n_{\vB{}}}$, respectively.
\endfixed
%%%%%%%%%%%%%%%%%%%%%%%%%%%%%%%%%%%%%%%%%%%%%%%%%%%%%%%%%%%%%%
\definition
Let $\U$ be a subset of $\vB{}$.
\begin{equation}
\vbSum{\vbundle{1}}{\vbundle{2}}{\U}:=
\Union{\point}{\U}{\directsum{\fibervecs{\vbundle{1}}{\point}}{\fibervecs{\vbundle{2}}{\point}}}.
\end{equation}
Specifically, $\vbSum{\vbundle{1}}{\vbundle{2}}{\vB{}}$ can alternatively be denoted by $\VBSum{\vbundle{1}}{\vbundle{2}}$.
\endef
%%%%%%%%%%%%%%%%%%%%%%%%%%%%%%%%%%%%%%%%%%%%%%%%%%%%%%%%%%%%%%%%%%%%%%%%%%%%%%%%%%%%%%%%%%%%%%%%%%%%%%%%%%%%%%%%%%%%%%%%%%%%%%%%
%%%%%%%%%%%%%%%%%%%%%%%%%%%%%%%%%%%%%%%%%%%%%%%%%%%%%%%%%%%%%%%%%%%%%%%%%%%%%%%%%%%%%%%%%%%%%%%%%%%%%%%%%%%%%%%%%%%%%%%%%%%%%%%%
\definition\label{defSumbundlelocaltrivializations}
We associate to each element $\btriple{\U}{\phi}{\binary{\psi_1}{\psi_2}}$ of
$\vbchartlocaltt{\vbundle{1}}{\vbundle{2}}$ the mapping
$\function{\vbSumchart{\phi}{\psi_1}{\psi_2}}{\vbSum{\vbundle{1}}{\vbundle{2}}{\U}}
{\Cprod{\func{\phi}{\U}}{\(\directsum{\vbfiber{1}}{\vbfiber{2}}\)}}$ defined as,
\begin{align}
&\begin{aligned}
\Foreach{\point}{\U}
\Foreach{\opair{v_1}{v_2}}{\directsum{\fibervecs{\vbundle{1}}{\point}}{\fibervecs{\vbundle{2}}{\point}}}
\end{aligned}\cr
&\begin{aligned}
\func{\vbSumchart{\phi}{\psi_1}{\psi_2}}{\binary{v_1}{v_2}}\eqdef
\opair{\func{\phi}{\point}}{\opair{\func{\pltfib{\psi_1}{\point}}{v_1}}{\func{\pltfib{\psi_2}{\point}}{v_2}}}.
\end{aligned}
\end{align}
Furthermore, we associate to each element $\btriple{\U}{\phi}{\binary{\psi_1}{\psi_2}}$ of
$\vbchartlocaltt{\vbundle{1}}{\vbundle{2}}$ the mapping\\
$\function{\vbSumlocalt{\phi}{\psi_1}{\psi_2}}{\vbHom{\vbundle{1}}{\vbundle{2}}{\U}}
{\Cprod{\U}{\(\directsum{\vbfiber{1}}{\vbfiber{2}}\)}}$ defined as,
\begin{align}
&\begin{aligned}
\Foreach{\point}{\U}
\Foreach{\opair{v_1}{v_2}}{\directsum{\fibervecs{\vbundle{1}}{\point}}{\fibervecs{\vbundle{2}}{\point}}}
\end{aligned}\cr
&\begin{aligned}
\func{\vbSumlocalt{\phi}{\psi_1}{\psi_2}}{\binary{v_1}{v_2}}\eqdef
\opair{\point}{\opair{\func{\pltfib{\psi_1}{\point}}{v_1}}{\func{\pltfib{\psi_2}{\point}}{v_2}}}.
\end{aligned}
\end{align}
\endef
%%%%%%%%%%%%%%%%%%%%%%%%%%%%%%%%%%%%%%%%%%%%%%%%%%%%%%%%%%%%%%%%%%%%%%%%%%%%%%%%%%%%%%%%%%%%%%%%%%%%%%%%%%%%%%%%%%%%%%%%%%%%%%%%
\textcolor{Blue}{\lemma
The set $\defSet{\vbSumchart{\phi}{\psi_1}{\psi_2}}{\triple{\phi}{\psi_1}{\psi_2}\in\vbchartlocaltt{\vbundle{1}}{\vbundle{2}}}$
is a $\difclass{\infty}$ atlas on $\VBSum{\vbundle{1}}{\vbundle{2}}$ modeled on
the Banach-space $\Cprod{\R^{n_{\vB{}}}}{\(\directsum{\vbfiber{1}}{\vbfiber{2}}\)}$.}
\endlem
%%%%%%%%%%%%%%%%%%%%%%%%%%%%%%%%%%%%%%%%%%%%%%%%%%%%%%%%%%%%%%%%%%%%%%%%%%%%%%%%%%%%%%%%%%%%%%%%%%%%%%%%%%%%%%%%%%%%%%%%%%%%%%%%
\definition
The $\difclass{\infty}$ maximal-atlas on $\VBSum{\vbundle{1}}{\vbundle{2}}$ modeled on
the Banach-space $\Cprod{\R^{n_{\vB{}}}}{\(\directsum{\vbfiber{1}}{\vbfiber{2}}\)}$ generated by the atlas
$\defSet{\vbSumchart{\phi}{\psi_1}{\psi_2}}{\triple{\phi}{\psi_1}{\psi_2}\in\vbchartlocaltt{\vbundle{1}}{\vbundle{2}}}$
will be denoted by
$\vbSummaxatlas{\vbundle{1}}{\vbundle{2}}$. That is,
\begin{equation}
\vbSummaxatlas{\vbundle{1}}{\vbundle{2}}:=
\maxatlasgen{\infty}{\VBSum{\vbundle{1}}{\vbundle{2}}}{\Cprod{\R^{n_{\vB{}}}}{\(\directsum{\vbfiber{1}}{\vbfiber{2}}\)}}{
\defSet{\vbSumchart{\phi}{\psi_1}{\psi_2}}{\triple{\phi}{\psi_1}{\psi_2}\in\vbchartlocaltt{\vbundle{1}}{\vbundle{2}}}}.
\end{equation}
\endef
%%%%%%%%%%%%%%%%%%%%%%%%%%%%%%%%%%%%%%%%%%%%%%%%%%%%%%%%%%%%%%%%%%%%%%%%%%%%%%%%%%%%%%%%%%%%%%%%%%%%%%%%%%%%%%%%%%%%%%%%%%%%%%%%
\textcolor{Blue}{\theorem
The differentiable structure $\opair{\VBSum{\vbundle{1}}{\vbundle{2}}}{\vbSummaxatlas{\vbundle{1}}{\vbundle{2}}}$
is a $\difclass{\infty}$ manifold, which means the topology induced by the maximal atlas
$\vbSummaxatlas{\vbundle{1}}{\vbundle{2}}$
on $\VBSum{\vbundle{1}}{\vbundle{2}}$ is Hausdorff and second-countable.}
\endthm
%%%%%%%%%%%%%%%%%%%%%%%%%%%%%%%%%%%%%%%%%%%%%%%%%%%%%%%%%%%%%%%%%%%%%%%%%%%%%%%%%%%%%%%%%%%%%%%%%%%%%%%%%%%%%%%%%%%%%%%%%%%%%%%%
\definition
The manifold $\opair{\VBSum{\vbundle{1}}{\vbundle{2}}}{\vbSummaxatlas{\vbundle{1}}{\vbundle{2}}}$
will be denoted by $\vbSumMan{\vbundle{1}}{\vbundle{2}}$.
\endef
%%%%%%%%%%%%%%%%%%%%%%%%%%%%%%%%%%%%%%%%%%%%%%%%%%%%%%%%%%%%%%%%%%%%%%%%%%%%%%%%%%%%%%%%%%%%%%%%%%%%%%%%%%%%%%%%%%%%%%%%%%%%%%%%
\definition
The mapping $\function{\vbSumprojection{\vbundle{1}}{\vbundle{2}}}{\VBSum{\vbundle{1}}{\vbundle{2}}}{\vB{}}$ is
defined as,
\begin{equation}
\Foreach{\point}{\vB{}}
\Foreach{\opair{v_1}{v_2}}{\directsum{\fibervecs{\vbundle{1}}{\point}}{\fibervecs{\vbundle{2}}{\point}}}
\func{\vbSumprojection{\vbundle{1}}{\vbundle{2}}}{\binary{v_1}{v_2}}\eqdef\point.
\end{equation}
\endef
%%%%%%%%%%%%%%%%%%%%%%%%%%%%%%%%%%%%%%%%%%%%%%%%%%%%%%%%%%%%%%%%%%%%%%%%%%%%%%%%%%%%%%%%%%%%%%%%%%%%%%%%%%%%%%%%%%%%%%%%%%%%%%%%
\lemma\label{lemSumbundleisafiberbundle}
The quadruple
$\quadruple{\vbSumMan{\vbundle{1}}{\vbundle{2}}}{\vbSumprojection{\vbundle{1}}{\vbundle{2}}}{\vbbase{}}
{\vecsmanifold{\directsum{\vbfiber{1}}{\vbfiber{2}}}}$
is a smooth fiber bundle.
\proof
It is similar to the proof of \reflem{lemHombundleisafiberbundle}.
\endlem
%%%%%%%%%%%%%%%%%%%%%%%%%%%%%%%%%%%%%%%%%%%%%%%%%%%%%%%%%%%%%%%%%%%%%%%%%%%%%%%%%%%%%%%%%%%%%%%%%%%%%%%%%%%%%%%%%%%%%%%%%%%%%%%%
\corollary
The set $\defSet{\vbSumlocalt{\phi}{\psi_1}{\psi_2}}{\triple{\phi}{\psi_1}{\psi_2}\in\vbchartlocaltt{\vbundle{1}}{\vbundle{2}}}$
is an atlas of the smooth fiber bundle
$\quadruple{\vbSumMan{\vbundle{1}}{\vbundle{2}}}{\vbSumprojection{\vbundle{1}}{\vbundle{2}}}{\vbbase{}}
{\vecsmanifold{\directsum{\vbfiber{1}}{\vbfiber{2}}}}$.
\endcor
%%%%%%%%%%%%%%%%%%%%%%%%%%%%%%%%%%%%%%%%%%%%%%%%%%%%%%%%%%%%%%%%%%%%%%%%%%%%%%%%%%%%%%%%%%%%%%%%%%%%%%%%%%%%%%%%%%%%%%%%%%%%%%%%
\proposition
Let $\btriple{\U}{\phi}{\binary{\psi_1}{\psi_2}}$ and $\btriple{\V}{\p{\phi}}{\binary{\eta_1}{\eta_2}}$ be a pair of elements of
$\vbchartlocaltt{\vbundle{1}}{\vbundle{2}}$.
$\function{\cmp{\vbSumlocalt{\phi}{\psi_1}{\psi_2}}{\finv{\(\vbSumlocalt{\p{\phi}}{\eta_1}{\eta_2}\)}}}
{\Cprod{\(\U\cap\V\)}{\directsum{\vbfiber{1}}{\vbfiber{2}}}}{\Cprod{\(\U\cap\V\)}{\directsum{\vbfiber{1}}{\vbfiber{2}}}}$, and
the value of the transition map of the local trivializations $\vbSumlocalt{\phi}{\psi_1}{\psi_2}$ and
$\vbSumlocalt{\p{\phi}}{\eta_1}{\eta_2}$
of the fiber bundle $\quadruple{\vbSumMan{\vbundle{1}}{\vbundle{2}}}{\vbSumprojection{\vbundle{1}}{\vbundle{2}}}{\vbbase{}}
{\vecsmanifold{\directsum{\vbfiber{1}}{\vbfiber{2}}}}$ at every point $\point$ of its domain $\U\cap\V$ is a linear isomorphism
from $\directsum{\vbfiber{1}}{\vbfiber{2}}$ to itself. That is,
\begin{align}
\Foreach{\point}{\U\cap\V}
\func{\[\transition{}{\vbSumlocalt{\phi}{\psi_1}{\psi_2}}{\vbSumlocalt{\p{\phi}}{\eta_1}{\eta_2}}\]}{\point}\in
\GL{\directsum{\vbfiber{1}}{\vbfiber{2}}}{}.
\end{align}
\proof
It is trivial.
\endlem
%%%%%%%%%%%%%%%%%%%%%%%%%%%%%%%%%%%%%%%%%%%%%%%%%%%%%%%%%%%%%%%%%%%%%%%%%%%%%%%%%%%%%%%%%%%%%%%%%%%%%%%%%%%%%%%%%%%%%%%%%%%%%%%%
\definition
We will denote by $\vbSumatlas{\vbundle{1}}{\vbundle{2}}$
the maximal atlas of the smooth fiber bundle\\
$\quadruple{\vbSumMan{\vbundle{1}}{\vbundle{2}}}{\vbSumprojection{\vbundle{1}}{\vbundle{2}}}{\vbbase{}}
{\vecsmanifold{\directsum{\vbfiber{1}}{\vbfiber{2}}}}$
including $\defSet{\vbSumlocalt{\phi}{\psi_1}{\psi_2}}{\triple{\phi}{\psi_1}{\psi_2}\in\vbchartlocaltt{\vbundle{1}}{\vbundle{2}}}$,
endowed with which, the fiber bundle
$\quadruple{\vbSumMan{\vbundle{1}}{\vbundle{2}}}{\vbSumprojection{\vbundle{1}}{\vbundle{2}}}{\vbbase{}}
{\vecsmanifold{\directsum{\vbfiber{1}}{\vbfiber{2}}}}$ becomes a smooth vector bundle.
\endef
%%%%%%%%%%%%%%%%%%%%%%%%%%%%%%%%%%%%%%%%%%%%%%%%%%%%%%%%%%%%%%%%%%%%%%%%%%%%%%%%%%%%%%%%%%%%%%%%%%%%%%%%%%%%%%%%%%%%%%%%%%%%%%%%
\corollary
The quintuple
$\quintuple{\vbSumMan{\vbundle{1}}{\vbundle{2}}}{\vbSumprojection{\vbundle{1}}{\vbundle{2}}}{\vbbase{}}
{\vecsmanifold{\directsum{\vbfiber{1}}{\vbfiber{2}}}}{\vbSumatlas{\vbundle{1}}{\vbundle{2}}}$ is a smooth vector bundle.
\endcor
%%%%%%%%%%%%%%%%%%%%%%%%%%%%%%%%%%%%%%%%%%%%%%%%%%%%%%%%%%%%%%%%%%%%%%%%%%%%%%%%%%%%%%%%%%%%%%%%%%%%%%%%%%%%%%%%%%%%%%%%%%%%%%%%
\definition
We will denote by $\vbSumbundle{\vbundle{1}}{\vbundle{2}}$ the smooth vector bundle\\
$\quintuple{\vbSumMan{\vbundle{1}}{\vbundle{2}}}{\vbSumprojection{\vbundle{1}}{\vbundle{2}}}{\vbbase{}}
{\vecsmanifold{\directsum{\vbfiber{1}}{\vbfiber{2}}}}{\vbSumatlas{\vbundle{1}}{\vbundle{2}}}$, which is referred to as the
$\quotl$Whitney-sum of the smooth vector bundles $\vbundle{1}$ and $\vbundle{2}$$\quotr$.
\endef
%%%%%%%%%%%%%%%%%%%%%%%%%%%%%%%%%%%%%%%%%%%%%%%%%%%%%%%%%%%%%%%%%%%%%%%%%%%%%%%%%%%%%%%%%%%%%%%%%%%%%%%%%%%%%%%%%%%%%%%%%%%%%%%%
\corollary
For every point $\point$ of $\vbbase{}$, the fiber space of $\vbSumbundle{\vbundle{1}}{\vbundle{2}}$ over $\point$ equals
the direct sum of the fiber space of $\vbundle{1}$ over $\point$ and
the fiber space of $\vbundle{2}$ over $\point$. That is,
\begin{align}
\Foreach{\point}{\vB{}}
\fibervecs{\vbSumbundle{\vbundle{1}}{\vbundle{2}}}{\point}=
\directsum{\fibervecs{\vbundle{1}}{\point}}{\fibervecs{\vbundle{2}}{\point}}.
\end{align}
\endcor
%%%%%%%%%%%%%%%%%%%%%%%%%%%%%%%%%%%%%%%%%%%%%%%%%%%%%%%%%%%%%%%%%%%%%%%%%%%%%%%%%%%%%%%%%%%%%%%%%%%%%%%%%%%%%%%%%%%%%%%%%%%%%%%%
\textcolor{Blue}{\theorem
Suppose that $\vbundle{1}$ and $\vbundle{2}$ are smooth vector sub-bundles of the vector bundle $\vbundle{}$. If
\begin{equation}
\Foreach{\point}{\vB{}}
\fibervecs{\vbundle{}}{\point}=
\directsum{\fibervecs{\vbundle{1}}{\point}}{\fibervecs{\vbundle{2}}{\point}},
\end{equation}
then the smooth vector bundles $\vbundle{}$ and $\vbSumbundle{\vbundle{1}}{\vbundle{2}}$ are isomorphic.
}
\endthm
%%%%%%%%%%%%%%%%%%%%%%%%%%%%%%%%%%%%%%%%%%%%%%%%%%%%%%%%%%%%%%%%%%%%%%%%%%%%%%%%%%%%%%%%%%%%%%%%%%%%%%%%%%%%%%%%%%%%%%%%%%%%%%%%
%%%%%%%%%%%%%%%%%%%%%%%%%%%%%%%%%%%%%%%%%%%%%%%%%%%%%%%%%%%%%%%%%%%%%%%%%%%%%%%%%%%%%%%%%%%%%%%%%%%%%%%%%%%%%%%%%%%%%%%%%%%%%%%%
%%%%%%%%%%%%%%%%%%%%%%%%%%%%%%%%%%%%%%%%%%%%%%%%%%%%%%%%%%%%%%%%%%%%%%%%%%%%%%%%%%%%%%%%%%%%%%%%%%%%%%%%%%%%%%%%%%%%%%%%%%%%%%%%
%%%%%%%%%%%%%%%%%%%%%%%%%%%%%%%%%%%%%%%%%%%%%%%%%%%%%%%%%%%%%%%%%%%%%%%%%%%%%%%%%%%%%%%%%%%%%%%%%%%%%%%%%%%%%%%%%%%%%%%%%%%%%%%%
%%%%%%%%%%%%%%%%%%%%%%%%%%%%%%%%%%%%%%%%%%%%%%%%%%%%%%%%%%%%%%%%%%%%%%%%%%%%%%%%%%%%%%%%%%%%%%%%%%%%%%%%%%%%%%%%%%%%%%%%%%%%%%%%
\subsection{Direct Product of Smooth Vector Bundles}
%%%%%%%%%%%%%%%%%%%%%%%%%%%%%%%%%%%%%%%%%%%%%%%%%%%%%%%%%%%%%%
\definition
Let $\U$ be a subset of $\vB{}$.
\begin{equation}
\vbPro{\vbundle{1}}{\vbundle{2}}{\U_1}{\U_2}:=
\Union{\opair{\point_1}{\point_2}}{\Cprod{\U_1}{\U_2}}{\directsum{\fibervecs{\vbundle{1}}{\point_1}}{\fibervecs{\vbundle{2}}{\point_2}}}.
\end{equation}
Specifically, $\vbPro{\vbundle{1}}{\vbundle{2}}{\vB{1}}{\vB{2}}$ can alternatively be denoted by $\VBPro{\vbundle{1}}{\vbundle{2}}$.
\endef
%%%%%%%%%%%%%%%%%%%%%%%%%%%%%%%%%%%%%%%%%%%%%%%%%%%%%%%%%%%%%%%%%%%%%%%%%%%%%%%%%%%%%%%%%%%%%%%%%%%%%%%%%%%%%%%%%%%%%%%%%%%%%%%%
\definition
\begin{align}
\provbchartlocaltt{\vbundle{1}}{\vbundle{2}}:=
\defSet{\quadruple{\phi_1}{\phi_2}{\psi_1}{\psi_2}}{\[\phi_i\in\maxatlas{\vB{i}},~
\psi_i\in\vbatlas{i},~\domain{\phi_i}=\func{\image{\vbprojection{i}}}{\domain{\psi_i}}~for~i=1,~2\]}.
\end{align}
An element $\quadruple{\phi_1}{\phi_2}{\psi_1}{\psi_2}$ of $\provbchartlocaltt{\vbundle{1}}{\vbundle{2}}$ can
alternatively be denoted by
$\btriple{\binary{\U_1}{\U_2}}{\binary{\phi_1}{\phi_2}}{\binary{\psi_1}{\psi_2}}$ where
$\U_i:=\domain{\phi_i}=\func{\image{\vbprojection{i}}}{\domain{\psi_i}}$ for $i=1,~2$.
\endef
%%%%%%%%%%%%%%%%%%%%%%%%%%%%%%%%%%%%%%%%%%%%%%%%%%%%%%%%%%%%%%%%%%%%%%%%%%%%%%%%%%%%%%%%%%%%%%%%%%%%%%%%%%%%%%%%%%%%%%%%%%%%%%%%
\lemma
The set\\
$\defsets{\Cprod{\U_1}{\U_2}}{\Cprod{\vB{1}}{\vB{2}}}
{\[\Exists{\quadruple{\phi_1}{\phi_2}{\psi_1}{\psi_2}}{\provbchartlocaltt{\vbundle{1}}{\vbundle{2}}}\U_i=\domain{\phi_i}=
\func{\image{\vbprojection{i}}}{\domain{\psi_i}}~(i=1,~2)\]}$ is an open covering of $\manprod{\vbbase{1}}{\vbbase{2}}$.
\proof
It is trivial.
\endlem
%%%%%%%%%%%%%%%%%%%%%%%%%%%%%%%%%%%%%%%%%%%%%%%%%%%%%%%%%%%%%%%%%%%%%%%%%%%%%%%%%%%%%%%%%%%%%%%%%%%%%%%%%%%%%%%%%%%%%%%%%%%%%%%%
\definition\label{defProbundlelocaltrivializations}
We associate to each element $\btriple{\opair{\U_1}{\U_2}}{\binary{\phi_1}{\phi_2}}{\binary{\psi_1}{\psi_2}}$ of
$\provbchartlocaltt{\vbundle{1}}{\vbundle{2}}$ the mapping
$\function{\vbProchart{\phi_1}{\phi_2}{\psi_1}{\psi_2}}{\vbPro{\vbundle{1}}{\vbundle{2}}{\U_1}{\U_2}}
{\Cprod{\(\Cprod{\func{\phi_1}{\U_1}}{\func{\phi_2}{\U_2}}\)}{\directsum{\vbfiber{1}}{\vbfiber{2}}}}$, defined as,
\begin{align}
&\begin{aligned}
\Foreach{\opair{\point_1}{\point_2}}{\Cprod{\U_1}{\U_2}}
\Foreach{\opair{v_1}{v_2}}{\directsum{\fibervecs{\vbundle{1}}{\point_1}}{\fibervecs{\vbundle{2}}{\point_2}}}
\end{aligned}\cr
&\begin{aligned}
\func{\vbProchart{\phi_1}{\phi_2}{\psi_1}{\psi_2}}{\binary{v_1}{v_2}}\eqdef
\opair{\opair{\func{\phi_1}{\point_1}}{\func{\phi_2}{\point_2}}}
{\opair{\func{\pltfib{\psi_1}{\point_1}}{v_1}}{\func{\pltfib{\psi_2}{\point_2}}{v_2}}}.
\end{aligned}
\end{align}
Furthermore, we associate to each element $\btriple{\opair{\U_1}{\U_2}}{\binary{\phi_1}{\phi_2}}{\binary{\psi_1}{\psi_2}}$ of
$\provbchartlocaltt{\vbundle{1}}{\vbundle{2}}$ the mapping
$\function{\vbProlocalt{\phi_1}{\phi_2}{\psi_1}{\psi_2}}{\vbPro{\vbundle{1}}{\vbundle{2}}{\U_1}{\U_2}}
{\Cprod{\(\Cprod{\U_1}{\U_2}\)}{\directsum{\vbfiber{1}}{\vbfiber{2}}}}$ defined as,
\begin{align}
&\begin{aligned}
\Foreach{\opair{\point_1}{\point_2}}{\Cprod{\U_1}{\U_2}}
\Foreach{\opair{v_1}{v_2}}{\directsum{\fibervecs{\vbundle{1}}{\point_1}}{\fibervecs{\vbundle{2}}{\point_2}}}
\end{aligned}\cr
&\begin{aligned}
\func{\vbProlocalt{\phi_1}{\phi_2}{\psi_1}{\psi_2}}{\binary{v_1}{v_2}}\eqdef
\opair{\opair{\point_1}{\point_2}}{\opair{\func{\pltfib{\psi_1}{\point_1}}{v_1}}{\func{\pltfib{\psi_2}{\point_2}}{v_2}}}.
\end{aligned}
\end{align}
\endef
%%%%%%%%%%%%%%%%%%%%%%%%%%%%%%%%%%%%%%%%%%%%%%%%%%%%%%%%%%%%%%%%%%%%%%%%%%%%%%%%%%%%%%%%%%%%%%%%%%%%%%%%%%%%%%%%%%%%%%%%%%%%%%%%
\textcolor{Blue}{\lemma
The set $\defSet{\vbProchart{\phi_1}{\phi_2}{\psi_1}{\psi_2}}{\quadruple{\phi_1}{\phi_2}{\psi_1}{\psi_2}
\in\provbchartlocaltt{\vbundle{1}}{\vbundle{2}}}$
is a $\difclass{\infty}$ atlas on the set $\VBPro{\vbundle{1}}{\vbundle{2}}$ modeled on
the Banach-space $\Cprod{\(\Cprod{\R^{n_{\vB{1}}}}{\R^{n_{\vB{2}}}}\)}{\(\directsum{\vbfiber{1}}{\vbfiber{2}}\)}$.}
\endlem
%%%%%%%%%%%%%%%%%%%%%%%%%%%%%%%%%%%%%%%%%%%%%%%%%%%%%%%%%%%%%%%%%%%%%%%%%%%%%%%%%%%%%%%%%%%%%%%%%%%%%%%%%%%%%%%%%%%%%%%%%%%%%%%%
\definition
The $\difclass{\infty}$ maximal-atlas on the set $\VBPro{\vbundle{1}}{\vbundle{2}}$ modeled on
the Banach-space $\Cprod{\(\Cprod{\R^{n_{\vB{1}}}}{\R^{n_{\vB{2}}}}\)}{\(\directsum{\vbfiber{1}}{\vbfiber{2}}\)}$
generated by the atlas
$\defSet{\vbProchart{\phi_1}{\phi_2}{\psi_1}{\psi_2}}{\quadruple{\phi_1}{\phi_2}{\psi_1}{\psi_2}
\in\provbchartlocaltt{\vbundle{1}}{\vbundle{2}}}$
will be denoted by
$\vbPromaxatlas{\vbundle{1}}{\vbundle{2}}$. That is,
\begin{align}
&~\vbPromaxatlas{\vbundle{1}}{\vbundle{2}}\cr
:=&~\maxatlasgen{\infty}{\VBPro{\vbundle{1}}{\vbundle{2}}}
{\Cprod{\(\Cprod{\R^{n_{\vB{1}}}}{\R^{n_{\vB{2}}}}\)}{\(\directsum{\vbfiber{1}}{\vbfiber{2}}\)}}{
\defSet{\vbProchart{\phi_1}{\phi_2}{\psi_1}{\psi_2}}{\quadruple{\phi_1}{\phi_2}{\psi_1}{\psi_2}
\in\provbchartlocaltt{\vbundle{1}}{\vbundle{2}}}}.\cr
&{}
\end{align}
\endef
%%%%%%%%%%%%%%%%%%%%%%%%%%%%%%%%%%%%%%%%%%%%%%%%%%%%%%%%%%%%%%%%%%%%%%%%%%%%%%%%%%%%%%%%%%%%%%%%%%%%%%%%%%%%%%%%%%%%%%%%%%%%%%%%
\textcolor{Blue}{\theorem
The differentiable structure $\opair{\VBPro{\vbundle{1}}{\vbundle{2}}}{\vbPromaxatlas{\vbundle{1}}{\vbundle{2}}}$
is a $\difclass{\infty}$ manifold, which means the topology induced by the maximal atlas
$\vbPromaxatlas{\vbundle{1}}{\vbundle{2}}$
on $\VBPro{\vbundle{1}}{\vbundle{2}}$ is Hausdorff and second-countable.}
\endthm
%%%%%%%%%%%%%%%%%%%%%%%%%%%%%%%%%%%%%%%%%%%%%%%%%%%%%%%%%%%%%%%%%%%%%%%%%%%%%%%%%%%%%%%%%%%%%%%%%%%%%%%%%%%%%%%%%%%%%%%%%%%%%%%%
\definition
The manifold $\opair{\VBPro{\vbundle{1}}{\vbundle{2}}}{\vbPromaxatlas{\vbundle{1}}{\vbundle{2}}}$
will be denoted by $\vbProMan{\vbundle{1}}{\vbundle{2}}$.
\endef
%%%%%%%%%%%%%%%%%%%%%%%%%%%%%%%%%%%%%%%%%%%%%%%%%%%%%%%%%%%%%%%%%%%%%%%%%%%%%%%%%%%%%%%%%%%%%%%%%%%%%%%%%%%%%%%%%%%%%%%%%%%%%%%%
\definition
The mapping $\function{\vbProprojection{\vbundle{1}}{\vbundle{2}}}{\VBPro{\vbundle{1}}{\vbundle{2}}}{\Cprod{\vB{1}}{\vB{2}}}$ is
defined as,
\begin{equation}
\Foreach{\opair{\point_1}{\point_2}}{\Cprod{\vB{1}}{\vB{2}}}
\Foreach{\opair{v_1}{v_2}}{\directsum{\fibervecs{\vbundle{1}}{\point_1}}{\fibervecs{\vbundle{2}}{\point_2}}}
\func{\vbProprojection{\vbundle{1}}{\vbundle{2}}}{\binary{v_1}{v_2}}\eqdef\opair{\point_1}{\point_2}.
\end{equation}
\endef
%%%%%%%%%%%%%%%%%%%%%%%%%%%%%%%%%%%%%%%%%%%%%%%%%%%%%%%%%%%%%%%%%%%%%%%%%%%%%%%%%%%%%%%%%%%%%%%%%%%%%%%%%%%%%%%%%%%%%%%%%%%%%%%%
\lemma\label{lemSumbundleisafiberbundle}
The quadruple
$\quadruple{\vbProMan{\vbundle{1}}{\vbundle{2}}}{\vbProprojection{\vbundle{1}}{\vbundle{2}}}{\manprod{\vbbase{1}}{\vbbase{2}}}
{\vecsmanifold{\directsum{\vbfiber{1}}{\vbfiber{2}}}}$
is a smooth fiber bundle.
\proof
It is similar to the proof of \reflem{lemHombundleisafiberbundle}.
\endlem
%%%%%%%%%%%%%%%%%%%%%%%%%%%%%%%%%%%%%%%%%%%%%%%%%%%%%%%%%%%%%%%%%%%%%%%%%%%%%%%%%%%%%%%%%%%%%%%%%%%%%%%%%%%%%%%%%%%%%%%%%%%%%%%%
\corollary
The set $\defSet{\vbProlocalt{\phi_1}{\phi_2}{\psi_1}{\psi_2}}{\quadruple{\phi_1}{\phi_2}{\psi_1}{\psi_2}
\in\provbchartlocaltt{\vbundle{1}}{\vbundle{2}}}$
is an atlas of the smooth fiber bundle
$\quadruple{\vbProMan{\vbundle{1}}{\vbundle{2}}}{\vbSumprojection{\vbundle{1}}{\vbundle{2}}}{\manprod{\vbbase{1}}{\vbbase{2}}}
{\vecsmanifold{\directsum{\vbfiber{1}}{\vbfiber{2}}}}$.
\endcor
%%%%%%%%%%%%%%%%%%%%%%%%%%%%%%%%%%%%%%%%%%%%%%%%%%%%%%%%%%%%%%%%%%%%%%%%%%%%%%%%%%%%%%%%%%%%%%%%%%%%%%%%%%%%%%%%%%%%%%%%%%%%%%%%
\proposition
Let $\btriple{\binary{\U_1}{\U_2}}{\binary{\phi_1}{\phi_2}}{\binary{\psi_1}{\psi_2}}$ and
$\btriple{\binary{\V_1}{\V_2}}{\binary{\p{\phi_1}}{\p{\phi_2}}}{\binary{\eta_1}{\eta_2}}$ be a pair of elements of
$\provbchartlocaltt{\vbundle{1}}{\vbundle{2}}$.
%such that there exists a pair $\opair{\U}{\phi}$ and $\opair{\V}{\p{\phi}}$ of charts of the manifold $\vbbase{}$,
$\function{\cmp{\vbProlocalt{\phi_1}{\phi_2}{\psi_1}{\psi_2}}{\finv{\(\vbProlocalt{\p{\phi_1}}{\p{\phi_2}}{\eta_1}{\eta_2}\)}}}
{\Cprod{\(\(\Cprod{\U_1}{\U_2}\)\cap\(\Cprod{\V_1}{\V_2}\)\)}
{\directsum{\vbfiber{1}}{\vbfiber{2}}}}{\Cprod{\(\(\Cprod{\U_1}{\U_2}\)\cap\(\Cprod{\V_1}{\V_2}\)\)}
{\directsum{\vbfiber{1}}{\vbfiber{2}}}}$, and
the value of the transition map of the local trivializations $\vbProlocalt{\phi_1}{\phi_2}{\psi_1}{\psi_2}$ and
$\vbProlocalt{\p{\phi_1}}{\p{\phi_2}}{\eta_1}{\eta_2}$
of the fiber bundle
$\quadruple{\vbProMan{\vbundle{1}}{\vbundle{2}}}{\vbProprojection{\vbundle{1}}{\vbundle{2}}}{\manprod{\vbbase{1}}{\vbbase{2}}}
{\vecsmanifold{\directsum{\vbfiber{1}}{\vbfiber{2}}}}$ at every point $\opair{\point_1}{\point_2}$ of its domain
$\(\Cprod{\U_1}{\U_2}\)\cap\(\Cprod{\V_1}{\V_2}\)$ is a linear isomorphism
from $\directsum{\vbfiber{1}}{\vbfiber{2}}$ to itself. That is,
\begin{align}
\Foreach{\opair{\point_1}{\point_2}}{\(\Cprod{\U_1}{\U_2}\)\cap\(\Cprod{\V_1}{\V_2}\)}
\func{\[\transition{}{\vbProlocalt{\phi_1}{\phi_2}{\psi_1}{\psi_2}}{\vbProlocalt{\p{\phi_1}}{\p{\phi_2}}{\eta_1}{\eta_2}}\]}{\point}\in
\GL{\directsum{\vbfiber{1}}{\vbfiber{2}}}{}.
\end{align}
\proof
It is trivial.
\endlem
%%%%%%%%%%%%%%%%%%%%%%%%%%%%%%%%%%%%%%%%%%%%%%%%%%%%%%%%%%%%%%%%%%%%%%%%%%%%%%%%%%%%%%%%%%%%%%%%%%%%%%%%%%%%%%%%%%%%%%%%%%%%%%%%
\definition
We will denote by $\vbProatlas{\vbundle{1}}{\vbundle{2}}$
the maximal atlas of the smooth fiber bundle\\
$\quadruple{\vbProMan{\vbundle{1}}{\vbundle{2}}}{\vbProprojection{\vbundle{1}}{\vbundle{2}}}{\vbbase{}}
{\vecsmanifold{\directsum{\vbfiber{1}}{\vbfiber{2}}}}$
including $\defSet{\vbProlocalt{\phi_1}{\phi_2}{\psi_1}{\psi_2}}{\quadruple{\phi_1}{\phi_2}{\psi_1}{\psi_2}\in
\provbchartlocaltt{\vbundle{1}}{\vbundle{2}}}$,
endowed with which, the fiber bundle
$\quadruple{\vbProMan{\vbundle{1}}{\vbundle{2}}}{\vbProprojection{\vbundle{1}}{\vbundle{2}}}{\manprod{\vbbase{1}}{\vbbase{2}}}
{\vecsmanifold{\directsum{\vbfiber{1}}{\vbfiber{2}}}}$ becomes a smooth vector bundle.
\endef
%%%%%%%%%%%%%%%%%%%%%%%%%%%%%%%%%%%%%%%%%%%%%%%%%%%%%%%%%%%%%%%%%%%%%%%%%%%%%%%%%%%%%%%%%%%%%%%%%%%%%%%%%%%%%%%%%%%%%%%%%%%%%%%%
\corollary
The quintuple
$\quintuple{\vbProMan{\vbundle{1}}{\vbundle{2}}}{\vbProprojection{\vbundle{1}}{\vbundle{2}}}{\manprod{\vbbase{1}}{\vbbase{2}}}
{\vecsmanifold{\directsum{\vbfiber{1}}{\vbfiber{2}}}}{\vbProatlas{\vbundle{1}}{\vbundle{2}}}$ is a smooth vector bundle.
\endcor
%%%%%%%%%%%%%%%%%%%%%%%%%%%%%%%%%%%%%%%%%%%%%%%%%%%%%%%%%%%%%%%%%%%%%%%%%%%%%%%%%%%%%%%%%%%%%%%%%%%%%%%%%%%%%%%%%%%%%%%%%%%%%%%%
\definition
We will denote by $\vbProbundle{\vbundle{1}}{\vbundle{2}}$ the smooth vector bundle\\
$\quintuple{\vbProMan{\vbundle{1}}{\vbundle{2}}}{\vbProprojection{\vbundle{1}}{\vbundle{2}}}{\manprod{\vbbase{1}}{\vbbase{2}}}
{\vecsmanifold{\directsum{\vbfiber{1}}{\vbfiber{2}}}}{\vbProatlas{\vbundle{1}}{\vbundle{2}}}$, which is referred to as the
$\quotl$direct product of the smooth vector bundles $\vbundle{1}$ and $\vbundle{2}$$\quotr$.
\endef
%%%%%%%%%%%%%%%%%%%%%%%%%%%%%%%%%%%%%%%%%%%%%%%%%%%%%%%%%%%%%%%%%%%%%%%%%%%%%%%%%%%%%%%%%%%%%%%%%%%%%%%%%%%%%%%%%%%%%%%%%%%%%%%%
\corollary
For every point $\point_1$ of $\vbbase{1}$ and every point $\point_2$ of $\vbbase{2}$,
the fiber space of $\vbProbundle{\vbundle{1}}{\vbundle{2}}$ over $\opair{\point_1}{\point_2}$ equals
the direct sum of the fiber space of $\vbundle{1}$ over $\point_1$ and
the fiber space of $\vbundle{2}$ over $\point_2$. That is,
\begin{align}
\Foreach{\opair{\point_1}{\point_2}}{\Cprod{\vB{1}}{\vB{2}}}
\fibervecs{\vbProbundle{\vbundle{1}}{\vbundle{2}}}{\binary{\point_1}{\point_2}}=
\directsum{\fibervecs{\vbundle{1}}{\point_1}}{\fibervecs{\vbundle{2}}{\point_2}}.
\end{align}
\endcor
%%%%%%%%%%%%%%%%%%%%%%%%%%%%%%%%%%%%%%%%%%%%%%%%%%%%%%%%%%%%%%%%%%%%%%%%%%%%%%%%%%%%%%%%%%%%%%%%%%%%%%%%%%%%%%%%%%%%%%%%%%%%%%%%
%%%%%%%%%%%%%%%%%%%%%%%%%%%%%%%%%%%%%%%%%%%%%%%%%%%%%%%%%%%%%%%%%%%%%%%%%%%%%%%%%%%%%%%%%%%%%%%%%%%%%%%%%%%%%%%%%%%%%%%%%%%%%%%%
%%%%%%%%%%%%%%%%%%%%%%%%%%%%%%%%%%%%%%%%%%%%%%%%%%%%%%%%%%%%%%%%%%%%%%%%%%%%%%%%%%%%%%%%%%%%%%%%%%%%%%%%%%%%%%%%%%%%%%%%%%%%%%%%
\subsection{Induced Bundle}
\fixed
$\Man{}=\opair{\M{}}{\maxatlas{M}}$ is fixed as an $n$-dimensional and $\difclass{\infty}$ manifold
modeled the Banach-space $\R^n$. $f$ is fixed as an element of $\mapdifclass{\infty}{\Man{}}{\vbbase{}}$.
\endfixed
%%%%%%%%%%%%%%%%%%%%%%%%%%%%%%%%%%%%%%%%%%%%%%%%%%%%%%%%%%%%%%%%%%%%%%%%%%%%%%%%%%%%%%%%%%%%%%%%%%%%%%%%%%%%%%%%%%%%%%%%%%%%%%%%
\definition
By definition,
\begin{align}
\inducedbundle{f}{\vbundle{}}{\Man{}}:&=
\defset{\opair{\point}{v}}{\Cprod{\M{}}{\vTot{}}}{\func{f}{\point}=\func{\vbprojection{}}{v}}\cr
&=\Union{\point}{\M{}}{\Cprod{\seta{\point}}{\func{\finv{\vbprojection{}}}{\seta{\func{f}{\point}}}}}.
\end{align}
\endef
%%%%%%%%%%%%%%%%%%%%%%%%%%%%%%%%%%%%%%%%%%%%%%%%%%%%%%%%%%%%%%%%%%%%%%%%%%%%%%%%%%%%%%%%%%%%%%%%%%%%%%%%%%%%%%%%%%%%%%%%%%%%%%%%
\definition
The mapping
$\function{\inducedvbprojection{\vbprojection{}}{}}{\inducedbundle{f}{\vbundle{}}{\Man{}}}{\Cprod{\M{}}{\vbfiber{}}}$
is defined as
\begin{align}
\Foreach{\opair{\point}{v}}{\inducedbundle{f}{\vbundle{}}{\Man{}}}
\func{\inducedvbprojection{\vbprojection{}}{}}{\binary{\point}{v}}\eqdef\point.
\end{align}
\endef
%%%%%%%%%%%%%%%%%%%%%%%%%%%%%%%%%%%%%%%%%%%%%%%%%%%%%%%%%%%%%%%%%%%%%%%%%%%%%%%%%%%%%%%%%%%%%%%%%%%%%%%%%%%%%%%%%%%%%%%%%%%%%%%%
\definition
We associate to each local trivialization $\opair{\U}{\phi}$ of the smooth vector bundle $\vbundle{}$
the mapping $\displaystyle\function{\inducedvblocalt{\phi}{f}}
{\Union{\point}{\func{\finv{f}}{\U}}{\Cprod{\seta{\point}}{\func{\finv{\vbprojection{}}}{\seta{\func{f}{\point}}}}}}
{\Cprod{\func{\finv{f}}{\U}}{\vbfiber{}}}$,
defined as
\begin{align}
\Foreach{\opair{q}{v}}{
\Union{\point}{\func{\finv{f}}{\U}}{\Cprod{\seta{\point}}{\func{\finv{\vbprojection{}}}{\seta{\func{f}{\point}}}}}}
\func{\[\inducedvblocalt{\phi}{f}\]}{\binary{q}{v}}\eqdef
\opair{q}{\func{\pltfib{\phi}{\func{f}{q}}}{v}}.
\end{align}
\endef
%%%%%%%%%%%%%%%%%%%%%%%%%%%%%%%%%%%%%%%%%%%%%%%%%%%%%%%%%%%%%%%%%%%%%%%%%%%%%%%%%%%%%%%%%%%%%%%%%%%%%%%%%%%%%%%%%%%%%%%%%%%%%%%%
\textcolor{Blue}{\lemma
\begin{itemize}
\item
$\inducedbundle{f}{\vbundle{}}{\Man{}}$ is a regular submanifold of the product manifold
$\manprod{\Man{}}{\vbtotal{}}$.
\item
The quadruple
$\quadruple{\manprod{\Man{}}{\vbtotal{}}}{\inducedvbprojection{\vbprojection{}}{}}{\Man{}}{\vecsmanifold{\vbfiber{}}}$
is a smooth fiber bundle, and
$\defSet{\inducedvblocalt{\phi}{f}}{\phi\in\vbatlas{}}$ is an atlas of this smooth fiber bundle.
\item
For every pair $\opair{\U}{\phi}$ and $\opair{\V}{\psi}$ of local trivializations of the smooth vector bundle
$\vbundle{}$,
\begin{equation}
\Foreach{\point}{\func{\finv{f}}{\U\cap\V}}
\func{\transition{}{\inducedvblocalt{\phi}{f}}{\inducedvblocalt{\psi}{f}}}{\point}\in\GL{\vbfiber{}}{}.
\end{equation}
\end{itemize}
}
\endlem
%%%%%%%%%%%%%%%%%%%%%%%%%%%%%%%%%%%%%%%%%%%%%%%%%%%%%%%%%%%%%%%%%%%%%%%%%%%%%%%%%%%%%%%%%%%%%%%%%%%%%%%%%%%%%%%%%%%%%%%%%%%%%%%%
\definition
We will denote by $\vbIndatlas{\vbundle{}}{\Man{}}$
the maximal atlas of the smooth fiber bundle\\
$\quadruple{\manprod{\Man{}}{\vbtotal{}}}{\inducedvbprojection{\vbprojection{}}{}}{\Man{}}{\vecsmanifold{\vbfiber{}}}$
including $\defSet{\inducedvblocalt{\phi}{f}}{\phi\in\vbatlas{}}$,
endowed with which, this fiber bundle becomes a smooth vector bundle.
\endef
%%%%%%%%%%%%%%%%%%%%%%%%%%%%%%%%%%%%%%%%%%%%%%%%%%%%%%%%%%%%%%%%%%%%%%%%%%%%%%%%%%%%%%%%%%%%%%%%%%%%%%%%%%%%%%%%%%%%%%%%%%%%%%%%
\corollary
The quintuple
$\quintuple{\manprod{\Man{}}{\vbtotal{}}}{\inducedvbprojection{\vbprojection{}}{}}{\Man{}}{\vecsmanifold{\vbfiber{}}}
{\vbIndatlas{\vbundle{1}}{\Man{}}}$ is a smooth vector bundle.
\endcor
%%%%%%%%%%%%%%%%%%%%%%%%%%%%%%%%%%%%%%%%%%%%%%%%%%%%%%%%%%%%%%%%%%%%%%%%%%%%%%%%%%%%%%%%%%%%%%%%%%%%%%%%%%%%%%%%%%%%%%%%%%%%%%%%
\definition
We will denote by $\vbIndbundle{\vbundle{}}{\Man{}}{f}$ the smooth vector bundle\\
$\quintuple{\manprod{\Man{}}{\vbtotal{}}}{\inducedvbprojection{\vbprojection{}}{}}{\Man{}}{\vecsmanifold{\vbfiber{}}}
{\vbIndatlas{\vbundle{1}}{\Man{}}}$, which is referred to as the
$\quotl$induced bundle from the smooth vector bundle $\vbundle{}$ on the manifold $\Man{}$ via the smooth map $f$$\quotr$.
\endef
%%%%%%%%%%%%%%%%%%%%%%%%%%%%%%%%%%%%%%%%%%%%%%%%%%%%%%%%%%%%%%%%%%%%%%%%%%%%%%%%%%%%%%%%%%%%%%%%%%%%%%%%%%%%%%%%%%%%%%%%%%%%%%%%
\corollary
For every point $\point$ of $\Man{}$,
the fiber space of $\vbIndbundle{\vbundle{}}{\Man{}}{f}$ over $\point$ has the following linear structure.
\begin{align}
&\Foreach{\opair{v_1}{v_2}}{\Cprod{\fibervecs{\vbundle{}}{\func{f}{\point}}}{\fibervecs{\vbundle{}}{\func{f}{\point}}}}
\Foreach{\opair{c_1}{c_2}}{\Cprod{\algfield{}}{\algfield{}}}\cr
&c_1\opair{\point}{v_1}+c_2\opair{\point}{v_2}=
\opair{\point}{c_1v_1+c_2v_2}.
\end{align}
\endcor
\newpage
\Bibliography{}
\renewcommand{\addcontentsline}[3]{}

\let\addcontentsline\oldaddcontentsline
%%%%%%%%%%%%%%%%%%%%%%%%%%%%%%%%%%%%%%%%%%%%%%%%%%%%%%%%%%%%%%%%%%%%%%%%%%%%%%%%%%%%%%%%%%%%%%%%%%%%%%%


\begin{thebibliography}{99}
\bibitem{Berger}
Marcel Berger and Bernard Gostiaux, $\quotl$Differential Geometry: Manifolds, Curves, and Surfaces$\quotr$;
Springer-Verlag,
%New York,
1988
\bibitem{JLee}
Jeffrey M. Lee, $\quotl$Manifolds and Differential Geometry$\quotr$; 
American Mathematical Society,
%Rhode Islands,
2009
\bibitem{Lang}
Serge Lang, $\quotl$Fundamentals of Differential Geometry$\quotr$; Springer-Verlag,
%New York,
1999
\bibitem{Milnor}
John W. Milnor and  James D. Stasheff, $\quotl$Characteristic Classes$\quotr$; Princeton University Press,
%Princeton, New Jersey,
1974
\bibitem{Steenrod}
Norman Steenrod, $\quotl$The Topology of Fibre Bundles$\quotr$; Princeton University Press,
%Princeton, New Jersey,
1951
\bibitem{Husemoller}
Dale Husemoller, $\quotl$Fibre Bundles$\quotr$, Third Edition; Springer-Verlag,
%New York,
1994
\bibitem{Hatcher}
Allen Hatcher, $\quotl$Vector Bundles and K-Theory$\quotr$, Version 2.2, 2017
\bibitem{Cartan}
Henri Cartan, $\quotl$Differential Calculus$\quotr$; Hermann,
%Paris,
1971
\bibitem{Cartandifforms}
Henri Cartan, $\quotl$Differential Forms$\quotr$; Hermann,
%Paris,
1967
\bibitem{Viro}
O. Ya. Viro, O. A. Ivanov, N. Yu. Netsvetaev, V. M. Kharlamov
$\quotl$Elementary Topology: Problem Textbook$\quotr$; American Mathematical Society,
%Rhode Islands,
2008
\bibitem{Munkres}
James R. Munkres, $\quotl$Topology$\quotr$, Second Edition; Prentice Hall, 2000
\bibitem{Hoffman}
Kenneth Hoffman and Ray Kunze, $\quotl$Linear Algebra$\quotr$, Second Edition;  Prentice-Hall, 1971
\bibitem{MacLane}
Saunders MacLane and Garrett Birkhoff, $\quotl$Algebra$\quotr$, Third Edition;
Chelsea Publishing Company,
%New York,
1988
\bibitem{Hungerford}
Thomas W. Hungerford, $\quotl$Algebra$\quotr$;
Springer-Verlog,
%New York,
1974
\end{thebibliography}
\end{document}